\documentclass[10pt]{amsart}
\usepackage[T1]{fontenc}
\usepackage[letterpaper,textheight=249mm,vcentering,left=2.5cm,right=2.5cm,headsep=10pt]{geometry}
\usepackage{amsmath,amssymb,amsthm,mathtools,mathrsfs,bm}
\usepackage{booktabs,longtable,array,enumitem,xcolor}
\usepackage{needspace}
\usepackage[hidelinks]{hyperref}
\hypersetup{unicode=true,
 pdftitle={Fibrewise Real-Analytic Landsberg Metrics in Dimension at Least Four Are Berwald},
 pdfauthor={Ming Li and Kefeng Liu},
 pdfsubject={English mathematical manuscript}}
\allowdisplaybreaks[2]
\makeatletter
\renewcommand{\@pnumwidth}{3em}
\makeatother
\newtheorem{theorem}{Theorem}[section]
\newtheorem{lemma}[theorem]{Lemma}
\newtheorem{proposition}[theorem]{Proposition}
\newtheorem{corollary}[theorem]{Corollary}
\theoremstyle{definition}
\newtheorem{definition}[theorem]{Definition}

\theoremstyle{remark}
\newtheorem{remark}[theorem]{Remark}
\newcommand{\R}{\mathbb R}\newcommand{\C}{\mathbb C}\newcommand{\J}{\mathcal J}
\newcommand{\dd}{\mathrm d}\newcommand{\Id}{\operatorname{Id}}\newcommand{\Sym}{\operatorname{Sym}}
\newcommand{\Hess}{\operatorname{Hess}}\newcommand{\tr}{\operatorname{tr}}\newcommand{\rank}{\operatorname{rank}}
\newcommand{\rad}{\operatorname{rad}}\newcommand{\Cent}{\operatorname{Cent}}\newcommand{\Ric}{\operatorname{Ric}}\newcommand{\Scal}{\operatorname{Scal}}
\newcommand{\ip}[2]{\langle #1,#2\rangle}
\numberwithin{equation}{section}
\title[Fibrewise real-analytic Landsberg metrics]{Fibrewise Real-Analytic Landsberg Metrics in Dimension at Least Four Are Berwald}
\author{Ming Li}
\author{Kefeng Liu}
\address{Mathematical Science Research Center, Chongqing University of Technology}
\email[Ming Li]{mingli@cqut.edu.cn}
\email[Kefeng Liu]{kefengliu@cqut.edu.cn}
\date{October 29, 2026}
\begin{document}
\raggedbottom
\begin{abstract}
Let $M$ be a smooth manifold of dimension $n\ge4$, and let $F$ be a smooth, strongly convex Landsberg metric, regular on the entire slit tangent bundle. We prove that $F$ is Berwald if, at each fixed base point, the Minkowski norm is real-analytic on the entire slit tangent space. The proof rests on two analytic rigidity theorems for a single Minkowski space. The first asserts that, in dimension at least $3$, a real-analytic Minkowski norm whose Hessian metric is flat is induced by an inner product. The second asserts that, in dimension at least $4$, the Hessian metric of a real-analytic Minkowski norm is flat if the entire indicatrix carries a nonzero smooth symmetric quartic form satisfying a twisted Codazzi equation. The quartic form need not be real-analytic, and the second theorem is proved independently of the first. Applying these results to the angular Berwald quartic form shows that the metric is Riemannian on the open set of base points on whose fibres the Berwald curvature is not identically zero. Since every Riemannian metric is Berwald, this set must be empty. No analytic structure, compactness or completeness of the base manifold is assumed, and the metric need not be reversible.
\end{abstract}
\maketitle
\tableofcontents
\section{Introduction and main results}\label{sec:introduction}

A Finsler metric replaces the inner-product norm of a Riemannian metric by a strongly convex Minkowski norm on each tangent space. Its canonical spray coefficients need not be quadratic in the fibre variables. A Berwald metric has spray coefficients quadratic in those variables. A Landsberg metric has vanishing Landsberg curvature; equivalently, its Berwald tensor is completely symmetric after its upper index is lowered with the fundamental tensor. Every Berwald metric is Landsberg. The converse is the Landsberg--Berwald problem, also called the unicorn problem. All metrics considered here are smooth and strongly convex on the entire slit tangent bundle. This regularity in all nonzero directions is an essential part of the problem. For the basic definitions we refer to \cite{bcs}; for the history, principal methods and related developments, see the surveys of Feng--Li \cite{feng-li2024} and Crampin \cite{crampin2024}.

Matsumoto's 1996 survey collected classical reduction results and asked for explicit examples of Landsberg metrics which are not Berwald \cite[p.~82]{matsumoto1996}. Bao called such metrics ``unicorns'' in his 2007 survey \cite{bao2007}. The existence question also appears as Problem~1 in Shen's 2009 list of open problems \cite{shen-problems2009}.

Rigidity results for regular Landsberg metrics have been obtained under additional curvature assumptions and within special classes of metrics. A classical result states that a Landsberg metric of dimension $n\ge3$ with vanishing Douglas curvature is Berwald \cite[Theorem 7]{matsumoto1996}. Numata proved that a Landsberg space of dimension $n\ge3$ with nowhere-vanishing scalar flag curvature is a Riemannian space of constant sectional curvature \cite{numata1975}. In his systematic study of $(\alpha,\beta)$-metrics, Shen proved that every regular Landsberg $(\alpha,\beta)$-metric in dimension $n\ge3$ is Berwald \cite[Theorem 1.1]{shen2009}. Starting from equivalence properties of the unit indicatrices of Minkowski norms, Feng--Han--Li extended this conclusion to general $(\alpha,\beta)$-metrics in dimension $n\ge4$ \cite{feng-han-li2021}. Xu--Deng proved that every Landsberg metric of $(\alpha_1,\alpha_2)$-type is Berwald, using a linear Landsberg equation associated with a fixed Minkowski norm and the action of its linear isometry group on the solution space \cite[Theorem 1.2 and Section 5]{xu-deng2021}.

Another approach concerns the relations of the Landsberg condition to mean curvatures and volume. Using centroaffine geometry, Li established equivalence theorems for Minkowski spaces and proved that a Landsberg metric with vanishing mean Berwald curvature is Berwald in dimension $n\ge3$ \cite{li2019,li-minkowski}. Crampin showed that it is enough for the trace of the mean Berwald curvature with respect to the fundamental tensor to vanish \cite{crampin2018}. Li--Zhang proved, for general Finsler metrics, that vanishing Berwald scalar curvature implies vanishing mean Berwald curvature, giving another derivation of this Berwald criterion for Landsberg metrics \cite[Theorems 1 and 2]{li-zhang2020}. In 2024 Crampin further proved that, on a compact manifold without boundary of dimension $n\ge3$, a regular Landsberg metric whose Berwald scalar curvature is invariant under the geodesic flow is Berwald \cite[Theorem 4]{crampin2024}. The noncompact version also requires forward or backward geodesic completeness and boundedness of the scalar $g^{ij}\partial_{y^i}\partial_{y^j}\log\det(g_{ab})$ on the unit indicatrix bundle. In 2025, using integral identities on that bundle, Li proved that a Landsberg metric on a closed connected manifold is Berwald if its Berwald scalar curvature does not change sign on the entire unit indicatrix bundle \cite[Theorem 3.7]{li2025}.

Nonregular examples and generalized structures exhibit further features of the Landsberg equations. Shen also constructs almost regular Landsberg $(\alpha,\beta)$-metrics which are not Berwald; his explicit example is singular in two extremal directions in each tangent space \cite[Theorem 1.2 and Example 1.1]{shen2009}. Further constructions and related results are surveyed in \cite{tayebi2021}. By Cartan--K\"ahler theory and path geometry, Sabau--Shibuya--Shimada established the existence of nontrivial generalized Landsberg structures \cite{sabau-shibuya-shimada2010}. Realizing such a structure as a regular Finsler metric on the entire slit tangent bundle remains a separate problem. A preprint submitted by Elgendi in October 2026 further studies nonregular, including almost regular, non-Berwald Landsberg metrics of general $(\alpha,\beta)$-type \cite{elgendi2026}. The domains or regularity of these constructions differ from the hypotheses of the present paper; they therefore do not furnish counterexamples to the regular unicorn problem. Szab\'o announced a proof in the general regular case in 2008, but his subsequent correction stated that the proof was incomplete. The conclusion retained there concerns the stronger dual Landsberg condition \cite{szabo2008,szabo2009}.

Recent work also approaches the problem through sprays and curvature distributions. Elgendi--Muzsnay studied the image, kernel and rank of the Berwald curvature of Landsberg metrics, together with their relations to the holonomy distribution and metrizability of the geodesic spray, and obtained corresponding structural restrictions \cite{elgendi-muzsnay2025}. A 2026 preprint of Mao--Zhou gives a proof that every smooth, strongly convex, regular Landsberg metric on a closed three-dimensional manifold is Berwald, without assuming reversibility or real analyticity \cite[Theorem 1.1]{mao-zhou}. Closedness enters the final compactness argument in that work, while dimension three is also used in the preceding analysis of the indicatrix and curvature algebra.

The unicorn problem is closely related to Hessian rigidity on a single Minkowski space. Laugwitz's early work connected the geometry of Minkowski indicatrices with the theory of centroaffine hypersurfaces \cite{laugwitz1957}. His subsequent conjecture asserts that, for a smooth, strongly convex Minkowski norm $F$ regular in all nonzero directions in dimension $N\ge3$, flatness of $g=\Hess(F^2/2)$ on the entire slit vector space implies that $F^2$ is a positive-definite quadratic form \cite{laugwitz1965,xu-matveev}. Only positive homogeneity is assumed, so that $F(-y)\ne F(y)$ is allowed. By the cone decomposition of the Hessian metric, this is equivalent to ellipsoid rigidity when the centroaffine metric of the indicatrix has constant sectional curvature $1$. The dimension assumption cannot be omitted: the Hessian metric of every two-dimensional Minkowski norm is locally flat, and the same conclusion is false in dimension two \cite[Section 1.1]{xu-matveev}.

Among the classical results on the Laugwitz conjecture, Brickell established rigidity in the absolutely homogeneous, or reversible, case \cite{brickell1967}. By affine-geometric methods, Schneider extended the conclusion to the case in which the centroid of the unit ball is the origin \cite{schneider1968}; see \cite{feng-li2024,li-zuo2023} for related discussions and new proofs. These assumptions concern the symmetry of the norm or the volume centroid of its unit ball and need not hold for a general nonreversible norm. Xu--Matveev proved Laugwitz rigidity for smooth, strongly convex Minkowski norms with the standard block-diagonal $SO(k)\times SO(N-k)$ symmetry. Xu subsequently extended flat-Hessian rigidity, in dimension $N\ge3$, to Minkowski norms induced by isoparametric foliations of the standard sphere \cite{xu-matveev,xu-isoparametric2022}. Li--Zuo also proved that a smooth, strongly convex, closed indicatrix of dimension at least $2$ is an ellipsoid centred at the origin if it has constant sectional curvature and its traceless Tchebychev operator is parallel with respect to the Levi-Civita connection of the centroaffine metric \cite{li-zuo2023}. These results use symmetry, an isoparametric structure or an additional condition on a centroaffine tensor.

Mao--Zhou proved the Laugwitz conjecture in ambient dimension $n=3$, or indicatrix dimension $m=2$, for smooth strongly convex Minkowski norms regular on the entire slit space, without reversibility or real analyticity \cite[Corollary 3.3]{mao-zhou}. Thus the three-dimensional case of Theorem A below is already covered by their result. Our proof treats all dimensions $n\ge3$ under the real-analyticity hypothesis.

The distinction between Hessian isometry and linear equivalence is central to the relation between these two problems. Li's equivalence theorem states that, in dimension at least $3$, a norm-preserving, positively homogeneous diffeomorphism of degree one is linear if it preserves both the Hessian metric and the Cartan trace form \cite{li2019,li-minkowski}. Preservation of the Cartan trace form cannot simply be omitted. Xu--Matveev's analysis of Hessian isometries between symmetric norms further shows the following: fix $1\le k\le n-1$ and suppose that, on a Finsler manifold of dimension $n\ge3$, every tangent-space norm has the standard $SO(k)\times SO(n-k)$ symmetry in suitable linear coordinates. Then the Landsberg condition implies the Berwald condition \cite[Corollary 1.10]{xu-matveev}. This smooth theorem includes several familiar special classes of metrics, but its fibrewise symmetry remains a substantive hypothesis.

We prove Berwald rigidity under fibrewise real analyticity in dimension $n\ge4$, without symmetry assumptions. Part I establishes two independent rigidity theorems for a single Minkowski space; Part II applies them to the angular Berwald quartic form and concludes by an argument on an open subset of the base manifold. We quote classical results and prove the lemmas specific to this argument.

\subsection{The main theorem and its hypotheses}
We say that $F$ is regular on the entire slit tangent bundle if it is smooth and strictly positive on all of $TM\setminus\{0\}$, satisfies $F(x,\lambda y)=\lambda F(x,y)$ for $\lambda>0$, and has positive-definite fundamental tensor
\[
 g_{ij}(x,y)=\frac12\frac{\partial^2 F^2}{\partial y^i\partial y^j}(x,y).
\]
We denote the canonical spray coefficients by $G^a$ and the Berwald curvature by $\mathcal B^a{}_{ijk}=\partial_{y^i}\partial_{y^j}\partial_{y^k}G^a$. Our main result is the following.

\begin{theorem}\label{thm:main}
Let $M$ be a smooth manifold of dimension $n\ge4$, and let $F$ be a smooth, strongly convex Landsberg metric regular on the entire slit tangent bundle $TM\setminus\{0\}$. Suppose that, for every $x\in M$, the function $F_x=F(x,\cdot)$ is real-analytic on all of $T_xM\setminus\{0\}$. Then $F$ is a Berwald metric, that is, $\mathcal B\equiv0$.
\end{theorem}

Analyticity is required only in the fibre variables at each fixed base point. The base manifold need not carry an analytic structure, and no uniform radius of analyticity over the fibres is assumed. We require neither compactness nor geodesic completeness of $M$, nor the identity $F(x,-y)=F(x,y)$. The compact space used in the proof is the unit indicatrix
\[
 \Sigma_x=\{y\in T_xM:F(x,y)=1\},\qquad \dim\Sigma_x=m=n-1.
\]
Regularity in all nonzero directions and strong convexity imply that it is compact and diffeomorphic to $S^m$.

The conclusion of the main theorem allows non-Riemannian Berwald metrics. For example, a nonquadratic Minkowski norm independent of the base point defines a locally Minkowski metric whose canonical spray and Berwald curvature both vanish. Quadratic fibres arise in the proof only under the temporary assumption that the Berwald quartic form is nonzero on a fibre; they are an intermediate step in excluding that assumption.

\subsection{Two analytic rigidity theorems for Minkowski spaces}
Both rigidity theorems are stated on a single vector space. Let $V$ be a real vector space of dimension $n\ge3$, let $D$ be its standard flat affine connection, and let $F$ be a Minkowski norm on $V$. Put
\[
 g=D^2\!\left(\tfrac12F^2\right),\qquad \mathcal C=\tfrac12Dg,\qquad \Sigma=\{y\in V:F(y)=1\},\qquad h=g|_{T\Sigma},\quad C=\mathcal C|_{T\Sigma}.
\]
We denote the Levi-Civita connection of $h$ by $\nabla^h$ and use
\[
 h(C_XY,Z)=C(X,Y,Z)
\]
to define the self-adjoint contraction operator $C_X$. The cone structure of the Hessian metric and the Gauss equation of the indicatrix take the form
\[
 g=\dd\rho^2+\rho^2h,\qquad R^h=R_0-[C,C],
\]
where $\rho=F$ and $R_0(X,Y)Z=h(Y,Z)X-h(X,Z)Y$. In particular, $g$ is flat if and only if $h$ has constant sectional curvature $1$.

The first result is analytic Laugwitz rigidity, which we shall also call Theorem A.
\begin{theorem}[Analytic Laugwitz rigidity]\label{thm:laugwitz-intro}
Let $V$ be a real vector space of dimension $n\ge3$, and let $F$ be a smooth, strongly convex Minkowski norm, regular and real-analytic on the entire slit space $V\setminus\{0\}$. If the Hessian metric
\[
 g=D^2\!\left(\tfrac12F^2\right)
\]
is flat on $V\setminus\{0\}$, then $F^2$ is a positive-definite quadratic form. The identity $F(-y)=F(y)$ is not assumed.
\end{theorem}

Theorem A is restated and proved as Theorem \ref{thm:laugwitz}. Flatness reduces the problem to a real-analytic cubic form on the entire standard sphere $(S^m,h_0)$ satisfying
\[
 \nabla^{h_0}C\in\Sym^4T^*S^m,\qquad [C_X,C_Y]=0
\]
where $m=n-1$ and $\nabla^{h_0}$ is the Levi-Civita connection of $h_0$. An extremum of the trace potential gives a point where the cubic tensor vanishes. Near this point, separation of the leading term of a homogeneous potential and the use of nested normal ideals give the local zero-subsphere theorem, Theorem \ref{thm:local-zero}. The first nonzero orders in different transverse directions may differ. The full collection of Gram relations bounds the total number of active axes by $2R\le m+1$, so that the induction terminates after finitely many steps. We then fix an ambient normal along a great-circle arc on which the cubic tensor vanishes and obtain a differential equation for an analytic periodic scalar function. Integrating this equation over one period gives a sum of nonnegative integrals equal to zero. The scalar function must therefore vanish, which proves global rigidity of the cubic form.

The second result deduces Hessian flatness from a twisted Codazzi equation on the indicatrix. We shall also call it Theorem B.
\begin{theorem}[Analytic rigidity of twisted Codazzi quartic forms on an indicatrix]\label{thm:quartic-rigidity-intro}
Let $V$ be a real vector space of dimension $n\ge4$, and let $F$ be a smooth, strongly convex Minkowski norm, regular and real-analytic on the entire slit space $V\setminus\{0\}$. Define $g,\Sigma,h,C$ as above. Suppose that the entire indicatrix $\Sigma$ carries a smooth, completely symmetric covariant tensor of order four
\[
 B\in\Gamma(\Sym^4T^*\Sigma),\qquad B\not\equiv0,
\]
such that, for all tangent vector fields $X,Y,Z,U,W$,
\begin{equation}\label{eq:twisted-intro}
 (\nabla^h_XB)(Y,Z,U,W)-(\nabla^h_YB)(X,Z,U,W) =B(X,C_YZ,U,W)-B(Y,C_XZ,U,W).
\end{equation}
Then $h$ has constant sectional curvature $1$. Equivalently, $g$ is flat on the entire slit space $V\setminus\{0\}$.
\end{theorem}

Theorem B requires only smoothness of $B$. It need not be real-analytic or arise from the Berwald curvature of a Finsler metric. Equation \eqref{eq:twisted-intro} and complete symmetry already imply the algebraic compatibility identity used below; that identity is not an additional hypothesis. The analyticity assumption concerns $F$ and the tensors $h,C$ determined by it. Reversibility of $F$ is not required. Theorem B is proved in Section \ref{sec:quartic-rigidity} without using Theorem A.

To describe the two cases in the proof of Theorem B, put $m=n-1$ and denote by
\[
 q=\max_{y\in\Sigma}\rank\bigl\{
 T_y\Sigma\longrightarrow\Sym^3T_y^*\Sigma,\quad
 v\longmapsto B_y(v,\cdot,\cdot,\cdot)
 \bigr\}
\]
the maximal rank of $B$. Since $B\not\equiv0$, we have $1\le q\le m$, and the points of maximal rank form a nonempty open set. If $0<q<m$, the radical distribution, the centroid of the restricted quartic form, Riccati transport along centroaffine geodesics, and the restrictions imposed by closed orbits imply that $C$ vanishes on a nonempty open set. Since $F$ is real-analytic, the vanishing of $C$ extends over the connected indicatrix. The Gauss equation and cone structure then give flatness of $g$. This branch proves the conclusion of Theorem B directly; it is unnecessary first to exclude a nonzero quartic form of deficient maximal rank.

If $q=m$, we first establish generic irreducibility of the contraction operators on the set where the radical is trivial. On a fixed positive-definite inner-product space, the angular algebraic compatibility identity is converted, by a full contraction identity, into curvature compatibility of each $B_{uv}$ with $K=[C,C]$. The theorem of Mantica--Molinari on curvature-compatible tensors gives Jordan closure \cite[arXiv version, Theorem 4.2]{mm}; we quote this result, specifying the index and curvature conventions needed here. An extremum of a nonzero quartic form on the unit sphere provides a contraction operator with a simple extreme eigenvalue, and hence a rank-one spectral projection. Irreducibility and this projection imply that the generated Jordan algebra consists of all self-adjoint operators. It follows that
\[
 K=sR_0
\]
(Theorem \ref{thm:algebra}). The Gauss equation gives $R^h=\kappa R_0$, where $\kappa=1-s$. The real-analytic identity theorem extends this curvature isotropy from a nonempty open set to the entire indicatrix. Since $m\ge3$, Schur's lemma makes $\kappa$ constant; an extremum of the trace potential and a deformation of the connection then show that $\kappa=1$. Thus this branch also gives flatness of $g$, without using a low-dimensional classification of Jordan algebras.

In Part II we apply the two theorems for a single fibre to actual Landsberg metrics. The Berwald tensor with its upper index lowered,
\[
 B_{aijk}=g_{ab}\mathcal B^b{}_{ijk}
\]
is completely symmetric precisely when the Landsberg condition holds. Under this condition it vanishes whenever a radial vector is inserted, and its tangential restriction to each indicatrix satisfies \eqref{eq:twisted-intro}. The set $\mathscr U=\{x\in M:\mathcal B(x,\cdot)\not\equiv0\}$ is open. At every $x\in\mathscr U$, the angular quartic form $B_x$ is nonzero, so Theorem B makes the fibre Hessian metric flat and Theorem A makes $F_x^2$ quadratic. Smoothness of the fundamental tensor then gives a Riemannian metric on $\mathscr U$. Its Berwald curvature must vanish, contradicting the definition of $\mathscr U$ unless this set is empty. This proves the main theorem.

The planarity and completeness of centroaffine geodesics are taken, respectively, from \cite[Proposition 5.1 and Theorem 5.4]{opozda-completeness}. Section \ref{sec:centro-periodic} supplies the explicit parameter, periodicity of position and velocity, and uniform bounds needed for the subsequent transport argument.

We use the Einstein summation convention: repeated upper and lower indices in a term are summed over $1,\ldots,n$.

\subsection{Scope and organization of the paper}
Theorem A requires $n\ge3$. The restriction $n\ge4$ in Theorem B and the main theorem enters through Schur's lemma on the indicatrix, whose dimension is $m=n-1\ge3$. If $n=3$, pointwise isotropy on the two-dimensional indicatrix does not control the variation of its Gaussian curvature; the present full-rank argument therefore lacks the corresponding constancy step. Extending these arguments to smooth, nonanalytic norms in higher dimensions would require substitutes for the finite-order leading terms and global continuation in analytic Laugwitz rigidity, and for the extension of Cartan-tensor vanishing or curvature isotropy from an open set to the entire indicatrix in Theorem B. We do not pursue these questions here, nor do we treat base dimension two.

Part I concerns a fixed Minkowski space. Section \ref{sec:paper-2} provides the geometry and constant-curvature normalization. Sections \ref{sec:paper-3} and \ref{sec:paper-4} prove Theorem A through local zero sets and continuation along great circles. Sections \ref{sec:paper-6}--\ref{sec:quartic-rigidity} prove Theorem B: radical transport treats deficient rank, while Sections \ref{sec:paper-7} and \ref{sec:paper-8} supply generic irreducibility and Jordan rigidity for full rank.

Part II applies these results to Finsler geometry. Section \ref{sec:finsler-foundations} gives the foundations, Section \ref{sec:paper-5} establishes the angular equations, and Section \ref{sec:paper-9} proves the main theorem.

\part{Geometry and analytic rigidity of Minkowski spaces}\label{part:minkowski}
We work on a fixed Minkowski space and prove Theorems A and B independently. The quartic tensor in Theorem B is arbitrary subject to the stated equation; its realization as Berwald curvature is considered only in Part II.

\section{Minkowski spaces, indicatrices, and centroaffine geometry}\label{sec:paper-2}
This section collects the classical geometry of a single Minkowski space and its indicatrix, and fixes the notation for the subsequent rigidity arguments. For Minkowski norms, Cartan tensors, and homogeneity identities, see \cite[Sections 1.1, 1.2, and 1.4]{bcs} and \cite[Sections 1.1 and 1.5]{chern-shen-2005}; for the curvature and cone decomposition of the Hessian metric and the Gauss equation of the indicatrix, see \cite[Sections 14.2, 14.3, 14.5, and 14.6]{bcs}. The centroaffine interpretation of the indicatrix follows the correspondences in \cite[Sections 1 and 2]{li-minkowski} and their relation to the classical theory in \cite{nomizu-sasaki,opozda-completeness}.

Unless otherwise stated, the norm is smooth, regular on the entire slit space, and strongly convex; neither reversibility nor real analyticity is required. We specify dimensions and conventions where needed, quote classical results, and prove the supplementary assertions used later.

\paragraph{Notation for connections and metrics.}
We use the following conventions throughout. The standard flat affine connection on the vector space $V$ is denoted by $D$; the Levi-Civita connection of the fibre Hessian metric $g$ by $\nabla^g$; the Levi-Civita connection of the induced indicatrix metric $h$ by $\nabla^h$; and the affine connection induced on the indicatrix by the centroaffine normalization by $\nabla$. The Riemannian metric of the standard unit sphere $S^m$ and its Levi-Civita connection are denoted by $h_0$ and $\nabla^{h_0}$, respectively. The corresponding curvatures are denoted by $R^g$, $R^h$, $R^{\nabla}$, and $R^{h_0}$.

For an induced connection on a quotient bundle or another vector bundle, a superscript on $\nabla$ specifies the bundle on which it acts; its inducing connection and the permitted directions of differentiation are stated at each use. For example, a connection on a quotient bundle $Q$ is denoted by $\nabla^Q$. Coordinate partial derivatives retain the notation $\partial_i$ and the like; the symbol $\nabla$ is not used for an unspecified Levi-Civita connection or gradient.

\subsection{Minkowski norms, homogeneity, and strong convexity}\label{subsec:paper-1-1}

\begin{definition}[Minkowski norm]\label{def:minkowski-norm}
Fix an $n$-dimensional real vector space $V$. A Minkowski norm on $V$ is a continuous function $F:V\to[0,\infty)$ such that $F(0)=0$, and such that on $V\setminus\{0\}$ it is smooth, strictly positive, and positively homogeneous of degree one:
\[
 F(\lambda y)=\lambda F(y),\qquad \lambda>0.
\]
In addition, the Hessian of the energy function $E(y):=\tfrac12F(y)^2$ is required to be positive definite at every nonzero vector. This condition is called \textbf{strong convexity}.

We do not require $F(-y)=F(y)$; when this equality holds, $F$ is called reversible or symmetric. Thus a Minkowski norm in this paper need not satisfy the symmetry condition in the usual definition of a norm.
\end{definition}

On $V$, choose linear coordinates $y^1,\ldots,y^n$, and put
\[
 E_i=\frac{\partial E}{\partial y^i},\qquad
 g_{ij}=\frac{\partial^2E}{\partial y^i\partial y^j}.
\]
For each fixed $y\ne0$, the positive definite matrix $(g_{ij}(y))$ defines an inner product on $V$ by
\[
 g_y(u,v)=g_{ij}(y)u^iv^j.
\]
We call $g_y$ the \textbf{fundamental tensor} of $F$ at the reference vector $y$.

For each fixed reference vector $y\ne0$, the function $g_y(u,v)$ is linear in each of $u,v$, although $g_y$ itself may vary with $y$. This variation is described by the \textbf{Cartan cubic tensor}
\[
 \mathcal C_{ijk}(y)
 :=\frac12\frac{\partial g_{ij}}{\partial y^k}
 =\frac12E_{ijk}(y).
\]
Since third partial derivatives commute, $\mathcal C$ is completely symmetric. Its coordinate-free expression is
\[
 \mathcal C_y(u,v,w)
 =\frac12\left.\frac{\dd}{\dd t}\right|_{t=0}g_{y+tw}(u,v).
\]
In what follows, $C$ denotes the restriction of $\mathcal C$ to the tangent spaces of the indicatrix.

\begin{remark}[Normalization of the Cartan tensor]\label{rem:cartan-normalization}
Our tensor $\mathcal C=\tfrac14\partial_y^3(F^2)$ agrees with the convention for Cartan torsion in \cite[Section 1.2]{shen-2001} and \cite[Section 1.5]{chern-shen-2005}. Bao--Chern--Shen use the normalized components $A_{ijk}=F\mathcal C_{ijk}$; see \cite[equation (1.4.2)]{bcs}. Thus $A_{ijk}$ and $\mathcal C_{ijk}$ have different degrees of homogeneity, whereas their tangential restrictions to the indicatrix $F=1$ coincide. All subsequent statements concerning homogeneity degrees refer to components in fixed linear coordinates.
\end{remark}

The following proposition is a direct consequence of Euler's homogeneity identity; see \cite[Theorem 1.2.1; equations (1.4.3)--(1.4.6)]{bcs} and \cite[equations (1.49)--(1.51)]{chern-shen-2005}.
\begin{proposition}[Homogeneity identities and radial vanishing]\label{prop:minkowski-euler}
Let $F$ be a Minkowski norm as above, and let $E=F^2/2$. Then
\begin{equation}\label{gC}
	g_y(y,u)=\dd E_y(u),\qquad g_y(y,y)=F(y)^2,
	\qquad \mathcal C_y(y,u,v)=0.
\end{equation}
 
The equivalent coordinate identities are
\begin{equation}\label{eq:euler}
 y^iE_i=2E,\qquad y^ig_{ij}=E_j,\qquad y^iE_{ijk}=0.
\end{equation}
The radial vanishing identity in \eqref{gC} holds in every slot of $\mathcal C$ upon insertion of $y$. Moreover, for $\lambda>0$ and fixed vectors $u,v,w\in V$,
\[
 g_{\lambda y}(u,v)=g_y(u,v),\qquad
 \mathcal C_{\lambda y}(u,v,w)=\lambda^{-1}\mathcal C_y(u,v,w).
\]

\end{proposition}

The following proposition is the classical characterization of a Euclidean structure by the vanishing of the Cartan tensor; see \cite[Proposition 1.2.1]{shen-2001} and \cite[Section 1.5]{chern-shen-2005}.
\begin{proposition}[Characterization of vanishing Cartan tensor]\label{prop:cartan-zero-quadratic}
Suppose that $\dim V=n\ge2$. Then $\mathcal C$ vanishes identically on the whole of $V\setminus\{0\}$ if and only if $F^2$ is a positive definite quadratic form.
\end{proposition}

\begin{remark}[Radial vanishing, pointwise vanishing, and vanishing identically]
\leavevmode
\begin{enumerate}[label=\arabic*.,leftmargin=2em,itemsep=3pt]
\item $\mathcal C_y(y,\cdot,\cdot)=0$: this holds for every Minkowski norm, as a consequence of homogeneity.
\item $\mathcal C_y=0$: the entire cubic tensor vanishes at a fixed reference vector $y$, that is, $\mathcal C_y(u,v,w)=0$ for all $u,v,w\in V$.
\item $\mathcal C\equiv0$: the cubic tensor vanishes at every nonzero reference vector; when $n\ge2$, this is equivalent to $F$ being induced by an inner product.
\end{enumerate}
Pointwise vanishing does not imply identical vanishing. In Section \ref{sec:quartic-rigidity}, analyticity is applied to vanishing on a nonempty open set.
\end{remark}

\subsection{Compactness and spherical topology of the indicatrix}\label{subsec:paper-1-2}
Continue to fix a Minkowski space $(V,F)$, and define its \textbf{indicatrix} by
\[
 \Sigma=\{y\in V:F(y)=1\}
       =\{y\in V:E(y)=\tfrac12\}.
\]
For the geometry of the indicatrix, see \cite[Section 1.1]{chern-shen-2005}, \cite[Section 1.4]{bcs}, and \cite[Lemmas 6 and 7]{li-minkowski}.

\begin{proposition}[Topology and induced metric of the indicatrix]\label{prop:indicatrix-geometry}
Suppose that $n\ge2$. The indicatrix $\Sigma=\{F=1\}$ is a compact smooth hypersurface without boundary, diffeomorphic to $S^{n-1}$. For every $y\in\Sigma$,
\[
 T_y\Sigma=\ker\dd E_y=\{X:g_y(y,X)=0\}.
\]
The restriction $h=g|_{T\Sigma}$ is a Riemannian metric, and
\[
 V=\R y\mathbin{\oplus^{\perp_g}}T_y\Sigma,\qquad
 u^\top=u-g_y(y,u)y.
\]

\end{proposition}

\begin{remark}
We shall see that the Riemannian metric $h$ can be regarded both as the first fundamental form of the indicatrix as a submanifold of $V$ equipped with its Hessian structure, and as the centroaffine metric of the indicatrix.
\end{remark}

\subsection{Connection and curvature of the Hessian metric}\label{subsec:paper-2-1}\label{sec:hess-curv}
Put $\partial_i=\partial/\partial y^i$, and denote by $D$ the standard flat affine connection of $V$. Under the natural identification $T_y(V\setminus\{0\})\simeq V$, the fundamental tensor $g_y$ defines a Riemannian metric on the slit vector space by
\[
 g=g_{ij}(y)\,\dd y^i\dd y^j.
\]
We denote its Levi-Civita connection by $\nabla^g$. All connections and curvatures in this section are computed within this fixed fibre. To fix the curvature sign, for any connection $\mathcal D$ we use the convention
\[
 R^{\mathcal D}(X,Y)Z=\mathcal D_X\mathcal D_YZ-\mathcal D_Y\mathcal D_XZ-\mathcal D_{[X,Y]}Z,
 \qquad [A,B]=AB-BA.
\]
With this convention, the sectional curvature of the plane spanned by orthonormal vectors $X,Y$ is $g(R^g(X,Y)Y,X)$. For the connection and curvature formulas below, see \cite[equations (14.2.2) and (14.2.6)]{bcs}; replacing $A/F$ there by our $\mathcal C$ gives the stated form.

\begin{proposition}[Hessian connection and curvature]\label{prop:hessian-connection-curvature}
Denote by $\mathcal C^\sharp$ the vector-valued bilinear form obtained by using $g$ to raise one index of $\mathcal C$, and put $\mathcal C_XY=\mathcal C^\sharp(X,Y)$. Then
\begin{equation}\label{eq:hessian-curv}
 \nabla^g=D+\mathcal C^\sharp,\qquad
 R^g(X,Y)=-[\mathcal C_X,\mathcal C_Y].
\end{equation}
Consequently, $g$ is flat if and only if all contraction operators $\mathcal C_X$ commute pairwise at each point.
\end{proposition}

\begin{lemma}[Codazzi property of the Hessian cubic form]\label{lem:hessian-codazzi}
For every smooth strongly convex Minkowski norm, $\nabla^g\mathcal C$ is a completely symmetric covariant four-tensor. This does not require $g$ to be flat.
\end{lemma}
For a Hessian metric, $(g,D)$ is a statistical structure. Its $g$-dual connection $D^*$ is flat as well: the curvature duality identity $g(R^D(X,Y)Z,W)=-g(Z,R^{D^*}(X,Y)W)$ and $R^D=0$ give $R^{D^*}=0$. Thus the structure is conjugate symmetric. Its difference tensor relative to $\nabla^g$ has cubic form $-\mathcal C$, so the assertion follows from \cite[Lemma 2.1]{opozda-completeness}.

\subsection{Cone structure and the Gauss equation of the indicatrix}\label{subsec:paper-2-2}\label{sec:cone}
For the material in this subsection, see \cite[equations (14.3.4)--(14.3.6)]{bcs}.
\begin{proposition}[Cone structure of the Hessian metric]\label{prop:hessian-cone}
The map $\Phi:(0,\infty)\times\Sigma\to V\setminus\{0\}$, $\Phi(r,u)=ru$, is a diffeomorphism, and
\begin{equation}\label{eq:cone}
 \Phi^*g=\dd r^2+r^2h.
\end{equation}
If $X,Y$ are tangential lifts from the indicatrix, independent of $r$, and $\nabla^h$ is the Levi-Civita connection of $h$, then
\[
 \nabla^g_{\partial_r}\partial_r=0,\quad \nabla^g_X\partial_r=\nabla^g_{\partial_r}X=r^{-1}X,\qquad \nabla^g_XY=\nabla^h_XY-rh(X,Y)\partial_r.
\]
Here vector fields on the two sides are identified by $\Phi$.
\end{proposition}

We use the classical Gauss formulas for the indicatrix in \cite[equations (14.6.5)--(14.6.6), Proposition 14.6.1]{bcs}, together with the Codazzi property of Lemma \ref{lem:hessian-codazzi}.
\begin{proposition}[Gauss--Codazzi identities of the indicatrix]\label{prop:indicatrix-gauss-codazzi}
Put $C=\mathcal C|_{T\Sigma}$, and use $h(C_XY,Z)=C(X,Y,Z)$ to define $C_X$. Then
\begin{equation}\label{eq:gauss}
 R^h=R_0-[C,C],\qquad R_0(X,Y)Z=h(Y,Z)X-h(X,Z)Y,\qquad \nabla^h C\in\Sym^4T^*\Sigma.
\end{equation}
Here $[C,C](X,Y)Z=[C_X,C_Y]Z$. The following three conditions are equivalent: $g$ is flat; $R^h=R_0$; and all $C_X$ commute pairwise at each point. When $n\ge3$, the second condition says that $h$ has constant sectional curvature $1$.
\end{proposition}

\subsection{Basic formulas for the centroaffine connection}\label{subsec:paper-2-3}
Regard $\Sigma$ again as a hypersurface in $V$. Since the position vector $y$ is transverse to $T_y\Sigma$, the ambient derivative has a unique decomposition
\[
 D_XY=\nabla_XY+\beta(X,Y)y,
\]
whose first term is tangent to $\Sigma$. This defines a tangential connection $\nabla$.

\begin{definition}[Centroaffine connection]\label{def:centroaffine-connection}
The connection $\nabla$ above is called the connection induced by the centroaffine normalization of the indicatrix with transversal vector field $-y$. The corresponding centroaffine metric is $-\beta$.
\end{definition}

For the centroaffine geometry of indicatrices, see \cite[Chapter I, Section 3]{nomizu-sasaki} and \cite[Section 1.1, Lemmas 6 and 7]{li-minkowski}.

The following proposition describes the correspondence between the metrics and cubic forms; see \cite[Lemmas 7 and 8]{li-minkowski}.

\begin{proposition}[Basic centroaffine identities]\label{prop:centroaffine-identities}
The centroaffine metric of this normalization is precisely $h=g|_{T\Sigma}$, and the induced connection $\nabla$ is torsion-free and satisfies
\begin{equation}\label{eq:centro}
 D_XY=\nabla_XY-h(X,Y)y,\qquad \nabla=\nabla^h-C,\qquad \nabla h=2C,\qquad R^{\nabla}=R_0.
\end{equation}
In the difference of connections, $C$ denotes the vector-valued bilinear form obtained by raising an index with $h$.
\end{proposition}

\begin{remark}[The two notions of completeness]
Since $\nabla$ need not preserve $h$, completeness of $(\Sigma,h)$ does not directly imply completeness of $\nabla$. We use the centroaffine completeness theorem in the next subsection.
\end{remark}

\subsection{Completeness and periodicity of centroaffine geodesics}\label{subsec:paper-7-3}\label{sec:centro-periodic}

We use the planarity, area-parameter characterization, and completeness results of \cite[Section 5]{opozda-completeness} to establish periodicity of position and velocity and the required uniform bounds.

Let $V$ be a real vector space of dimension at least $2$, and let $F$ be a Minkowski norm on $V$, smooth and strongly convex on $V\setminus\{0\}$. Put
\[
 E=\tfrac12F^2,\qquad
 \Sigma=\{y\in V:E(y)=\tfrac12\},\qquad
 g_y=\Hess_yE,\qquad h_y=g_y|_{T_y\Sigma}.
\]
Denote by $D$ the ordinary affine derivative in $V$, and by $\nabla$ the centroaffine connection of $\Sigma$, so that
\[
 D_XY=\nabla_XY-h(X,Y)y.
\]

\begin{lemma}[Ambient equation and fixed plane]\label{lem:centro-geodesic-plane}
Given $y_0\in\Sigma$ and $0\ne v_0\in T_{y_0}\Sigma$, the corresponding local $\nabla$-geodesic has nowhere vanishing velocity and satisfies
\begin{equation}\label{eq:plane}
 \ddot y=-a(t)y,\qquad
 a(t)=h_{y(t)}(\dot y,\dot y)>0.
\end{equation}
Moreover, $y(t)\wedge\dot y(t)=y_0\wedge v_0\ne0$. Thus the curve lies in the fixed two-dimensional plane $\Pi=\operatorname{span}\{y_0,v_0\}$, along its section $\Gamma=\Sigma\cap\Pi$.
\end{lemma}
For the planarity of geodesics, see \cite[Proposition 5.1]{opozda-completeness} and \cite[Chapter I, Section 3; Chapter II, Section 3, Remark 3.2]{nomizu-sasaki}; equation \eqref{eq:plane} is our centroaffine Gauss formula along a geodesic.

\begin{definition}[Area parameter of a planar section]\label{def:centro-area-parameter}
On the two-dimensional plane $\Pi$, choose a positive area form $\det$. If a regular curve $z(t)\in\Gamma$ satisfies $\det(z,\dot z)=c$ for a nonzero constant, we call $t$ a \textbf{parameter of constant areal velocity}. The oriented area $\mathcal S$ swept out by the radial segment satisfies
\[
 \frac{\dd\mathcal S}{\dd t}=\frac12\det(z,\dot z)=\frac c2.
\]
Thus $c$ is the constant determinant value, while the areal velocity is $c/2$. Choosing the orientation to agree with the direction of motion, we may take $c>0$.
\end{definition}

\begin{lemma}[Equivalence of affine and area parameters]\label{lem:centro-affine-area}
For the section $\Gamma$, any smooth regular parametrization $z(t)$ satisfies
\[
 \nabla_{\dot z}\dot z=0
 \quad\Longleftrightarrow\quad
 \det(z,\dot z)\text{ is constant}.
\]
Thus the parameters of constant areal velocity on this planar section are precisely the affine geodesic parameters of $\nabla$.
\end{lemma}
This correspondence of parameters is given in \cite[equation (32) in Lemma 5.2 and the proof of Theorem 5.3]{opozda-completeness}; here we use the areal-velocity normalization of Definition \ref{def:centro-area-parameter}.

By \cite[Theorem 5.4]{opozda-completeness}, the induced and dual connections of a centroaffine ovaloid with transversal $-y$ are complete. We prove the additional periodicity and uniform bounds below.

\begin{theorem}[Completeness and periodicity of centroaffine geodesics]\label{thm:centro-periodic}
Under the smooth strongly convex Minkowski assumptions above, $\nabla$ is geodesically complete. More precisely, for every $y_0\in\Sigma$ and $0\ne v_0\in T_{y_0}\Sigma$, there exists a unique affinely parametrized geodesic
\[
 y:\R\longrightarrow\Sigma,\qquad
 \nabla_{\dot y}\dot y=0,\qquad
 y(0)=y_0,\quad \dot y(0)=v_0.
\]
Its image is the planar section
\[
 \Sigma\cap\Pi,\qquad \Pi=\operatorname{span}\{y_0,v_0\},
\]
and there exists $\tau_{\mathrm{per}}>0$ such that, for every $t\in\R$,
\[
 y(t+\tau_{\mathrm{per}})=y(t),\qquad \dot y(t+\tau_{\mathrm{per}})=\dot y(t).
\]
The function $a(t)=h_{y(t)}(\dot y(t),\dot y(t))$ along the curve is smooth, positive, and $\tau_{\mathrm{per}}$-periodic. Consequently, there are constants
\[
 0<a_0\le a(t)\le a_1<\infty\qquad(t\in\R).
\]
Here $\tau_{\mathrm{per}},a_0,a_1$ may depend on the initial data.
\end{theorem}
\begin{proof}[Supplementary proof of periodicity and uniform bounds]
The indicatrix $\Sigma$ is a centroaffine ovaloid with transversal vector field $-y$. Hence $\nabla$ is geodesically complete by \cite[Theorem 5.4]{opozda-completeness}. It remains to verify the periodicity and uniform bounds required below.

Fix nonzero initial data $(y_0,v_0)$, and let $y(t)$ be the corresponding geodesic defined on the whole real axis. By Lemma \ref{lem:centro-geodesic-plane}, its image lies in the fixed plane $\Pi=\operatorname{span}\{y_0,v_0\}$ and its velocity never vanishes. Choose an auxiliary Euclidean inner product and an orientation on $\Pi$ so that
\[
 c=\det(y_0,v_0)>0.
\]
The ambient equation \eqref{eq:plane} gives $\det(y(t),\dot y(t))=c$.

Choose an oriented orthonormal basis $e_1,e_2$, and put
\[
 e(\theta)=\cos\theta\,e_1+\sin\theta\,e_2,\qquad
 r(\theta)=\frac1{F(e(\theta))},\qquad
 \eta(\theta)=r(\theta)e(\theta).
\]
Positive homogeneity ensures that $\eta$ parametrizes the entire section $\Gamma=\Sigma\cap\Pi$. The function $r$ is smooth, strictly positive, and $2\pi$-periodic, so there exist
\[
 0<r_{\min}\le r(\theta)\le r_{\max}<\infty.
\]
A direct computation gives
\[
 \det(\eta,\eta_\theta)=r(\theta)^2.
\]
For $y(t)$, choose a continuous lift $\theta(t)$ of its angular parameter such that $\theta(0)=\theta_0$ and $y_0=\eta(\theta_0)$. Since the determinant velocity is constant,
\[
 r(\theta(t))^2\dot\theta(t)=c,\qquad
 \dot\theta(t)=\frac{c}{r(\theta(t))^2}>0.
\]
In particular, $\dot\theta\ge c/r_{\max}^2$. Hence $\theta(t)$, as $t\to\pm\infty$, tends respectively to $\pm\infty$, and the image of the geodesic is the whole of $\Gamma$.

Integration of the parameter equation gives
\[
 t=G(\theta(t)),\qquad
 G(\theta)=\frac1c\int_{\theta_0}^{\theta}r(s)^2\,\dd s.
\]
Put
\[
 \tau_{\mathrm{per}}=\frac1c\int_{\theta_0}^{\theta_0+2\pi}r(s)^2\,\dd s>0.
\]
Since $r^2$ is $2\pi$-periodic,
\[
 G(\theta+2\pi)-G(\theta)
 =\frac1c\int_{\theta}^{\theta+2\pi}r(s)^2\,\dd s=\tau_{\mathrm{per}}.
\]
The function $G$ is strictly increasing, whence
\[
 \theta(t+\tau_{\mathrm{per}})=\theta(t)+2\pi.
\]
It follows that $y(t+\tau_{\mathrm{per}})=\eta(\theta(t)+2\pi)=y(t)$. Differentiating this identity also gives $\dot y(t+\tau_{\mathrm{per}})=\dot y(t)$.

Finally, $\dot y$ is nowhere zero and $h$ is positive definite, so
\[
 a(t)=h_{y(t)}(\dot y(t),\dot y(t))>0.
\]
The smooth periodicity of position and velocity implies that $a$ is also smooth and $\tau_{\mathrm{per}}$-periodic. Taking its extrema over a full period, we obtain
\[
 a_0:=\min_{0\le t\le \tau_{\mathrm{per}}}a(t)>0,\qquad
 a_1:=\max_{0\le t\le \tau_{\mathrm{per}}}a(t)<\infty.
\]
Periodicity then gives the uniform estimates for all $t\in\R$.
\end{proof}

\subsection{Trace potentials, constant-curvature normalization, and reduction to the sphere}\label{subsec:paper-2-4}\label{subsec:paper-12-1}\label{sec:k-one}
The Cartan trace potential will give constant-curvature normalization and a zero of a commuting cubic tensor on the sphere. The trace-potential identity uses smoothness and strong convexity; the subsequent global conclusions require $n\ge3$, or $m=n-1\ge2$.

The following proposition is the standard trace-potential formula for the mean Cartan torsion; see \cite[equation (1.48)]{chern-shen-2005} and \cite[Definition 4 in Section 2.1, Lemma 8]{li-minkowski}.

\begin{proposition}[Cartan trace potential]\label{prop:cartan-trace-potential}
In any fixed system of linear coordinates, put
\[
 \varphi:=\left.\tfrac12\log\det(g_{ij})\right|_\Sigma,
 \qquad I(X):=\tr_hC(X,\cdot,\cdot).
\]
Then $I=\dd\varphi$. A change of linear coordinates adds only a constant to $\varphi$, so $I$ is independent of the coordinate choice.
\end{proposition}

The next proposition normalizes an already constant curvature. Its constant-curvature hypothesis will be verified separately in the full-rank branch of Section \ref{sec:quartic-rigidity}.

\begin{proposition}[Constant-curvature normalization for a compact Minkowski indicatrix]\label{prop:compact-indicatrix-normalization}
Let $F$ be a smooth strongly convex Minkowski norm on a finite-dimensional real vector space $V$, regular on the entire slit space. Let $\Sigma=\{y\in V:F(y)=1\}$, and let $h$ be its indicatrix metric. If $m=\dim\Sigma\ge2$ and $h$ has constant sectional curvature $\kappa$, then $\kappa=1$.
\end{proposition}
\begin{proof}

We use the Gauss equation $R^h=R_0-[C_X,C_Y]$. At a point of the indicatrix, choose an $h$-orthonormal basis $e_1,\ldots,e_m$, and put
\[
 C_{ijk}=C(e_i,e_j,e_k),\qquad
 C_i=C_{e_i},\qquad
 I_r=\sum_i C_{iir}.
\]
Thus $|I|^2=\sum_r I_r^2$ and $|C|^2=\sum_{i,j,r}C_{ijr}^2$. By the definition of the commutator and the complete symmetry of $C$,
\[
 \sum_{i,j}h([C_i,C_j]e_j,e_i) =\sum_{i,j,r}\bigl(C_{iir}C_{jjr}-C_{ijr}^2\bigr) =|I|^2-|C|^2.
\]
On the other hand, the constant-sectional-curvature hypothesis is $R^h=\kappa R_0$. Substitution in the Gauss equation gives
\[
 [C_X,C_Y]=(1-\kappa)R_0(X,Y).
\]
In the orthonormal basis, $h(R_0(e_i,e_j)e_j,e_i)$ equals, for $i\ne j$, the value $1$, and for $i=j$, the value $0$. Its double contraction is therefore $m(m-1)$. Comparing the two calculations, we obtain
\begin{equation}\label{eq:trace-k}
 m(m-1)(1-\kappa)=|I|^2-|C|^2.
\end{equation}
\Needspace{4\baselineskip}

By Proposition \ref{prop:cartan-trace-potential}, the one-form $I$ has a smooth global potential $\varphi$ on $\Sigma$. Since $\Sigma$ is compact, $\varphi$ attains a maximum at some point $p$, where $I_p=\dd\varphi_p=0$. Applying \eqref{eq:trace-k} at this point gives
\[
 m(m-1)(1-\kappa)=-|C|^2\le0.
\]
Since $m\ge2$, we have $m(m-1)>0$, and hence $\kappa\ge1$.

We must exclude strict inequality. Suppose, to the contrary, that $\kappa>1$, and put
\[
 \lambda=\sqrt{\frac{\kappa}{\kappa-1}},\qquad
 \nabla^{(\lambda)}_XY=\nabla^h_XY+\lambda C_XY.
\]
Since $\nabla^h$ is torsion-free and $C_XY=C_YX$, so is $\nabla^{(\lambda)}$.

The curvature deformation formula for connections gives
\[
 R^{\nabla^{(\lambda)}}(X,Y) =R^h(X,Y) +\lambda\bigl((\nabla^h_XC)_Y-(\nabla^h_YC)_X\bigr) +\lambda^2[C_X,C_Y].
\]
The Codazzi equation gives $(\nabla^h_XC)_Y=(\nabla^h_YC)_X$, so the linear term vanishes. Using $R^h=\kappa R_0$ and $[C_X,C_Y]=(1-\kappa)R_0(X,Y)$, we obtain
\[
 R^{\nabla^{(\lambda)}}(X,Y)
 =\{\kappa+\lambda^2(1-\kappa)\}R_0(X,Y)
 =\left(\kappa+\frac{\kappa}{\kappa-1}(1-\kappa)\right)R_0(X,Y)
 =0.
\]
Thus $\nabla^{(\lambda)}$ is flat.

The Minkowski indicatrix $\Sigma$ is diffeomorphic to $S^m$; for $m\ge2$, it is therefore compact, connected, and simply connected. The connection induced by $\nabla^{(\lambda)}$ on $T^*\Sigma$ is also flat. Parallel translation for this cotangent connection depends only on the homotopy class of the path. Thus any nonzero covector at one point extends by parallel translation along paths to a global parallel $1$-form $\alpha$. Simple connectedness makes the extension independent of the path, and the fact that parallel translation is a linear isomorphism ensures that $\alpha$ is nowhere zero. In the formulas below, $\nabla^{(\lambda)}\alpha$ denotes differentiation with this induced cotangent connection.

Torsion-freeness implies that $\alpha$ is closed. Indeed,
\[
 \dd\alpha(X,Y) =(\nabla^{(\lambda)}_X\alpha)(Y)-(\nabla^{(\lambda)}_Y\alpha)(X) +\alpha\bigl(\nabla^{(\lambda)}_XY-\nabla^{(\lambda)}_YX-[X,Y]\bigr)=0.
\]
The first two covariant-derivative terms vanish because $\nabla^{(\lambda)}\alpha=0$, and the torsion term vanishes because $\nabla^{(\lambda)}$ is torsion-free. Simple connectedness also makes a closed $1$-form globally exact, so there is a smooth function $f$ with $\alpha=\dd f$. By compactness, $f$ attains a maximum, where $\dd f=0$, contradicting the fact that $\alpha$ is nowhere zero. Hence $\kappa>1$ is impossible. Together with $\kappa\ge1$, this gives $\kappa=1$.

\end{proof}

The argument uses compactness and simple connectivity, but neither metric compatibility nor geodesic completeness of $\nabla^{(\lambda)}$.

\begin{remark}
When $\kappa=1$, the Gauss equation gives $[C_X,C_Y]=0$, and for every finite $\lambda$ the curvature deformation above gives $R^{\nabla^{(\lambda)}}=R_0$, rather than a flat connection. Thus the normalization argument does not permit one to deduce from $\kappa=1$ that $C=0$, nor does it alone imply that the Minkowski norm is quadratic. In the real-analytic setting of this paper, quadraticity also requires the Laugwitz rigidity result of Section \ref{sec:paper-4}.
\end{remark}

\begin{proposition}[Reduction of a flat Hessian metric to the sphere]\label{prop:flat-hessian-sphere}
Let $n\ge3$ and put $m=n-1$. Suppose that the fibre Hessian metric $g$ is flat on the whole of $V\setminus\{0\}$. Then $(\Sigma,h)$ is isometric to the standard unit sphere $(S^m,h_0)$. The cubic form $C$ transferred by this isometry satisfies the spherical Codazzi--commutation system
\begin{equation}\label{eq:sphere-system}
 \nabla^{h_0}C\in\Sym^4T^*S^m,\qquad [C_X,C_Y]=0,
\end{equation}
where $\nabla^{h_0}$ is the Levi-Civita connection of the standard spherical metric $h_0$. If, in addition, $F$ is real-analytic on the entire slit vector space, the transferred tensor $C$ is real-analytic as well.
\end{proposition}
\begin{proof}
The smooth isometric reduction and its global nature are given by \cite[Proposition 14.6.1 and the subsequent remark; Theorem 13.4.1]{bcs}. We give a short verification of analyticity. If $F$ is real-analytic on the entire slit vector space, then $\Sigma,h,C$ are real-analytic. Let $\Phi:(\Sigma,h)\to(S^{n-1},h_0)$ be the isometry and fix an arbitrary point $p\in\Sigma$. On sufficiently small normal neighbourhoods of $p$ and $\Phi(p)$, it satisfies $\Phi=\exp^{h_0}_{\Phi(p)}\circ\dd\Phi_p\circ(\exp^h_p)^{-1}$. The exponential maps of the analytic metrics are analytic, and the local inverse of $\exp^h_p$ is analytic. Hence $\Phi$ is analytic near $p$. Since $p$ is arbitrary, the isometry and the transferred cubic tensor are real-analytic everywhere.
\end{proof}

In Sections \ref{sec:paper-3} and \ref{sec:paper-4} we prove that every global analytic solution of the spherical Codazzi--commutation system \eqref{eq:sphere-system} vanishes. Radial vanishing, homogeneity, and Proposition \ref{prop:cartan-zero-quadratic} then give the quadraticity of $F^2$.

\begin{lemma}[Trace--norm identity for a commuting cubic form]\label{lem:commuting-cubic-trace}
Let $(W,h)$ be an $m$-dimensional positive definite inner product space, and let $C\in\Sym^3W^*$ be such that all self-adjoint operators $C_X$ commute pairwise. Then there is an orthonormal basis whose dual basis $\theta^1,\ldots,\theta^m$ satisfies
\[
 C=\sum_{i=1}^m c_i(\theta^i)^{\otimes3},\qquad
 |\tr_hC|^2=|C|^2=\sum_i c_i^2.
\]
Here tensor norms are defined by full contraction with $h$.
\end{lemma}
\begin{proof}

The proof uses simultaneous orthogonal diagonalization and complete symmetry. In the inner product space $W$, complete symmetry of $C$ makes each $C_X$ self-adjoint. Since these operators commute, they can be simultaneously orthogonally diagonalized. It is enough to diagonalize the finite family $C_{e_1},\ldots,C_{e_m}$ simultaneously, because every $C_X$ is a linear combination of these operators.

Simultaneous diagonalization follows successively from the spectral theorem. First diagonalize one operator orthogonally. Since the remaining operators commute with it, each of its eigenspaces is invariant under them. Diagonalize the next operator on each eigenspace, and continue. After finitely many steps we obtain a common orthonormal eigenbasis.

Choose such an orthonormal basis $e_1,\ldots,e_m$. If $i\ne j$,
\[
 C(e_a,e_i,e_j)=h(C_{e_a}e_i,e_j)=0.
\]
By complete symmetry, whenever at least two of the three indices differ, those two indices can be placed in the last two slots. The corresponding component therefore vanishes. The only remaining components are
\[
 c_i:=C(e_i,e_i,e_i).
\]
If $\theta^i$ is the dual basis, then $C=\sum_i c_i(\theta^i)^{\otimes3}$, and
\[
 I(e_i)=\sum_jC(e_i,e_j,e_j)=c_i.
\]
It follows that
\begin{equation}\label{eq:trace-norm}
 |\operatorname{tr}_h C|^2=\sum_i c_i^2=|C|^2.
\end{equation}
This identity holds at each point and does not require the diagonalizing frame to depend smoothly on the point.
\end{proof}

\begin{lemma}[Trace potential and a zero of the cubic tensor on the sphere]\label{lem:spherical-cubic-zero-point}
Let $m\ge2$, and let $C$ be a smooth completely symmetric cubic form on the standard sphere $(S^m,h_0)$ satisfying the spherical Codazzi--commutation system \eqref{eq:sphere-system}. Then its trace one-form $I=\tr_{h_0}C$ has a global smooth potential. Moreover, there exists at least one point $p\in S^m$ at which the entire cubic tensor satisfies $C_p=0$.

\end{lemma}
This conclusion does not require analyticity and does not assert that $C$ vanishes identically.
\begin{proof}
We first use the Codazzi condition and simple connectedness of the sphere to construct a global trace potential, then use compactness and Lemma \ref{lem:commuting-cubic-trace} to obtain a zero of the cubic tensor.

Since $\nabla^{h_0}$ preserves $h_0$, choosing a $\nabla^{h_0}$-normal frame at a point gives
\[
 (\dd I)(X,Y) =(\nabla^{h_0}_XI)(Y)-(\nabla^{h_0}_YI)(X) =\sum_i\bigl((\nabla^{h_0}_XC)(Y,e_i,e_i) -(\nabla^{h_0}_YC)(X,e_i,e_i)\bigr)=0.
\]
The last equality is precisely the Codazzi condition. Thus $I$ is a closed one-form. Since $S^m$ is simply connected for $m\ge2$, there is a global function $\tau$ such that $I=\dd\tau$.

By compactness, $\tau$ attains a maximum at some point $p$. Since the sphere has no boundary, $\dd\tau_p=0$, that is, $I_p=0$. Equation \eqref{eq:trace-norm} gives
\[
 |C_p|^2=|I_p|^2=0.
\]
Hence $C_p=0$. Analyticity will be used only in the subsequent local expansion and continuation argument.
\end{proof}

Section \ref{sec:paper-3} extends this zero to a great-circle arc; Section \ref{sec:paper-4} uses periodicity to obtain global vanishing.

\section{Local structure of the zero set of an analytic commuting cubic form}\label{sec:paper-3}
In this section we prove the local zero subsphere theorem \ref{thm:local-zero}. Fix the standard unit sphere $(S^m,h_0)$, $m\ge2$, and denote its Levi-Civita connection by $\nabla^{h_0}$. Let $C$ be a real-analytic symmetric cubic tensor on an open neighbourhood $U$ of $p$, satisfying the spherical Codazzi--commutation system \eqref{eq:sphere-system} on $U$, with $C_p=0$. Here the contraction operators are defined by $h_0(C_XY,Z)=C(X,Y,Z)$. We shall prove that the local zero set of $C$ contains a portion of a positive-dimensional great subsphere through $p$, and that the entire cubic tensor vanishes on this portion, including its components with normal inputs.

The proof passes from a homogeneous potential and the full affine commutation equation to fixed orthogonal separation and same-layer Gram relations. Nested normal ideals then control delays between successive leading terms; cross-layer Gram relations give a dimension bound that terminates the induction.

All Taylor expansions and orders of vanishing below are taken at the origin of the specified local coordinates. For an analytic function or tensor component, $O(|w|^d)$ means that every Taylor term of degree less than $d$ vanishes; the corresponding estimates after differentiation follow by termwise differentiation of the convergent Taylor series. We state in each subsection the hypotheses actually used: the construction of the potential requires only the Codazzi condition, commutativity is used in the subsequent coefficient comparisons, and vanishing of the entire tensor at the point ensures that its leading term has degree $k\ge1$. The conclusions of this section are local throughout.

\subsection{A homogeneous potential of degree two for a spherical cubic form}\label{subsec:paper-3-1}\label{sec:potential}

The next two lemmas construct the radial extension $\bar C$ and a homogeneous potential $P$ of degree two, then recover $C$ from the affine slice $f(w)=P(w,1)$ and normalize its lower-order terms. Only the Codazzi condition is used here; commutativity enters in the next subsection.

\begin{lemma}[Radial extension and its homogeneous potential of degree two]\label{lem:local-cubic-potential}
Let $C$ be a real-analytic symmetric cubic tensor on an open subset $U$ of the standard unit sphere $(S^m,h_0)$, and let $\nabla^{h_0}$ be the Levi-Civita connection of $h_0$. Suppose that
\[
 \nabla^{h_0}C\in\Sym^4T^*U.
\]
For any $p\in U$, choose orthogonal coordinates such that $p=(0,1)$. Then there exists $\varepsilon>0$ such that
\[
 \mathcal K_\varepsilon=\{(w,t):t>0,\ |w|<\varepsilon t\},
 \qquad U_\varepsilon=\mathcal K_\varepsilon\cap S^m\subset U,
\]
and the following assertions hold on this open cone.
\begin{enumerate}[label=\textup{(\arabic*)},leftmargin=2em,itemsep=4pt]
\item Properties of the extended tensor.
Put $\Pi_uA=A-\langle A,u\rangle u\in T_uS^m$ and define
\[
 \bar C_{\rho u}(A,B,Z)
 :=\rho^{-1}C_u(\Pi_uA,\Pi_uB,\Pi_uZ),
 \qquad \rho>0,\quad u\in U_\varepsilon.
\]
Then $\bar C$ is a real-analytic symmetric cubic tensor whose tangential restriction to the sphere is $C$, and
\[
 \bar C_q(q,A,B)=0,\qquad
 \bar C_{\lambda q}(A,B,Z)=\lambda^{-1}\bar C_q(A,B,Z),
 \qquad\lambda>0.
\]
Moreover, $D\bar C$ is symmetric in all four arguments, where $D$ is the standard flat connection of the ambient space.
\item Existence of a homogeneous potential.
There exist a real-analytic symmetric matrix field $M$ and a real-analytic function $P$ satisfying
\[
 \partial_\alpha M_{\beta\gamma}=\bar C_{\alpha\beta\gamma},
 \qquad D^2P=M,\qquad D^3P=\bar C,
\]
and
\[
 M(\lambda q)=M(q),\qquad
 P(\lambda q)=\lambda^2P(q)\quad(\lambda>0).
\]
\end{enumerate}
\end{lemma}

\begin{proof}

Regard $S^m$ as the unit sphere in the ambient Euclidean space $\R^{m+1}$, with inner product $\langle\cdot,\cdot\rangle$ and standard flat connection $D$. Its Levi-Civita connection remains denoted by $\nabla^{h_0}$.

After an orthogonal change of ambient coordinates, we may assume that
\[
 p=(0,\ldots,0,1),\qquad
 q=(w,t)=(w^1,\ldots,w^m,t).
\]
In the ambient component calculations of this subsection, Greek indices $\alpha,\beta,\gamma,\delta$ range over $1,\ldots,m+1$, whereas Roman indices $a,b,c$ range only over $1,\ldots,m$ and denote the $w$ directions. The last direction is written separately as $t$, as in $\bar C_{abt}$. Thus
\[
 \partial_\alpha=\frac{\partial}{\partial q^\alpha},\qquad
 \partial_a=\frac{\partial}{\partial w^a},\qquad
 \partial_t=\frac{\partial}{\partial t}.
\]
Here $t$ is the last coordinate, whereas the radial length is
\[
 \rho=|q|=\sqrt{|w|^2+t^2}.
\]
These are generally different.

Choose a spherical open neighbourhood $U$ of $p$ and define the open cone generated by it:
\[
 \mathcal K_U
 :=\{\rho u:\rho>0,\ u\in U\}
 =\left\{q\ne0:\frac{q}{|q|}\in U\right\}.
\]
The map $(\rho,u)\mapsto\rho u$ is a diffeomorphism from $(0,\infty)\times S^m$ onto $\R^{m+1}\setminus\{0\}$, with inverse $q\mapsto(|q|,q/|q|)$. Hence $\mathcal K_U$ is open, and $\rho,u$ are uniquely determined by $q$.

Since we shall integrate along line segments, we also need a convex open cone. Choose $\varepsilon>0$ sufficiently small that
\[
 U_\varepsilon
 :=\{(v,s)\in S^m:s>0,\ |v|<\varepsilon s\}
 \subset U.
\]
The corresponding open cone is precisely
\[
 \mathcal K_\varepsilon
 =\{(w,t)\in\R^m\times\R:t>0,\ |w|<\varepsilon t\}.
\]
Indeed, $q/|q|\in U_\varepsilon$ is equivalent to $t/|q|>0$ and $|w|/|q|<\varepsilon t/|q|$, which are the two inequalities above. If $q_i=(w_i,t_i)\in\mathcal K_\varepsilon$ and $0\le\theta\le1$, then
\[
 |(1-\theta)w_1+\theta w_2| \le(1-\theta)|w_1|+\theta|w_2| <\varepsilon\bigl((1-\theta)t_1+\theta t_2\bigr).
\]
The last coordinate remains positive, so $\mathcal K_\varepsilon$ is convex and contains $p$. Thus every segment $[p,q]$, $q\in\mathcal K_\varepsilon$, stays in this cone, on which we work below.

For $u\in U_\varepsilon$, orthogonal projection onto the tangent space of the sphere is given by
\[
 \Pi_uA=A-\langle A,u\rangle u.
\]
For $q=\rho u\in\mathcal K_\varepsilon$ and arbitrary ambient vectors $A,B,Z$, define
\[
 \bar C_{\rho u}(A,B,Z)
 :=\rho^{-1}C_u(\Pi_uA,\Pi_uB,\Pi_uZ).
\]
This is a real-analytic symmetric cubic tensor. When $\rho=1$ and all three arguments are tangent, $\Pi_u$ leaves them unchanged, so the tangential restriction of $\bar C$ is $C$.

Since $\Pi_uq=\Pi_u(\rho u)=0$,
\[
 \bar C_q(q,A,B)=0.
\]
By complete symmetry, inserting the position vector $q$ into any slot gives zero. On the other hand, for fixed ambient vectors $A,B,Z$ and $\lambda>0$,
\[
 \bar C_{\lambda q}(A,B,Z) =(\lambda\rho)^{-1}C_u(\Pi_uA,\Pi_uB,\Pi_uZ)=\lambda^{-1}\bar C_q(A,B,Z).
\]
Thus, in fixed orthogonal coordinates,
\[
 \bar C_{\alpha\beta\gamma}(\lambda q)
 =\lambda^{-1}\bar C_{\alpha\beta\gamma}(q).
\]
Here only the point $q$ is rescaled; the three arguments are not rescaled with it.

First calculate at $q=u\in U_\varepsilon$, where $\rho=1$. Since
\[
 \R^{m+1}=\R u\mathbin{\oplus^\perp}T_uS^m,
\]
each argument splits into a radial and a tangential part. As $D\bar C$ is multilinear in its four arguments, it suffices to consider separately the possible numbers of radial arguments.

First suppose that all four arguments are tangent.
We begin with the relation between the spherical and ambient connections. If $X,Y$ are tangent vector fields on the sphere, then
\[
 0=X\langle Y,u\rangle
  =\langle D_XY,u\rangle+\langle Y,X\rangle,
\]
where we used $D_Xu=X$ for the position vector field. Thus the normal component of $D_XY$ is $-h_0(X,Y)u$, and its tangential component is $\nabla^{h_0}_XY$. Hence
\[
 D_XY=\nabla^{h_0}_XY-h_0(X,Y)u.
\]
For tangent fields $X,Y,Z,W$, expand the covariant derivative of the tensor:
\[
 (D_X\bar C)(Y,Z,W) ={}X[\bar C(Y,Z,W)]-\bar C(D_XY,Z,W) -\bar C(Y,D_XZ,W)-\bar C(Y,Z,D_XW).
\]
Substitute the connection formula term by term. All normal terms containing $u$ vanish by radial vanishing, leaving
\[
 \begin{aligned}
 (D_X\bar C)(Y,Z,W)
 ={}&X[C(Y,Z,W)]-C(\nabla^{h_0}_XY,Z,W) -C(Y,\nabla^{h_0}_XZ,W)-C(Y,Z,\nabla^{h_0}_XW)\\
 ={}&(\nabla^{h_0}_XC)(Y,Z,W).
 \end{aligned}
\]
The right-hand side is symmetric in all four arguments by hypothesis.

Next suppose that exactly one argument is radial.
If it occupies the derivative slot, fix $u$ and $X,Y,Z\in T_uS^m$ and differentiate along the line $(1+s)u$. This gives
\[
 (D_u\bar C)_u(X,Y,Z) =\left.\frac{\dd}{\dd s}\right|_{s=0} \bar C_{(1+s)u}(X,Y,Z) =\left.\frac{\dd}{\dd s}\right|_{s=0} (1+s)^{-1}C_u(X,Y,Z) =-C_u(X,Y,Z).
\]
If the radial argument occupies a slot of the cubic tensor, extend $Y,Z$ as constant ambient vectors and differentiate $\bar C_q(q,Y,Z)=0$ in the direction $X$. Since $D_Xq=X$, the product rule gives
\[
 (D_X\bar C)_q(q,Y,Z)+\bar C_q(X,Y,Z)=0.
\]
At $q=u$, this yields
\[
 (D_X\bar C)_u(u,Y,Z)=-C_u(X,Y,Z).
\]
The tensor $D_X\bar C$ remains symmetric in its last three arguments, so placing the radial argument in any of these three slots gives the same result. Together with the complete symmetry of $C$, this permits all four slots to be interchanged.

Finally, suppose that at least two arguments are radial.
If the derivative argument is tangent, differentiate $\bar C_q(q,q,Z)=0$:
\[
 (D_X\bar C)_q(q,q,Z)
 +\bar C_q(X,q,Z)+\bar C_q(q,X,Z)=0.
\]
The last two terms vanish, so $(D_X\bar C)_u(u,u,Z)=0$. If the derivative argument is radial, homogeneity of degree minus one gives
\[
 (D_u\bar C)_u(u,Y,Z)
 =-\bar C_u(u,Y,Z)=0.
\]
An additional radial argument also gives zero. Thus these values are unchanged by permutations of the four arguments.

This proves complete symmetry when $\rho=1$. To extend it to arbitrary $\rho>0$, differentiate the homogeneity formula with respect to $q^\delta$:
\[
 \lambda(\partial_\delta\bar C_{\alpha\beta\gamma})(\lambda q)
 =\lambda^{-1}(\partial_\delta\bar C_{\alpha\beta\gamma})(q).
\]
Setting $q=u$ and $\lambda=\rho$, we obtain
\[
 (\partial_\delta\bar C_{\alpha\beta\gamma})(\rho u)
 =\rho^{-2}(\partial_\delta\bar C_{\alpha\beta\gamma})(u).
\]
The common factor $\rho^{-2}$ preserves symmetry, proving assertion \textup{(1)} throughout $\mathcal K_\varepsilon$.

For assertion \textup{(2)}, we construct $M$ by line integration.
For fixed $\beta,\gamma$, define the one-form
\[
 \omega_{\beta\gamma}
 :=\bar C_{\alpha\beta\gamma}\,\dd q^\alpha.
\]
Since
\[
 \partial_\delta\bar C_{\alpha\beta\gamma}
 =\partial_\alpha\bar C_{\delta\beta\gamma},
\]
the corresponding terms in the exterior derivative cancel, and therefore $\dd\omega_{\beta\gamma}=0$.

\begin{samepage}
Choose a constant symmetric matrix $H_0=((H_0)_{\beta\gamma})$. Along the segment
\[
 q_\sigma=p+\sigma(q-p),\qquad 0\le\sigma\le1,
\]
define
\[
 M_{\beta\gamma}(q) :=(H_0)_{\beta\gamma}+\int_{[p,q]}\omega_{\beta\gamma} =(H_0)_{\beta\gamma}+\int_0^1 \bar C_{\alpha\beta\gamma}(q_\sigma) (q^\alpha-p^\alpha)\,\dd\sigma.
\]
\end{samepage}
Then $M_{\beta\gamma}(p)=(H_0)_{\beta\gamma}$. Since both $H_0$ and $\bar C$ are symmetric in $\beta,\gamma$, we also have $M_{\beta\gamma}=M_{\gamma\beta}$.

We use Stokes' theorem to verify that this is indeed a primitive of $\omega_{\beta\gamma}$. Let $e_\delta$ be the $\delta$th constant coordinate vector. For sufficiently small real $h$, the triangle with vertices $p,q,q+he_\delta$ lies in the convex set $\mathcal K_\varepsilon$. Orient its boundary as $p\to q\to q+he_\delta\to p$. Then
\[
 \int_{[p,q]}\omega_{\beta\gamma}
 +\int_{[q,q+he_\delta]}\omega_{\beta\gamma}
 -\int_{[p,q+he_\delta]}\omega_{\beta\gamma}
 =\int_{\triangle}\dd\omega_{\beta\gamma}=0.
\]
Hence
\[
 M_{\beta\gamma}(q+he_\delta)-M_{\beta\gamma}(q) =\int_{[q,q+he_\delta]}\omega_{\beta\gamma} =\int_0^h\bar C_{\delta\beta\gamma}(q+re_\delta)\,\dd r.
\]
Dividing by $h$ and letting $h\to0$, we obtain
\[
 \partial_\delta M_{\beta\gamma}=\bar C_{\delta\beta\gamma},
 \qquad
 \dd M_{\beta\gamma}
 =\bar C_{\alpha\beta\gamma}\,\dd q^\alpha.
\]
Since the one-form is real-analytic, so is $M$: near any point, termwise integration of its convergent power series along short segments gives a convergent power series for $M$ minus its value at that point.

Radial vanishing gives
\[
 q^\alpha\partial_\alpha M_{\beta\gamma}
 =q^\alpha\bar C_{\alpha\beta\gamma}=0,
 \qquad
 (\partial_\alpha M_{\beta\gamma})q^\beta
 =\bar C_{\alpha\beta\gamma}q^\beta=0.
\]
The first identity says that $M$ is constant along each ray. More explicitly,
\[
 \frac{\dd}{\dd\lambda}M_{\beta\gamma}(\lambda q) =q^\alpha(\partial_\alpha M_{\beta\gamma})(\lambda q) =\lambda^{-1}(\lambda q)^\alpha \bar C_{\alpha\beta\gamma}(\lambda q)=0.
\]
Thus $M(\lambda q)=M(q)$, so $M$ is homogeneous of degree zero. Define
\[
 P(q):=\tfrac12q^\beta M_{\beta\gamma}(q)q^\gamma.
\]
Apply the product rule to compute the first derivative:
\[
 \partial_\alpha P ={}\tfrac12M_{\alpha\gamma}q^\gamma +\tfrac12q^\beta(\partial_\alpha M_{\beta\gamma})q^\gamma +\tfrac12q^\beta M_{\beta\alpha}.
\]
The middle term is $\tfrac12q^\beta q^\gamma\bar C_{\alpha\beta\gamma}=0$. The other two terms combine by symmetry of $M$, giving
\[
 \partial_\alpha P=M_{\alpha\beta}q^\beta.
\]
Differentiate again, retaining both terms of the product rule:
\[
 \partial_\delta\partial_\alpha P =(\partial_\delta M_{\alpha\beta})q^\beta+M_{\alpha\delta} =\bar C_{\delta\alpha\beta}q^\beta+M_{\alpha\delta} =M_{\alpha\delta}.
\]
One further differentiation gives
\[
 \partial_\gamma\partial_\delta\partial_\alpha P
 =\partial_\gamma M_{\alpha\delta}
 =\bar C_{\gamma\alpha\delta}.
\]
Hence $D^2P=M$ and $D^3P=\bar C$. On the other hand,
\[
 P(\lambda q) =\tfrac12(\lambda q)^\beta M_{\beta\gamma}(\lambda q)(\lambda q)^\gamma =\lambda^2P(q),
\]
so $P$ is homogeneous of degree two. This does not mean that $P$ is a quadratic polynomial; nor do we require the Hessian of the auxiliary function $P$ to be positive definite.

\end{proof}

\Needspace{4\baselineskip}
We now restrict $P$ to the affine hyperplane $t=1$, reducing the study of the spherical tensor to that of a function $f$ of $m$ variables.

\begin{lemma}[The affine slice function and recovery of the original tensor]\label{lem:affine-slice-potential}
Use the notation of Lemma \ref{lem:local-cubic-potential}, and let $P$ be the homogeneous potential of degree two constructed there. Put
\[
 B_\varepsilon=\{w\in\R^m:|w|<\varepsilon\},\qquad f(w)=P(w,1),\qquad R(w)=\sqrt{1+|w|^2},\qquad u(w)=\frac{(w,1)}{R(w)}\in U_\varepsilon.
\]
Then $f$ is real-analytic on $B_\varepsilon$, and the following assertions hold.
\begin{enumerate}[label=\textup{(\arabic*)},leftmargin=2em,itemsep=4pt]
\item Recovery of the original tensor from third derivatives.
The potential on the open cone is determined by its affine slice:
\[
 P(w,t)=t^2f(w/t).
\]
For $a,b,c\in\{1,\ldots,m\}$, write $f_{abc}=\partial_a\partial_b\partial_cf$. Then
\begin{equation}\label{eq:pullback}
 C_{u(w)}(\partial_a u,\partial_bu,\partial_cu)
 =\frac{f_{abc}(w)}{1+|w|^2}.
\end{equation}
Thus the third derivatives of $f$ determine all values of $C$ on $U_\varepsilon$. In particular,
\[
 C_{u(w)}=0
 \quad\Longleftrightarrow\quad
 f_{abc}(w)=0\quad\text{for all }a,b,c\in\{1,\ldots,m\}.
\]
\item Normalization at a zero.
We may subtract a homogeneous quadratic polynomial from $P$, and correspondingly a polynomial of degree at most $2$ from $f$, so that
\begin{equation}\label{eq:affine-normalization}
 f(0)=0,\qquad f_a(0)=0,\qquad f_{ab}(0)=0.
\end{equation}
This adjustment leaves the third derivatives of $f$ unchanged. If $C_p=0$, then also $f_{abc}(0)=0$. Thus, if the normalized $f$ is not identically zero, its first nonzero homogeneous Taylor term has degree at least $4$.
\end{enumerate}
\end{lemma}

\begin{proof}

By its definition above, the affine slice $f$ is the restriction of $P$ to the hyperplane $t=1$; the point $(w,1)$ generally does not lie on the unit sphere. Homogeneity of degree two gives
\[
 P(w,t)=P\!\left(t\left(\frac wt,1\right)\right)
       =t^2f(w/t),\qquad (w,t)\in\mathcal K_\varepsilon.
\]
Thus $f$ determines $P$ on the open cone.

From now on, $f_a,f_{ab},f_{abc}$ denote partial derivatives with respect to the $m$ variables $w^1,\ldots,w^m$. Accordingly, $D$ also denotes the standard flat connection on the affine slice and its linear coordinate subspaces, and $D^3f$ denotes the covariant cubic tensor formed by all third partial derivatives. Since $t$ is fixed at $1$ when taking the restriction, the chain rule gives
\[
 f_a(w)=P_a(w,1),\quad
 f_{ab}(w)=P_{ab}(w,1),\quad
 f_{abc}(w)=P_{abc}(w,1)=\bar C_{abc}(w,1).
\]
The function $f$ has no $t$ variable, so notation such as $f_{abt}$ is not used; components containing $t$ belong to the ambient tensor $\bar C$.

Put $T_{abc}(w):=f_{abc}(w)$. Euler's identity for $P$ and its next two partial derivatives give the radial identity as follows:
\[
 q^\alpha P_\alpha=2P,\qquad P_\beta+q^\alpha P_{\alpha\beta}=2P_\beta, \qquad\text{that is, }q^\alpha P_{\alpha\beta}=P_\beta,\qquad P_{\gamma\beta}+q^\alpha P_{\alpha\beta\gamma} =P_{\beta\gamma}.
\]
The second derivatives in the last identity cancel, giving $q^\alpha\bar C_{\alpha\beta\gamma}=0$. Now restrict to $q=(w,1)$ and split the Greek-index sum into the $w$ directions and the $t$ direction.

First take $\beta=a,\gamma=b$. We obtain
\[
 0=w^c\bar C_{cab}(w,1)+\bar C_{tab}(w,1).
\]
Thus $\bar C_{abt}(w,1)=-w^cT_{abc}(w)$. Next take $\beta=a,\gamma=t$:
\[
 \bar C_{att}(w,1) =-w^b\bar C_{bat}(w,1) =-w^b\bigl(-w^cT_{abc}(w)\bigr) =w^bw^cT_{abc}(w).
\]
Finally, taking $\beta=\gamma=t$ gives
\[
 \bar C_{ttt}(w,1)=-w^a\bar C_{att}(w,1)
                 =-w^aw^bw^cT_{abc}(w).
\]
Collecting these three identities, with all right-hand sides evaluated at $w$ and all left-hand sides at $(w,1)$, we have
\begin{equation}\label{eq:radial-components}
 \bar C_{abt}=-w^cT_{abc},\qquad \bar C_{att}=w^bw^cT_{abc},\qquad \bar C_{ttt}=-w^aw^bw^cT_{abc}.
\end{equation}
Thus $T_{abc}$ determines every cubic tensor component on $t=1$, and homogeneity of degree minus one then determines $\bar C$ on the entire open cone.

The vector $e_a=\partial/\partial w^a$ is tangent to $t=1$ but generally not to the sphere, so $T_{abc}=\bar C_{(w,1)}(e_a,e_b,e_c)$ cannot in general be identified directly with a spherical coordinate component of $C$. Their relation is obtained as follows.

Put $q=(w,1)$ and $R=\sqrt{1+|w|^2}$, so that $q=Ru(w)$. Since
\[
 \partial_aR=\frac{w^a}{R},\qquad
 \partial_a(R^{-1})=-R^{-3}w^a,
\]
the product rule gives
\[
 \partial_a u =\partial_a\bigl(R^{-1}(w,1)\bigr) =R^{-1}e_a-R^{-3}w^a q.
\]
The second term is a multiple of the position vector $q$. Substituting three such expressions into $\bar C_q$ and expanding by multilinearity, all terms containing $q$ vanish by radial vanishing. Hence
\[
 \bar C_q(\partial_a u,\partial_bu,\partial_cu) =R^{-3}\bar C_q(e_a,e_b,e_c) =R^{-3}T_{abc}(w).
\]
On the other hand, differentiation of $|u(w)|^2=1$ gives $\langle\partial_a u,u\rangle=0$. Thus $\partial_a u\in T_uS^m$, and $\Pi_u$ leaves it unchanged. The definition of the radial extension therefore gives
\[
 \bar C_q(\partial_a u,\partial_bu,\partial_cu)
 =R^{-1}C_u(\partial_a u,\partial_bu,\partial_cu).
\]
Comparing the two identities and multiplying by $R$, we obtain
\[
 C_{u(w)}(\partial_a u,\partial_bu,\partial_cu)
 =R^{-2}T_{abc}(w)
 =\frac{f_{abc}(w)}{1+|w|^2}.
\]
The inverse of $u:B_\varepsilon\to U_\varepsilon$ is $(v,s)\mapsto v/s$, so $\partial_1u,\ldots,\partial_mu$ form a spherical tangent basis. The formula gives all components of $C$; since $1+|w|^2>0$, they vanish exactly when all $f_{abc}(w)$ vanish. This proves assertion \textup{(1)}.

We next normalize the potential at a point where the entire tensor vanishes.

At present $M(p)=H_0$. Put
\[
 Q(q)=\tfrac12q^\alpha (H_0)_{\alpha\beta}q^\beta.
\]
Since $H_0$ is a constant symmetric matrix, $D^2Q=H_0$ and $D^3Q=0$. Replacing $P,M$ by $P-Q$ and $M-H_0$, respectively, we still have
\[
 D^2P=M,\qquad D^3P=\bar C,
\]
and $P$ remains homogeneous of degree two, with $M(p)=0$. Replace $f$ by the new $P(w,1)$. Since $Q(w,1)$ has degree at most $2$, the third derivatives of $f$ and the recovery formula are unchanged. Equivalently, take $H_0=0$ in the integral formula, as we do henceforth. The preceding formulas give
\[
 P(p)=\tfrac12p^\alpha M_{\alpha\beta}(p)p^\beta=0,\qquad \partial_\alpha P(p)=M_{\alpha\beta}(p)p^\beta=0,\qquad \partial_{\alpha\beta}P(p)=M_{\alpha\beta}(p)=0.
\]

Since $p=(0,1)$, the preceding lower-order vanishing identities give the normalization conditions \eqref{eq:affine-normalization}.
If we also assume $C_p=0$, then the definition of the extension gives, for arbitrary ambient vectors $A,B,Z$,
\[
 \bar C_p(A,B,Z)=C_p(\Pi_pA,\Pi_pB,\Pi_pZ)=0.
\]
Hence
\[
 f_{abc}(0)=\bar C_{abc}(p)=0.
\]
Thus all derivatives of orders zero through three vanish. Real analyticity gives a convergent Taylor expansion
\[
 f(w)=\sum_{|\nu|\ge4}
       \frac{\partial^\nu f(0)}{\nu!}w^\nu,
\]
where $\nu=(\nu_1,\ldots,\nu_m)\in\mathbb N_0^m$, $|\nu|=\nu_1+\cdots+\nu_m$, $\nu!=\nu_1!\cdots\nu_m!$, and $w^\nu=(w^1)^{\nu_1}\cdots(w^m)^{\nu_m}$. If $f$ is not identically zero, its first nonzero homogeneous term has degree at least $4$; if $f$ is identically zero, there is no such term. The condition $C_p=0$ is used only here to eliminate the cubic term; it was not needed in the construction of the potential.

\end{proof}

The convention that Greek indices range over $1,\ldots,m+1$ was used for the ambient calculations in this subsection. In the next subsection and thereafter, when $\alpha,\beta$ occur as indices of $T$, they range only over $1,\ldots,m$, just as the other indices of $T$ do. Greek letters will then distinguish the angular indices under consideration, rather than denote the additional $t$ direction.

\subsection{The full affine-coordinate form of the commutation equations}\label{subsec:paper-3-2}
We now express $[C_X,C_Y]=0$ in terms of $D^3f$. We compute the contractions using the inverse spherical metric and give a brief ambient explanation of the positive radial term after the proof.

Using the notation of Lemma \ref{lem:affine-slice-potential}, put
\[
 R=\sqrt{1+|w|^2},\qquad q=(w,1),\qquad
 u(w)=\frac qR,\qquad E_a=\partial_a u,\qquad T_{abc}=f_{abc}.
\]
In this subsection the Roman indices $a,b,c,e,r,s$ all range over $1,\ldots,m$. We use $r,s$ as summation indices to avoid confusion with the spherical point $u(w)$. The symbol $w^a$ denotes the $a$th Euclidean coordinate of $w$; later we shall also write it as $w_a$ in Euclidean coordinates. This placement of indices does not involve raising or lowering them with the spherical metric.

\begin{proposition}[Affine-coordinate form of commutativity]\label{prop:affine-commutation}
Let $f$ be the function constructed in the preceding subsection. The contraction operators of the spherical tensor $C$ satisfy $[C_X,C_Y]=0$ if and only if, for all $a,b,c,e\in\{1,\ldots,m\}$,
\begin{equation}\label{eq:wdvv}
 \sum_{r,s=1}^m(\delta^{rs}+w^rw^s)
 (T_{abr}T_{ces}-T_{acr}T_{bes})=0,
\end{equation}
where $\delta^{rs}$ is the Kronecker symbol.
\end{proposition}

\begin{proof}

The calculation in the preceding subsection gives
\[
 E_a=R^{-1}e_a-R^{-3}w^a q,
\]
where $e_a$ is the $a$th constant unit coordinate vector of the ambient space. We also have
\[
 \langle e_a,e_b\rangle=\delta_{ab},\qquad
 \langle e_a,q\rangle=w^a,\qquad
 \langle q,q\rangle=R^2.
\]
Expanding the coordinate components of the spherical metric term by term, we obtain
\[
 \begin{aligned}
 (h_0)_{ab}:=h_0(E_a,E_b)
 &=\langle R^{-1}e_a-R^{-3}w^a q,
              R^{-1}e_b-R^{-3}w^b q\rangle\\
 &=R^{-2}\delta_{ab}-R^{-4}w^aw^b
                       -R^{-4}w^aw^b+R^{-6}w^aw^bR^2 =R^{-2}\delta_{ab}-R^{-4}w^aw^b.
 \end{aligned}
\]
Regarding $w$ as a column vector, we may write this matrix as
\[
 ((h_0)_{ab})=R^{-2}I-R^{-4}ww^{\mathsf T}.
\]
Its inverse is
\[
 ((h_0)^{rs})=R^2(I+ww^{\mathsf T}),\qquad
 (h_0)^{rs}=R^2(\delta^{rs}+w^rw^s).
\]
To verify this, use $(ww^{\mathsf T})^2=|w|^2ww^{\mathsf T}$ and multiply directly:
\[
 (R^{-2}I-R^{-4}ww^{\mathsf T})R^2(I+ww^{\mathsf T}) =I+ww^{\mathsf T}-R^{-2}ww^{\mathsf T} -R^{-2}|w|^2ww^{\mathsf T} =I+\left(1-\frac{1+|w|^2}{R^2}\right)ww^{\mathsf T} =I.
\]
Thus the rank-one term in the inverse metric does indeed have a positive sign.

To distinguish the spherical coordinate components from $T_{abc}$, temporarily write
\[
 \widetilde C_{abc}:=C(E_a,E_b,E_c).
\]
Lemma \ref{lem:affine-slice-potential} gives
\[
 \widetilde C_{abc}=R^{-2}T_{abc}.
\]
For fixed $a,b$, put $V_{ab}:=C_{E_a}E_b$ and write
\[
 V_{ab}=\sum_s V_{ab}^{\,s}E_s.
\]
The contraction operators are defined by $h_0(C_XY,Z)=C(X,Y,Z)$, so
\[
 \sum_s (h_0)_{rs}V_{ab}^{\,s}=\widetilde C_{abr},\qquad
 V_{ab}^{\,s}=\sum_r (h_0)^{sr}\widetilde C_{abr}.
\]
In particular, $V_{ab}=V_{ba}$. Taking the inner product of two such vectors gives
\[
 h_0(V_{ab},V_{ce}) =\sum_s V_{ab}^{\,s}h_0(E_s,V_{ce}) =\sum_{r,s}(h_0)^{sr}\widetilde C_{abr}\widetilde C_{ces}.
\]
Since $C$ is completely symmetric, each $C_X$ is self-adjoint. Hence
\[
 \begin{aligned}
 h_0(C_{E_b}C_{E_c}E_e,E_a)
 &=h_0(C_{E_c}E_e,C_{E_b}E_a)=h_0(V_{ce},V_{ab}),\\
 h_0(C_{E_c}C_{E_b}E_e,E_a)
 &=h_0(C_{E_b}E_e,C_{E_c}E_a)=h_0(V_{be},V_{ac}).
 \end{aligned}
\]
Subtracting according to the convention $[A,B]=AB-BA$, we obtain
\[
 h_0([C_{E_b},C_{E_c}]E_e,E_a) =\sum_{r,s}(h_0)^{rs} \bigl(\widetilde C_{abr}\widetilde C_{ces} -\widetilde C_{acr}\widetilde C_{bes}\bigr).
\]
Substitute the inverse metric computed above and $\widetilde C=R^{-2}T$. The right-hand side becomes
\[
 R^{-2}\sum_{r,s}(\delta^{rs}+w^rw^s)
 (T_{abr}T_{ces}-T_{acr}T_{bes}),
\]
where the common factor is $R^2\cdot R^{-2}\cdot R^{-2}=R^{-2}$.

If the contraction operators commute pairwise, the left-hand side vanishes, and cancellation of the nonzero factor $R^{-2}$ gives \eqref{eq:wdvv}. Conversely, if \eqref{eq:wdvv} holds for every choice of indices, then $[C_{E_b},C_{E_c}]E_e$ is orthogonal to every $E_a$. Since $E_1,\ldots,E_m$ form a basis of the tangent space and $h_0$ is positive definite, this vector is zero. Taking linear combinations of $E_e$, $E_b$, and $E_c$, respectively, gives $[C_X,C_Y]=0$ for arbitrary tangent vectors $X,Y$.
\end{proof}

The positive sign of the radial term can also be seen from the ambient components. At $q=(w,1)=\sum_r w^re_r+e_t$, we have $\bar C_{abr}=T_{abr}$, while radial vanishing gives
\[
 0=\bar C_q(e_a,e_b,q)=\sum_r w^rT_{abr}+\bar C_{abt},
 \qquad \bar C_{abt}=-\sum_r w^rT_{abr}.
\]
Thus the contribution of the last coordinate direction to the full ambient contraction is
\[
 \bar C_{abt}\bar C_{cet}-\bar C_{act}\bar C_{bet}
 =\sum_{r,s}w^rw^s(T_{abr}T_{ces}-T_{acr}T_{bes}),
\]
which is the radial contraction term in \eqref{eq:wdvv}; the positive sign is the product of two negative signs. The fixed coordinate vector $e_t=(0,1)$ is generally not parallel to the position vector $q=(w,1)$. Radial vanishing means that inserting $q$ into the tensor gives zero; it does not mean that components containing $t$ vanish. The contraction in the last coordinate direction therefore cannot be omitted.

Taking $a=b=i$, $c=\alpha$, and $e=\beta$ in \eqref{eq:wdvv}, we obtain
\begin{equation}\label{eq:Q}
 Q_{i\alpha\beta}:={}\sum_{r,s=1}^m(\delta^{rs}+w^rw^s)\cdot(T_{iir}T_{\alpha\beta s} -T_{i\alpha r}T_{i\beta s})=0.
\end{equation}
Here $\alpha,\beta$, like all the other indices of $T$, range over $1,\ldots,m$. In subsequent applications we take $\alpha,\beta\ne i$, but allow $\alpha=\beta$.

In \eqref{eq:wdvv}, $\delta^{rs}$ gives the Euclidean contraction, while $w^rw^s$ comes from the additional term in the inverse spherical metric. We shall also call the latter the radial contraction term, since in the ambient calculation it arises from the components containing $t$.

\begin{lemma}[Comparison of the lowest orders in the commutation equation]\label{lem:commutation-orders}
Let $f$ be the nonzero analytic potential normalized at a point where the entire tensor vanishes as in Lemma \ref{lem:affine-slice-potential}, and suppose that \eqref{eq:wdvv} holds. If its leading term is $f_{k+3}$, then $k\ge1$. Put $T^{(j)}=D^3f_{j+3}$. The part of degree $2k$ in the commutation equation contains only the Euclidean self-contraction of $T^{(k)}$, and the part of degree $2k+1$ contains only the Euclidean cross-contractions of $T^{(k)}$ and $T^{(k+1)}$. At degree $2k+2$, one must also include the self-contraction of $T^{(k+1)}$, the cross-contractions of $T^{(k)}$ and $T^{(k+2)}$, and the contraction of two copies of $T^{(k)}$ by $w^rw^s$. Higher homogeneous terms do not enter these three degrees.
\end{lemma}
\begin{proof}
Assume that the entire tensor $C$ vanishes at $p$, and normalize $f$ as in Lemma \ref{lem:affine-slice-potential}. If $f$ is not identically zero near $0$, real analyticity ensures that its Taylor expansion has a first nonzero term. We may write
\[
 f=f_{k+3}+f_{k+4}+f_{k+5}+\cdots,\qquad k\ge1,
\]
where $f_d$ is a homogeneous polynomial of degree $d$ and $f_{k+3}\ne0$. Taking third derivatives term by term gives
\[
 T=T^{(k)}+T^{(k+1)}+T^{(k+2)}+\cdots,\qquad
 T^{(j)}=D^3f_{j+3}.
\]
At least one component of $T^{(k)}$ is nonzero. Otherwise every third derivative of $f_{k+3}$ would vanish, so it would be a polynomial of degree at most $2$, contrary to its degree being at least $4$.

\begin{samepage}
The product of tensor components of degrees $j,\ell$ has degree $j+\ell$; multiplication by $w^rw^s$ raises this to $j+\ell+2$. Consequently, comparison of homogeneous parts in \eqref{eq:wdvv} gives:
\begin{itemize}[leftmargin=2em,itemsep=3pt]
\item degree $2k$ contains only the contraction of $T^{(k)}$ with $T^{(k)}$ by $\delta^{rs}$;
\item degree $2k+1$ contains only the cross-contractions of $T^{(k)}$ and $T^{(k+1)}$ by $\delta^{rs}$;
\item at degree $2k+2$, the Euclidean contractions come from the three degree pairs $(k,k+2)$, $(k+1,k+1)$, and $(k+2,k)$, together with the contraction of two copies of $T^{(k)}$ by $w^rw^s$.
\end{itemize}
\end{samepage}
In particular, the equation of lowest degree is
\[
 \sum_r\bigl(T^{(k)}_{abr}T^{(k)}_{cer}
              -T^{(k)}_{acr}T^{(k)}_{ber}\bigr)=0.\qedhere
\]
\end{proof}
The lowest-degree relation is studied next; the radial term enters at degree $2k+2$ in the Gram relations of Section \ref{sec:next-layer}.

We shall often restrict \eqref{eq:Q} to a coordinate subspace. All third partial derivatives must first be computed, and only then may the coordinate equations of the subspace be substituted. Differentiating the restricted function alone does not recover all the components $T_{abc}$.

For example, in two variables take $f(x,z)=x^3z$. Although
\[
 f|_{z=0}=0,
\]
we have
\[
 \partial_zf=x^3,\qquad
 \partial_x\partial_zf=3x^2,\qquad
 \left.\partial_x^2\partial_zf\right|_{z=0}=6x.
\]
Thus vanishing of the restricted function does not imply that its ambient third derivatives involving normal directions vanish on the subspace. The internal sum in the commutation equation still ranges over all $r,s\in\{1,\ldots,m\}$; setting a coordinate equal to zero does not justify deleting derivatives in that direction. Tangential differentiation along the subspace does commute with restriction. What must be retained are the normal derivatives, which the restricted function itself does not determine.

\begin{proposition}[All third derivatives and a local zero subsphere]\label{prop:affine-zero-subsphere}
Use the coordinate map $u(w)$ of Lemma \ref{lem:affine-slice-potential}. Let $L\subset\R^m$ be a linear subspace through the origin. Within the coordinate domain, all components of $D^3f$ vanish along $L$ if and only if $C$, as an entire spherical cubic tensor, vanishes along $u(L)$. Here $u(L)$ is the portion of the great subsphere $(L\times\R)\cap S^m$ lying in this coordinate neighbourhood.
\end{proposition}
\begin{proof}
By the recovery formula \eqref{eq:pullback} and $R^{-2}>0$, vanishing of all $T_{abc}$ at a point is equivalent to $C=0$ at the corresponding spherical point. This concerns all third derivatives, not merely the tangential components along a subspace.

Moreover, if $L\subset\R^m$ is a linear subspace through $0$, then
\[
 u(L\cap B_\varepsilon)
 =U_\varepsilon\cap\bigl((L\times\R)\cap S^m\bigr).
\]
The right-hand side is the portion of a great subsphere through $p$ contained in the coordinate neighbourhood. Thus proving that all components of $T$ vanish on a coordinate subspace gives local vanishing of the entire tensor on the corresponding great subsphere.

\end{proof}

\subsection{Fixed orthogonal separation of the leading homogeneous term}\label{subsec:paper-3-3}\label{sec:leading}
For the leading term $g=f_{k+3}$, the contraction matrices of $D^3g$ commute by the preceding subsection. Pointwise simultaneous diagonalization is insufficient for the expansion: using homogeneity, we obtain fixed orthogonal coordinates and separation into pure powers.

\begin{lemma}[Orthogonal separation of the leading homogeneous term]\label{lem:homogeneous-separation}
Let $g\ne0$ be a homogeneous polynomial of degree $k+3$ on $\R^m$, where $k\ge1$. In standard orthogonal coordinates $w=(w^1,\ldots,w^m)$, define the symmetric matrices
\[
 A_i(w):=\bigl(g_{iab}(w)\bigr)_{a,b=1}^m,\qquad
 g_{iab}=\partial_i\partial_a\partial_b g.
\]
Suppose that $[A_i(w),A_j(w)]=0$ for all $w\in\R^m$ and $i,j\in\{1,\ldots,m\}$. Then there exists a constant orthogonal matrix $Q$ such that, in the new coordinates $w=Qx$, writing $g(x)$ for $g(Qx)$, we have
\begin{equation}\label{eq:leading}
 g(x)=\sum_{i=1}^r\frac{a_i x_i^{k+3}}{(k+1)(k+2)(k+3)},
 \qquad a_i\ne0,\quad 1\le r\le m.
\end{equation}
The corresponding third derivatives are
\[
 D^3g=\sum_{i=1}^r a_i x_i^k(\dd x_i)^{\otimes3}.
\]
\end{lemma}

\begin{proof}
Put
\[
 \mathsf{H}=\Hess g=(g_{ab}).
\]
Each $g_{ab}$ is a homogeneous polynomial of degree $k+1$. Applying Euler's identity componentwise gives
\[
 \sum_i w^i\partial_i g_{ab}=(k+1)g_{ab},
 \qquad\text{that is,}\qquad
 (k+1)\mathsf{H}=\sum_i w^i\partial_i\mathsf{H}.
\]
Since $\partial_i\mathsf{H}=A_i$, we have
\[
 \mathsf{H}=\frac1{k+1}\sum_i w^i A_i.
\]
At each fixed point the coefficients $w^i$ are real numbers. Thus, for every $j$,
\[
 [\mathsf{H},\partial_j\mathsf{H}]
 =[\mathsf{H},A_j]
 =\frac1{k+1}\sum_i w^i[A_i,A_j]=0.
\]
Hence $\partial_j\mathsf{H}$ preserves every eigenspace of $\mathsf{H}$. We must also prove that these eigenspaces themselves are independent of the point.

To construct spectral projections, first choose an open set on which the eigenvalue multiplicities are constant.
\begin{samepage}
Let $s(w)$ be the number of distinct eigenvalues of the real symmetric matrix $\mathsf{H}(w)$. Since $s(w)$ takes values in the finite set $1,\ldots,m$, we may choose $w_0$ where it attains its maximum $s$. List the distinct eigenvalues of $\mathsf{H}(w_0)$ in increasing order as
\[
 \lambda_1^0<\cdots<\lambda_s^0,
\]
with multiplicities $n_1,\ldots,n_s$, where $n_1+\cdots+n_s=m$.
\end{samepage}

Choose pairwise disjoint open intervals around the $\lambda_\alpha^0$. The ordered eigenvalues of a real symmetric matrix depend continuously on the matrix. Hence, on a sufficiently small open ball $\Omega$ about $w_0$, the $\alpha$th interval contains exactly $n_\alpha$ eigenvalues, counted with multiplicity.

The eigenvalues within each interval must be equal. Indeed, each interval contains at least one distinct eigenvalue; if one interval contained two or more distinct eigenvalues, their total number would exceed $s$, contrary to maximality. Thus, on $\Omega$, we may write
\[
 \lambda_1(w)<\cdots<\lambda_s(w),
\]
and the multiplicity of $\lambda_\alpha(w)$ is constantly $n_\alpha$. Continuity of the eigenvalues also shows that these functions are continuous.

We next prove that they are real-analytic, so that the spectral projections can be differentiated. For $j=1,\ldots,s$, the spectral theorem gives
\[
 \tr(\mathsf{H}^j)=\sum_{\alpha=1}^s n_\alpha\lambda_\alpha^j.
\]
Consider the map
\[
 \Phi(\ell_1,\ldots,\ell_s)
 =\left(\sum_{\alpha=1}^s n_\alpha\ell_\alpha^j\right)_{j=1}^s.
\]
Its Jacobian matrix is $J_{j\alpha}=j n_\alpha\ell_\alpha^{j-1}$. At $\ell_\alpha=\lambda_\alpha^0$,
\[
 \det J
 =\left(\prod_{j=1}^s j\right)
  \left(\prod_{\alpha=1}^s n_\alpha\right)
  \prod_{\alpha<\beta}(\lambda_\beta^0-\lambda_\alpha^0)
 \ne0.
\]
The last product is a Vandermonde determinant. By the analytic inverse function theorem, $\Phi$ has an analytic inverse near this point. Since the entries of $\mathsf{H}$ are polynomials, each $\tr(\mathsf{H}^j)$ is real-analytic. After shrinking $\Omega$, the continuous tuple of eigenvalues must satisfy
\[
 (\lambda_1(w),\ldots,\lambda_s(w))
 =\Phi^{-1}\bigl(\tr \mathsf{H}(w),\ldots,\tr(\mathsf{H}(w)^s)\bigr).
\]
Thus every $\lambda_\alpha$ is real-analytic.

Define the matrix fields
\[
 P_\alpha
 =\prod_{\beta\ne\alpha}
  \frac{\mathsf{H}-\lambda_\beta I}{\lambda_\alpha-\lambda_\beta}.
\]
If $s=1$, the empty product is understood as $I$. On the $\lambda_\gamma$ eigenspace of $\mathsf{H}$, this polynomial equals $1$ when $\gamma=\alpha$ and $0$ when $\gamma\ne\alpha$. Hence $P_\alpha$ is precisely the orthogonal projection onto the $\lambda_\alpha$ eigenspace. It is real-analytic and satisfies
\[
 \mathsf{H}P_\alpha=P_\alpha \mathsf{H}=\lambda_\alpha P_\alpha,\qquad P_\alpha^{\mathsf T}=P_\alpha,\qquad P_\alpha^2=P_\alpha,\qquad P_\alpha P_\beta=0\quad(\alpha\ne\beta),\qquad \sum_{\alpha=1}^sP_\alpha=I.
\]

Fix a coordinate direction $w^i$. Multiplying the identity $[\mathsf{H},\partial_i\mathsf{H}]=0$ on the left by $P_\beta$ and on the right by $P_\alpha$ gives
\[
 0=P_\beta[\mathsf{H},\partial_i\mathsf{H}]P_\alpha =(\lambda_\beta-\lambda_\alpha) P_\beta(\partial_i\mathsf{H})P_\alpha.
\]
Therefore, when $\beta\ne\alpha$,
\[
 P_\beta(\partial_i\mathsf{H})P_\alpha=0.
\]
Now differentiate $\mathsf{H}P_\alpha=\lambda_\alpha P_\alpha$:
\[
 (\partial_i\mathsf{H})P_\alpha+\mathsf{H}(\partial_iP_\alpha)
 =(\partial_i\lambda_\alpha)P_\alpha
                   +\lambda_\alpha(\partial_iP_\alpha).
\]
Multiply on the left by $P_\beta$, with $\beta\ne\alpha$, and use $P_\beta P_\alpha=0$. This gives
\[
 P_\beta(\partial_i\mathsf{H})P_\alpha
 =(\lambda_\alpha-\lambda_\beta)P_\beta(\partial_iP_\alpha).
\]
The left-hand side is zero and the two eigenvalues are distinct, so
\[
 P_\beta(\partial_iP_\alpha)=0\qquad(\beta\ne\alpha).
\]
Summing over all $\beta\ne\alpha$ gives
\[
 (I-P_\alpha)(\partial_iP_\alpha)=0.
\]
Since both $P_\alpha$ and $\partial_iP_\alpha$ are symmetric, transposition also gives
\[
 (\partial_iP_\alpha)(I-P_\alpha)=0.
\]
These two identities state, respectively, that
\[
 \partial_iP_\alpha=P_\alpha(\partial_iP_\alpha),\qquad
 \partial_iP_\alpha=(\partial_iP_\alpha)P_\alpha.
\]
Together they give
\[
 \partial_iP_\alpha=P_\alpha(\partial_iP_\alpha)P_\alpha.
\]
Finally, differentiation of $P_\alpha^2=P_\alpha$ yields
\[
 (\partial_iP_\alpha)P_\alpha
       +P_\alpha(\partial_iP_\alpha)=\partial_iP_\alpha.
\]
Multiplying on both sides by $P_\alpha$ and using $P_\alpha^2=P_\alpha$, we obtain
\[
 2P_\alpha(\partial_iP_\alpha)P_\alpha
   =P_\alpha(\partial_iP_\alpha)P_\alpha.
\]
Thus $P_\alpha(\partial_iP_\alpha)P_\alpha=0$, and hence
\[
 \partial_iP_\alpha=0.
\]
Every coordinate partial derivative vanishes. Since $\Omega$ is a connected open ball, $P_\alpha$ is a constant matrix. Thus each $\operatorname{im}P_\alpha$ is a fixed linear subspace.

Choose a fixed orthonormal basis in each fixed subspace $\operatorname{im}P_\alpha$. These subspaces are mutually orthogonal and their direct sum is $\R^m$, so the combined bases form a fixed orthonormal basis of $\R^m$. Let $Q$ be the constant orthogonal matrix with these basis vectors as columns.

By the spectral decomposition
\[
 \mathsf{H}(w)=\sum_{\alpha=1}^s\lambda_\alpha(w)P_\alpha,
\]
the matrix $\mathsf{H}(w)$ is diagonal in this basis for every $w\in\Omega$. If an eigenvalue has multiplicity greater than one, $\mathsf{H}$ acts as a scalar on its eigenspace, so any fixed orthonormal basis of that space will do.

Put $w=Qx$ and temporarily write $\widehat g(x)=g(Qx)$. Since $Q$ is constant, the chain rule gives
\[
 \partial_{x_i}\partial_{x_j}\widehat g(x)
 =\sum_{a,b}Q_{ai}g_{ab}(Qx)Q_{bj}
 =\bigl(Q^{\mathsf T}\mathsf{H}(Qx)Q\bigr)_{ij}.
\]
For $x\in Q^{\mathsf T}\Omega$, the right-hand side is zero when $i\ne j$.

\begin{samepage}
Each $\partial_{x_i}\partial_{x_j}\widehat g$ is a polynomial, and $Q^{\mathsf T}\Omega$ is a nonempty open set. A polynomial vanishing on a nonempty open set vanishes identically: at a point in that set, all its partial derivatives vanish because it is zero in a neighbourhood; the polynomial equals its finite Taylor expansion at that point and is therefore zero everywhere. Consequently,
\[
 \partial_{x_i}\partial_{x_j}\widehat g\equiv0
 \qquad(i\ne j)
\]
holds on all of $\R^m$. We henceforth write $g$ for $\widehat g$.
\end{samepage}

There is no requirement that $\Omega$ contain the origin. Indeed, the entries of $\mathsf{H}$ have degree $k+1\ge2$, so $\mathsf{H}(0)=0$ and all eigenvalues coincide at the origin. Finding fixed coordinates first on the above open set and then extending by the polynomial identity principle avoids selecting eigenvectors pointwise near the origin.

Finally, we deduce separation into pure powers from the vanishing of mixed second derivatives.
We know that $\partial_j(\partial_i g)=0$ for every $j\ne i$. Hence $\partial_i g$ depends only on $x_i$, and we may write
\[
 \partial_i g=p_i(x_i).
\]
A linear change of coordinates preserves homogeneous degree, so $\partial_i g$ is a homogeneous polynomial of degree $k+2$. Expanding the one-variable polynomial $p_i$ and comparing degrees gives
\[
 p_i(x_i)=b_i x_i^{k+2}
\]
where $b_i$ is constant and may be zero. Applying Euler's identity to $g$ once more gives
\[
 (k+3)g=\sum_i x_i\partial_i g
        =\sum_i b_i x_i^{k+3}.
\]
Put $a_i=(k+1)(k+2)b_i$ and reorder the coordinates so that the nonzero coefficients occupy the first $r$ positions. This gives \eqref{eq:leading}. Since $g\ne0$, at least one coefficient is nonzero, and therefore $1\le r\le m$.

Finally, take third derivatives term by term. For each nonzero term,
\[
 \frac{\partial^3}{\partial x_i^3}
 \left(\frac{a_i x_i^{k+3}}{(k+1)(k+2)(k+3)}\right)
 =a_i x_i^k.
\]
All mixed third derivatives involving distinct coordinates vanish. Therefore
\[
 D^3g=\sum_{i=1}^r a_i x_i^k(\dd x_i)^{\otimes3},
\]

\end{proof}

To apply Lemma \ref{lem:homogeneous-separation} to the first nonzero homogeneous term $g=f_{k+3}$, we first verify the matrix commutation hypothesis. The degree $2k$ part of \eqref{eq:wdvv} gives
\[
 \sum_r(g_{abr}g_{cer}-g_{acr}g_{ber})=0.
\]
Complete symmetry of the third partial derivatives gives
\[
 (A_bA_c)_{ae} =\sum_r g_{bar}g_{cre}=\sum_r g_{abr}g_{cer},\qquad (A_cA_b)_{ae} =\sum_r g_{car}g_{bre}=\sum_r g_{acr}g_{ber}.
\]
Subtracting yields $[A_b,A_c]=0$. Since all components are homogeneous polynomials, this is a polynomial identity for all $w\in\R^m$. The lemma therefore applies to the leading term of the normalized potential.
After applying the lemma to the first nonzero homogeneous term of $f$, apply the same orthogonal change $w=Qx$ to the entire function $f$, retaining the notation $f$. Since $Q$ is a constant orthogonal matrix, $|Qx|=|x|$; the third partial derivatives transform as tensor components, and internal contractions are preserved. Hence \eqref{eq:wdvv} remains valid in the new coordinates and has the same form.

The first $r$ coordinate directions are called the \textbf{active axes} of this layer. Denote the remaining variables by
\[
 z=(x_{r+1},\ldots,x_m).
\]
These are called the \textbf{remaining variables} or \textbf{inactive variables}. Here fixed means that $Q$ is independent of the point; inactive means only that these variables do not occur in the current lowest homogeneous term. They may still occur in higher-order terms of $f$. If $r=m$, there are no remaining variables in this layer.

\subsection{The next homogeneous term and Gram relations within one layer}\label{subsec:paper-3-4}\label{sec:next-layer}
Keeping the fixed coordinates of the preceding subsection, we compare degrees $2k+1$ and $2k+2$ to determine the next homogeneous term and the same-layer Gram relations, respectively.

Suppose that the normalized $f$ is not identically zero near $0$. In the fixed orthogonal coordinates of Lemma \ref{lem:homogeneous-separation}, write
\[
 \begin{gathered}
 w=(x,z)=(x_1,\ldots,x_r,z_1,\ldots,z_{m-r}),\qquad
 f=g+\mathsf{h}+p+\sum_{d\ge k+6}f_d,\\
 g=f_{k+3},\quad \mathsf{h}=f_{k+4},\quad p=f_{k+5}.
 \end{gathered}
\]
Here $g$ is given by \eqref{eq:leading}, with $a_i\ne0$, and $f_d$ denotes a homogeneous polynomial of degree $d$, which may be zero. Put
\[
 L=D^3g,\qquad U=D^3\mathsf{h},\qquad \mathsf{V}=D^3p.
\]
Their degrees are $k,k+1,k+2$, respectively, and $T=D^3f=L+U+\mathsf{V}+\cdots$.

Active indices $i,j$ range over $1,\ldots,r$, and remaining indices $\alpha,\beta,\gamma$ over $1,\ldots,m-r$. Components with Greek indices are written according to the names of the variables; for example,
\[
 U_{ii\alpha}
 :=\frac{\partial^3\mathsf{h}}{\partial x_i^2\partial z_\alpha},
 \qquad
 U_{i\alpha\beta}
 :=\frac{\partial^3\mathsf{h}}{\partial x_i\partial z_\alpha\partial z_\beta}.
\]
If all coordinates are numbered consecutively from $1,\ldots,m$, the $z_\alpha$ direction has index $r+\alpha$. In the matrix calculations below, $a,b,\mu,\nu,s$ all refer to this common numbering; in block calculations we return to active and remaining indices. If $r=m$, all sums over remaining indices are empty.

\begin{lemma}[Form of the next homogeneous term]\label{lem:next-homogeneous-form}
With the preceding notation, if $f$ satisfies the full commutation equation \eqref{eq:wdvv}, then its next homogeneous term has the form
\begin{equation}\label{eq:next}
 \mathsf{h}=f_{k+4}={}\sum_{i=1}^r\beta_i x_i^{k+4} +\sum_{i=1}^r\sum_{\alpha=1}^{m-r} \frac{d_{i\alpha}x_i^{k+3}z_\alpha}{(k+2)(k+3)} +G_{k+4}(z),
\end{equation}
where $\beta_i,d_{i\alpha}$ are constants and $G_{k+4}$ is a homogeneous polynomial of degree $k+4$ in the remaining variables alone.
\end{lemma}

\begin{proof}

Put $c_i=a_i x_i^k$. In the common numbering, set $c_a=a_a x_a^k$ for active directions and $c_a=0$ for remaining directions. The contraction matrices of the leading term and of the next term are, respectively,
\[
 L_a=(L_{a\mu\nu})=c_a e_a e_a^{\mathsf T},\qquad
 U_a=(U_{a\mu\nu}),
\]
where $e_a$ is the constant unit coordinate vector. Thus the only possibly nonzero components of $L$ are $L_{iii}=c_i$.

At degree $2k+1$, \eqref{eq:wdvv} contains only the Euclidean cross-products of $L$ and $U$; radial contractions first occur at degree $2k+2$ or higher. Using the matrix notation of Section \ref{sec:leading}, we obtain
\[
 [L_a,U_b]+[U_a,L_b]=0.
\]
To use this identity componentwise, first compute the components of its four matrix products:
\[
 (L_aU_b)_{\mu\nu}=c_a\delta_{\mu a}U_{ba\nu},\quad (U_bL_a)_{\mu\nu}=c_aU_{b\mu a}\delta_{\nu a},\quad (U_aL_b)_{\mu\nu}=c_bU_{a\mu b}\delta_{\nu b},\quad (L_bU_a)_{\mu\nu}=c_b\delta_{\mu b}U_{ab\nu}.
\]
The linear constraint is therefore explicitly
\[
 0={}c_a(\delta_{\mu a}U_{ba\nu}-U_{b\mu a}\delta_{\nu a}) +c_b(U_{a\mu b}\delta_{\nu b}-\delta_{\mu b}U_{ab\nu}).
\]

We first eliminate components involving three distinct coordinate directions.
Take the active direction $a=i$ and directions $b,s$ such that $i,b,s$ are pairwise distinct, and set $\mu=i,\nu=s$. Only the first term in the preceding equation remains, giving
\[
 c_iU_{bis}=0.
\]
The polynomial $c_i=a_ix_i^k$ is nonzero, and the polynomial ring has no zero divisors. Hence $U_{bis}=0$. Thus any component of $U$ involving three distinct coordinate directions vanishes if at least one of them is active. This is a polynomial identity and also holds at points where $x_i=0$.

We next eliminate the remaining mixed components between distinct active axes.
Fix distinct active indices $i,j$ and set $a=i,b=j,\mu=i,\nu=j$ in the matrix equation above. This gives
\[
 c_iU_{ijj}+c_jU_{iij}=0.
\]
Since $c_i,c_j$ involve distinct variables, they are relatively prime. The divisibility $c_i\mid c_jU_{iij}$ therefore implies $c_i\mid U_{iij}$, so there is a polynomial $\ell$ such that
\[
 U_{iij}=c_i\ell,\qquad U_{ijj}=-c_j\ell.
\]
Since $U$ has degree $k+1$ and $c_i$ has degree $k$, the polynomial $\ell$ is homogeneous of degree one and may be zero.

As $U=D^3\mathsf{h}$, symmetry of the fourth partial derivatives gives
\[
 \partial_jU_{iij}=\partial_iU_{ijj}.
\]
Substituting the preceding expressions and using $\partial_jc_i=\partial_ic_j=0$, we obtain
\[
 a_ix_i^k\partial_j\ell=-a_jx_j^k\partial_i\ell.
\]
Both $\partial_i\ell,\partial_j\ell$ are constants. Since $k\ge1$, the monomials $x_i^k,x_j^k$ are distinct, and comparison of coefficients gives
\[
 \partial_i\ell=\partial_j\ell=0.
\]
For any direction $s$ distinct from $i,j$, the vanishing result for three distinct coordinate directions gives $U_{ijs}=0$. Hence
\[
 \partial_sU_{iij}=\partial_iU_{ijs}=0,
 \qquad c_i\partial_s\ell=0.
\]
Again using the absence of zero divisors in the polynomial ring, we obtain $\partial_s\ell=0$. Every partial derivative of $\ell$ therefore vanishes, so it is constant; since it is homogeneous of degree one, it is zero. Thus
\[
 U_{iij}=U_{ijj}=0\qquad(i\ne j\text{ both active}).
\]

We next determine the mixed components involving remaining directions.
Fix an active index $i$ and a remaining index $\alpha$. In the matrix equation with common numbering, take
\[
 a=i,\qquad b=r+\alpha,\qquad\mu=i,\qquad\nu=r+\alpha.
\]
Since $c_{r+\alpha}=0$, we obtain $c_iU_{i\alpha\alpha}=0$, and hence $U_{i\alpha\alpha}=0$. Together with the preceding vanishing result for three distinct directions, this gives
\[
 U_{i\alpha\beta}=0
 \qquad\text{for all remaining indices }\alpha,\beta.
\]

Now consider $U_{ii\alpha}$. For any coordinate direction $s$ other than $x_i$, symmetry of partial derivatives gives
\[
 \partial_sU_{ii\alpha}=\partial_iU_{i\alpha s}.
\]
If $s$ is a remaining direction, the right-hand side vanishes by the result just proved. If $s$ is another active direction, it involves three distinct directions and vanishes by the preceding result for such components. Hence $U_{ii\alpha}$ depends only on $x_i$. Since it is homogeneous of degree $k+1$, there is a constant $d_{i\alpha}$ such that
\[
 U_{ii\alpha}=d_{i\alpha}x_i^{k+1}.
\]

For $U_{iii}$, the vanishing of mixed components between distinct active axes and symmetry of partial derivatives give
\[
 \partial_jU_{iii}=\partial_iU_{iij}=0\qquad(j\ne i\text{ active}),
\]
while for every remaining direction,
\[
 \partial_{z_\alpha}U_{iii}
 =\partial_iU_{ii\alpha}
 =(k+1)d_{i\alpha}x_i^k.
\]
Thus subtracting $(k+1)\sum_\alpha d_{i\alpha}x_i^kz_\alpha$ from $U_{iii}$ leaves a function depending only on $x_i$. Homogeneity then gives
\[
 U_{iii}=b_ix_i^{k+1}
        +(k+1)\sum_\alpha d_{i\alpha}x_i^kz_\alpha,
\]
where $b_i$ is constant.

Finally, we determine the purely remaining components and recover the potential.
For a component all of whose indices are remaining directions,
\[
 \partial_iU_{\alpha\beta\gamma}
 =\partial_{z_\alpha}U_{i\beta\gamma}=0,
\]
so $U_{\alpha\beta\gamma}$ depends only on $z$. Put
\[
 G_{k+4}(z)=\mathsf{h}(0,z),\qquad
 \beta_i=\frac{b_i}{(k+2)(k+3)(k+4)},
\]
and define a candidate polynomial $\widehat{\mathsf{h}}$ by the right-hand side of \eqref{eq:next}. Its purely active components are given by
\[
 \partial_i^3(\beta_i x_i^{k+4})=b_ix_i^{k+1}
\]
The two possibly nonzero components of the mixed terms are
\[
 \partial_i^2\partial_{z_\alpha} \left(\frac{d_{i\alpha}x_i^{k+3}z_\alpha}{(k+2)(k+3)}\right) =d_{i\alpha}x_i^{k+1},\qquad \partial_i^3 \left(\frac{d_{i\alpha}x_i^{k+3}z_\alpha}{(k+2)(k+3)}\right) =(k+1)d_{i\alpha}x_i^kz_\alpha.
\]
The third derivatives in remaining directions alone are supplied by $G_{k+4}$ and agree with $U_{\alpha\beta\gamma}$: since the latter is independent of the active variables, it suffices to compare at $x=0$. All other mixed third derivatives vanish, in agreement with the calculations of each type of mixed component above.

Thus $D^3(\mathsf{h}-\widehat{\mathsf{h}})=0$: all third partial derivatives of the difference vanish identically, so the difference is a polynomial of degree at most $2$. On the other hand, $\mathsf{h}-\widehat{\mathsf{h}}$ is homogeneous of degree $k+4\ge5$, and therefore is zero. This proves \eqref{eq:next}.
\end{proof}

In \eqref{eq:next}, the mixed term associated with the $i$th active axis is
\[
 \frac{x_i^{k+3}}{(k+2)(k+3)}
       \sum_{\alpha=1}^{m-r} d_{i\alpha}z_\alpha,
\]
where $\sum_\alpha d_{i\alpha}z_\alpha$ is a linear function on the space of remaining variables. Under the standard Euclidean inner product it is represented by a constant vector:
\[
 d_i:=(d_{i1},\ldots,d_{i,m-r})\in\R^{m-r},\qquad
 \sum_\alpha d_{i\alpha}z_\alpha=\ip{d_i}{z}.
\]
The coefficient of this active axis in the leading term is $a_i\ne0$. Dividing the coefficient vector of the mixed term by the leading coefficient, define
\[
 \xi_i:=\frac{d_i}{a_i}
       =\left(\frac{d_{i\alpha}}{a_i}\right)_{\alpha=1}^{m-r}
       \in\R^{m-r}.
\]
Each $\xi_i$ is the constant vector of mixed-to-leading coefficient ratios associated with the $i$th active axis. It lies in the remaining-variable space and need not be a unit vector. The mixed term and tensor component are, respectively,
\[
 \frac{a_i x_i^{k+3}}{(k+2)(k+3)}\ip{\xi_i}{z},\qquad
 U_{ii\alpha}=a_i x_i^{k+1}\xi_{i\alpha}.
\]
For example, if there are two remaining variables, $a_i=2$, and the mixed term is
\[
 \frac{x_i^{k+3}(3z_1-2z_2)}{(k+2)(k+3)},
\]
then $d_i=(3,-2)$ and $\xi_i=(3/2,-1)$. This example only illustrates the correspondence between coefficients and vectors; vectors associated with distinct $i$ must also satisfy the constraints below.

\begin{lemma}[Gram relations within one layer]\label{lem:same-layer-gram}
Under the hypotheses of Lemma \ref{lem:next-homogeneous-form}, put $\xi_i=d_i/a_i\in\R^{m-r}$. Then, for any two distinct active indices $i,j$,
\begin{equation}\label{eq:gram}
 \ip{\xi_i}{\xi_j}=-1\qquad(i\ne j).
\end{equation}
\end{lemma}

\begin{proof}
Fix distinct active indices $i,j$. Use the particular component of \eqref{eq:Q}
\[
 Q_{ijj}=\sum_{s,v=1}^m(\delta^{sv}+w^sw^v)
        (T_{iis}T_{jjv}-T_{ijs}T_{ijv})=0,
\]
where $s,v$ range over all $m$ angular coordinates, including every remaining direction.

We first list all contributions of degree $2k+2$.
In the Euclidean contraction, this degree comes from the cross-products of $L$ and $\mathsf{V}$ and the product of $U$ with itself:
\[
 \sum_s\bigl(L_{iis}\mathsf{V}_{jjs}+\mathsf{V}_{iis}L_{jjs}+U_{iis}U_{jjs} -L_{ijs}\mathsf{V}_{ijs}-\mathsf{V}_{ijs}L_{ijs}-U_{ijs}^2\bigr).
\]
Since $i\ne j$, we have $L_{ijs}=0$. Moreover, $L_{iis}$ equals $c_i$ only when $s=i$, and $L_{jjs}$ equals $c_j$ only when $s=j$. Thus the expression reduces to
\[
 c_i\mathsf{V}_{ijj}+c_j\mathsf{V}_{iij}
 +\sum_s(U_{iis}U_{jjs}-U_{ijs}^2).
\]

The radial contraction introduces two additional coordinate factors, so only two copies of the leading term $L$ can contribute at degree $2k+2$. Its contribution is
\[
 \sum_{s,v}w^sw^v(L_{iis}L_{jjv}-L_{ijs}L_{ijv}) =x_ix_jc_ic_j =a_ia_jx_i^{k+1}x_j^{k+1}.
\]
Combining the two contributions gives
\[
 0={}c_i\mathsf{V}_{ijj}+c_j\mathsf{V}_{iij} +\sum_s(U_{iis}U_{jjs}-U_{ijs}^2) +a_ia_jx_i^{k+1}x_j^{k+1}.
\]

After computing all derivatives, restrict to the plane spanned by the two active axes.
Let $\Pi_{ij}$ be the coordinate plane on which only $x_i,x_j$ may be nonzero; all other active coordinates and all remaining coordinates are set equal to zero. By the preceding lemma, distinct active indices $i,j$ satisfy
\[
 U_{ijs}=0\qquad\text{for every coordinate direction }s.
\]
For the product $U_{iis}U_{jjs}$, examine each type of internal summation direction:
\begin{itemize}[leftmargin=2em,itemsep=3pt]
\item if $s=i$, then $U_{jji}=0$, so the product vanishes;
\item if $s=j$, then $U_{iij}=0$, so the product vanishes;
\item if $s$ is another active direction, both components vanish;
\item if $s$ is a remaining direction $z_\alpha$, then
\[
 U_{ii\alpha}U_{jj\alpha}
 =d_{i\alpha}d_{j\alpha}x_i^{k+1}x_j^{k+1}.
\]
\end{itemize}
Therefore
\[
 \left.\sum_s(U_{iis}U_{jjs}-U_{ijs}^2)\right|_{\Pi_{ij}}
 =\left(\sum_\alpha d_{i\alpha}d_{j\alpha}\right)
       x_i^{k+1}x_j^{k+1}.
\]
Although $z=0$, derivatives in remaining directions still enter the contraction: these components are obtained by differentiating first and then restricting.

Put $q(x_i,x_j):=p|_{\Pi_{ij}}$. This is a homogeneous polynomial of degree $k+5$ in two variables, possibly zero. The derivatives $\mathsf{V}_{ijj}$ and $\mathsf{V}_{iij}$ involve only directions tangent to this plane and hence commute with restriction:
\[
 \mathsf{V}_{ijj}|_{\Pi_{ij}}=q_{ijj},\qquad
 \mathsf{V}_{iij}|_{\Pi_{ij}}=q_{iij}.
\]
The degree $2k+2$ equation therefore becomes
\[
 0={}a_ix_i^kq_{ijj}+a_jx_j^kq_{iij} +\left(a_ia_j+\sum_\alpha d_{i\alpha}d_{j\alpha}\right) x_i^{k+1}x_j^{k+1}.
\]

First compare the coefficients of two pure powers.
Denote the coefficient of the indicated monomial in $q$ by $[x_i^ux_j^v]q$, and put
\[
 A=[x_i^{k+3}x_j^2]q,\qquad
 B=[x_i^2x_j^{k+3}]q.
\]
First compare the coefficient of $x_i^{2k+2}$. The term $a_jx_j^kq_{iij}$ contains the factor $x_j^k$ and, since $k\ge1$, cannot contribute a pure $x_i$ term; the last term also contains both $x_i,x_j$. Thus the contribution can only come from $a_ix_i^kq_{ijj}$.

On a monomial $x_i^ux_j^v$, the operator $x_i^k\partial_i\partial_j^2$ changes the exponents to $(u+k-1,v-2)$, provided that the derivative is nonzero. To obtain $(2k+2,0)$, we must have $u=k+3,v=2$. Thus the coefficient $A$ is the only possible source, and
\[
 a_ix_i^k\partial_i\partial_j^2
       (Ax_i^{k+3}x_j^2)
 =2(k+3)a_iA\,x_i^{2k+2}.
\]
The coefficient equation gives $2(k+3)a_iA=0$. Since $a_i\ne0$, we have $A=0$. Interchanging $i,j$ and comparing the pure power $x_j^{2k+2}$ gives similarly
\[
 2(k+3)a_jB=0,\qquad B=0.
\]

Now compare the coefficient of $x_i^{k+1}x_j^{k+1}$. In $a_ix_i^kq_{ijj}$, the equality $(u+k-1,v-2)=(k+1,k+1)$ requires
\[
 (u,v)=(2,k+3).
\]
Thus this contribution can only come from $Bx_i^2x_j^{k+3}$, with coefficient
\[
 2(k+2)(k+3)a_iB.
\]
Similarly, in $a_jx_j^kq_{iij}$ the only possible source is $Ax_i^{k+3}x_j^2$, contributing
\[
 2(k+2)(k+3)a_jA.
\]
The full equation for the desired coefficient is therefore
\[
 2(k+2)(k+3)(a_iB+a_jA)
       +a_ia_j+\sum_\alpha d_{i\alpha}d_{j\alpha}=0.
\]
Comparison of the two pure-power coefficients has already given $A=B=0$. Hence
\[
 a_ia_j+\sum_\alpha d_{i\alpha}d_{j\alpha}=0.
\]
Finally, dividing by the nonzero product $a_ia_j$, we obtain
\[
 \ip{\xi_i}{\xi_j}
 =\sum_\alpha\frac{d_{i\alpha}}{a_i}\frac{d_{j\alpha}}{a_j}
 =-1.
\]
This step requires only that $a_i,a_j$ be nonzero, not that they have the same sign. The lemma is proved.
\end{proof}

\begin{corollary}[A bound on the number of active axes in one layer]\label{cor:single-layer-rank-bound}
Under the hypotheses of Lemma \ref{lem:same-layer-gram}, the number $r$ of active axes in the first layer satisfies $2r\le m+1$.
\end{corollary}
\begin{proof}
The pairwise inner products of $\xi_1,\ldots,\xi_r$ form their Gram matrix. Equation \eqref{eq:gram} determines its off-diagonal entries. Add a last coordinate to each vector and put
\[
 \widehat\xi_i=(\xi_i,1)\in\R^{m-r+1}.
\]
For $i\ne j$,
\[
 \ip{\widehat\xi_i}{\widehat\xi_j}
 =\ip{\xi_i}{\xi_j}+1=0.
\]
Each vector has last coordinate $1$, and hence is nonzero. Thus $\R^{m-r+1}$ contains $r$ pairwise orthogonal nonzero vectors, and necessarily
\[
 r\le m-r+1,\qquad\text{that is,}\qquad 2r\le m+1.
\]
\end{proof}

This is a single-layer bound. Applying it to all active axes in Theorem \ref{thm:local-zero} also requires the cross-layer Gram relations of Section \ref{sec:cross-gram}.

The $-1$ in \eqref{eq:gram} comes from the positive radial contraction term $a_ia_jx_i^{k+1}x_j^{k+1}$. If the radial term were omitted from the full commutation equation, the same coefficient comparison would give only $\ip{\xi_i}{\xi_j}=0$. This would allow all $\xi_i=0$ and would not imply the preceding bound on the number of active axes.

\subsection{Nested normal ideals and the induction conditions}\label{subsec:paper-4-1}\label{sec:ideals}

Let $m\ge2$, let $U\subset S^m$ be a sufficiently small open neighbourhood of $p$, and let $C$ be a real-analytic symmetric cubic tensor on $U$ satisfying the spherical Codazzi--commutation system \eqref{eq:sphere-system} on $U$, with $C_p=0$.
Use the affine coordinates and analytic potential $f$ on a neighbourhood $W$ of $0$ from Sections \ref{subsec:paper-3-1}--\ref{subsec:paper-3-4}, with $p$ at $w=0$. All homogeneous terms, orders and germs are taken at this origin. The full commutation equation retains every contraction direction and the radial term; $C$ need not be defined on the entire sphere.

The subspaces $H_a$ below are defined by setting active variables to zero, and all assertions on them concern $W\cap H_a$ near $0$. Under $u(w)=(w,1)/\sqrt{1+|w|^2}$, this becomes a local portion of a great subsphere $\Sigma_a$ through $p$. Vanishing there means vanishing of all components of $C$, including those with normal arguments. We may shrink $U,W$ but make no assertion on the entire subspace or subsphere.

Lemma \ref{lem:homogeneous-separation} identifies the first active axes in the leading term, but the remaining variables may still occur in higher terms of $f$. Throughout what follows, $T=D^3f$ denotes the third derivative tensor of the original potential in all $m$ affine variables. Let $X$ denote all active variables already selected and $z$ the remaining variables, and put $H=\{X=0\}$ and $F(z)=f(0,z)$. After the first layer, $X$ consists of the active variables in \eqref{eq:leading}. We study $F$ to find the active directions that first appear at higher degrees.

If $F$ is a nonzero analytic germ, denote its first nonzero homogeneous term by $F_{k_{\mathrm{new}}+3}$; the leading term of its third derivative $D^3F$ in the variables $z$ then has degree $k_{\mathrm{new}}$. We seek to repeat the analysis of Lemmas \ref{lem:homogeneous-separation}--\ref{lem:same-layer-gram} on this term and select the next active axes. If $F=0$, we seek to end the selection process and prove that the whole tensor $T$ vanishes on $H$. Both steps require control of derivatives in old active directions, beyond information about $F$ alone. Below, $d$ in $f_d$ always denotes the degree of a homogeneous term of the potential.

Denote the sets of $X$ and $z$ indices by $\mathcal N$ and $\mathcal S$, respectively. They are disjoint and together comprise all $m$ coordinate indices. For $\alpha,\beta,\gamma\in\mathcal S$, purely tangential differentiation commutes with restriction, so
\[
 F_{\alpha\beta\gamma}=T_{\alpha\beta\gamma}|_H.
\]
This identity, however, contains no information about normal derivatives. Third-order components with normal indices may be nonzero even when the restricted function vanishes identically. Moreover, the internal indices of the original commutation equation range over all directions; restricting the point of evaluation to $H$ does not remove normal contractions. We therefore cannot simply assert that $F$ satisfies the full commutation equation in the remaining variables. The next proposition identifies the defect left by restriction and gives the order estimate sufficient for the next layer.

\begin{proposition}[The restricted commutation equation and its normal defect]\label{prop:restricted-commutation-defect}
For the decomposition $w=(X,z)$ above, suppose that $f$ satisfies the full commutation equation \eqref{eq:wdvv}, and put $F=f|_H$. For arbitrary $\alpha,\beta,\gamma,\delta\in\mathcal S$, its tangential equation satisfies
\[
 \sum_{\rho,\sigma\in\mathcal S} (\delta^{\rho\sigma}+z^\rho z^\sigma) (F_{\alpha\beta\rho}F_{\gamma\delta\sigma} -F_{\alpha\gamma\rho}F_{\beta\delta\sigma})=-\sum_{i\in\mathcal N} \left.(T_{\alpha\beta i}T_{\gamma\delta i} -T_{\alpha\gamma i}T_{\beta\delta i})\right|_H.
\]
If every component $T_{i\alpha\beta}|_H$ with one normal index vanishes to order at least $k_{\mathrm{new}}+2$, then the right-hand side vanishes to order at least $2k_{\mathrm{new}}+4$. Hence it does not affect coefficient comparisons in degrees $2k_{\mathrm{new}},2k_{\mathrm{new}}+1,2k_{\mathrm{new}}+2$. This conclusion does not require the right-hand side to vanish identically.
\end{proposition}
\begin{proof}
Take the four external indices of the full commutation equation \eqref{eq:wdvv} to be $\alpha,\beta,\gamma,\delta\in\mathcal S$, and then set $X=0$. The normal coordinates in the radial contraction vanish, but the normal indices in the Euclidean contraction must still be retained. Thus
\[
 0={}\sum_{\rho,\sigma\in\mathcal S} (\delta^{\rho\sigma}+z^\rho z^\sigma) (F_{\alpha\beta\rho}F_{\gamma\delta\sigma} -F_{\alpha\gamma\rho}F_{\beta\delta\sigma}) +\sum_{i\in\mathcal N} \left.(T_{\alpha\beta i}T_{\gamma\delta i} -T_{\alpha\gamma i}T_{\beta\delta i})\right|_H.
\]
The first sum is exactly the left-hand side of the desired commutation equation for $F$, whereas the second sum is the normal contraction remaining after restriction. Moving the second sum to the right gives the defect term of the restricted equation. Each $T_{\alpha\beta i}|_H$ must be obtained by differentiating first and then restricting; it cannot be replaced by a derivative of $F$.

If the leading term of $D^3F$ has degree $k_{\mathrm{new}}$, repeating the arguments of Lemmas \ref{lem:homogeneous-separation}, \ref{lem:next-homogeneous-form} and \ref{lem:same-layer-gram} for leading-term separation, the next homogeneous term and the same-layer Gram relations requires coefficient comparisons in degrees $2k_{\mathrm{new}}$, $2k_{\mathrm{new}}+1$ and $2k_{\mathrm{new}}+2$, respectively. It is not necessary first to prove that the defect vanishes identically. Provided that every component with one normal index satisfies
\[
 T_{i\alpha\beta}|_H=O(|z|^{k_{\mathrm{new}}+2}),
 \qquad i\in\mathcal N,\quad\alpha,\beta\in\mathcal S,
\]
the products in the normal contraction begin at degree $2k_{\mathrm{new}}+4$ and do not affect those comparisons. Thus the essential induction step is to raise the order bounds for the old normal components sufficiently far before studying a new leading term. This is also why it is not enough to keep track of the restricted function $F$ alone.

\end{proof}

Proposition \ref{prop:restricted-commutation-defect} specifies the estimate needed to continue the induction; it does not yet establish that estimate. Lemma \ref{lem:local-arbitrary-delay} will control mixed terms before the next leading term appears, and Lemma \ref{lem:next-mixed-normal} will treat the following degree and establish the required bound for components with one normal index. Only then can Lemma \ref{lem:inherited-layer} construct the new layer. For this purpose, we introduce ideals that record how many selected-variable factors each monomial retains. The same tool will handle termination: when the remaining potential is zero, sufficiently many such factors ensure that every third derivative retains one, so that the full tensor vanishes on the remaining subspace.

Write the active variables selected in the nonzero layers already constructed as successive nonempty blocks $X_1,\ldots,X_t$. Put $H_0=\mathbb R^m$ and, for $1\le a\le t$, define
\[
H_a=\{X_1=\cdots=X_a=0\}.
\]
Denote by $z$ the coordinates not yet selected as active variables. At layer $a$, the old normal variables are $X_1,\ldots,X_{a-1}$, and the full tuple of remaining coordinates on $H_{a-1}$ is $w^{(a)}=(X_a,\ldots,X_t,z)$. The notation $T|_{H_{a-1}}$ means that third partial derivatives are taken in all variables before the point of evaluation is restricted to $H_{a-1}$; each of the three arguments may still be any ambient coordinate direction. For this full component array, define
\[
 k_a:=\min_{1\le i,j,\ell\le m}
 \operatorname{ord}_0\!\left(
 \left.\partial_i\partial_j\partial_\ell f\right|_{H_{a-1}}
 \right).
\]
Here $\operatorname{ord}_0$ is the lowest nonzero total Taylor degree in all variables $w^{(a)}$, with order $+\infty$ assigned to the zero germ. The full restricted array at each constructed layer is nonzero, so $k_a$ is finite. It is the degree of the first nonzero homogeneous part of that array; the analytic tensor need not itself be homogeneous, and its covariant order is always $3$. The constructed layers satisfy
\[
1\le k_1<k_2<\cdots<k_t.
\]

The third derivatives of the restricted potential $f|_{H_{a-1}}$ in the variables $w^{(a)}$ record only those components whose three indices are tangent to $H_{a-1}$. Induction condition (1) below will require the degree-$k_a$ part of the full array to be nonzero and to have the stated purely tangential form. Under this condition, the leading term of the restricted potential has degree $k_a+3$. This correspondence must be established at each layer; it does not follow from restriction alone.

\begin{definition}[Nested normal ideals]\label{def:nested-normal-ideals}
For these active blocks and leading degrees $k_a$ of the full restricted arrays, define, in the polynomial ring $\mathbb R[X_1,\ldots,X_t,z]$,
\begin{equation}\label{eq:ideals}
 I_a=(X_1,\ldots,X_a),\quad
 \J_t=\sum_{a=1}^t I_a^{k_a+3},\quad
 H_a=\{X_1=\cdots=X_a=0\},
\end{equation}
where $(X_1,\ldots,X_a)$ is the ideal generated by all coordinate variables in the first $a$ blocks. A polynomial belongs to $I_a^q$ if and only if each of its nonzero monomials has total degree at least $q$ in these variables. This does not require degree at least $q$ in each variable, nor does it mean that only the $q$th powers of the individual variables are used as generators. For example, $x^2y^2\in(x,y)^4$, but $x^2y^2\notin(x^4,y^4)$. Total degree is used because three differentiations may remove degrees from different normal variables.

\end{definition}

\begin{lemma}[Ideal membership and vanishing of all third derivatives]\label{lem:ideal-full-third-vanishing}
Use the notation of \eqref{eq:ideals}. If a polynomial $P$ belongs to $\J_t$, then each of its third partial derivatives vanishes on $H_t$, including components involving normal differentiation. For a monomial, denote its cumulative degree in the variables of the first $a$ blocks by $P_a$. It lies outside $\J_t$ if and only if $P_a\le k_a+2$ for every $a$.
\end{lemma}
\begin{proof}
The degree $k_a+3$ also has a direct consequence. Applying any three partial derivatives to a monomial in $I_a^{k_a+3}$ decreases its total degree in the first $a$ blocks by at most $3$. If the result is nonzero, that degree is therefore still at least $k_a\ge1$. Thus all third derivatives vanish on $H_a$, and hence on the smaller subspace $H_t$. Applying this term by term to a polynomial in $\J_t$ gives
\[
 P\in\J_t\quad\Longrightarrow\quad
 (D^3P)|_{H_t}=0.
\]
This applies also to normal derivatives. Section \ref{subsec:paper-4-2} will establish the ideal control first for homogeneous terms and then for convergent analytic series.

If a monomial has cumulative degree $P_a$ in the first $a$ blocks, then it lies outside $I_a^{k_a+3}$ if and only if $P_a\le k_a+2$. A sum of monomial ideals is generated by the union of the monomials in its summands. Hence
\begin{equation}\label{eq:outside}
 \text{the monomial lies outside }\J_t\quad\Longleftrightarrow\quad
 P_a\le k_a+2\quad\text{for every }a.
\end{equation}
The sum is essential here: a monomial belongs to $\J_t$ as soon as it reaches the cumulative degree threshold of any one layer; to lie outside $\J_t$, it must lie below the threshold of every layer simultaneously. The degree comparison in Section \ref{subsec:paper-4-2} uses all these inequalities, not only $P_t\le k_t+2$.

\end{proof}

For example, take two layers $X_1=(x)$ and $X_2=(y)$, and let $k_1=1,k_2=3$. Then
\[
\J_2=(x)^4+(x,y)^6.
\]
The monomial $x^4$ belongs to the first ideal although its total degree is only $4$; $y^6$ belongs to the second ideal; and $x^3y^2$ has cumulative degrees $3$ and $5$, below the respective thresholds $4$ and $6$, so it lies outside $\J_2$. In contrast, $x^3y^3$ has total degree $6$ and belongs to $\J_2$. Below we eliminate, one by one, monomials outside $\J_t$ that contain selected variables. Terms in the remaining variables alone must be retained: they are precisely the Taylor terms of $F=f|_{H_t}$.

With the conventions above for $T|_{H_{a-1}}$, $w^{(a)}$ and $k_a$, we denote the conjunction of the following five conditions at stage $t$ by $(\mathrm{IH}_t)$, with the assertions about layer $a$ required for every $1\le a\le t$. These conditions are verified initially and restored after each new block; they are not assumed for an arbitrary division into blocks.
\begin{enumerate}[leftmargin=2em,itemsep=3pt]
\item\begin{samepage} The form of the leading term of the full restriction. The homogeneous part of degree $k_a$ in the full array $T|_{H_{a-1}}$ is
\[
 \sum_{i\in X_a}a_i x_i^{k_a}(\mathrm dx_i)^{\otimes3},\qquad a_i\ne0.
\]
Here $i\in X_a$ means that $i$ ranges over the coordinate indices of block $a$. All components are obtained by taking third derivatives in the full coordinate space and then restricting the point of evaluation, including indices normal to $H_{a-1}$. This is stronger than knowing only the leading term of the tangential third derivatives of the restricted function.
\end{samepage}
\item The order of components with one normal index. If $i$ is an old normal index in the first $a-1$ blocks and $\alpha,\beta$ are tangential indices of $H_{a-1}$, then
\[
 T_{i\alpha\beta}|_{H_{a-1}}
 =O(|w^{(a)}|^{k_a+2}).
\]
The order here is computed in all variables $w^{(a)}$ on $H_{a-1}$, not merely in the final variables $z$. By Proposition \ref{prop:restricted-commutation-defect}, the defect of the restricted commutation equation begins no earlier than degree $2k_a+4$. This suffices for the comparisons in the three degrees required by Lemmas \ref{lem:homogeneous-separation}--\ref{lem:same-layer-gram}.
\item The form of the next homogeneous term. The homogeneous term of degree $k_a+4$ in the restricted potential $f|_{H_{a-1}}$ has the form \eqref{eq:next}, with active variables $X_a$ and remaining variables $(X_{a+1},\ldots,X_t,z)$. Its mixed coefficients define the vectors $\xi_i=d_i/a_i$ for the active axes of this layer.
\item Ideal membership of the lower-degree potential terms. For the current last layer, $f_d\in\J_t$ for $d\le k_t+3$. This controls degrees at the level of the potential; Section \ref{subsec:paper-4-2} extends it through every degree before the next leading term in the remaining variables alone appears.
\item Gram relations in the common remaining space. A vector $\xi_i$ from layer $a$ initially lies in the space of all subsequent coordinate directions. Its components in the already selected later blocks $X_{a+1},\ldots,X_t$ are zero, so all selected vectors can be regarded as vectors in the same final $z$ space. In this space, vectors associated with distinct active axes have inner product $-1$. After a new block is selected, Section \ref{subsec:paper-4-4} proves that the old vectors have zero components in it, and Section \ref{subsec:paper-4-5} proves the inner-product relations between old and new vectors.
\end{enumerate}

The initial case $t=1$ was established in Sections \ref{subsec:paper-3-1}--\ref{subsec:paper-3-4}: the leading and next terms are given by Lemmas \ref{lem:homogeneous-separation} and \ref{lem:next-homogeneous-form}, respectively; the leading term of the potential belongs to $I_1^{k_1+3}$; and the relations between vectors of the same layer follow from Lemma \ref{lem:same-layer-gram}. There are no old normal indices at this stage, so condition (2) is vacuous; in condition (5) there are as yet no components in subsequent blocks to check.

In the present local analytic setting, start with $(\mathrm{IH}_t)$ and a nonzero remaining germ $F=f|_{H_t}$. Analyticity first supplies a leading term $F_{k_{\mathrm{new}}+3}$, and condition (4) gives $k_{\mathrm{new}}>k_t$, with no bound on the delay. At this point $k_{\mathrm{new}}$ is determined by the restricted potential; the old normal components still require control before it can be identified as the leading degree of the full array $T|_{H_t}$. Section \ref{subsec:paper-4-2} extends ideal control across that delay. Section \ref{subsec:paper-4-3} identifies the possible degree-$k_{\mathrm{new}}+4$ boundary correction and bounds the components with one normal index, putting the defect beyond degree $2k_{\mathrm{new}}+2$. Section \ref{subsec:paper-4-4} can then construct the new layer and eliminate the correction by showing that old vectors have zero components along the new axes; it also confirms that the full restricted array has leading degree $k_{\mathrm{new}}$. Section \ref{subsec:paper-4-5} proves the cross-layer Gram relations. Setting $k_{t+1}=k_{\mathrm{new}}$ gives $(\mathrm{IH}_{t+1})$. The proof of Theorem \ref{thm:local-zero} checks the restoration of all five conditions and the preservation of earlier layers, and uses the dimension bound to prove finite termination.

\subsection{Ideal control for an arbitrary delay}\label{subsec:paper-4-2}\label{sec:nested-lemma}
After the first $t$ layers have been selected, we prove ideal control of all lower-degree terms involving selected variables before the next leading term of $F=f|_{H_t}$. The argument uses the full commutation equation, differentiating before restriction; it does not assume that the tangential restriction satisfies that equation without a defect.

\begin{lemma}[Local ideal control for an arbitrary delay]\label{lem:local-arbitrary-delay}
Use the local coordinates and block notation of Section \ref{sec:ideals}. Suppose that the real-analytic potential $f$ satisfies the full commutation equation \eqref{eq:wdvv} near $0$ and conditions (1) and (4) of the induction in that subsection. Expand $f$ and its restriction as
\[
 f=\sum_{d\ge4}f_d,\qquad
 F(z):=f(0,\ldots,0,z)=\sum_{d\ge4}F_d(z).
\]
Here $F_d$ is also regarded as a polynomial in all variables, independent of the selected variables.
\begin{enumerate}[label=\textup{(\arabic*)},leftmargin=2em,itemsep=3pt]
\item If the first nonzero homogeneous term of $F$ is $F_{k_{\mathrm{new}}+3}$, with $k_{\mathrm{new}}>k_t$, then
\begin{equation}\label{eq:nested}
 f_d-F_d\in\J_t\qquad(d\le k_{\mathrm{new}}+3).
\end{equation}
\item If $F=0$ as an analytic function germ at $0$, then $f\in\J_t$ in the ring of convergent analytic function germs. Moreover, there is a neighbourhood $W'$ of $0$ such that all third partial derivatives of $f$ vanish on $W'\cap H_t$. Consequently, all components of the original spherical tensor $C$ vanish on the corresponding portion of the great subsphere near $p$.
\end{enumerate}
\end{lemma}

\begin{proof}
We first prove the degree-by-degree assertion for a finite delay, and then pass from convergent series to analytic germs when $F=0$. Two notions of degree must be distinguished: total degree selects the part of the commutation equation containing the current unknown term; cumulative degree in the first several blocks rules out cancellation by older terms.

Condition (4) in Section \ref{sec:ideals} already gives
\[
 f_e\in\J_t\qquad(e\le k_t+3).
\]
At these degrees $F_e=0$, so \eqref{eq:nested} holds. Now fix
\[
 k_t+4\le d\le k_{\mathrm{new}}+3,
\]
and suppose that the conclusion holds for every $e<d$. Since $F$ has no nonzero term below degree $k_{\mathrm{new}}+3$, we have $F_e=0$ for every $e<d$, and hence
\[
 f_e\in\J_t\qquad(e<d).
\]
We call the contributions to the commutation equation from pairwise products of these lower-degree potential terms the old source terms. The homogeneous term to be controlled at this stage is $f_d$. At the endpoint $d=k_{\mathrm{new}}+3$, $F_d$ may be nonzero, but it does not enter the old source terms formed from degrees below $d$.

Fix a block $X_a$, take $i\in X_a$, and choose indices $\alpha,\beta$ in the final remaining directions. Use the following component of \eqref{eq:Q}:
\[
 Q_{i\alpha\beta}
 =\sum_{r,s=1}^m(\delta^{rs}+w^rw^s)
   (T_{iir}T_{\alpha\beta s}-T_{i\alpha r}T_{i\beta s})=0.
\]
Here $r,s$ range over all coordinate directions, including those normal to $H_{a-1}$. After differentiation, restrict the point of evaluation to $H_{a-1}$.

Condition (1) of the induction says that the full restricted array has lowest homogeneous degree $k_a$, with leading term
\[
 L=\sum_{j\in X_a}a_jx_j^{k_a}(\dd x_j)^{\otimes3}.
\]
Thus $L_{iir}$ equals $a_ix_i^{k_a}$ only when $r=i$. Since $\alpha,\beta$ are final remaining indices,
\[
 L_{\alpha\beta s}=L_{i\alpha r}=L_{i\beta s}=0
 \qquad\text{for all }r,s.
\]
On the other hand, $D^3f_d$ has degree $d-3$. Its first product with the leading term occurs at total degree $k_a+d-3$. In the Euclidean contraction, only the following term remains:
\begin{equation}\label{eq:principal}
 a_i x_i^{k_a}\partial_{i\alpha\beta}(f_d|_{H_{a-1}}).
\end{equation}
All three derivative directions here are tangent to $H_{a-1}$, so differentiation commutes with this restriction.

No other contribution containing $f_d$ enters this degree. Pairing it with a tensor component of degree greater than $k_a$ gives a higher total degree; the product of two copies of $D^3f_d$ also has higher degree, since $d-3>k_a$. A radial contraction introduces two further coordinate factors. Potential terms of degree greater than $d$ have third derivatives of degree greater than $d-3$, and even pairing them with the lowest degree $k_a$ part of the full restricted array exceeds the degree being compared. Thus the equation at this degree is
\[
 0=a_i x_i^{k_a}\partial_{i\alpha\beta}(f_d|_{H_{a-1}})
       +S_{a,i,\alpha,\beta},
\]
where $S_{a,i,\alpha,\beta}$ is formed entirely from potential terms of degree less than $d$.

We next establish a lower bound for the cumulative degrees of the old source terms.
Every lower-degree potential term $f_e$ belongs to the monomial ideal $\J_t$. Hence each of its nonzero monomials belongs to at least one $I_b^{k_b+3}$. If it belongs to several such ideals, choose one and assign the monomial to it.

If the assigned block satisfies $b<a$, the monomial has cumulative degree at least $k_b+3\ge4$ in the first $b$ blocks. After any three differentiations, a nonzero result still has degree at least one in these variables and therefore vanishes on restriction to $H_{a-1}$. Thus, for a pair of potential monomials surviving in an old source term, their assigned blocks may be ordered as
\[
 a\le b\le c.
\]
They belong to $I_b^{k_b+3}$ and $I_c^{k_c+3}$, respectively. Since the first $b$ blocks are contained in the first $c$ blocks, the undifferentiated product has cumulative degree in the first $c$ blocks at least
\[
 (k_b+3)+(k_c+3)=k_b+k_c+6.
\]

The two products in an ordinary contraction have the forms
\[
 T_{iir}T_{\alpha\beta r},
 \qquad T_{i\alpha r}T_{i\beta r}.
\]
Taken together, the two factors involve two derivatives in direction $i$, one in direction $\alpha$, one in direction $\beta$, and two in the internal contraction directions. Since $i\in X_a$ and $a\le c$, the two $i$ derivatives decrease the cumulative degree by two in total. The final remaining directions $\alpha,\beta$ do not decrease the degree in the first $c$ blocks, and the two internal derivatives decrease it by at most two more. Thus every nonzero ordinary old source term has cumulative degree in the first $c$ blocks at least
\begin{equation}\label{eq:source-bound}
 k_b+k_c+6-2-2=k_b+k_c+2.
\end{equation}
This estimate already includes internal indices in old normal directions.

For the radial contraction $w^rw^s$, if an internal derivative is in a variable of the first $c$ blocks, its coordinate multiplier restores the one degree removed; if it is in any other variable, neither operation changes the degree in the first $c$ blocks. Thus only the two external $i$ derivatives cause a net decrease, giving the lower bound
\[
 k_b+k_c+6-2=k_b+k_c+4.
\]
If a coordinate multiplier lies in a direction already set equal to zero, the term vanishes outright. Hence the weaker lower bound in \eqref{eq:source-bound} applies to every old source term.

We now select a monomial to be excluded and take a nonzero derivative of it.
Choose a monomial $\lambda M$ in $f_d$ that contains selected variables but does not belong to $\J_t$. Let $X_a$ be the first block of selected variables occurring in $M$, and choose a variable $x_i$ in that block with positive exponent. Since $M$ contains no variables from the first $a-1$ blocks, it survives restriction to $H_{a-1}$.

Let $P_c$ be the cumulative degree of $M$ in the first $c$ blocks. By \eqref{eq:outside},
\[
 P_c\le k_c+2\qquad(1\le c\le t).
\]
Since $M$ has total degree $d$, its degree in the final remaining variables is
\[
 d-P_t\ge d-(k_t+2)\ge2.
\]
Thus we can choose $\alpha,\beta$ so that $\partial_{\alpha\beta}M\ne0$: if some remaining variable has exponent at least $2$, differentiate twice in that direction; otherwise choose two distinct remaining directions with positive exponents. Using also the positive exponent of $x_i$, we obtain
\[
 \partial_{i\alpha\beta}M\ne0.
\]
If there are no remaining variables, the preceding degree inequality already rules out such an $M$, so no further exclusion is needed.

We next rule out cancellation by old source terms and by other unknown monomials.
For any $c\ge a$, the principal operation changes the cumulative degree of $M$ in the first $c$ blocks to
\[
 \deg_{I_c}\bigl(x_i^{k_a}\partial_{i\alpha\beta}M\bigr)
 =k_a+P_c-1.
\]
Indeed, $\partial_i$ decreases it by one, the two remaining derivatives leave it unchanged, and multiplication by $x_i^{k_a}$ increases it by $k_a$.

Now examine the old source terms one by one. Suppose a surviving term comes from assigned blocks $a\le b\le c$. The target monomial has cumulative degree satisfying
\[
 k_a+P_c-1
 \le k_a+k_c+1
 < k_b+k_c+2,
\]
where the strict inequality uses $k_b\ge k_a$. The cumulative degree estimate above, however, shows that the degree of this old source term in the first $c$ blocks is at least $k_b+k_c+2$. Thus they cannot be the same monomial, and this old source term cannot cancel the target. Applying this argument with the corresponding $c$ for each old source term rules out all of them. This step requires \eqref{eq:outside} for every $c$.

The degree comparisons are summarized in Table \ref{tab:delay-source-bounds}.
The total degree being compared is $k_a+d-3$; let $e_1,e_2<d$ be the degrees of the two old potential terms.
For each surviving old source, order its assigned blocks as $a\le b\le c$ and use that source's own prefix $I_c$; the same prefix is used for the target in the first row.
Sources with an assigned block below $a$ vanish on $H_{a-1}$.
\begin{center}
\small
\refstepcounter{table}\label{tab:delay-source-bounds}
\textsc{Table~\thetable}. Source bounds in the arbitrary-delay comparison.
\par\smallskip
\renewcommand{\arraystretch}{1.18}
\begin{tabular}{@{}p{0.40\linewidth}p{0.27\linewidth}p{0.27\linewidth}@{}}
\toprule
Contribution & Total degree & Degree in the prefix $I_c$ \\
\midrule
Target $x_i^{k_a}\partial_{i\alpha\beta}M$ & $k_a+d-3$ & $\le k_a+k_c+1$ \\
Old ordinary contraction & $e_1+e_2-6$ & $\ge k_b+k_c+2$ \\
Old radial contraction & $e_1+e_2-4$ & $\ge k_b+k_c+4$ \\
Other nonzero terms involving $f_d$, or terms of higher potential degree & $>k_a+d-3$ & No contribution at the degree compared \\
\bottomrule
\end{tabular}
\end{center}
Only old sources whose total degree equals $k_a+d-3$ can enter the current equation; the last column excludes their cancellation of the target. Distinct current monomials are dealt with separately by the injective exponent map below.

It remains to rule out cancellation between different monomials of $f_d$. Number all variables together, write $M=w^\nu$, and let $\mu$ be the multi-index of the fixed derivative $\partial_{i\alpha\beta}$. For each monomial not annihilated by this derivative,
\[
 x_i^{k_a}\partial^\mu w^\nu
 =\left(\prod_{\ell=1}^m
       \frac{\nu_\ell!}{(\nu_\ell-\mu_\ell)!}\right)
       w^{\nu-\mu+k_ae_i},
\]
where $e_i$ is the unit vector in the $i$th exponent direction. The exponent map
$\nu\mapsto\nu-\mu+k_ae_i$ is injective, and the integer coefficient in front is nonzero. Hence no other unknown monomial can produce the same output term. Terms containing variables of the first $a-1$ blocks vanish after restriction and do not affect this conclusion.

The coefficient of the target monomial in the equation is therefore $a_i\lambda$ times a nonzero integer. Since $a_i\ne0$, we obtain $\lambda=0$. Every monomial outside the ideal that contains selected variables has such a first block $X_a$, so all are excluded. The only monomials outside the ideal that remain are those in the remaining variables alone, and they form precisely $F_d$. Thus $f_d-F_d\in\J_t$, completing the induction and proving assertion \textup{(1)}.

If $F=0$ as an analytic germ, then every $F_d=0$, and the preceding induction applies to every total degree. Thus every nonzero monomial in the convergent Taylor series of $f$ belongs to $\J_t$.

Since the ideal has finitely many monomial generators, the following direct convergence estimate turns coefficientwise membership into analytic germ membership. Choose finitely many monomial generators $G_j=w^{\eta^{(j)}}$, $1\le j\le N$, of $\J_t$. Write
\[
 f(w)=\sum_\nu c_\nu w^\nu.
\]
Assign each nonzero monomial to a generator dividing it, obtaining disjoint index sets $E_1,\ldots,E_N$. Define
\[
 A_j(w)=\sum_{\nu\in E_j}c_\nu w^{\nu-\eta^{(j)}}.
\]
Since $f$ is real-analytic, we may choose $r>0$ sufficiently small that
$\sum_\nu|c_\nu|r^{|\nu|}<\infty$. Then
\[
 \sum_{\nu\in E_j}|c_\nu|r^{|\nu|-|\eta^{(j)}|}
 \le r^{-|\eta^{(j)}|}
       \sum_\nu|c_\nu|r^{|\nu|}<\infty.
\]
Each $A_j$ therefore defines a convergent analytic function in a smaller neighbourhood, and
\[
 f=\sum_{j=1}^NG_jA_j.
\]
This proves $f\in\J_t$ in the ring of convergent analytic function germs. We have used the convergent coefficients of the given analytic function; we have not constructed an infinite sequence of coordinate changes or a different formal solution.

Every generator has total degree at least $4$ in the selected variables. In any third partial derivative of $G_jA_j$, each term in the Leibniz formula differentiates the generator at most three times, so at least one selected-variable factor remains. After shrinking the neighbourhood, setting these variables equal to zero on $H_t$ gives
\[
 (D^3f)|_{W'\cap H_t}=0.
\]
Any of the three derivative directions may be normal. By the recovery formula in Lemma \ref{lem:affine-slice-potential}, the entire cubic tensor $C$ also vanishes at the corresponding spherical points. This conclusion concerns only the portion near $p$. Assertion \textup{(2)} follows.
\end{proof}

\subsection{Mixed terms at the next order and estimates for components with one normal index}\label{subsec:paper-4-3}\label{sec:boundary}
Section \ref{sec:nested-lemma} gives $f_d-F_d\in\J_t$ for $d\le k_{\mathrm{new}}+3$ when $F_{k_{\mathrm{new}}+3}$ is the leading term of $F=f|_{H_t}$. At the next degree, its pairing with an earlier layer's next term introduces new sources. We determine their mixed terms, exclude other mixed terms outside the ideal, and estimate the components with one old normal index.

We continue to work only near $p$. Let $\mathcal N$ denote the indices of the old active variables in the first $t$ blocks, and let $\alpha,\beta,\gamma$ denote indices of the final remaining variables $z$. Write
\[
 G(z):=F_{k_{\mathrm{new}}+3}(z),\qquad
 L_i(z):=\ip{\xi_i}{\operatorname{grad}_zG(z)}\quad(i\in\mathcal N).
\]
Here $\operatorname{grad}_zG$ denotes the Euclidean gradient in the remaining variables $z$. Inductive condition (5) ensures that every old vector $\xi_i$ already lies in the final $z$-space. Thus $L_i$ depends only on $z$ and is homogeneous of degree $k_{\mathrm{new}}+2$. All restrictions and order estimates below are taken near the origin of $H_t$; $f$ need not be defined on the entire coordinate subspace.

\Needspace{4\baselineskip}
\begin{lemma}[Mixed terms at the next order and estimates for components with one normal index]\label{lem:next-mixed-normal}
Use the local notation of Section \ref{sec:ideals}. Suppose that its inductive conditions hold and that the first nonzero homogeneous term of the remaining potential $F$ is $G=F_{k_{\mathrm{new}}+3}$, where $k_{\mathrm{new}}>k_t$. Then the following assertions hold.
\begin{enumerate}[label=\textup{(\arabic*)},leftmargin=2em,itemsep=3pt]
\item There exists $P_{k_{\mathrm{new}}+4}\in\J_t$ such that
\begin{equation}\label{eq:boundary}
 f_{k_{\mathrm{new}}+4}=F_{k_{\mathrm{new}}+4}
 -\frac12\sum_{i\in\mathcal N}x_i^2\ip{\xi_i}{\operatorname{grad}_zG}
 +P_{k_{\mathrm{new}}+4}.
\end{equation}
Thus, after subtracting the term in the remaining variables alone, the mixed part is uniquely determined modulo $\J_t$.
\item For every old normal index $i\in\mathcal N$ and remaining indices $\alpha,\beta$, we have
\begin{equation}\label{eq:old-normal}
 \left.f_{x_i}\right|_{H_t}=O(|z|^{k_{\mathrm{new}}+4}),\qquad
 \left.T_{i\alpha\beta}\right|_{H_t}=O(|z|^{k_{\mathrm{new}}+2}).
\end{equation}
Here $T=D^3f$, and both estimates are local estimates as $z\to0$.
\end{enumerate}
\end{lemma}

The term quadratic in the old variables in \eqref{eq:boundary} will also be called the \emph{boundary correction} of this induction step: it occurs at the first degree beyond the range controlled in the preceding section. This name does not refer to the boundary of a spatial domain. We do not yet prove that the correction vanishes; Section \ref{sec:inherited} will eliminate it by using the leading term of the new layer.

\begin{proof}

Fix an old active index $i\in X_a$, and write $k=k_a$. In the full commutativity equation, use
\[
 Q_{i\alpha\beta}
 =\sum_{r,s=1}^m(\delta^{rs}+w^rw^s)
 \bigl(T_{iir}T_{\alpha\beta s}-T_{i\alpha r}T_{i\beta s}\bigr)=0.
\]
As before, we take all third derivatives before restricting the base point to $H_{a-1}$. Inductive condition (1) says that the leading term of the full restricted array is
\[
 \sum_{j\in X_a}a_jx_j^k(\dd x_j)^{\otimes3}.
\]
The third derivatives of the unknown potential term $f_{k_{\mathrm{new}}+4}$ have degree $k_{\mathrm{new}}+1$. Hence, in the part of $Q_{i\alpha\beta}|_{H_{a-1}}$ of degree $k+k_{\mathrm{new}}+1$, the only contribution involving this unknown term is
\[
 a_ix_i^k\partial_{i\alpha\beta}
       \bigl(f_{k_{\mathrm{new}}+4}|_{H_{a-1}}\bigr).
\]
This is the comparison in Section \ref{sec:nested-lemma} that determines the contribution of the current unknown term, now with $d=k_{\mathrm{new}}+4$. All three differentiation directions are tangent to $H_{a-1}$, so differentiation commutes with restriction here. A radial contraction involving the current term has two additional coordinate factors. Its pairings with higher old terms, its self-product, and contributions from higher potential terms also have degree greater than $k+k_{\mathrm{new}}+1$.

Lemma \ref{lem:local-arbitrary-delay} shows that all potential terms of degree less than $k_{\mathrm{new}}+4$ belong to $\J_t$, with the sole possible exception of the term $G=F_{k_{\mathrm{new}}+3}$ in the remaining variables alone. Thus every new source term must contain a third derivative of $G$, of degree $k_{\mathrm{new}}$.

If the other tensor factor has degree $j$, an ordinary contraction must satisfy
\[
 j+k_{\mathrm{new}}=k+k_{\mathrm{new}}+1.
\]
Hence $j=k+1$. Inductive condition (3) determines the next homogeneous potential term of the old layer. Its part involving the $i$-axis and the final remaining variables is
\[
 \frac{x_i^{k+3}}{(k+2)(k+3)}
       \sum_\gamma d_{i\gamma}z_\gamma.
\]
Differentiating twice in $x_i$ and once in $z_\gamma$, we obtain
\[
 \partial_{ii\gamma}
 \left(\frac{x_i^{k+3}}{(k+2)(k+3)}
       \sum_\delta d_{i\delta}z_\delta\right)
 =d_{i\gamma}x_i^{k+1}.
\]
The corresponding component of a pure power of $x_i$ can have only the internal index $i$, and therefore cannot pair with a derivative of $G$. The other old active axes provide no component with two $i$ indices. Thus the contribution to $T_{iir}T_{\alpha\beta r}$ is precisely
\[
 x_i^{k+1}\sum_\gamma d_{i\gamma}G_{\alpha\beta\gamma}.
\]
Since $G$ does not contain $x_i$, it cannot occur in $T_{iir}$ or $T_{i\alpha r}$. The other product, $T_{i\alpha r}T_{i\beta r}$, therefore has no contribution involving $G$.

We must also check the radial contractions and the pairings at the lowest degree. The ordinary contraction of the old leading term with $G$ vanishes: in the leading term, $T_{iir}$ is nonzero only for $r=i$, whereas $G_{\alpha\beta i}=0$. Their radial contraction has degree $k+k_{\mathrm{new}}+2$, one greater than the target degree. More generally, a radial pairing requires $j+k_{\mathrm{new}}+2=k+k_{\mathrm{new}}+1$, or $j=k-1$, and no term of such low degree occurs in the full restricted array. Two factors involving only the remaining variables cannot contribute either, since at least one factor in each product is differentiated in the $i$-direction. In particular, even when $k_{\mathrm{new}}=k+1$ and the degrees of two $G$ factors sum to the target degree, their product is still zero.

Thus the equation at the current degree can be written as
\[
 0={}a_ix_i^k\partial_{i\alpha\beta} \bigl(f_{k_{\mathrm{new}}+4}|_{H_{a-1}}\bigr) +x_i^{k+1}\sum_\gamma d_{i\gamma}G_{\alpha\beta\gamma} +S_{a,i,\alpha,\beta}.
\]
Here $S_{a,i,\alpha,\beta}$ arises entirely from two lower-degree terms in the old ideal, including both their ordinary and their radial contractions.

We cancel the new source by the following correction.
Define
\[
 \beta_{\mathrm{bd}}=-\frac12\sum_{j\in\mathcal N}x_j^2L_j(z).
\]
Since $L_j$ depends only on $z$, differentiation by $\partial_i$ leaves only the term with $j=i$. Therefore
\[
 \partial_i\beta_{\mathrm{bd}}=-x_iL_i,\qquad
 \partial_{i\alpha\beta}\beta_{\mathrm{bd}}
 =-x_i\sum_\gamma\xi_{i\gamma}G_{\alpha\beta\gamma}.
\]
Using $a_i\xi_{i\gamma}=d_{i\gamma}$, we obtain
\[
 a_ix_i^k\partial_{i\alpha\beta}\beta_{\mathrm{bd}}
 =-x_i^{k+1}\sum_\gamma d_{i\gamma}G_{\alpha\beta\gamma},
\]
which cancels exactly the new source term above. This also explains the coefficient $-\tfrac12$. Since $\operatorname{grad}_zG$ has degree $k_{\mathrm{new}}+2$, $\beta_{\mathrm{bd}}$ indeed has degree $k_{\mathrm{new}}+4$.

We now prove uniqueness of this correction modulo the old ideal.
Put
\[
 R=f_{k_{\mathrm{new}}+4}-F_{k_{\mathrm{new}}+4}-\beta_{\mathrm{bd}}.
\]
Since $\partial_iF_{k_{\mathrm{new}}+4}=0$, subtracting the correction just constructed gives
\[
 0=a_ix_i^k\partial_{i\alpha\beta}
       \bigl(R|_{H_{a-1}}\bigr)+S_{a,i,\alpha,\beta}.
\]
This only rewrites the linear comparison at the current degree. We do not require $R$ itself to satisfy the full nonlinear commutativity equation. The remaining exclusion uses the same linear mechanism as in Lemma \ref{lem:local-arbitrary-delay}: the cumulative-degree bounds for the old sources, followed by injectivity of the exponent map. We give the coefficient comparison here for $R$, after the new source involving $G$ has been removed.

Table \ref{tab:boundary-source-bounds} records the separate comparison at total degree $k+k_{\mathrm{new}}+1$, where $k=k_a$. Let $e_1,e_2<k_{\mathrm{new}}+4$ be the degrees of old ideal terms; in the rows involving $G$, the other tensor factor has degree $j$. For each old source, use its assigned blocks $a\le b\le c$ and its own prefix $I_c$. The last row concerns a candidate monomial $M$ of $R$ outside $\J_t$, with prefix degree $P_c$.
\begin{center}
\small
\refstepcounter{table}\label{tab:boundary-source-bounds}
\textsc{Table~\thetable}. Sources at the boundary degree and after its correction.
\par\smallskip
\renewcommand{\arraystretch}{1.18}
\begin{tabular}{@{}p{0.36\linewidth}p{0.28\linewidth}p{0.30\linewidth}@{}}
\toprule
Contribution & Total degree & Outcome or prefix bound \\
\midrule
Ordinary pairing of $D^3G$ with an old ideal term & $j+k_{\mathrm{new}}$ & Requires $j=k+1$: the new source canceled by $\beta_{\mathrm{bd}}$ \\
Radial pairing of $D^3G$ with an old ideal term & $j+k_{\mathrm{new}}+2$ & Requires $j=k-1$: absent by the leading-degree bound \\
Two factors involving only remaining variables & Any degree & Zero: each product has an $i$ derivative \\
Old ordinary contraction & $e_1+e_2-6$ & Prefix degree $\ge k_b+k_c+2$ \\
Old radial contraction & $e_1+e_2-4$ & Prefix degree $\ge k_b+k_c+4$ \\
Target $x_i^k\partial_{i\alpha\beta}M$ from $R$ & $k+k_{\mathrm{new}}+1$ & Prefix degree $\le k+k_c+1$ \\
\bottomrule
\end{tabular}
\end{center}
Only sources of total degree $k+k_{\mathrm{new}}+1$ enter this comparison. All other nonzero terms involving the current unknown $f_{k_{\mathrm{new}}+4}$, and all terms of higher potential degree, lie above that degree. Thus the prefix comparison for $R$ applies only after the displayed source involving $G$ has been canceled; it does not assert that $R$ satisfies the full nonlinear equation.

Suppose that $R$ contains a monomial $\lambda M$ involving an old active variable but not belonging to $\J_t$. Take the first block $X_a$ occurring in it, and then a variable $x_i$ in that block with positive exponent. Let $P_c$ be the cumulative degree of $M$ in the first $c$ blocks. By \eqref{eq:outside},
\[
 P_c\le k_c+2\quad(1\le c\le t).
\]
Since the total degree of $M$ is $k_{\mathrm{new}}+4$, its degree in the remaining variables satisfies
\[
 k_{\mathrm{new}}+4-P_t\ge k_{\mathrm{new}}-k_t+2\ge3.
\]
We can therefore choose remaining indices $\alpha,\beta$ for which $\partial_{i\alpha\beta}M\ne0$. If there are no remaining variables, the preceding inequality already rules out such an $M$.

For every $c\ge a$, the cumulative degree of the target monomial in the first $c$ blocks is
\[
 \deg_{I_c}\bigl(x_i^k\partial_{i\alpha\beta}M\bigr)
 =k+P_c-1\le k+k_c+1.
\]
On the other hand, the cumulative-degree estimates for the old source terms in Section \ref{sec:nested-lemma} still apply. If a surviving old source term comes from assigned blocks $a\le b\le c$, the cumulative degree of its ordinary contraction in the first $c$ blocks is at least $k_b+k_c+2$, and that of its radial contraction is at least $k_b+k_c+4$. Since $k_b\ge k$,
\[
 k+P_c-1\le k+k_c+1<k_b+k_c+2.
\]
Thus this old source term cannot produce the target monomial. Applying this argument to the corresponding $c$ for each old source term excludes cancellation by any old source term.

Finally, for fixed $i,\alpha,\beta$, the operation $x_i^k\partial_{i\alpha\beta}$ sends the exponent $\nu$ of every monomial not annihilated by differentiation to
\[
 \nu-e_i-e_\alpha-e_\beta+ke_i.
\]
This map is injective, and the integer coefficient produced by differentiation is nonzero. Hence distinct unknown monomials cannot cancel one another either. Terms containing variables from the first $a-1$ blocks vanish after tangential differentiation and restriction to $H_{a-1}$. Since $a_i\ne0$, the target coefficient must satisfy $\lambda=0$.

This excludes all mixed terms of $R$ outside the ideal. Moreover, $f_{k_{\mathrm{new}}+4}|_{H_t}=F_{k_{\mathrm{new}}+4}$ and $\beta_{\mathrm{bd}}|_{H_t}=0$ give $R|_{H_t}=0$, so $R$ has no term in the remaining variables alone. Thus $R\in\J_t$. Taking $P_{k_{\mathrm{new}}+4}=R$ proves \eqref{eq:boundary}. In particular, among terms whose total degree in the old variables is less than $4$, the only mixed terms are the diagonal quadratic terms displayed in $\beta_{\mathrm{bd}}$. There are no terms linear in the old variables or of the form $x_ix_jA(z)$ with $i\ne j$.

We finally obtain the normal estimates from these degree bounds for the potential.
For $d\le k_{\mathrm{new}}+3$, Lemma \ref{lem:local-arbitrary-delay} gives $f_d-F_d\in\J_t$. For $d=k_{\mathrm{new}}+4$, equation \eqref{eq:boundary}, just proved, gives the additional term $\beta_{\mathrm{bd}}$. Differentiate each of these terms once in an old normal direction $\partial_i$ and then restrict to $H_t$. The derivative of a term involving only the remaining variables is zero. Every monomial of the old ideal contains at least four old-variable factors and still contains an old variable after one differentiation. Also, $\partial_i\beta_{\mathrm{bd}}=-x_iL_i$ vanishes on $H_t$.

Consequently, all terms of degree at most $k_{\mathrm{new}}+3$ in the Taylor expansion of $f_{x_i}|_{H_t}$ vanish. The first possible contribution comes from $f_{k_{\mathrm{new}}+5}$ and has degree $k_{\mathrm{new}}+4$ after one differentiation. Hence
\[
 \left.f_{x_i}\right|_{H_t}=O(|z|^{k_{\mathrm{new}}+4}).
\]
Differentiation in two remaining directions commutes with restriction, so
\[
 T_{i\alpha\beta}|_{H_t}
 =\partial_{\alpha\beta}\bigl(f_{x_i}|_{H_t}\bigr)
 =O(|z|^{k_{\mathrm{new}}+2}).
\]
The derivative estimate here follows by termwise differentiation of the analytic Taylor series; we are not differentiating an arbitrary $O$ estimate. This proves \eqref{eq:old-normal}.
\end{proof}

\begin{corollary}[The leading boundary term of components with two normal indices]\label{cor:boundary-two-normal}
Under the hypotheses of Lemma \ref{lem:next-mixed-normal}, for old normal indices $i,j$ and a remaining index $\alpha$, we have
\[
 T_{ii\alpha}|_{H_t}
 =-\partial_\alpha\langle\xi_i,\operatorname{grad}_zG\rangle+O(|z|^{k_{\mathrm{new}}+2}).
\]
If $i\ne j$, then $T_{ij\alpha}|_{H_t}=O(|z|^{k_{\mathrm{new}}+2})$. The same lower bound for the order holds when all three indices are old normal indices. The homogeneous term displayed in the first formula has degree $k_{\mathrm{new}}+1$ and cannot yet be discarded on the strength of this estimate alone.
\end{corollary}
\begin{proof}
Equation \eqref{eq:old-normal} controls the components with exactly one old normal index. Differentiating the correction once more in an old normal direction gives $\partial_{ii}\beta_{\mathrm{bd}}=-L_i$, and hence
\[
 \left.T_{ii\alpha}\right|_{H_t}
 =-\partial_\alpha\ip{\xi_i}{\operatorname{grad}_zG}
   +O(|z|^{k_{\mathrm{new}}+2}).
\]
The displayed homogeneous term has degree $k_{\mathrm{new}}+1$; we cannot yet assert that it is zero. If $i\ne j$ are both old normal indices, then $\partial_{ij}\beta_{\mathrm{bd}}=0$, so $T_{ij\alpha}|_{H_t}=O(|z|^{k_{\mathrm{new}}+2})$. Components with three old normal indices have the same lower order bound, since $\beta_{\mathrm{bd}}$ is only quadratic in the old variables. All third derivatives of old ideal terms vanish on $H_t$ and therefore do not alter these conclusions.

\end{proof}

Section \ref{sec:inherited} will use the possible degree-$k_{\mathrm{new}}+1$ components with two normal indices in a further comparison of the commutativity equation. On the other hand, Section \ref{sec:ideals} has shown that the commutativity defect of the tangential restriction contains only products of two components with one normal index. Thus \eqref{eq:old-normal} already gives
\[
 O(|z|^{k_{\mathrm{new}}+2})\,O(|z|^{k_{\mathrm{new}}+2})=O(|z|^{2k_{\mathrm{new}}+4}),
\]
which does not affect the comparisons at degrees $2k_{\mathrm{new}}$, $2k_{\mathrm{new}}+1$, and $2k_{\mathrm{new}}+2$ needed for the new layer. This allows us first to construct the new layer in the next section, and then to eliminate the boundary correction found here.

\subsection{The inherited defect and construction of the new layer}\label{subsec:paper-4-4}\label{sec:inherited}
The estimates in Section \ref{sec:boundary} put the restricted commutativity defect above the degrees used to construct a new layer, although $T_{ii\alpha}$ may still have a degree-$k_{\mathrm{new}}+1$ term. After constructing the layer, we eliminate that term by a coefficient of the full equation.

All the following conclusions concern a neighbourhood of $p$. Let $\Sigma_t$ be the great subsphere corresponding to the coordinate subspace $H_t$, with local portion $U\cap\Sigma_t$, and let $C_t^\top$ denote the tangential restriction of $C$ to this portion. Since $\Sigma_t$ is locally totally geodesic, $C_t^\top$ still satisfies the Codazzi condition. For vectors $U,V,W,Z$ tangent to $\Sigma_t$, define the tangential commutativity expression
\[
 \mathcal G_{C_t^\top}(U,V,W,Z)
 :=\sum_\rho\bigl\{
 C_t^\top(U,V,e_\rho)C_t^\top(W,Z,e_\rho)
 -C_t^\top(U,W,e_\rho)C_t^\top(V,Z,e_\rho)
 \bigr\},
\]
where $\{e_\rho\}$ is any orthonormal basis of $T\Sigma_t$. Completing it to an orthonormal basis of $TS^m$, the full commutativity relation gives the normal contraction defect
\begin{equation}\label{eq:defect}
 \mathcal G_{C_t^\top}(U,V,W,Z)
 =-\sum_i\{C(n_i,U,V)C(n_i,W,Z)-C(n_i,U,W)C(n_i,V,Z)\}.
\end{equation}
Here the $n_i$ form a local orthonormal frame of the normal bundle of this subsphere portion, and the sum runs over all old normal directions. By \eqref{eq:old-normal}, every component with one normal index is $O(|z|^{k_{\mathrm{new}}+2})$, so the defect is $O(|z|^{2k_{\mathrm{new}}+4})$. These orders are computed from the Taylor expansion at $p$. To check the orders of the components with one normal index, put $R=\sqrt{1+|z|^2}$ on $H_t$, and let $e_i$ be the fixed ambient unit vector corresponding to an old coordinate axis. The old coordinates satisfy $w^i=0$, so $\partial_i u=R^{-1}e_i$, while $e_i$ itself is tangent to the sphere and orthogonal to $\Sigma_t$. The recovery formula gives
\[
 C(e_i,\partial_\alpha u,\partial_\beta u)
 =R^{-1}T_{i\alpha\beta}|_{H_t}.
\]
Replacing $\partial_\alpha u$ by an analytic orthonormal frame of the subsphere amounts to an analytic invertible tangential change of frame, with no contribution from components with two normal indices. Thus \eqref{eq:old-normal} indeed gives the lower order bound for every component with one normal index in \eqref{eq:defect}.

Consequently, the restricted equation has the same homogeneous coefficients at degrees $2k_{\mathrm{new}},2k_{\mathrm{new}}+1,2k_{\mathrm{new}}+2$ as the commutativity equation without a defect. The new-layer analysis in Sections \ref{subsec:paper-3-1}--\ref{subsec:paper-3-4} uses only these degrees: degree $2k_{\mathrm{new}}$ determines the separated form of the leading term, degree $2k_{\mathrm{new}}+1$ determines the next homogeneous term, and degree $2k_{\mathrm{new}}+2$ gives the Gram relations between the new axes. We do not assert that $C_t^\top$ satisfies the commutativity equation without a defect throughout a neighbourhood.

\Needspace{4\baselineskip}
\begin{lemma}[The new layer to finite order and elimination of old-vector components]\label{lem:inherited-layer}
Use the local notation of Section \ref{sec:ideals} and the conclusions of Section \ref{sec:boundary}. Suppose that the leading term of the remaining potential is $F_{k_{\mathrm{new}}+3}\ne0$, where $k_{\mathrm{new}}>k_t$. Then a constant orthogonal change in the final remaining space divides the coordinates into new active variables $y_1,\ldots,y_s$ and further remaining variables $v$, so that
\begin{align}
 F_{k_{\mathrm{new}}+3}&=\sum_{j=1}^s
 \frac{b_jy_j^{k_{\mathrm{new}}+3}}{(k_{\mathrm{new}}+1)(k_{\mathrm{new}}+2)(k_{\mathrm{new}}+3)},\label{eq:new-leading}\\
 F_{k_{\mathrm{new}}+4}&=\sum_{j=1}^s
 \frac{\beta_jy_j^{k_{\mathrm{new}}+4}}{(k_{\mathrm{new}}+1)(k_{\mathrm{new}}+2)(k_{\mathrm{new}}+3)(k_{\mathrm{new}}+4)}
 +\sum_{j,\delta}
 \frac{e_{j\delta}y_j^{k_{\mathrm{new}}+3}v_\delta}{(k_{\mathrm{new}}+2)(k_{\mathrm{new}}+3)}
 +\text{terms in }v\text{ alone}.\label{eq:new-next}
\end{align}
Here $b_j\ne0$. Put $\eta_j=(e_{j\delta}/b_j)_\delta$. Distinct new axes satisfy $\ip{\eta_j}{\eta_\ell}=-1$. Moreover, for every old axis $i$ and every new axis $y_j$,
\[
 \xi_{ij}=0.
\]
Thus the boundary correction of Section \ref{sec:boundary} vanishes, and every component with an old normal index satisfies
\[
 T_{A B C}|_{H_t}=O(|z|^{k_{\mathrm{new}}+2})
 \quad\text{whenever at least one of }A,B,C\text{ is an old normal index}.
\]
Taking the new variables as active block number $t+1$, denoted by $X_{t+1}$ and setting $k_{t+1}=k_{\mathrm{new}}$, inductive conditions (1)--(4) hold for the new layer. The components of the old vectors along the new axes vanish; the inner-product relations between new and old vectors will be proved in Section \ref{sec:cross-gram}.
\end{lemma}

\begin{proof}

Write $F_{k_{\mathrm{new}}+3}$ as $G$. Every product in the defect \eqref{eq:defect} contains two components with one normal index and therefore has degree at least $2k_{\mathrm{new}}+4$. Thus no additional terms occur in the comparisons of the restricted commutativity equation at degrees $2k_{\mathrm{new}}$, $2k_{\mathrm{new}}+1$, and $2k_{\mathrm{new}}+2$. More precisely, the equation of degree $2k_{\mathrm{new}}$ is exactly the matrix commutativity condition in Lemma \ref{lem:homogeneous-separation}, which we first apply to $G$. The equation of degree $2k_{\mathrm{new}}+1$ gives the same linear constraints as the expansion of the matrix equation in the proof of Lemma \ref{lem:next-homogeneous-form}. Its polynomial argument eliminating mixed components and recovering the potential therefore applies term by term. The equation of degree $2k_{\mathrm{new}}+2$ is likewise the one used to compare coefficients of pure powers and mixed terms in the proof of Lemma \ref{lem:same-layer-gram}, so the corresponding Gram relations hold. This gives the separated leading term \eqref{eq:new-leading}, the next term \eqref{eq:new-next}, and the same-layer relations for the vectors $\eta_j$, namely $\ip{\eta_j}{\eta_\ell}=-1$. These conclusions depend only on the indicated Taylor coefficients of finite order; the local tangential restriction has not been treated as a tensor satisfying the full commutativity equation exactly.

We next compute the components with two normal indices contributed by the boundary correction.
By \eqref{eq:boundary}, the boundary correction is
\[
 \beta_{\mathrm{bd}}=-\frac12\sum_{i\text{ old}}x_i^2
          \ip{\xi_i}{\operatorname{grad}_zG}.
\]
Section \ref{subsec:paper-4-3} has shown that this is the only source of the homogeneous term of degree $k_{\mathrm{new}}+1$ in $T_{ii\alpha}|_{H_t}$: third derivatives of old ideal terms vanish on $H_t$, and terms involving only the remaining variables have no derivatives in old normal directions. Hence
\[
 T_{ii\alpha}^{[k_{\mathrm{new}}+1]}|_{H_t}
 =-\partial_\alpha\ip{\xi_i}{\operatorname{grad}_zG}.
\]
By \eqref{eq:new-leading},
\[
 G_{y_jy_jy_j}=b_jy_j^{k_{\mathrm{new}}},
 \qquad
 \partial_{y_j}\ip{\xi_i}{\operatorname{grad}_zG}
 =\frac{\xi_{ij}b_j}{k_{\mathrm{new}}+1}y_j^{k_{\mathrm{new}}+1},
\]
so
\[
 T_{iij}^{[k_{\mathrm{new}}+1]}|_{H_t}
 =-\frac{\xi_{ij}b_j}{k_{\mathrm{new}}+1}y_j^{k_{\mathrm{new}}+1},
 \qquad (j\text{ denotes the }y_j\text{-direction}).
\]

We now compare $Q_{ijj}$ in the full equation.
Fix an old axis $i$ and a new axis $y_j$. Use the full commutativity equation $Q_{ijj}=0$, take all derivatives before restricting to $H_t$, and compare the coefficient of $y_j^{2k_{\mathrm{new}}+1}$. In the first Euclidean contraction $\sum_rT_{iir}T_{jjr}$, the factors with $r=j$ give the only possible term of degree $2k_{\mathrm{new}}+1$:
\[
 T_{iij}^{[k_{\mathrm{new}}+1]}T_{jjj}^{[k_{\mathrm{new}}]}
 =-\frac{\xi_{ij}b_j^2}{k_{\mathrm{new}}+1}y_j^{2k_{\mathrm{new}}+1}.
\]
Here $T_{jjj}^{[k_{\mathrm{new}}]}=b_jy_j^{k_{\mathrm{new}}}$. If $r$ is another new axis or a remaining direction, the separated form of the pure leading term $G$ gives $T_{jjr}^{[k_{\mathrm{new}}]}=0$, and the contribution has at least one higher degree. If $r$ is an old normal direction, both factors contain an old normal index and again have higher order.

The second Euclidean contraction is $-\sum_rT_{ijr}T_{ijr}$. Its potentially lowest-order component, $T_{iji}=T_{iij}$, starts at order $k_{\mathrm{new}}+1$, so its square has order at least $2k_{\mathrm{new}}+2$. Terms with other internal indices have no lower order. Radial contractions carry two additional coordinate factors and have order at least two greater than the corresponding Euclidean terms. Thus the coefficient equation for $y_j^{2k_{\mathrm{new}}+1}$ is simply
\[
 -\frac{\xi_{ij}b_j^2}{k_{\mathrm{new}}+1}=0.
\]
Since $b_j\ne0$, we obtain $\xi_{ij}=0$. This holds for every old axis and every new axis.

We can now eliminate the correction and restore the inductive conditions.
By the separated leading-term formula, $\operatorname{grad}_zG$ has components only in the new axis directions, whereas we have just proved that $\xi_i$ has zero component in each $y_j$-direction. Hence
\[
 \ip{\xi_i}{\operatorname{grad}_zG}=0\quad\text{for every old axis }i,
 \qquad \beta_{\mathrm{bd}}=0.
\]
Equation \eqref{eq:boundary} reduces to $f_{k_{\mathrm{new}}+4}-F_{k_{\mathrm{new}}+4}\in\J_t$. Put
\[
 I_{t+1}=(X_1,\ldots,X_t,y_1,\ldots,y_s),
 \qquad
 \J_{t+1}=\J_t+I_{t+1}^{k_{\mathrm{new}}+3}.
\]
The leading term $F_{k_{\mathrm{new}}+3}$ belongs to $I_{t+1}^{k_{\mathrm{new}}+3}$. Combining this with the conclusion of Section \ref{sec:nested-lemma} for $d\le k_{\mathrm{new}}+3$, we obtain
\[
 f_d\in\J_{t+1}\qquad(d\le k_{\mathrm{new}}+3).
\]
This is precisely inductive condition (4) for the new layer. At degree $k_{\mathrm{new}}+4$, the pure powers of $y_j$ and the terms $y_j^{k_{\mathrm{new}}+3}v_\delta$ in \eqref{eq:new-next} all belong to $I_{t+1}^{k_{\mathrm{new}}+3}$, but the terms in $v$ alone need not belong to this ideal. Writing $V_{k_{\mathrm{new}}+4}(v)=F_{k_{\mathrm{new}}+4}(0,v)$, the precise conclusion is
\[
 f_{k_{\mathrm{new}}+4}-V_{k_{\mathrm{new}}+4}\in\J_{t+1}.
\]
The retained term $V_{k_{\mathrm{new}}+4}$ is a possible leading term of the remaining potential in the next step and cannot be removed here.

Consider now the components with old normal indices. Each monomial of $\J_t$ has degree at least $k_a+3\ge4$ in the old normal variables of some prefix of blocks. After any three partial derivatives, every nonzero result still contains an old-variable factor and therefore vanishes on $H_t$. The terms $F_d$ in the remaining variables alone have zero derivatives in old normal directions, and the boundary correction has already vanished. Thus the part of the potential of degree at most $k_{\mathrm{new}}+4$ contributes nothing to any third-order component with an old normal index. The first possible contribution comes from $f_{k_{\mathrm{new}}+5}$, whose third derivatives have degree $k_{\mathrm{new}}+2$. Hence
\[
 T_{A B C}|_{H_t}=O(|z|^{k_{\mathrm{new}}+2})
 \quad\text{whenever at least one index is an old normal index}.
\]
The leading term, estimates for old normal components, ideal membership and next homogeneous term are now established for the new layer. Old vectors have zero components along the new axes and therefore lie in the remaining $v$-space with their previous inner products unchanged. The new vectors $\eta_j$ lie there too and satisfy the same-layer Gram relations. The cross-layer relations remain to be proved in Section \ref{sec:cross-gram}.
\end{proof}

\subsection{Gram relations between new and old layers}\label{subsec:paper-4-5}\label{sec:cross-gram}
By Section \ref{sec:inherited}, the old vectors $\xi_i$ and new vectors $\eta_j$ lie in the same remaining $v$-space. Their respective same-layer inner products are known. We prove that every inner product between an old and a new vector is $-1$, allowing the dimension bound to apply to all active axes simultaneously.

The leading degrees $k<k_{\mathrm{new}}$ differ. We compare two pure powers to eliminate the potential coefficients that could enter the desired identity, then compare a mixed term. Only finite Taylor data and the preceding local conclusions are used; neither the entire sphere nor analytic continuation enters.

\begin{lemma}[Gram relations between new and old layers]\label{lem:cross-layer-gram}
Use the notation and conclusions of Lemma \ref{lem:inherited-layer}. For any old active axis $i\in X_a$ and any new active axis $j$, regard the corresponding vectors $\xi_i,\eta_j$ as vectors in the remaining $v$-space. Then
\[
 \ip{\xi_i}{\eta_j}=-1.
\]
The old axis may come from any layer among $X_1,\ldots,X_t$, not only from $X_t$, the layer immediately preceding the new one.
\end{lemma}

\begin{proof}

Fix $i\in X_a$ and a new axis $j$, and write
\[
 k=k_a,\qquad x=x_i,\qquad y=y_j,\qquad
 a=a_i,\qquad b=b_j.
\]
Then $1\le k<k_{\mathrm{new}}$, and both $a,b$ are nonzero. Let $\delta$ denote the remaining $v$-directions after the new block has been chosen. The relevant coefficients of the old and new layers satisfy
\[
 d_{i\delta}=a\xi_{i\delta},\qquad
 e_{j\delta}=b\eta_{j\delta},\qquad
 d_{ij}=a\xi_{ij}=0.
\]
The last equality follows from Lemma \ref{lem:inherited-layer}; the old vector also has zero components along all other new axes.

Take the following component of the full commutativity equation:
\[
 Q_{ijj}={}\sum_r(T_{iir}T_{jjr}-T_{ijr}^{\,2}) +\sum_{r,s}w^rw^s (T_{iir}T_{jjs}-T_{ijr}T_{ijs})=0.
\]
The internal indices $r,s$ always run through all coordinate directions. Let $\Pi_{ij}$ denote the coordinate plane on which only $x,y$ may be nonzero. In each calculation below we take the third derivatives of $f$ first, and only then set the other coordinates to zero. In particular, setting $v=0$ in $f$ before differentiating would lose the internal contractions in the $v$-directions.

For a polynomial $P$, write $[x^\mu y^\nu]P$ for the coefficient of $x^\mu y^\nu$, with zero exponents for all other variables. Define
\[
 p_{ij}=[x^{k+3}y^2]f_{k+5},\qquad
 q_{ij}=[x^2y^{k_{\mathrm{new}}+3}]f_{k_{\mathrm{new}}+5}.
\]
These are the original monomial coefficients in the potential, without division by any additional factorial. All differentiation factors below are computed directly from this definition.

First compare the pure power $x^{2k+2}$ to prove $p_{ij}=0$.
On the $x$-axis, the leading term of the old layer gives
\[
 T_{iii}=ax^k+O(x^{k+1}).
\]
In the next homogeneous term of the old layer, the relevant monomial involving $y$ is
\[
 \frac{d_{ij}}{(k+2)(k+3)}x^{k+3}y.
\]
It contributes $d_{ij}x^{k+1}$ to $T_{iji}$, but contributes zero to $T_{jji}$, which involves two derivatives in $y$. Thus, on the $x$-axis,
\[
 T_{iji}=d_{ij}x^{k+1}+O(x^{k+2}),\qquad T_{jji}=2(k+3)p_{ij}x^{k+2}+O(x^{k+3}).
\]
The coefficient in the second formula comes from
\[
 \partial_i\partial_j^2
       (p_{ij}x^{k+3}y^2)
 =2(k+3)p_{ij}x^{k+2}.
\]
There are no lower-degree contributions: all potential terms of degree at most $k+5\le k_{\mathrm{new}}+4$, apart from those in the remaining variables alone, belong to $\J_t$. Terms in the remaining variables alone cannot produce these low-degree components on the $x$-axis, and the forms of the leading and next terms of the old layer exclude the corresponding terms.

We must also examine the internal contractions with $r\ne i$. On the $x$-axis the full leading term of the old layer has only the $iii$ component, so $T_{iir}=O(x^{k+1})$. In $T_{jjr}$, the two $j$-derivatives do not lower the degree in the old variables, and the internal derivative lowers it by at most one. Every surviving ideal term therefore leaves degree at least $k+2$ in $x$, giving $T_{jjr}=O(x^{k+2})$. The term of degree $k+1$ in $T_{ijr}$ also vanishes: when $r$ is tangent to $H_{a-1}$, use the form of the next term of the old layer; when $r$ is normal to an earlier layer, use the order estimate for components with one normal index. Thus
\[
 T_{iir}T_{jjr}=O(x^{2k+3}),\qquad
 T_{ijr}^{\,2}=O(x^{2k+4})\qquad(r\ne i).
\]
On the $x$-axis the radial contraction retains only $r=s=i$, with the additional factor $x^2$, and cannot contribute to $x^{2k+2}$ either.

Hence the coefficient of this pure power comes only from $T_{iii}T_{jji}-T_{iji}^{\,2}$, giving
\[
 2a(k+3)p_{ij}-d_{ij}^{\,2}=0.
\]
Substituting $d_{ij}=0$ and using $a\ne0$, we obtain $p_{ij}=0$.

Next compare the pure power $y^{2k_{\mathrm{new}}+2}$ to prove $q_{ij}=0$.
Restrict the base point to the $y$-axis, which lies in $H_t$. Lemma \ref{lem:inherited-layer} has bounded all components with an old normal index by $O(|z|^{k_{\mathrm{new}}+2})$. Thus on this axis
\[
 T_{iir}=O(y^{k_{\mathrm{new}}+2}),\qquad
 T_{ijr}=O(y^{k_{\mathrm{new}}+2})\qquad\text{for all }r.
\]
The full leading term of the new layer has only its $jjj$ component nonzero on this axis; that is,
\[
 T_{jjj}=by^{k_{\mathrm{new}}}+O(y^{k_{\mathrm{new}}+1}),\qquad
 T_{jjr}=O(y^{k_{\mathrm{new}}+1})\quad(r\ne j).
\]
Thus the positive products with $r\ne j$ have order at least $2k_{\mathrm{new}}+3$, and all negative squares have order at least $2k_{\mathrm{new}}+4$. The radial terms carry two additional coordinate factors and have order at least $2k_{\mathrm{new}}+4$. Only $T_{iij}T_{jjj}$ can contribute to $y^{2k_{\mathrm{new}}+2}$.

For $T_{iij}$ to contain $y^{k_{\mathrm{new}}+2}$ on the $y$-axis, the only possible monomial of the potential is $q_{ij}x^2y^{k_{\mathrm{new}}+3}$, since
\[
 \partial_i^2\partial_j
       (q_{ij}x^2y^{k_{\mathrm{new}}+3})
 =2(k_{\mathrm{new}}+3)q_{ij}y^{k_{\mathrm{new}}+2}.
\]
Thus the pure-power coefficient equation is
\[
 2b(k_{\mathrm{new}}+3)q_{ij}=0.
\]
Since $b\ne0$, we obtain $q_{ij}=0$. The total degree compared here is $2k_{\mathrm{new}}+2$, whereas the preceding comparison on the $x$-axis had degree $2k+2$. The new and old layers cannot be treated as having leading terms of the same order.

Return to the plane $\Pi_{ij}$ and compare the coefficient of $x^{k+1}y^{k_{\mathrm{new}}+1}$, of total degree $k+k_{\mathrm{new}}+2$. We first compute the following four kinds of contributions, and then show, case by case, that there are no others.

\begin{enumerate}[label=\textup{(\arabic*)},leftmargin=2em,itemsep=5pt]
\item The old leading term and the $q_{ij}$ term. Take the internal index $r=i$, and use the old leading term of $T_{iii}$ and the corresponding term in $T_{jji}$. Since
\[
 \partial_i\partial_j^2(q_{ij}x^2y^{k_{\mathrm{new}}+3})
 =2(k_{\mathrm{new}}+2)(k_{\mathrm{new}}+3)q_{ij}xy^{k_{\mathrm{new}}+1},
\]
multiplication by $ax^k$ gives the target coefficient $2a(k_{\mathrm{new}}+2)(k_{\mathrm{new}}+3)q_{ij}$.

\item The new leading term and the $p_{ij}$ term. Take the internal index $r=j$, and use the new leading term of $T_{jjj}$ and the corresponding term in $T_{iij}$. From
\[
 \partial_i^2\partial_j(p_{ij}x^{k+3}y^2)
 =2(k+2)(k+3)p_{ij}x^{k+1}y,
\]
multiplication by $by^{k_{\mathrm{new}}}$ gives the target coefficient $2b(k+2)(k+3)p_{ij}$.

\item Ordinary contractions in the common remaining directions. The next potential terms of the old and new layers contain, respectively,
\[
 \frac{d_{i\delta}x^{k+3}v_\delta}{(k+2)(k+3)},
 \qquad
 \frac{e_{j\delta}y^{k_{\mathrm{new}}+3}v_\delta}{(k_{\mathrm{new}}+2)(k_{\mathrm{new}}+3)}.
\]
Applying $\partial_i^2\partial_\delta$ and $\partial_j^2\partial_\delta$, respectively, gives
\[
 d_{i\delta}x^{k+1},\qquad
 e_{j\delta}y^{k_{\mathrm{new}}+1}.
\]
Thus summing over all $v$-directions contributes $\sum_\delta d_{i\delta}e_{j\delta}$.

\item The radial contraction of the two leading terms. The old leading term has only $T_{iii}=ax^k$, and the new one has only $T_{jjj}=by^{k_{\mathrm{new}}}$. In the radial sum take $r=i,s=j$, with the coordinate factor $xy$. Then
\[
 xy\,(ax^k)(by^{k_{\mathrm{new}}})=abx^{k+1}y^{k_{\mathrm{new}}+1}.
\]
This contributes the term $+ab$; no sign is assumed for $ab$. Both leading terms have zero $T_{ijr}$ components, so there is no corresponding negative term.
\end{enumerate}

Collecting these contributions, the target coefficient equation should be
\begin{equation}\label{eq:cross}
 0={}2a_i(k_{\mathrm{new}}+2)(k_{\mathrm{new}}+3)q_{ij} +2b_j(k+2)(k+3)p_{ij} +\sum_\delta d_{i\delta}e_{j\delta}+a_ib_j.
\end{equation}

We now exclude all other contributions, case by case.
Lemma \ref{lem:inherited-layer} gives
\[
 f_d-F_d\in\J_t\qquad(d\le k_{\mathrm{new}}+4).
\]
We accordingly divide the Taylor monomials of the potential into three classes: those belonging to the old ideal $\J_t$; those involving only the remaining variables $z$ (the new active variables and $v$); and those involving old variables but not belonging to $\J_t$. The third class can start only at degree $k_{\mathrm{new}}+5$. Terms in the remaining variables alone start at degree $k_{\mathrm{new}}+3$, and Lemma \ref{lem:inherited-layer} has determined their first two terms.

We first give a common lower degree bound for terms in the old ideal. Assign such a monomial to some $I_c^{k_c+3}$. If $c<a$, three differentiations still leave a variable from the first $a-1$ blocks, so the term vanishes on $\Pi_{ij}$. For a surviving term we may therefore take $c\ge a$, and its total degree in the old variables is at least $k_c+3\ge k+3$. After restriction to $\Pi_{ij}$, only $x$ can remain among the old variables; an internal derivative may also remove one occurrence of another old variable. This possible decrease of degree under the internal derivative must be included.

\begin{enumerate}[label=\textup{(\arabic*)},leftmargin=2em,itemsep=5pt]
\item Pairing two terms of the old ideal. In $T_{iir}T_{jjr}$, the three derivatives of the first factor remove at most three old-variable factors, and those of the second remove at most one. Thus a surviving product has degree in $x$ at least
\[
 (k+3-3)+(k+3-1)=2k+2>k+1.
\]
In $T_{ijr}^{\,2}$, each factor loses at most two old-variable factors, so the product again retains degree at least $2(k+1)=2k+2$ in $x$. Neither can produce the target mixed term.

In a radial contraction, if the internal derivative is in an old variable, the corresponding coordinate multiplier restores the lost degree. If that multiplier is in a direction outside the plane, the term vanishes outright. Thus the net loss of old-variable degree comes only from the two external $i$-derivatives, and the lower bound increases to $2k+4$. There is again no contribution.

\item Pairing an old ideal term with a term in the remaining variables alone. In an ordinary contraction, the latter term cannot be differentiated in the $i$-direction. Hence the only possibly nonzero positive product has $T_{iir}$ from the old term and $T_{jjr}$ from the term in the remaining variables, with $r$ necessarily a remaining direction.

Let the degrees of the two potential terms be $D,E$. Three differentiations lower each degree by three, so a contribution of total degree $k+k_{\mathrm{new}}+2$ requires
\[
 D+E=k+k_{\mathrm{new}}+8,\qquad D\ge k+3,\qquad E\ge k_{\mathrm{new}}+3.
\]
There are exactly three possible degree pairings:
\[
 (D,E)=(k+3,k_{\mathrm{new}}+5),\ (k+4,k_{\mathrm{new}}+4),\ (k+5,k_{\mathrm{new}}+3).
\]
In the first pair, the third derivatives of the old term can give only the old leading $iii$ component, with internal direction $i$. But the derivative in $i$ of the term involving only the remaining variables is zero, so this pair contributes nothing.

In the second pair, the next term of the old layer gives $T_{ii\delta}=d_{i\delta}x^{k+1}$, and the next term of the new layer gives $T_{jj\delta}=e_{j\delta}y^{k_{\mathrm{new}}+1}$. This is precisely contribution \textup{(3)} above. If the internal index is the new axis $y_j$, a possible pure power of $y$ in the new layer is multiplied by $d_{ij}=0$. The other new axis directions also contribute nothing, since the old vectors have zero components along every new axis.

In the third pair, the new leading term forces the internal index to be $j$ and gives the factor $by^{k_{\mathrm{new}}}$. To leave the target exponents, the other factor must contain $x^{k+1}y$. The unique monomial before differentiation is $p_{ij}x^{k+3}y^2$, which is contribution \textup{(2)} above.

A radial contraction adds two to the total degree, so the potential degrees must satisfy $D+E=k+k_{\mathrm{new}}+6$. The two lower bounds force $(D,E)=(k+3,k_{\mathrm{new}}+3)$, pairing the old and new leading terms and giving only contribution \textup{(4)} above.

\item Pairing two terms in the remaining variables alone. Each product contains a component with an $i$ index, while each such potential term has zero derivative in the $i$-direction. Both ordinary and radial contractions therefore vanish.

\item A pairing involving at least one higher mixed term of the third class. Its potential degree is at least $k_{\mathrm{new}}+5$. On $\Pi_{ij}\subset H_{a-1}$ the full third-order array has no part of degree less than $k$, so the potential degree of the other surviving factor is at least $k+3$. The ordinary contraction requires the sum of the potential degrees to be $k+k_{\mathrm{new}}+8$, leaving only
\[
 k_{\mathrm{new}}+5\quad\text{and}\quad k+3.
\]
The old leading term has only the $iii$ component. Hence the higher term can pair with $ax^k$ only through $T_{jji}$. The target exponents force it to come from $q_{ij}x^2y^{k_{\mathrm{new}}+3}$, which is contribution \textup{(1)} above. Mixed terms of higher degree exceed the total degree under comparison.

This also excludes pairing it with a term in the remaining variables alone or with another higher mixed term: in the former case, for example, the sum of the potential degrees is at least $2k_{\mathrm{new}}+8>k+k_{\mathrm{new}}+8$. A radial contraction requires the still smaller sum $k+k_{\mathrm{new}}+6$, so no term of the third class contributes radially.
\end{enumerate}
This accounts for all internal directions and the negative-square terms, proving that \eqref{eq:cross} omits no contribution.

Finally, we normalize to obtain all Gram relations between new and old layers.
The two pure-power comparisons have already given $p_{ij}=q_{ij}=0$. Substitution in \eqref{eq:cross} yields
\[
 \sum_\delta d_{i\delta}e_{j\delta}+a_ib_j=0.
\]
Dividing by the nonzero number $a_ib_j$, we obtain
\[
 \ip{\xi_i}{\eta_j}
 =\sum_\delta\frac{d_{i\delta}}{a_i}
                  \frac{e_{j\delta}}{b_j}
 =-1.
\]
Only the nonvanishing of $a_i,b_j$ is used; they need not have the same sign. The old layer $a$ and the axes $i,j$ were arbitrary, so the relations hold between every old layer and the new one. This completes inductive condition (5).
\end{proof}

\subsection{The local zero subsphere theorem}\label{subsec:paper-4-6}\label{sec:rank-budget}
We now use the joint Gram dimension bound to terminate the local induction near $p$, then convert the resulting ideal membership into vanishing of the entire spherical tensor.

\begin{theorem}[Local zero subsphere]\label{thm:local-zero}
Let $m\ge2$, let $U\subset S^m$ be an open neighbourhood of $p$, and let $C$ be a real-analytic symmetric cubic tensor on $U$ satisfying the spherical Codazzi--commutation system \eqref{eq:sphere-system} on $U$, with $C_p=0$.
If the germ of $C$ at $p$ is nonzero, then there exist an integer
\[
 1\le R\le\left\lfloor\frac{m+1}{2}\right\rfloor,
\]
a great subsphere $\Sigma\cong S^{m-R}$ through $p$, and an open neighbourhood $U'\subset U$ of $p$, such that
\[
 C_q=0\in\Sym^3T_q^*S^m
 \qquad\text{for every }q\in U'\cap\Sigma.
\]
In particular, there is a nontrivial great-circle arc through $p$ on which the entire cubic tensor $C$ vanishes.
\end{theorem}

\begin{proof}
Choose the normalized analytic potential $f$ of Section \ref{subsec:paper-3-1}, with $T=D^3f$. Normalization gives $f_0=f_1=f_2=0$, while $C_p=0$ and Lemma \ref{lem:affine-slice-potential} give $f_3=0$. The recovery formula in that lemma also shows that the germ of $f$ is nonzero when the germ of $C$ is nonzero. Thus Sections \ref{subsec:paper-3-1}--\ref{subsec:paper-3-4} construct the first layer from its leading term.

Suppose that several layers have been completed and that a total of $R$ active axes have been selected. The remaining coordinate space has dimension $m-R$. Number the vectors from all layers as
\[
 \xi_1,\ldots,\xi_R\in\R^{m-R}.
\]
The old vectors have zero components along active axes selected at later steps. Thus their previous inner products are unchanged when they are regarded as vectors in this common remaining space. The same-layer Gram relations together with Lemma \ref{lem:cross-layer-gram} give
\[
 \ip{\xi_i}{\xi_j}=-1\qquad(i\ne j).
\]
As in Section \ref{subsec:paper-3-4}, append one coordinate to each vector and put
\[
 \widehat\xi_i=(\xi_i,1)\in\R^{m-R+1}.
\]
For $i\ne j$,
\[
 \ip{\widehat\xi_i}{\widehat\xi_j}
 =\ip{\xi_i}{\xi_j}+1=0,
 \qquad
 |\widehat\xi_i|^2=|\xi_i|^2+1>0.
\]
Thus $\widehat\xi_1,\ldots,\widehat\xi_R$ are $R$ pairwise orthogonal nonzero vectors in a space of dimension $m-R+1$. They are linearly independent, so $R\le m-R+1$, or
\begin{equation}\label{eq:rank-budget}
 2R\le m+1.
\end{equation}
This bound holds after every completed step. The original vectors $\xi_i$ need not be linearly independent. The argument also allows one of the $\xi_i$ to be zero, since each augmented vector is nonzero.

We next check each inductive condition after the new block has been chosen.
Suppose that the first $t$ layers have been completed and that the leading term of $F=f|_{H_t}$ is $F_{k_{\mathrm{new}}+3}\ne0$, where $k_{\mathrm{new}}>k_t$. Sections \ref{subsec:paper-4-2}--\ref{subsec:paper-4-5} allow us to choose a nonempty new active block $X_{t+1}$, divide the remaining variables into new active variables $y$ and further remaining variables $v$, and put $k_{t+1}=k_{\mathrm{new}}$. The five conditions of Section \ref{sec:ideals} are restored as follows.
\begin{enumerate}[label=\textup{(\arabic*)},leftmargin=2em,itemsep=5pt]
\item The leading term of the full tensor. Lemma \ref{lem:inherited-layer} gives the purely tangential leading term
\[
 \sum_{j\in X_{t+1}}b_jy_j^{k_{\mathrm{new}}}(\dd y_j)^{\otimes3},\qquad b_j\ne0,
\]
and shows that every component with an old normal index has order at least $k_{\mathrm{new}}+2$. Hence the lowest-degree part of the full array restricted to $H_t$ is exactly the displayed expression, with no additional normal leading terms.
\item The order of components with one normal index. The same lemma gives
\[
 T_{i\alpha\beta}|_{H_t}=O(|z|^{k_{\mathrm{new}}+2}),
\]
where $i$ is an old normal index, $\alpha,\beta$ are tangent indices of $H_t$, and $z=(y,v)$. This is precisely condition (2) for the new layer.
\item The form of the next homogeneous term. Equation \eqref{eq:new-next} gives the next term of the new layer; its mixed coefficients determine the new vectors $\eta_j$.
\item Ideal membership of the lower-degree potential terms. Put
\[
 I_{t+1}=(X_1,\ldots,X_t,y),\qquad
 \J_{t+1}=\J_t+I_{t+1}^{k_{\mathrm{new}}+3}.
\]
Using $f_d-F_d\in\J_t$ for $d\le k_{\mathrm{new}}+3$, $F_d=0$ for $d<k_{\mathrm{new}}+3$, and $F_{k_{\mathrm{new}}+3}\in I_{t+1}^{k_{\mathrm{new}}+3}$, we obtain at each degree
\[
 f_d\in\J_{t+1}\qquad(d\le k_{\mathrm{new}}+3).
\]
Only degrees through $k_{\mathrm{new}}+3$ need to be controlled here. If $V_{k_{\mathrm{new}}+4}(v)=F_{k_{\mathrm{new}}+4}(0,v)$, the precise conclusion at the next degree is
\[
 f_{k_{\mathrm{new}}+4}-V_{k_{\mathrm{new}}+4}\in\J_{t+1};
\]
The term in $v$ alone must be retained, since it may become the leading term at the next step.
\item All Gram relations in the common remaining space. Lemma \ref{lem:inherited-layer} shows that the old vectors have zero components along every new axis. They can therefore all be regarded as vectors in the $v$-space, with their old inner products unchanged. The same lemma gives the same-layer relations between the new vectors; Lemma \ref{lem:cross-layer-gram} supplies the inner product $-1$ between any old and any new vector.
\end{enumerate}
Thus all five conditions are available at the next step. Choosing the new axes involves only a constant orthogonal change in the current remaining variables. It fixes the old active variables and orthogonally transforms the remaining components of the old vectors at the same time. Hence it preserves the old leading terms, old ideals, and established inner-product relations. The conditions already established for earlier layers are therefore retained.

To distinguish the steps, write
\[
 F^{(t)}=f|_{H_t},\qquad
 R_t=\sum_{a=1}^t\dim X_a,\qquad
 N=\left\lfloor\frac{m+1}{2}\right\rfloor.
\]
Each active block is nonempty, so $t\le R_t$; the preceding Gram dimension estimate gives $R_t\le N$.

If $F^{(t)}$ is a nonzero analytic germ, its convergent Taylor expansion has a leading term of finite degree. By ideal membership of the lower-degree potential terms, restricting $f_d\in\J_t$ to $H_t$ gives
\[
 F_d^{(t)}=0\qquad(d\le k_t+3).
\]
Thus the degree of this leading term is $k_{\mathrm{new}}+3$ with $k_{\mathrm{new}}>k_t$. The new-layer construction above then gives a nonempty new block, and
\[
 R_{t+1}=R_t+\dim X_{t+1}\ge R_t+1.
\]
Consequently, as long as the remaining analytic germ is nonzero, the total number of active axes strictly increases. But it cannot exceed the fixed integer $N$. The induction must therefore stop after at most $N$ nonempty active blocks, and at termination we must have
\[
 F^{(t)}=0\qquad\text{as an analytic function germ at }0.
\]
Otherwise another nonempty block could be added, contradicting the bound on the number of active axes. No bound has been imposed on the differences $k_{a+1}-k_a$ between consecutive layer orders; only the number of nonempty layers is bounded. The whole procedure uses only finitely many coordinate changes and neighbourhood shrinkings, so an open neighbourhood of $p$ remains at the end.

Fix the last layer obtained at termination, and let $R$ denote its total number of active axes. The identity $f|_{H_t}=0$ alone does not eliminate third derivatives with normal indices. We use part \textup{(2)} of Lemma \ref{lem:local-arbitrary-delay}, which gives
\[
 f\in\J_t
 \qquad\text{in the ring of convergent analytic function germs}.
\]
More explicitly, in some neighbourhood $W'$ of the origin we can write a finite sum
\[
 f=\sum_{\ell=1}^L G_\ell A_\ell,
\]
where the $G_\ell$ are monomial generators of $\J_t$ and the $A_\ell$ are analytic functions. Every $G_\ell$ has total degree at least $4$ in the active variables, since the generating degree $k_a+3$ of each ideal is at least $4$.

Take any differentiation multi-index $\mu$ of length $3$. Leibniz's formula gives
\[
 \partial^\mu(G_\ell A_\ell)
 =\sum_{\nu\le\mu}\binom\mu\nu
       (\partial^\nu G_\ell)
       (\partial^{\mu-\nu}A_\ell).
\]
Here $|\nu|\le3$. Thus, whenever $\partial^\nu G_\ell$ is nonzero, it still has active-variable degree at least $4-|\nu|\ge1$. Setting all active variables to zero on $H_t$ makes every summand vanish. Hence
\[
 \left.\partial^\mu f\right|_{W'\cap H_t}=0
 \qquad(|\mu|=3).
\]
The differentiation directions are arbitrary and may include normal directions. Thus all components of $T=D^3f$ vanish. By the recovery formula in Lemma \ref{lem:affine-slice-potential}, the entire original spherical tensor $C$ also vanishes at the corresponding points.

Let $Z=H_t$ be the final remaining coordinate subspace, of dimension $m-R$. The original affine coordinates naturally identify $\R^m$ with $T_pS^m=p^\perp$, so $Z$ may be regarded as a subspace of the ambient space $\R^{m+1}$ orthogonal to $p$. Put
\[
 \Sigma=S^m\cap(Z\oplus\R p).
\]
Since $p\notin Z$, the space $Z\oplus\R p$ has dimension $m-R+1$. Its unit sphere $\Sigma$ is a great subsphere through $p$, with $\Sigma\cong S^{m-R}$.

The affine coordinate map can be written as
\[
 u(w)=\frac{p+w}{\sqrt{1+|w|^2}}.
\]
It maps a neighbourhood of the origin in $Z$ to a neighbourhood of $p$ in $\Sigma$. Shrinking $W'$ so that $U'=u(W')\subset U$, we have
\[
 u(W'\cap Z)=U'\cap\Sigma.
\]
Thus $C_q=0$ for every $q\in U'\cap\Sigma$, as asserted.

Finally, \eqref{eq:rank-budget} and $m\ge2$ give
\[
 \dim Z=m-R
 \ge m-\left\lfloor\frac{m+1}{2}\right\rfloor
 =\left\lfloor\frac m2\right\rfloor\ge1.
\]
We may therefore choose a unit vector $e\in Z$ and define
\[
 \gamma(s)=\cos s\,p+\sin s\,e.
\]
Since $p,e$ are orthogonal unit vectors, $\gamma$ is a unit-speed great circle in $\Sigma$ through $p$, and $\gamma'(0)=e\ne0$. For $|s|<\pi/2$,
\[
 u((\tan s)e)=\cos s\,p+\sin s\,e=\gamma(s).
\]
Choose $\varepsilon>0$ small enough that $(\tan s)e\in W'$ for $|s|<\varepsilon$. Then
\[
 C_{\gamma(s)}=0\qquad(|s|<\varepsilon).
\]
This is a nontrivial great-circle arc, and the theorem is proved.
\end{proof}

\begin{remark}[The local conclusion and the role of analyticity]\label{rem:local-zero-scope}
The vanishing is of the entire tensor, including normal components, on the local portion $U'\cap\Sigma$; it asserts neither that $C$ is defined on all of $\Sigma$ nor that it vanishes on the open set $U'$. If the germ at $p$ is zero, $C$ already vanishes on a neighbourhood, and a great-circle arc through $p$ may be chosen without selecting active axes.

Analyticity enters twice: a nonzero remaining germ has a finite leading degree, allowing the induction to continue; and Section \ref{subsec:paper-4-2} passes from monomial membership of the convergent Taylor series to membership of the analytic germ in $\J_t$. The resulting analytic decomposition gives vanishing of the entire tensor on the local portion of the final subsphere. A nonzero smooth germ may have zero Taylor series, and formal ideal membership alone does not supply the required convergent decomposition.

\end{remark}

Starting from the zero supplied by Section \ref{subsec:paper-2-4}, we have obtained a local great-circle arc of zeros. Section \ref{sec:paper-4} uses an analytic identity along the entire circle and periodicity to prove global vanishing.

\section{Continuation along great circles and analytic Laugwitz rigidity}\label{sec:paper-4}
We apply Theorem \ref{thm:local-zero} to tensors on the entire standard sphere. Fixed normal axes along a great-circle arc of zeros, the Codazzi equations on a corresponding two-sphere, and estimates of all remaining normal contractions will give a scalar equation for the leading coefficient.

The same fixed axis defines an analytic periodic coefficient on the entire great circle. The division-free scalar equation extends by the analytic identity theorem, and integration over a period gives global rigidity of analytic cubic forms and analytic Laugwitz rigidity.

\subsection{Fixed normal axes along a great-circle arc of zeros}\label{subsec:paper-5-1}\label{sec:fixed-axis}
For the arc of zeros supplied by Sections \ref{subsec:paper-4-1}--\ref{subsec:paper-4-6}, we first prove that the active normal axes of the lowest nonzero term are fixed ambient vectors in $\R^{m+1}$.

A normal axis here lies in $T_{\gamma(s)}S^m$ perpendicular to the arc, not in the radial direction $\gamma(s)$. Before constancy is proved, it is a vector field $n_i(s)$ satisfying $h_0(n_i(s),\dot\gamma(s))=0$.

\begin{lemma}[Fixed normal axes along a great-circle arc of zeros]\label{lem:fixed-normal-axes}
Let $m\ge2$, let $U\subset S^m$ be connected and open, and let $C$ be a real-analytic symmetric cubic tensor on $U$, not identically zero, satisfying the spherical Codazzi--commutation system \eqref{eq:sphere-system} on $U$.
Let $\gamma:I\to U$ be a unit-speed great-circle arc defined on a nonempty open interval, with $C_{\gamma(s)}=0$ for all $s\in I$. Then there exist an integer $k\ge1$, a nonempty open subinterval $J\subset I$, pairwise orthogonal fixed unit vectors $n_1,\ldots,n_r\in\R^{m+1}$, and analytic functions $A_i$ nowhere zero on $J$, such that
\[
 (\nabla^{h_0})^jC|_{\gamma(J)}=0\qquad(0\le j<k),\qquad ((\nabla^{h_0})^kC)_{\gamma(s)} =k!\sum_{i=1}^r A_i(s) (n_i^\flat)^{\otimes(k+3)}\qquad(s\in J).
\]
Here $r\ge1$, each $n_i$ is perpendicular to the two-dimensional linear plane containing the great circle, and $n_i^\flat$ at $\gamma(s)$ denotes $h_0(n_i,\cdot)$.
\end{lemma}

\begin{proof}

The tensor $(\nabla^{h_0})^jC$ has $j+3$ covariant indices and includes covariant derivatives in all directions, not merely repeated differentiation along $\gamma$. Set
\[
 k=\min\{j\ge1:(\nabla^{h_0})^jC|_{\gamma(I)}\not\equiv0\}.
\]
This set is nonempty. Otherwise all covariant derivatives of $C$ would vanish at any point $q$ of the arc. Comparing covariant derivatives with coordinate partial derivatives successively shows that all Taylor coefficients at that point also vanish. Analyticity then gives $C=0$ in an open neighbourhood of $q$, and the analytic identity theorem on the connected open set $U$ gives $C\equiv0$, contrary to the hypothesis.

Thus $k$ is finite, and
\[
 (\nabla^{h_0})^jC|_{\gamma(I)}=0\qquad(0\le j<k).
\]
Write $X=\dot\gamma$ and $\Theta(s):=((\nabla^{h_0})^kC)_{\gamma(s)}$. Minimality of $k$ concerns the whole arc; $\Theta$ may vanish at individual points. We shall choose a subinterval near a point of maximal rank on which $\Theta$ is nowhere zero and has constant rank.

Fix a point $q$ on the arc and choose an orthonormal basis $e_1,\ldots,e_m$ of $T_qS^m$. Use the affine coordinates of Section \ref{subsec:paper-3-1} centred at $q$,
\[
 u_q(w)=\frac{q+\sum_a w^a e_a}{\sqrt{1+|w|^2}},
 \qquad u_q(0)=q,
\]
and take the corresponding local potential $f$. The recovery formula and the coordinate expression for the spherical metric are, respectively,
\[
 C_{abc}(w)=(1+|w|^2)^{-1}f_{abc}(w),\qquad (h_0)_{ab}(w)=\frac{\delta_{ab}}{1+|w|^2} -\frac{w^aw^b}{(1+|w|^2)^2}.
\]
Thus $(h_0)_{ab}(0)=\delta_{ab}$ and $\partial_c(h_0)_{ab}(0)=0$, so the connection coefficients satisfy $\Gamma^a_{bc}(0)=0$.

For $r\ge2$, repeated expansion of covariant derivatives gives, at the origin,
\[
 ((\nabla^{h_0})^rC)(q)=(\partial^rC)(0)
       +\text{linear terms formed from }(\partial^jC)(0),\ 0\le j\le r-2
          .
\]
Indeed, every correction containing partial derivatives of order $r-1$ is multiplied by the undifferentiated $\Gamma(0)$ and therefore vanishes. Whenever a derivative falls on $\Gamma$, at most $r-2$ derivatives fall on $C$. For $r=1$ we directly have $\nabla^{h_0}C(q)=\partial C(0)$.

From $(\nabla^{h_0})^jC|_{\gamma(I)}=0$ for $j<k$, successive comparison gives $\partial^jC(0)=0$ for $j<k$. The relation $D^3f=(1+|w|^2)C$ likewise gives $\partial^j(D^3f)(0)=0$ for $j<k$. The recovery factor equals $1$ at the origin and has zero first derivative there; its higher derivatives multiply only the derivatives of order less than $k$ already shown to vanish. Thus, for $r=k,k+1$, we have the componentwise identity
\[
 ((\nabla^{h_0})^rC)_{a_1\cdots a_r bcd}(q)
 =\partial_{a_1}\cdots\partial_{a_r}
       \partial_b\partial_c\partial_d f(0).
\]
The right-hand side is completely symmetric, so these two covariant derivative tensors are also completely symmetric at $q$. After normalization, the first possibly nonzero homogeneous term of $f$ is $f_{k+3}$; if $\Theta(q)\ne0$, this term is indeed nonzero.

At such a $q$, Lemma \ref{lem:homogeneous-separation} gives orthogonal active axes $n_i(q)$ such that
\[
 f_{k+3}=\sum_i
 \frac{A_i(q)x_i^{k+3}}{(k+1)(k+2)(k+3)}.
\]
Differentiating each term $k+3$ times along its axis changes its coefficient to
\[
 \frac{(k+3)!}{(k+1)(k+2)(k+3)}A_i(q)=k!A_i(q).
\]
Hence, pointwise,
\[
 \Theta(q)=k!\sum_i A_i(q)(n_i(q)^\flat)^{\otimes(k+3)},
 \qquad A_i(q)\ne0.
\]
Taking a covariant derivative along the arc of the zero tensor $(\nabla^{h_0})^{k-1}C$, we obtain
\[
 \Theta(X,\cdot,\ldots,\cdot)=0.
\]
Inserting $n_i$ into each of the remaining $k+2$ slots, orthogonality gives
\[
 0=\Theta(X,n_i,\ldots,n_i)=k!A_i h_0(X,n_i).
\]
Since $A_i\ne0$, we obtain $h_0(X,n_i)=0$. The pointwise axes are therefore perpendicular to the arc, but have not yet been chosen differentiably.

To choose the axes analytically on a shorter arc, define the linear map
\[
 \mathcal L_s:T_{\gamma(s)}S^m
       \longrightarrow\Sym^{k+2}T_{\gamma(s)}^*S^m,
 \qquad
 \mathcal L_s(Z)=\Theta(s)(Z,\cdot,\ldots,\cdot).
\]
In the pure-power decomposition above, $\mathcal L_s(n_i)=k!A_i(n_i^\flat)^{\otimes(k+2)}$. These images are linearly independent, while every vector orthogonal to all active axes belongs to the kernel. Thus $\rank\mathcal L_s$ is exactly the number of active axes.

Choose a point $s_0$ in $I$ at which this rank attains its maximum $r\ge1$. A nonzero minor of order $r$ remains nonzero nearby, so the rank is at least $r$ after shortening the interval; maximality ensures that it is no greater than $r$. The rank is therefore constant and equal to $r$ on this subinterval. Here $r$ counts the active axes of the first nonzero, order-$k$ jet at each point of this subinterval; $R$ in Theorem \ref{thm:local-zero} counts the axes accumulated over all layers at a fixed point.

Choose an analytic orthonormal parallel frame along the arc, identifying the tangent spaces with a fixed $\R^m$. Such a frame is obtained from the parallel translation equation with analytic coefficients. At $s_0$, choose the unit active axes $n_i^0$ from the pure-power decomposition of the leading term, and choose a tangent vector
\[
 v_0=\sum_{i=1}^r c_i n_i^0,
\]
so that the real numbers $k!A_i(s_0)c_i^{k+1}$ are all nonzero and pairwise distinct. This is possible: choose $c_i\ne0$ successively, avoiding at each step the finitely many values that give a previously obtained number. Let $v(s)$ be the parallel extension of $v_0$ along the arc.

Define an analytic symmetric operator field $\mathsf{B}(s)$ by
\[
 h_0(\mathsf{B}(s)Y,Z)=\Theta(s)(\underbrace{v(s),\ldots,v(s)}_{k+1\text{ entries}},Y,Z).
\]
The pointwise pure-power decomposition shows that its active eigenvalues are
\[
 k!A_i(s)h_0(v(s),n_i(s))^{k+1}.
\]
The corresponding eigenlines are exactly the active axes, and all remaining directions have eigenvalue zero. At this stage the expression is understood through the pointwise decomposition; analyticity of $n_i(s)$ has not been assumed.

At $s_0$, $\mathsf{B}$ has $r$ pairwise distinct nonzero simple eigenvalues. By the analytic implicit function theorem and continuity of eigenvalues, they extend on a smaller interval to analytic functions $\lambda_1(s),\ldots,\lambda_r(s)$, remaining nonzero and pairwise distinct. On the other hand, $\rank \mathsf{B}(s)\le\rank\mathcal L_s=r$, so no other nonzero eigenvalues can occur. The orthogonal spectral projections onto the corresponding eigenlines can be written explicitly as
\[
 P_i(s)=\frac{\mathsf{B}(s)}{\lambda_i(s)}
       \prod_{j\ne i}\frac{\mathsf{B}(s)-\lambda_j(s)\Id}
                              {\lambda_i(s)-\lambda_j(s)}.
\]
This polynomial equals $1$ at $\lambda_i$ and zero at every other $\lambda_j$ and at zero. It is therefore the required projection and depends analytically on $s$.

Extend $n_i^0$ by parallel translation to $E_i(s)$. Since $P_i(s_0)E_i(s_0)=n_i^0\ne0$, after shortening the interval we can define
\[
 n_i(s)=\frac{P_i(s)E_i(s)}{|P_i(s)E_i(s)|}.
\]
These vector fields are analytic, unit, and pairwise orthogonal. They remain precisely all the active axes of the pure-power decomposition, since the nonzero simple eigenlines must be those axes. The corresponding coefficients are determined by
\[
 A_i(s)=\frac1{k!}\Theta(s)(n_i(s),\ldots,n_i(s))
\]
and are analytic and nonzero. Fix such a nonempty open subinterval $J$. The leading-term decomposition and $h_0(X,n_i)=0$ are now analytic tensor identities on $J$.

Fix $s\in J$. Continue to use affine coordinates centred at $q=\gamma(s)$, with $n_1(s),\ldots,n_r(s)$ as active axes. Lemma \ref{lem:next-homogeneous-form} gives
\[
 f_{k+4}=\sum_j\beta_jx_j^{k+4}
       +\sum_{j,\alpha}c_{j\alpha}x_j^{k+3}z_\alpha
       +G_{k+4}(z).
\]
Fix an active axis $n_i$ and take any $W\in T_qS^m$ with $W\perp n_i$. The contraction of the leading tensor along the arc has already given $X\perp n_i$. Consider the derivative of total order $k+4$
\[
 \partial_X\partial_{n_i}^{\,k+2}\partial_W f_{k+4}.
\]
Every type of monomial contributes zero. The pure power $x_i^{k+4}$ is annihilated by any derivative transverse to its axis. The mixed term $x_i^{k+3}z_\alpha$ has degree only one in variables transverse to the axis and is annihilated by the two transverse derivatives in $X,W$. Terms associated with other active axes, and terms in $z$ alone, are annihilated by the $k+2$ derivatives along $n_i$. Even if $W=X$, there are still two transverse derivatives.

Comparing the full differential of the potential of order $k+4$ with $(\nabla^{h_0})^{k+1}C$ gives
\[
 ((\nabla^{h_0})^{k+1}C)(X,\underbrace{n_i,\ldots,n_i}_{k+2\text{ entries}},W)=0
 \qquad(W\perp n_i).
\]
All slots of $(\nabla^{h_0})^{k+1}C$ are completely symmetric here, and $\nabla^{h_0}_X\Theta$ is obtained by inserting $X$ into one covariant differentiation slot. Thus also
\[
 (\nabla^{h_0}_X\Theta)(n_i,\ldots,n_i,W)=0\qquad(W\perp n_i).
\]

Put $N=k+3$. Taking the covariant derivative along the arc of
\[
 \Theta=k!\sum_j A_j(n_j^\flat)^{\otimes N}
\]
and applying the product rule, we obtain
\[
 \nabla^{h_0}_X\Theta=k!\sum_j\left\{
 X(A_j)(n_j^\flat)^{\otimes N}
 +A_j\sum_{\ell=1}^N
 (n_j^\flat)^{\otimes(\ell-1)}\otimes \nabla^{h_0}_Xn_j^\flat
       \otimes(n_j^\flat)^{\otimes(N-\ell)}
 \right\}.
\]
Insert $n_i$ into the first $N-1=k+2$ slots and $W\perp n_i$ into the last slot, and examine each term.
\begin{enumerate}[label=\textup{(\arabic*)},leftmargin=2em,itemsep=4pt]
\item For $j=i$, the coefficient-derivative term contains the last factor $h_0(n_i,W)=0$ and therefore vanishes.
\item For $j\ne i$, the coefficient-derivative term vanishes. Each product-derivative term also vanishes, since at most one factor is differentiated and one of the remaining factors still pairs $n_j^\flat$ with $n_i$, giving $h_0(n_j,n_i)=0$.
\item For $j=i$, if the differentiated factor occurs in the first $N-1$ slots, the undifferentiated last factor still gives $h_0(n_i,W)=0$. Only the term with the differentiated factor in the last slot remains; all its other factors equal $h_0(n_i,n_i)=1$.
\end{enumerate}
Since $\nabla^{h_0}$ preserves the metric, $\nabla^{h_0}_Xn_i^\flat=h_0(\nabla^{h_0}_Xn_i,\cdot)$, and the sole remaining term is
\[
 (\nabla^{h_0}_X\Theta)(n_i,\ldots,n_i,W)=k!A_i h_0(\nabla^{h_0}_Xn_i,W).
\]
There is no additional factor $N$, since only the last of the $N$ product-rule terms is nonzero. The left-hand side has already been shown to vanish. Using $A_i\ne0$, we obtain
\[
 h_0(\nabla^{h_0}_Xn_i,W)=0\qquad\text{for all }W\perp n_i.
\]
Thus $\nabla^{h_0}_Xn_i$ is parallel to $n_i$. The unit-length condition also gives
\[
 0=Xh_0(n_i,n_i)=2h_0(\nabla^{h_0}_Xn_i,n_i).
\]
Hence $\nabla^{h_0}_Xn_i=0$, so every active axis is parallel along the arc.

The Gauss formula for the unit sphere in $\R^{m+1}$ gives
\[
 \frac{\dd n_i}{\dd s}
 =\nabla^{h_0}_Xn_i-h_0(X,n_i)\gamma(s).
\]
Both terms on the right vanish by $\nabla^{h_0}_Xn_i=0$ and $h_0(X,n_i)=0$. Therefore $\dd n_i/\dd s=0$, and each $n_i(s)$ is a fixed vector in the ambient space.

Let $L$ be the two-dimensional linear plane containing the great circle. For any $s_0\in J$, the orthogonal unit vectors $\gamma(s_0)$ and $X(s_0)$ span $L$. Since $n_i$ is tangent to the sphere, $\langle n_i,\gamma(s_0)\rangle=0$; also $\langle n_i,X(s_0)\rangle=0$. Hence $n_i\perp L$. These fixed vectors are therefore tangent to the sphere at every point of the entire great circle and are parallel along that circle. This proves the lemma.
\end{proof}

This construction is local on $J$ and uses neither compactness nor periodicity. For a tensor on the entire sphere, choose a fixed axis $n=n_i$. Section \ref{sec:fermi} uses its two-sphere, and Section \ref{sec:periodic} will define the scalar $\widehat A$ along the entire great circle by \eqref{eq:hatA}. It equals $A_i$ on $J$; no eigenframe need be extended through possible degeneracies outside the arc.

\subsection{The Codazzi equations in Fermi coordinates}\label{subsec:paper-5-2}\label{sec:fermi}
We compute the Codazzi equations on the two-sphere determined by the fixed axis $n$. Only spherical geometry and the Codazzi property are needed, not analyticity, commutativity, or periodicity.

\begin{lemma}[A Fermi frame near a great circle]\label{lem:fermi-frame}
Let $\gamma:J\to S^m$ be a unit-speed great-circle arc, let $L$ be its two-dimensional linear plane, and let $n\perp L$ be a fixed unit vector. Put $W=L\oplus\R n$. After shortening $J$ if necessary, the map
\begin{equation}\label{eq:fermi-map}
 \Phi(x,s)=\cos x\,\gamma(s)+\sin x\,n,\qquad |x|<\varepsilon<\pi/2,
\end{equation}
gives local coordinates on the totally geodesic two-sphere $S_W^2=S^m\cap W$. Its metric, orthonormal frame, and induced connection satisfy
\[
\begin{gathered}
 h_0=\dd x^2+\cos^2x\,\dd s^2,\qquad
 e_1=\partial_x,\qquad e_2=\sec x\,\partial_s,\\
 \nabla^{h_0}_{e_1}e_1=\nabla^{h_0}_{e_1}e_2=0,\qquad
 \nabla^{h_0}_{e_2}e_1=-\tan x\,e_2,\qquad \nabla^{h_0}_{e_2}e_2=\tan x\,e_1.
\end{gathered}
\]
Along $x=0$, we have $e_1=n$ and $e_2=\gamma'$. The restriction of $\nabla^{h_0}$ to directions tangent to the two-sphere is its Levi-Civita connection.
\end{lemma}
\begin{proof}
Since $n\perp L$,
\[
 W=L\oplus\R n,\qquad S_W^2=W\cap S^m
\]
are a three-dimensional linear space and its unit two-sphere, respectively. The arc $\gamma$ lies in $S_W^2$.

Let $D$ denote the Euclidean connection of the ambient space $\R^{m+1}$. At a point $p$ of the sphere, the Gauss formula is
\[
 \nabla^{h_0}_XY=D_XY+h_0(X,Y)p.
\]
If $X,Y$ are tangent to $S_W^2$, then $Y$ takes values in the fixed linear space $W$, so $D_XY\in W$. Since also $p\in W$, the right-hand side remains in $W$. At the same time, $\nabla^{h_0}_XY$ is tangent to $S^m$, so
\[
 \nabla^{h_0}_XY\in W\cap T_pS^m=T_pS_W^2.
\]
This proves that $S_W^2$ is totally geodesic in $S^m$ and that its induced Levi-Civita connection is the restriction of $\nabla^{h_0}$ to these tangent directions.

Total geodesy makes tangential restriction commute with covariant differentiation, without second-fundamental-form terms. Retaining the notation $C,\nabla^{h_0},h_0$ for the restrictions, complete symmetry of $\nabla^{h_0}C$ therefore gives the Codazzi condition on the two-sphere.

Shorten $J$ suitably and choose $0<\varepsilon<\pi/2$ so that the Fermi map $\Phi:(-\varepsilon,\varepsilon)\times J\to S_W^2$ defined by \eqref{eq:fermi-map} parametrizes a coordinate neighbourhood in the domain of $C$, where all subsequent calculations take place.

For fixed $s$, the curve $x\mapsto\Phi(x,s)$ is the unit-speed great circle starting at $\gamma(s)$ with initial velocity $n$. Thus $x$ is signed arc length in the direction perpendicular to the original great circle, while $s$ is the unit-speed parameter of the original circle.

From the unit-speed parametrization of the great circle and $n\perp L$, we have
\[
 |\gamma|=|\gamma'|=1,\quad
 \langle\gamma,\gamma'\rangle=0,\quad
 \gamma''=-\gamma,\quad
 \langle n,\gamma\rangle=\langle n,\gamma'\rangle=0.
\]
Taking partial derivatives gives
\[
 \Phi_x=-\sin x\,\gamma+\cos x\,n,
 \qquad \Phi_s=\cos x\,\gamma'.
\]
Hence
\[
 \langle\Phi_x,\Phi_x\rangle=1,\qquad
 \langle\Phi_x,\Phi_s\rangle=0,\qquad
 \langle\Phi_s,\Phi_s\rangle=\cos^2x,
\]
and the metric is
\[
 h_0=\dd x^2+\cos^2x\,\dd s^2.
\]
We make the usual identification of coordinate vectors with their images under the differential of $\Phi$, and take the orthonormal frame
\[
 e_1=\partial_x=-\sin x\,\gamma+\cos x\,n,
 \qquad
 e_2=\sec x\,\partial_s=\gamma'.
\]
In particular, on the original great circle $x=0$, we have $e_1=n$ and $e_2=\gamma'$. For any scalar function $u(x,s)$,
\[
 e_1(u)=u_x,\qquad e_2(u)=\sec x\,u_s.
\]
This orthonormal frame is not a coordinate frame, since
\[
 [e_1,e_2]=\partial_x(\sec x)\partial_s
           =\tan x\,e_2.
\]
The variation of the frame itself must be included when computing covariant derivatives of tensor components.

First differentiate along $e_1=\partial_x$. The ambient vector expressions above give
\[
 D_{e_1}e_1
 =\partial_x(-\sin x\,\gamma+\cos x\,n)=-\Phi,
 \qquad
 D_{e_1}e_2=\partial_x\gamma'=0.
\]
The Gauss formula therefore yields
\[
 \nabla^{h_0}_{e_1}e_1=-\Phi+h_0(e_1,e_1)\Phi=0,
 \qquad
 \nabla^{h_0}_{e_1}e_2=0+h_0(e_1,e_2)\Phi=0.
\]
Next differentiate along $e_2=\sec x\,\partial_s$. Since $n$ is fixed,
\[
 D_{e_2}e_1=\sec x\,\partial_s(-\sin x\,\gamma+\cos x\,n)
 =-\tan x\,\gamma'=-\tan x\,e_2.
\]
Here $h_0(e_2,e_1)=0$, so $\nabla^{h_0}_{e_2}e_1=-\tan x\,e_2$. Finally, since $\gamma''=-\gamma$,
\[
 D_{e_2}e_2=\sec x\,\gamma''=-\sec x\,\gamma,
\]
and hence
\[
 \nabla^{h_0}_{e_2}e_2=-\sec x\,\gamma+\Phi
 =-\frac{\sin^2x}{\cos x}\gamma+\sin x\,n
 =\tan x(-\sin x\,\gamma+\cos x\,n)=\tan x\,e_1.
\]
This gives all the connection formulas in the lemma.
Their antisymmetric part also satisfies $\nabla^{h_0}_{e_1}e_2-\nabla^{h_0}_{e_2}e_1=[e_1,e_2]$, as required by the absence of torsion.

\end{proof}

\begin{proposition}[The Codazzi equations in the Fermi frame]\label{prop:fermi-codazzi}
Let $C$ be a smooth symmetric cubic tensor on an open subset of the sphere containing the coordinate neighbourhood above, and suppose that $\nabla^{h_0}C$ is completely symmetric. Put
\[
 a=C(e_1,e_1,e_1),\quad b=C(e_1,e_1,e_2),\quad
 c=C(e_1,e_2,e_2),\quad d=C(e_2,e_2,e_2).
\]
Then on this coordinate neighbourhood the Fermi--Codazzi equations hold:
\begin{align}
 b_x&=\sec x\,a_s+3\tan x\,b,\label{eq:fermi1}\\
 c_x&=\sec x\,b_s-\tan x\,a+2\tan x\,c,\label{eq:fermi2}\\
 d_x&=\sec x\,c_s+\tan x\,d-2\tan x\,b.\label{eq:fermi3}
\end{align}
Neither commutativity of the contraction operators nor vanishing of $C$ on the original great circle is required. If, in addition, $C|_{\gamma(J)}=0$, then all four components vanish at $x=0$.
\end{proposition}
\begin{proof}
By complete symmetry of $C$, these four functions determine its restriction to directions tangent to the two-sphere. Subscripts $x,s$ denote partial derivatives of the scalar functions.

For any vector fields $X,Y,Z,V$, the covariant derivative is defined by
\[
 (\nabla^{h_0}_XC)(Y,Z,V) ={}X\bigl(C(Y,Z,V)\bigr)-C(\nabla^{h_0}_XY,Z,V) -C(Y,\nabla^{h_0}_XZ,V)-C(Y,Z,\nabla^{h_0}_XV).
\]
The last three terms correct for the changes in the vectors inserted into the three slots. Since $\nabla^{h_0}_{e_1}e_1=\nabla^{h_0}_{e_1}e_2=0$, these corrections vanish when differentiating along $e_1$, giving
\[
 (\nabla^{h_0}_{e_1}C)(e_2,e_1,e_1)=b_x,\qquad
 (\nabla^{h_0}_{e_1}C)(e_2,e_1,e_2)=c_x,\qquad
 (\nabla^{h_0}_{e_1}C)(e_2,e_2,e_2)=d_x.
\]

First insert $e_1$ into all three slots. Then
\[
 \begin{aligned}
 (\nabla^{h_0}_{e_2}C)(e_1,e_1,e_1)
 ={}&\sec x\,a_s-C(-\tan x\,e_2,e_1,e_1)\ -C(e_1,-\tan x\,e_2,e_1)-C(e_1,e_1,-\tan x\,e_2)\\
 ={}&\sec x\,a_s+3\tan x\,b.
 \end{aligned}
\]
Each correction is $+\tan x\,b$. Complete symmetry of $\nabla^{h_0}C$ allows us to interchange the differentiation slot and the first tensor slot, giving
\[
 (\nabla^{h_0}_{e_2}C)(e_1,e_1,e_1)
 =(\nabla^{h_0}_{e_1}C)(e_2,e_1,e_1)=b_x.
\]

Next differentiate $C(e_1,e_1,e_2)=b$:
\[
 \begin{aligned}
 (\nabla^{h_0}_{e_2}C)(e_1,e_1,e_2)
 ={}&\sec x\,b_s-C(-\tan x\,e_2,e_1,e_2)\ -C(e_1,-\tan x\,e_2,e_2)-C(e_1,e_1,\tan x\,e_1)\\
 ={}&\sec x\,b_s+2\tan x\,c-\tan x\,a.
 \end{aligned}
\]
The first two slots each contribute $+\tan x\,c$, and the last contributes $-\tan x\,a$. Again, the Codazzi condition gives
\[
 (\nabla^{h_0}_{e_2}C)(e_1,e_1,e_2)
 =(\nabla^{h_0}_{e_1}C)(e_2,e_1,e_2)=c_x.
\]

Finally differentiate $C(e_1,e_2,e_2)=c$:
\[
 \begin{aligned}
 (\nabla^{h_0}_{e_2}C)(e_1,e_2,e_2)
 ={}&\sec x\,c_s-C(-\tan x\,e_2,e_2,e_2)\ -C(e_1,\tan x\,e_1,e_2)-C(e_1,e_2,\tan x\,e_1)\\
 ={}&\sec x\,c_s+\tan x\,d-2\tan x\,b.
 \end{aligned}
\]
The first slot contributes $+\tan x\,d$, and the last two each contribute $-\tan x\,b$. Using
\[
 (\nabla^{h_0}_{e_2}C)(e_1,e_2,e_2)
 =(\nabla^{h_0}_{e_1}C)(e_2,e_2,e_2)=d_x,
\]
we obtain the Fermi--Codazzi equations \eqref{eq:fermi1}--\eqref{eq:fermi3} from these three slot comparisons.

If $C|_{\gamma(J)}=0$, the component definitions give their zero initial values at $x=0$.
\end{proof}

\begin{remark}[Restriction of the Codazzi equation and the full commutativity contraction]\label{rem:fermi-scope}
The first-order equations \eqref{eq:fermi1}--\eqref{eq:fermi3} require only the smooth Codazzi condition. Applied to the arc of zeros, they also have the initial values
\[
 a(0,s)=b(0,s)=c(0,s)=d(0,s)=0.
\]
Expansion in $x$ will give an ordinary differential equation in $s$.

The term $-\tan x\,a$ comes from the correction to the third tensor slot in the connection formula $\nabla^{h_0}_{e_2}e_2=\tan x\,e_1$. If
\[
 a=A(s)x^k+O(x^{k+1}),
\]
then $\tan x=x+O(x^3)$ gives
\[
 -\tan x\,a=-A(s)x^{k+1}+O(x^{k+2}).
\]
This contribution to the leading term of $c$ produces $-(k+1)^2A^2$ in \eqref{eq:arc-ode}. The negative sign ensures that both squared terms have positive coefficients after integration over a period.

Although the Codazzi equation restricts without correction by total geodesy, the commutativity contraction still includes all of $T_pS^m$, not only $e_1,e_2$. Section \ref{sec:arc-orders} estimates the remaining normal contractions.

\end{remark}

\subsection{Order estimates for normal contractions and the scalar equation}\label{subsec:paper-5-3}\label{sec:arc-orders}
We combine the tangential Codazzi equations with the full commutativity contraction to obtain an equation for $A(s)$. The tensor $C$ is still defined on $S^m$; its components $a,b,c,d$ do not include the normal contractions, which must first be estimated.

Retain the interval $J$, the lowest nonzero order $k\ge1$, and the fixed active axes from Section \ref{sec:fixed-axis}. Choose $n=n_1$ and put $A=A_1$. All computations take place in the local Fermi coordinate neighbourhood of Section \ref{sec:fermi}. For an analytic component, $F=O(x^\ell)$ means that, after shrinking the neighbourhood if necessary, $F(x,s)=x^\ell R(x,s)$ with $R$ analytic. Thus differentiation in $s$ does not lower the order of vanishing in $x$.

\begin{lemma}[Orders of normal contractions]\label{lem:fermi-normal-orders}
Let $C$, $k$, $J$, and the fixed active axes satisfy the hypotheses and conclusions of Lemma \ref{lem:fixed-normal-axes}. Take $n=n_1$ and use the coordinates and frame of Lemma \ref{lem:fermi-frame}. Choose a fixed orthonormal basis $w_1,\ldots,w_{m-2}$ of $W^\perp$, and put
\begin{equation}\label{eq:fermi-normal-matrix}
 B_\nu=
 \begin{pmatrix}
 C(w_\nu,e_1,e_1)&C(w_\nu,e_1,e_2)\\
 C(w_\nu,e_1,e_2)&C(w_\nu,e_2,e_2)
 \end{pmatrix}.
\end{equation}
On this coordinate neighbourhood, the full commutativity relation gives
\begin{equation}\label{eq:arc-gauss}
 ac+bd-b^2-c^2=-\sum_\nu\det B_\nu,
\end{equation}
and
\[
 (B_\nu)_{11}=O(x^{k+1}),\qquad
 (B_\nu)_{12},(B_\nu)_{22}=O(x^{k+2}),\qquad
 \det B_\nu=O(x^{2k+3}).
\]
If $m=2$, the normal sum is empty. These estimates do not assert that $B_\nu$ vanishes identically, nor do they raise the order of $(B_\nu)_{11}$ to $k+2$.
\end{lemma}
\begin{proof}
With $W=L\oplus\R n$ and $S_W^2=W\cap S^m$, choose an orthonormal basis of $W^\perp$,
\[
 w_1,\ldots,w_{m-2}.
\]
At any $p\in S_W^2$, we have $w_\nu\perp p$, so $w_\nu$ is tangent to $S^m$. Since $T_pS_W^2\subset W$, it is perpendicular to the tangent plane of the two-sphere. These are therefore \emph{normal vectors to the two-sphere within $S^m$}; they do not include the radial normal $p$ of the sphere. Thus
\[
 \mathcal E=\{e_1,e_2,w_1,\ldots,w_{m-2}\}
\]
is a full orthonormal basis of $T_pS^m$, with the normal part absent when $m=2$.

Consider the normal contraction matrices $B_\nu$ defined by \eqref{eq:fermi-normal-matrix}. Since $[C_{e_1},C_{e_2}]=0$, we have
\[
 0=h_0([C_{e_1},C_{e_2}]e_2,e_1).
\]
Each contraction operator is self-adjoint, so the two products satisfy, respectively,
\[
 \begin{aligned}
 h_0(C_{e_1}C_{e_2}e_2,e_1)
 &=h_0(C_{e_2}e_2,C_{e_1}e_1)
 =\sum_{E\in\mathcal E}C(e_2,e_2,E)C(e_1,e_1,E),\\
 h_0(C_{e_2}C_{e_1}e_2,e_1)
 &=h_0(C_{e_1}e_2,C_{e_2}e_1)
 =\sum_{E\in\mathcal E}C(e_1,e_2,E)^2.
 \end{aligned}
\]
Thus the full contraction identity is
\[
 0=\sum_{E\in\mathcal E}
 \bigl(C(e_1,e_1,E)C(e_2,e_2,E)-C(e_1,e_2,E)^2\bigr).
\]
The contribution from $E=e_1$ is $ac-b^2$, that from $E=e_2$ is $bd-c^2$, and that from $E=w_\nu$ is exactly $\det B_\nu$. Rearranging gives the full commutativity contraction identity \eqref{eq:arc-gauss}.

Fix $s$ and put $\eta_s(x)=\Phi(x,s)$. This is a geodesic with arc-length parameter $x$. The preceding section showed that $\nabla^{h_0}_{e_1}e_1=\nabla^{h_0}_{e_1}e_2=0$. For the fixed ambient vector $w_\nu$, the Gauss formula also gives
\[
 \nabla^{h_0}_{e_1}w_\nu
 =D_{e_1}w_\nu+h_0(e_1,w_\nu)\Phi=0.
\]
Thus $e_1,e_2,w_\nu$ are all parallel along $\eta_s$. Repeated differentiation of components with respect to $x$ has no corrections from variation of the frame.

Write $q=\gamma(s)$ and $X=\gamma'(s)$, and let $E_1=n$, $E_2=X$ denote the initial values of $e_1,e_2$ at $x=0$. For $\alpha,\beta\in\{1,2\}$,
\[
 \partial_x^j(B_\nu)_{\alpha\beta}(0,s)
 =((\nabla^{h_0})^jC)_q(\underbrace{n,\ldots,n}_{j\text{ entries}},w_\nu,E_\alpha,E_\beta).
\]
Since $(\nabla^{h_0})^jC|_{\gamma(J)}=0$ for $j<k$, the first $k$ Taylor coefficients of the three normal components vanish.

Consider now order $k$. Lemma \ref{lem:fixed-normal-axes} gives
\[
 ((\nabla^{h_0})^kC)_q=k!\sum_{i=1}^rA_i(s)(n_i^\flat)^{\otimes(k+3)}.
\]
At least one slot in the contraction above contains $n$, since $k\ge1$. For $i\ne1$, the corresponding pure-power term vanishes because $h_0(n_i,n)=0$. For $i=1$, it vanishes because $h_0(n,w_\nu)=0$. Therefore
\[
 \partial_x^k(B_\nu)_{11}(0,s)
 =\partial_x^k(B_\nu)_{12}(0,s)
 =\partial_x^k(B_\nu)_{22}(0,s)=0.
\]
Thus all three components are $O(x^{k+1})$.

In the affine coordinates centred at $q$, let the active variable $u_1$ correspond to the chosen axis $n=n_1$. By Lemma \ref{lem:next-homogeneous-form}, the next homogeneous term of the local potential is
\[
 f_{k+4}=\sum_i\beta_i u_i^{k+4}
       +\sum_{i,\alpha}c_{i\alpha}u_i^{k+3}z_\alpha
       +G_{k+4}(z).
\]
The only monomials here containing $u_1$ are $u_1^{k+4}$ and $u_1^{k+3}z_\alpha$, whose degree in the other variables is at most one. All other terms are independent of $u_1$.

Section \ref{sec:fixed-axis} showed that, at the origin of these centred coordinates, $(\nabla^{h_0})^{k+1}C$ equals the full differential of the potential of order $k+4$. Combining this with the Taylor-coefficient comparison in the parallel frame along the transverse geodesic gives
\[
 \partial_x^{k+1}(B_\nu)_{12}(0,s) =\partial_n^{\,k+2}\partial_{w_\nu}\partial_X f_{k+4}(0),\qquad \partial_x^{k+1}(B_\nu)_{22}(0,s) =\partial_n^{\,k+1}\partial_{w_\nu}\partial_X^2 f_{k+4}(0).
\]
Here $\partial_n,\partial_{w_\nu},\partial_X$ denote directional derivatives determined by the corresponding constant tangent-space vectors. Both $w_\nu,X$ are perpendicular to $n$. The first expression has two derivatives perpendicular to $n$, and the second has three. Both annihilate the displayed monomials containing $u_1$, since those monomials have degree at most one in the other variables. Terms not containing $u_1$ are annihilated by differentiation along $n$. Thus both derivatives of order $k+1$ vanish.

On the other hand,
\[
 \partial_x^{k+1}(B_\nu)_{11}(0,s)
 =\partial_n^{\,k+3}\partial_{w_\nu}f_{k+4}(0).
\]
This involves only one derivative perpendicular to $n$ and allows a nonzero contribution from $u_1^{k+3}z_\alpha$. We cannot assert that it vanishes. The precise estimates are therefore
\[
 (B_\nu)_{11}=O(x^{k+1}),\qquad
 (B_\nu)_{12}=O(x^{k+2}),\qquad
 (B_\nu)_{22}=O(x^{k+2}).
\]
It follows that
\[
 \det B_\nu=(B_\nu)_{11}(B_\nu)_{22}-(B_\nu)_{12}^{\,2}
 =O(x^{2k+3})-O(x^{2k+4})=O(x^{2k+3}).
\]
The full normal matrix need not vanish, but its determinant does not enter the comparison at degree $2k+2$.

\end{proof}

\begin{proposition}[The scalar equation for the leading coefficient]\label{prop:arc-scalar-equation}
Under the hypotheses of Lemma \ref{lem:fermi-normal-orders}, put $A=A_1$, and let a prime denote differentiation in the great-circle arc-length parameter $s$. The four tangential components satisfy
\[
\begin{aligned}
 a&=Ax^k+O(x^{k+1}),&\qquad
 b&=\frac{A'}{k+1}x^{k+1}+O(x^{k+2}),\\
 c&=\frac{A''-(k+1)A}{(k+1)(k+2)}x^{k+2}+O(x^{k+3}),&
 d&=O(x^{k+3}).
\end{aligned}
\]
Their leading coefficient satisfies, on $J$, the scalar equation
\begin{equation}\label{eq:arc-ode}
 (k+1)AA''-(k+2)(A')^2-(k+1)^2A^2=0.
\end{equation}
This is a local scalar identity. Its derivation uses the full commutativity contraction and does not require the tangential restriction to the two-sphere to satisfy the commutativity relation on its own.
\end{proposition}
\begin{proof}
The frame is parallel along $\eta_s$, so the lower derivatives of $a=C(e_1,e_1,e_1)$ are given by the corresponding covariant derivatives of $C$. The choice of the lowest nonzero order and the leading-term decomposition give
\[
 \partial_x^j a(0,s)=0\quad(j<k),\qquad
 \partial_x^k a(0,s)
 =((\nabla^{h_0})^kC)_q(n,\ldots,n)=k!A(s).
\]
Therefore
\[
 a(x,s)=A(s)x^k+O(x^{k+1}),\qquad A(s)\ne0\quad(s\in J).
\]
The arc of zeros also gives $b(0,s)=c(0,s)=d(0,s)=0$. Primes on functions of $s$ denote differentiation with respect to $s$.

Write \eqref{eq:fermi1} as
\[
 b_x-3\tan x\,b=\sec x\,a_s.
\]
Since $(\cos^3x)'=-3\tan x\cos^3x$, multiplication by the integrating factor $\cos^3x$ gives
\[
 (\cos^3x\,b)_x=\cos^2x\,a_s.
\]
Integrating from $0$ to $x$ and using $b(0,s)=0$, we obtain the exact formula
\[
 b(x,s)=\sec^3x\int_0^x\cos^2t\,a_s(t,s)\,\dd t.
\]
From $a_s(t,s)=A'(s)t^k+O(t^{k+1})$ and $\cos^2t=1+O(t^2)$, the integrand is
\[
 \cos^2t\,a_s(t,s)=A'(s)t^k+O(t^{k+1}).
\]
After integration and multiplication by $\sec^3x=1+O(x^2)$, this gives
\[
 b(x,s)=\frac{A'}{k+1}x^{k+1}+O(x^{k+2}).
\]

Equation \eqref{eq:fermi2} can be written as
\[
 c_x-2\tan x\,c=\sec x\,b_s-\tan x\,a.
\]
Multiplying by the integrating factor $\cos^2x$, we obtain
\[
 (\cos^2x\,c)_x=\cos x\,b_s-\sin x\cos x\,a.
\]
Since $c(0,s)=0$, the exact integral formula is
\[
 c(x,s)=\sec^2x\int_0^x
 \bigl(\cos t\,b_s(t,s)-\sin t\cos t\,a(t,s)\bigr)\,\dd t.
\]
Differentiating the expansion already obtained for $b$ with respect to $s$ gives
\[
 b_s(t,s)=\frac{A''}{k+1}t^{k+1}+O(t^{k+2}).
\]
Expand the two terms in the integrand separately:
\[
 \cos t\,b_s(t,s)=\frac{A''}{k+1}t^{k+1}+O(t^{k+2}),\qquad
 \sin t\cos t\,a(t,s)=A\,t^{k+1}+O(t^{k+2}).
\]
Thus the integrand is
\[
 \left(\frac{A''}{k+1}-A\right)t^{k+1}+O(t^{k+2}).
\]
Integration followed by multiplication by $\sec^2x=1+O(x^2)$ gives
\[
 c(x,s)=\frac{A''-(k+1)A}{(k+1)(k+2)}
          x^{k+2}+O(x^{k+3}).
\]
Here $-(k+1)A$ comes precisely from the term $-\tan x\,a$ in the Codazzi equation.

Finally, write \eqref{eq:fermi3} as
\[
 d_x-\tan x\,d=\sec x\,c_s-2\tan x\,b.
\]
Multiplication by the integrating factor $\cos x$ gives
\[
 (\cos x\,d)_x=c_s-2\sin x\,b.
\]
Using $d(0,s)=0$, we have
\[
 d(x,s)=\sec x\int_0^x
       \bigl(c_s(t,s)-2\sin t\,b(t,s)\bigr)\,\dd t.
\]
Since $c_s(t,s)=O(t^{k+2})$ and $\sin t\,b(t,s)=O(t^{k+2})$, integration gives $d=O(x^{k+3})$. The four tangential components therefore have the expansions
\[
 \begin{aligned}
 a&=A\,x^k+O(x^{k+1}),&\qquad
 b&=\frac{A'}{k+1}x^{k+1}+O(x^{k+2}),\\
 c&=\frac{A''-(k+1)A}{(k+1)(k+2)}x^{k+2}+O(x^{k+3}),&
 d&=O(x^{k+3}).
 \end{aligned}
\]

\begin{samepage}
Lemma \ref{lem:fermi-normal-orders} and the expansions of $a,b,c,d$ give the following orders:
\begin{center}
\begin{tabular}{cc}
\toprule
Term & Lower bound for the vanishing order in $x$\\
\midrule
$ac$ & $2k+2$\\
$b^2$ & $2k+2$\\
$bd$ & $2k+4$\\
$c^2$ & $2k+4$\\
$\det B_\nu$ & $2k+3$\\
\bottomrule
\end{tabular}
\end{center}
\end{samepage}
Thus only $ac-b^2$ contributes to the coefficient of degree $2k+2$ in the full commutativity equation \eqref{eq:arc-gauss}. Multiplying the corresponding expansions gives
\[
 ac=\frac{A\bigl(A''-(k+1)A\bigr)}{(k+1)(k+2)}x^{2k+2}+O(x^{2k+3}),\qquad
 b^2=\frac{(A')^2}{(k+1)^2}x^{2k+2}+O(x^{2k+3}).
\]
The coefficient of degree $2k+2$ must therefore satisfy
\[
 \frac{A\bigl(A''-(k+1)A\bigr)}{(k+1)(k+2)}
 -\frac{(A')^2}{(k+1)^2}=0.
\]
Multiplying by $(k+1)^2(k+2)$ and expanding gives the scalar equation \eqref{eq:arc-ode} for the leading coefficient.
No division by $A$ was made, so the equation remains meaningful at its zeros, as required for analytic continuation.

\end{proof}

\begin{remark}[The local scope of the scalar equation]\label{rem:arc-ode-scope}
Equation \eqref{eq:arc-ode} is local on $J$; neither compactness nor periodicity has been used. The leading and next homogeneous terms were obtained analytically, but once these forms and the vanishing of lower covariant derivatives are known, the calculations require only smoothness to the corresponding orders. We next extend the identity to the entire great circle.

\end{remark}

\subsection{Periodicity and global real-analytic rigidity}\label{subsec:paper-5-4}\label{sec:periodic}
For $C$ real-analytic on the entire sphere, we extend the local coefficient $A$ and equation \eqref{eq:arc-ode} from $J$ to the whole great circle. Global definition and analyticity justify this continuation; the subsequent integration over a period requires only $C^2$ regularity.

\begin{proposition}[Analytic periodic continuation of the local scalar equation]\label{prop:arc-analytic-periodic-extension}
Let $C$ be a real-analytic symmetric cubic tensor on the entire standard sphere $S^m$, satisfying the spherical Codazzi--commutation system \eqref{eq:sphere-system}. Suppose that a great-circle arc on which the entire tensor vanishes, together with its shortened interval $J$, satisfies the conclusions of Lemma \ref{lem:fixed-normal-axes}. Choose one of its fixed unit axes $n=n_1$ and the corresponding leading coefficient $A=A_1$. Extend $\gamma$ with unit arc-length parameter to the entire great circle. Then
\begin{equation}\label{eq:hatA}
 \widehat A(s):=\frac1{k!}((\nabla^{h_0})^kC)_{\gamma(s)}
       (\underbrace{n,\ldots,n}_{k+3\text{ entries}}),\qquad s\in\R,
\end{equation}
is a real-analytic $2\pi$-periodic function, equals $A$ on $J$, and satisfies \eqref{eq:arc-ode} on the whole real axis with $A$ replaced by $\widehat A$.
Only this scalar function and its equation are continued. The active-axis decomposition need not remain valid outside $J$.
\end{proposition}
\begin{proof}
Choose $s_0\in J$ and put
\[
 p=\gamma(s_0),\qquad v=\gamma'(s_0).
\]
Since $\gamma$ is a unit-speed great-circle arc, $p,v$ form an orthonormal basis of its plane $L$. Its full extension is
\[
 \gamma(s)=\cos(s-s_0)p+\sin(s-s_0)v,
 \qquad s\in\R.
\]
This agrees with $\gamma$ on the original arc and satisfies $\gamma(s+2\pi)=\gamma(s)$.

The vector $n$ obtained in Section \ref{sec:fixed-axis} is a fixed ambient unit vector with $n\perp L$. Hence, for all $s\in\R$,
\[
 \langle n,\gamma(s)\rangle=0,\qquad
 \langle n,\gamma'(s)\rangle=0.
\]
The first equality says that $n\in T_{\gamma(s)}S^m$. The second, together with the Gauss formula, gives
\[
 \nabla^{h_0}_{\gamma'}n
 =D_{\gamma'}n+\langle\gamma',n\rangle\gamma=0.
\]

The scalar $\widehat A$ defined by \eqref{eq:hatA} is globally defined because $(\nabla^{h_0})^kC$ exists on the whole sphere and $n$ is tangent along the circle. No active-axis decomposition outside $J$ is required.

The standard spherical connection is analytic, so $(\nabla^{h_0})^kC$ is analytic. Since $\gamma$ is analytic and $n$ is fixed, $\widehat A$ is real-analytic. Moreover, $\gamma(s+2\pi)=\gamma(s)$ and the inserted vector is always the same $n$, so
\[
 \widehat A(s+2\pi)=\widehat A(s).
\]
On the initial interval $J$, the leading-term decomposition gives
\[
 ((\nabla^{h_0})^kC)_{\gamma(s)}(n,\ldots,n)
 =k!\sum_iA_i(s)h_0(n_i,n)^{k+3}=k!A(s).
\]
Thus $\widehat A|_J=A$. In particular, since $A$ is nonzero on $J$, $\widehat A$ is not identically zero.

Define the function given by its left-hand side,
\[
 R(s)=(k+1)\widehat A(s)\widehat A''(s)
       -(k+2)\bigl(\widehat A'(s)\bigr)^2
       -(k+1)^2\widehat A(s)^2.
\]
Since $\widehat A$ is real-analytic on the whole real axis, so is $R$. The preceding section proved that $A$ satisfies \eqref{eq:arc-ode}. Since $\widehat A=A$ on the open interval $J$, all their derivatives there agree, and hence $R|_J=0$. The real axis is connected, so the real-analytic identity theorem gives
\[
 R(s)=0\qquad\text{for all }s\in\R.
\]

This uses the analytic identity theorem, not an ODE continuation theorem: the polynomial in $\widehat A,\widehat A',\widehat A''$ involves no division by $\widehat A$, so zeros outside $J$ do not affect analyticity of $R$.

\end{proof}

The next lemma requires neither spherical geometry nor real analyticity.

\begin{lemma}[Vanishing for the periodic scalar equation]\label{lem:periodic-scalar-vanishing}
Let $k\ge1$ be an integer and let $u:\R\to\R$ be a $C^2$ function with period $T>0$. If
\[
 (k+1)uu''-(k+2)(u')^2-(k+1)^2u^2=0,
\]
then $u\equiv0$.
\end{lemma}
\begin{proof}
Integrate the equation over $[0,T]$. Periodicity gives $u(T)=u(0)$ and $u'(T)=u'(0)$, so the boundary term in integration by parts vanishes:
\[
 \int_0^Tuu''\,\dd s
 =[uu']_0^T-\int_0^T(u')^2\,\dd s
 =-\int_0^T(u')^2\,\dd s.
\]
Substituting into the integrated equation and rearranging gives
\[
 (2k+3)\int_0^T(u')^2\,\dd s
 +(k+1)^2\int_0^Tu^2\,\dd s=0.
\]
Both coefficients are positive and both integrands are nonnegative, so $\int_0^Tu^2\,\dd s=0$. By continuity, $u$ vanishes throughout one period, and periodicity gives $u\equiv0$.

\end{proof}

Proposition \ref{prop:arc-analytic-periodic-extension} and Lemma \ref{lem:periodic-scalar-vanishing} give $\widehat A\equiv0$, contradicting $\widehat A|_J=A\ne0$. Together with the local zero subsphere theorem, this yields the global result.

\begin{theorem}[Global rigidity of real-analytic cubic forms]\label{thm:global-cubic}
Let $m\ge2$, let $h_0$ be the metric on the entire standard unit sphere $S^m$, and let $\nabla^{h_0}$ be its Levi-Civita connection. If a globally defined real-analytic symmetric cubic tensor $C$ satisfies the spherical Codazzi--commutation system \eqref{eq:sphere-system}, where $h_0(C_XY,Z)=C(X,Y,Z)$, then $C\equiv0$.
\end{theorem}
\begin{proof}
By Lemma \ref{lem:spherical-cubic-zero-point}, there is a point $p\in S^m$ at which $C_p=0$.

If $C$ vanishes in an open neighbourhood of $p$, connectedness of the sphere and the identity theorem for real-analytic tensors immediately give $C\equiv0$.

If $C$ is not identically zero in any neighbourhood of $p$, its analytic germ at $p$ is nonzero. Theorem \ref{thm:local-zero} gives a great-circle arc on which the entire cubic tensor $C$ vanishes. Since $C$ is also not identically zero on the sphere, Lemma \ref{lem:fixed-normal-axes} applies and gives a finite lowest nonzero order $k$, a shortened interval $J$, a fixed axis $n$, and an analytic coefficient $A$ that is nowhere zero.

Proposition \ref{prop:fermi-codazzi} and Lemma \ref{lem:fermi-normal-orders} supply the tangential equations and normal estimates for Proposition \ref{prop:arc-scalar-equation}, giving \eqref{eq:arc-ode} on $J$. Proposition \ref{prop:arc-analytic-periodic-extension} continues the scalar coefficient and its equation around the great circle, where Lemma \ref{lem:periodic-scalar-vanishing} forces it to vanish, contradicting $A\ne0$ on $J$. Thus the nonzero-germ case is impossible: $C$ vanishes near $p$, and the identity theorem gives $C\equiv0$.
\end{proof}

Simple connectivity and compactness supply the global trace potential and its extremum in the zero-point lemma. The entire great circle supplies the periodic function $\widehat A$. Analyticity is needed for the local zero-set and leading-term analysis and for continuation of the identity; the final $C^2$ periodic argument does not extend those preceding steps to the smooth category.

\begin{theorem}[Analytic Laugwitz rigidity]\label{thm:laugwitz}
Let $n\ge3$, and let $F$ be a Minkowski norm on $\R^n$, regular, strongly convex, and real-analytic on the entire $\R^n\setminus\{0\}$. Reversibility is not assumed. If the Hessian metric $g=\Hess(F^2/2)$ is flat on the entire slit space, then there is a positive-definite symmetric bilinear form $q$ such that
\[
 F(y)^2=q(y,y)\qquad(y\in\R^n).
\]
\end{theorem}
\begin{proof}
Put $E=F^2/2$, $\Sigma=\{F=1\}$, and $h=g|_{T\Sigma}$. By the spherical reduction cited in Proposition \ref{prop:flat-hessian-sphere} and its discussion of analyticity, flatness of $g$ implies that the compact indicatrix $(\Sigma,h)$ is analytically isometric to the entire standard sphere $(S^{n-1},h_0)$. Under this isometry the angular Cartan tensor becomes a global analytic cubic tensor satisfying the spherical Codazzi--commutation system \eqref{eq:sphere-system}. Since $n-1\ge2$, Theorem \ref{thm:global-cubic} applies and gives $C=0$.

By radial vanishing and homogeneity in Proposition \ref{prop:minkowski-euler}, the vanishing of $C$ on $\Sigma$ implies $\mathcal C=0$ on the entire slit space $\R^n\setminus\{0\}$. Proposition \ref{prop:cartan-zero-quadratic} therefore gives a positive-definite symmetric bilinear form $q$ such that $F(y)^2=q(y,y)$ for all $y\in\R^n$.
\end{proof}

\section{Twisted Codazzi quartic forms and maximal-rank branches}\label{sec:paper-6}
Let $(V,F)$ be a smooth, strongly convex Minkowski space regular on the entire slit space, and put
\[
 E=\frac12F^2,\qquad g=\Hess E,\qquad
 \Sigma=\{y\in V:F(y)=1\},\qquad m=\dim\Sigma.
\]
We retain the notation of Section \ref{sec:paper-2}: $h=g|_{T\Sigma}$ is the indicatrix metric, $C$ the angular Cartan cubic form, $\nabla^h$ the Levi-Civita connection of $h$, and $\nabla$ the induced centroaffine connection. Let $B$ be a smooth, fully symmetric quartic tensor field on all of $\Sigma$. This section and the next require only smoothness of $F$ and $B$.

\begin{definition}[Twisted Codazzi quartic form]\label{def:twisted-quartic}
Let $B\in\Gamma^\infty(\Sym^4T^*\Sigma)$. If, in any local coordinates on the indicatrix,
\begin{equation}\label{eq:twisted}
 (\nabla^h_eB)_{aijk}-(\nabla^h_aB)_{eijk}
 =C_{ai}{}^\ell B_{e\ell jk}-C_{ei}{}^\ell B_{a\ell jk},
\end{equation}
then we call $B$ a \textbf{twisted Codazzi quartic form} with respect to $(h,C)$. Here $e,a,i,j,k$ are tangential indices on the indicatrix, $C_{ai}{}^\ell=h^{\ell p}C_{aip}$, and all internal contractions are taken with respect to $h$. Equivalently, for any tangent vector fields $e,a,i,j,k$,
\begin{equation}\label{eq:twisted-vector}
 (\nabla^h_eB)(a,i,j,k)-(\nabla^h_aB)(e,i,j,k) =B\bigl(e,C^\sharp(a,i),j,k\bigr) -B\bigl(a,C^\sharp(e,i),j,k\bigr),
\end{equation}
where $h(C^\sharp(X,Y),Z)=C(X,Y,Z)$.
\end{definition}

\begin{lemma}[Algebraic compatibility of the twisted Codazzi equation]\label{lem:quartic-algebraic-compatibility}
If a fully symmetric quartic form $B$ satisfies \eqref{eq:twisted}, then it satisfies
\begin{equation}\label{eq:F}
 \mathcal F_{eaijk}:={}C_{ei}{}^\ell B_{a\ell jk} -C_{ai}{}^\ell B_{e\ell jk} -C_{ej}{}^\ell B_{a\ell ik} +C_{aj}{}^\ell B_{e\ell ik}=0.
\end{equation}
Thus \eqref{eq:F} is an algebraic consequence of the twisted Codazzi equation and requires no separate hypothesis.
\end{lemma}
\begin{proof}
Covariant differentiation preserves the symmetry in the arguments of the tensor being differentiated. In particular, since $B_{aijk}=B_{ajik}$ and $B_{eijk}=B_{ejik}$, we have
\[
 (\nabla^h_eB)_{aijk}=(\nabla^h_eB)_{ajik},\qquad
 (\nabla^h_aB)_{eijk}=(\nabla^h_aB)_{ejik}.
\]
The left-hand side of \eqref{eq:twisted} is therefore unchanged when $i,j$ are interchanged. Comparing the right-hand sides gives
\[
 C_{ai}{}^\ell B_{e\ell jk}-C_{ei}{}^\ell B_{a\ell jk}
 =
 C_{aj}{}^\ell B_{e\ell ik}-C_{ej}{}^\ell B_{a\ell ik}.
\]
Moving the right-hand side to the left and multiplying the resulting expression by $-1$ gives precisely \eqref{eq:F}.
\end{proof}

Assume henceforth that $B$ satisfies Definition \ref{def:twisted-quartic}. For $B\not\equiv0$, put $q=\max_{\Sigma}\rank B$. We treat the degenerate branch $0<q<m$ by radical transport and closed orbits; the full-rank branch $q=m$ is treated in the next section.

\subsection{The radical, active space, and centroid}\label{subsec:paper-7-1}\label{sec:quartic-rank-centroid}
We first introduce the radical, the active space, and the centroid, on which the operators $A_v,T_v,S_v$ will be defined.

On a fixed Minkowski indicatrix $\Sigma$, for each $y\in\Sigma$, let
\[
 E_y=T_y\Sigma,\qquad h_y=g_y|_{E_y},\qquad
 B_y\in\Sym^4(E_y^*).
\]
The following decomposition and centroid arguments use only positive definiteness of $h$, full symmetry of $B$, and smoothness.

\begin{definition}[Insertion map, radical, and rank of a quartic form]\label{def:quartic-radical-rank}
The \textbf{insertion map} of the quartic form $B_y$ is the linear map
\[
 \iota_y:E_y\longrightarrow\Sym^3(E_y^*),\qquad
 \iota_y(v):=B_y(v,\cdot,\cdot,\cdot).
\]
Its \textbf{radical} is defined by
\[
 N_y:=\rad B_y=\ker\iota_y
 =\{v\in E_y:B_y(v,a,b,c)=0\ \text{for all }a,b,c\in E_y\}.
\]
The \textbf{rank} of a quartic form will always mean the rank of its insertion map:
\[
 \rank B_y:=\rank\iota_y=m-\dim N_y.
\]
We say that $B_y$ \textbf{has trivial radical} if $N_y=\{0\}$. This is equivalent to the injectivity of $\iota_y$, and also to $\rank B_y=m$.
\end{definition}

By rank--nullity, $\rank\iota_y=m-\dim\ker\iota_y$. The condition $B_y(v,v,v,v)=0$ alone does not imply $v\in N_y$: all other arguments must be arbitrary. By symmetry, a radical vector in any slot annihilates $B_y$.

\begin{remark}
This rank is not the matrix rank obtained by regarding $B_y$ as a map $\Sym^2E_y\to\Sym^2E_y^*$. The dimensions of the radical and active space, and the rank branches considered below, all refer to the insertion-map rank in Definition \ref{def:quartic-radical-rank}.
\end{remark}

The rank takes values in the finite set $\{0,1,\ldots,m\}$, and $\Sigma$ is nonempty. We may therefore define
\begin{equation}\label{eq:maximal-rank-set}
 q:=\max_{y\in\Sigma}\rank B_y,\qquad
 U:=\{y\in\Sigma:\rank B_y=q\}.
\end{equation}
We call $U$ the \textbf{maximal-rank set}. Its rank is locally constant, giving the following decomposition.

We use the lower semicontinuity of rank and the kernel and image theorem for constant-rank vector bundle homomorphisms \cite[Theorem 10.34]{lee2013}.
\begin{proposition}[Smooth decomposition on the maximal-rank open set]\label{prop:maximal-rank-bundles}
The maximal-rank set $U$ is nonempty and open. On $U$, the spaces $N_y=\ker\iota_y$ form a smooth vector subbundle $N\subset T\Sigma|_U$ of rank $m-q$. Their $h$-orthogonal complements form a smooth subbundle of rank $q$, giving the orthogonal direct-sum decomposition
\[
 T\Sigma|_U=N\oplus N^\perp.
\]
\end{proposition}

\begin{definition}[Active space and restricted quartic form]\label{def:quartic-active-space}
On the maximal-rank open set $U$, put
\[
 P:=N^\perp,\qquad E_y=N_y\oplus P_y,\qquad \dim P_y=q.
\]
We call $P_y$ the \textbf{active space} of $B_y$, and $P$ the \textbf{active subbundle}. Let $\pi_{P_y}:E_y\to P_y$ denote the orthogonal projection, and define the \textbf{restricted quartic form}
\[
 b_y:=B_y|_{P_y\times P_y\times P_y\times P_y}\in\Sym^4(P_y^*).
\]
\end{definition}

The following lemma identifies the active restriction with the form induced on the quotient by the radical.

\begin{lemma}[Triviality of the radical after reduction]\label{lem:quartic-active-reduction}
For every $y\in U$ and $X_1,X_2,X_3,X_4\in E_y$, we have
\[
 B_y(X_1,X_2,X_3,X_4)
 =b_y(\pi_{P_y}X_1,\pi_{P_y}X_2,\pi_{P_y}X_3,\pi_{P_y}X_4).
\]
Moreover, $\rad b_y=\{0\}$. The inner product $h_y$ identifies $P_y$ linearly with the quotient space $E_y/N_y$, and $b_y$ is the representation on $P_y$ of the quartic form with trivial radical induced by $B_y$ on this quotient.
\end{lemma}
\begin{proof}
Write each argument uniquely as
\[
 X_r=n_r+p_r,\qquad n_r\in N_y,\quad p_r\in P_y,
 \qquad r=1,2,3,4.
\]
In the multilinear expansion, all terms containing any $n_r$ vanish. The sole remaining term is $b_y(p_1,p_2,p_3,p_4)$, proving the identity.

If $v\in\rad b_y$, then $v\in P_y$ and $b_y(v,p_2,p_3,p_4)=0$ for all $p_2,p_3,p_4\in P_y$. By the identity just proved, for any $X_2,X_3,X_4\in E_y$,
\[
 B_y(v,X_2,X_3,X_4)
 =b_y(v,\pi_{P_y}X_2,\pi_{P_y}X_3,\pi_{P_y}X_4)=0.
\]
Thus $v\in N_y$. But $v$ also belongs to $P_y=N_y^\perp$, so $v=0$. This proves $\rad b_y=\{0\}$.

To explain the quotient description, define
\[
 \overline B_y([X_1],[X_2],[X_3],[X_4])
 :=B_y(X_1,X_2,X_3,X_4),\qquad [X_r]\in E_y/N_y.
\]
Adding a radical vector to a representative does not change this value. The map $P_y\to E_y/N_y$, $p\mapsto[p]$, is injective since $P_y\cap N_y=0$ and surjective since $X=n+p$ has class $[p]$. It pulls $\overline B_y$ back to $b_y$, so the induced form has trivial radical. Only the choice of orthogonal representatives depends on $h_y$.
\end{proof}

For $q=0$, $B\equiv0$ and $P=0$, and the restricted quartic form on the zero-dimensional space $P$ has trivial radical; for $q=m$, $N=0$ and $P=T\Sigma$. The degenerate branch is $0<q<m$. In the centroid discussion below we assume $q>0$.

\begin{remark}
In this section, $N$ denotes the radical distribution of the quartic form, and $P=N^\perp$ denotes the active subbundle. The operator $T_v$ in the next subsection is the splitting operator of the radical distribution; it is unrelated to $T_{abc}=f_{abc}$ in Section \ref{sec:paper-3}.
\end{remark}

Fix $y\in U$ and write $P=P_y$, $b=b_y$; the next arguments are pointwise algebraic.

\begin{definition}[Centroid of a quartic form]\label{def:quartic-centroid}
The \textbf{centroid} of the quartic form $b$ is defined by
\[
 \Cent(b):=\{A\in\operatorname{End}_{\R}(P):\; b(Au_1,u_2,u_3,u_4)=b(u_1,Au_2,u_3,u_4)\qquad \text{for all }u_1,u_2,u_3,u_4\in P\}.
\]
Here $\operatorname{End}_{\R}(P)$ denotes the space of real linear endomorphisms of $P$. The definition requires that $A$ can be moved from the first argument slot to the second. It involves only $b$ and imposes no self-adjointness condition on $A$ with respect to $h|_P$.
\end{definition}

Full symmetry permits movement between all slots. Closure under composition remains to be proved.

\begin{lemma}[Moving centroid operators between slots]\label{lem:quartic-centroid-slots}
If $A\in\Cent(b)$, then for any $u_1,u_2,u_3,u_4\in P$,
\[
 b(Au_1,u_2,u_3,u_4) =b(u_1,Au_2,u_3,u_4) =b(u_1,u_2,Au_3,u_4) =b(u_1,u_2,u_3,Au_4).
\]
This property requires only full symmetry of $b$, not triviality of its radical.
\end{lemma}
\begin{proof}
The first equality is the definition. Interchanging the second and third arguments and then applying the definition gives
\[
 b(Au_1,u_2,u_3,u_4) =b(Au_1,u_3,u_2,u_4) =b(u_1,Au_3,u_2,u_4) =b(u_1,u_2,Au_3,u_4).
\]
The outer equalities use symmetry and the middle one the definition. Interchanging the second and fourth slots instead gives
\[
 b(Au_1,u_2,u_3,u_4)=b(u_1,u_2,u_3,Au_4).
\]
Thus the four expressions are equal. Reading these equalities in reverse permits $A$ to be moved between any two slots.
\end{proof}

\begin{proposition}[The centroid of a quartic form with trivial radical]\label{prop:quartic-centroid-algebra}
Let $P\ne0$ be a finite-dimensional real vector space, and let $b\in\Sym^4(P^*)$ satisfy $\rad b=0$. Then $\Cent(b)$ is a commutative real subalgebra of $\operatorname{End}_{\R}(P)$, with operator composition as multiplication and $\Id$ as identity. In particular, if $A\in\Cent(b)$ and $p\in\R[t]$, then $p(A)\in\Cent(b)$.
\end{proposition}
\begin{proof}
We first verify the linear structure. If $A,D\in\Cent(b)$ and $\alpha,\beta\in\R$, multilinearity and the definition give
\[
 \begin{aligned}
 b((\alpha A+\beta D)u,v,w,z) &=\alpha b(Au,v,w,z)+\beta b(Du,v,w,z)\\
 &=\alpha b(u,Av,w,z)+\beta b(u,Dv,w,z) =b(u,(\alpha A+\beta D)v,w,z).
 \end{aligned}
\]
Hence $\alpha A+\beta D\in\Cent(b)$. The identity $\Id$ satisfies the defining equality, so $\Id\in\Cent(b)$.

We prove commutativity before closure under multiplication. For arbitrary $u,v,w,z\in P$, Lemma \ref{lem:quartic-centroid-slots} allows successive movements of $A,D$, giving
\[
 b(ADu,v,w,z) =b(Du,Av,w,z) =b(u,Av,Dw,z) =b(Au,v,Dw,z) =b(DAu,v,w,z).
\]
Only individual movements of $A,D$ between slots were used, without assuming closure under multiplication. Consequently,
\[
 b((AD-DA)u,v,w,z)=0\qquad\text{for all }v,w,z\in P.
\]
Since $\rad b=0$, we obtain $(AD-DA)u=0$. As $u$ is arbitrary, $AD=DA$. This is where triviality of the radical is used.

Using the commutativity just proved, we now verify closure under multiplication:
\[
 b(ADu,v,w,z) =b(DAu,v,w,z) =b(Au,Dv,w,z) =b(u,ADv,w,z).
\]
Thus $AD\in\Cent(b)$. Associativity is inherited from composition. Since every power $A^j$ and $A^0=\Id$ belongs to this real algebra, so does every $p(A)$ with $p\in\R[t]$.
\end{proof}

The proposition gives commutativity and polynomial closure, but no spectral conclusion relative to an inner product.

\subsection{Geometric equations for the radical}\label{subsec:paper-7-2}\label{sec:radical-eq}
In the preceding subsection we obtained, on the maximal-rank open set $U$, the smooth decomposition
\[
 T\Sigma|_U=N\oplus P,\qquad N=\operatorname{rad}B,\qquad P=N^\perp,
\]
with $b=B|_{P^4}$ of trivial radical. We now insert radical and active arguments into \eqref{eq:F} and \eqref{eq:twisted} to derive the geometric equations.

The proofs apply on any constant-rank open set, where $N,P$ are smooth, without analyticity or compactness. We use the maximal-rank set $U$ for subsequent applications.

We retain $C^\sharp$ from Definition \ref{def:twisted-quartic} and the vector form \eqref{eq:twisted-vector} of the twisted Codazzi equation. The algebraic compatibility identity \eqref{eq:F} has the vector form
\begin{equation}\label{eq:quartic-compatibility-vector}
 0={}B\bigl(a,C^\sharp(e,i),j,k\bigr) -B\bigl(e,C^\sharp(a,i),j,k\bigr) -B\bigl(a,C^\sharp(e,j),i,k\bigr) +B\bigl(e,C^\sharp(a,j),i,k\bigr).
\end{equation}
Here $a,e,i,j,k$ are vector arguments. We use radical vanishing and injectivity of $W\mapsto b(W,\cdot,\cdot,\cdot)$ on $P$.

\begin{lemma}[Mixed Cartan components in radical directions]\label{lem:root-cartan-mixed}
On a constant-rank open set,
\begin{equation}\label{eq:CNNP}
 C(N,N,P)=0,
\end{equation}
and $C^\sharp(N,N)\subset N$, $C^\sharp(N,P)\subset P$.
\end{lemma}
\begin{proof}
In \eqref{eq:F}, take
\[
 e=u\in N,\qquad j=v\in N,\qquad a,i,k\in P.
\]
The first, second and fourth terms contain a radical argument in $B$ and vanish, leaving
\[
 0=-B\bigl(a,C^\sharp(u,v),i,k\bigr)
   =-b\bigl(a,(C^\sharp(u,v))_P,i,k\bigr).
\]
Since $a,i,k$ may vary arbitrarily in $P$, triviality of the radical gives
\[
 (C^\sharp(u,v))_P=0.
\]
For every $X\in P$, we consequently have $C(u,v,X)=h(C^\sharp(u,v),X)=0$. This proves \eqref{eq:CNNP}.
The identity may also be expressed by the two inclusions
\[
 C^\sharp(N,N)\subset N,\qquad C^\sharp(N,P)\subset P.
\]
The first is the conclusion just obtained; the second follows from full symmetry of $C$. Indeed, for $u,v\in N$ and $X\in P$,
\[
 h\bigl(C^\sharp(v,X),u\bigr)
 =C(v,X,u)=C(v,u,X)=0.
\]
Notice that \eqref{eq:CNNP} does not assert $C(N,P,P)=0$. These components are precisely those recorded by the operator $A_v$ below.
\end{proof}

\begin{definition}[Cartan contraction in a radical direction]\label{def:root-cartan-contraction}
For $y\in U$ and $v\in N_y$, define the \textbf{Cartan contraction operator} on $P_y$ by
\[
 A_vX:=(C^\sharp(v,X))_P,\qquad X\in P_y.
\]
By Lemma \ref{lem:root-cartan-mixed}, $C^\sharp(v,X)$ already belongs to $P_y$. The operator depends smoothly on $y$ and is real linear in $v$; it records the components $C(N,P,P)$ that have not yet been eliminated.
\end{definition}

\begin{proposition}[Centroid property and self-adjointness of Cartan contractions]\label{prop:root-cartan-centroid}
For every $y\in U$ and $v\in N_y$, we have $A_v\in\Cent(b_y)$, and $A_v$ is self-adjoint with respect to $h_y|_{P_y}$.
\end{proposition}
\begin{proof}
In \eqref{eq:F}, take $e=v\in N$ and $a,i,j,k\in P$. The second and fourth terms vanish because they contain the radical argument $v$, giving
\[
 b(a,A_vi,j,k) =b(a,A_vj,i,k) =b(a,i,A_vj,k).
\]
Thus $A_v$ may be moved between the second and third argument slots. Full symmetry of $b$ permits movement between the other slots as well. More explicitly,
\[
 b(A_va,i,j,k) =b(j,A_va,i,k) =b(j,a,A_vi,k) =b(a,A_vi,j,k).
\]
This is the defining condition of the centroid, so $A_v\in\Cent(b)$.

On the other hand, for any $X,Y\in P$,
\[
 h(A_vX,Y) =C(v,X,Y) =C(v,Y,X) =h(X,A_vY).
\]
Thus $A_v$ is self-adjoint by symmetry of $C$, whereas its centroid property follows from \eqref{eq:F}.
\end{proof}

\begin{proposition}[Autoparallelness of the radical distribution and local radical leaves]\label{prop:root-autoparallel}
The radical distribution is autoparallel with respect to the Levi-Civita connection $\nabla^h$:
\begin{equation}\label{eq:autoparallel}
 \nabla^h_NN\subset N.
\end{equation}
Consequently, $N$ is integrable, its local integral leaves are totally geodesic with respect to $h$, and $\nabla^h$-parallel translation along radical directions preserves $N$ and $P$.
\end{proposition}
\begin{proof}
We now turn to covariant derivatives. Choose smooth local vector fields
\[
 u,v\in\Gamma(N),\qquad a,j,k\in\Gamma(P).
\]
Putting $e=u$ and $i=v$ in \eqref{eq:twisted}, we obtain
\[
 (\nabla^h_uB)(a,v,j,k)-(\nabla^h_aB)(u,v,j,k) =B\bigl(u,C^\sharp(a,v),j,k\bigr) -B\bigl(a,C^\sharp(u,v),j,k\bigr)=0.
\]
The first term on the right contains $u\in N$, and the second contains $C^\sharp(u,v)\in N$ by Lemma \ref{lem:root-cartan-mixed}; both vanish.

We expand the two covariant derivatives on the left separately. For the term with two radical arguments,
\[
 (\nabla^h_aB)(u,v,j,k) ={}a\bigl(B(u,v,j,k)\bigr) -B(\nabla^h_au,v,j,k)-B(u,\nabla^h_av,j,k) -B(u,v,\nabla^h_aj,k)-B(u,v,j,\nabla^h_ak).
\]
The scalar function is identically zero, and every correction term retains a radical argument. Hence
\[
 (\nabla^h_aB)(u,v,j,k)=0.
\]
For the term with only one radical argument,
\[
 \begin{aligned}
 (\nabla^h_uB)(a,v,j,k) ={}&u\bigl(B(a,v,j,k)\bigr) -B(\nabla^h_ua,v,j,k)-B(a,\nabla^h_uv,j,k)\\
 &-B(a,v,\nabla^h_uj,k)-B(a,v,j,\nabla^h_uk) ={}-B(a,\nabla^h_uv,j,k).
 \end{aligned}
\]
All terms except that involving $\nabla^h_uv$ retain the radical argument $v$; no radical property of $\nabla^h_uv$ has yet been assumed. Thus
\[
 b\bigl(a,(\nabla^h_uv)_P,j,k\bigr)=0.
\]
Since $a,j,k$ are arbitrary and $b$ has trivial radical, $(\nabla^h_uv)_P=0$. This is \eqref{eq:autoparallel}.
This property is called autoparallelness of the radical distribution with respect to $\nabla^h$.

Since $\nabla^h$ is torsion-free, for $u,v\in\Gamma(N)$ we have
\[
 [u,v]=\nabla^h_uv-\nabla^h_vu\in\Gamma(N).
\]
The Frobenius theorem implies that $N$ is integrable. If $L$ is a local integral leaf, then $TL=N|_L$ and its normal bundle is $P|_L$. The second fundamental form of the leaf satisfies
\[
 \mathrm{II}(u,v)=(\nabla^h_uv)_P=0.
\]
Thus $L$ is totally geodesic.

The orthogonal complement $P$ is also preserved by differentiation along radical directions. Indeed, for $u,v\in\Gamma(N)$ and $X\in\Gamma(P)$, the identity $\nabla^h h=0$ gives
\[
 0=u\bigl(h(X,v)\bigr) =h(\nabla^h_uX,v)+h(X,\nabla^h_uv) =h(\nabla^h_uX,v).
\]
Hence $\nabla^h_uX\in\Gamma(P)$. In other words, $\nabla^h$-parallel translation along a radical leaf preserves the decomposition $N\oplus P$.
\end{proof}

\begin{remark}
Total geodesicity of local radical leaves does not imply their completeness or that a geodesic stays in the constant-rank set. Extension within that set will be proved separately.
\end{remark}

\begin{definition}[Two splitting operators of the radical distribution]\label{def:root-splitting-operators}
For a local radical vector field $v\in\Gamma(N)$ and an active vector field $X\in\Gamma(P)$, define
\[
 T_vX:=-(\nabla^h_Xv)_P,\qquad S_v:=T_v+A_v.
\]
We call $T_v$ the \textbf{splitting operator} with respect to $\nabla^h$. We shall prove below that $S_vX=-(\nabla_Xv)_P$, so $S_v$ is the splitting operator with respect to the centroaffine connection $\nabla$. These operators describe the active components of the transverse variation of radical directions. The minus sign is our convention.
\end{definition}

\begin{lemma}[Tensoriality of the splitting operators]\label{lem:root-splitting-tensor}
The preceding definitions are independent of the local extensions of the radical and active vectors and define smooth tensors
\[
 T,S\in\Gamma\bigl(N^*\otimes\operatorname{End}(P)\bigr).
\]
In particular, $v\mapsto T_v$ and $v\mapsto S_v$ are real linear maps on each radical space.
\end{lemma}
\begin{proof}
To check dependence only on pointwise values of $v,X$, let $f$ be smooth. Then
\[
 T_{fv}X =-\bigl(\nabla^h_X(fv)\bigr)_P =-\bigl(X(f)v+f\nabla^h_Xv\bigr)_P =fT_vX.
\]
Since $v\in N$, the $P$-projection of $X(f)v$ is zero. On the other hand, a connection is linear over functions in the differentiation direction, so
\[
 T_v(fX)=fT_vX.
\]
Expanding in smooth local frames of $N$ and $P$ now shows that $T$ is a smooth tensor:
\[
 T\in\Gamma\bigl(N^*\otimes\operatorname{End}(P)\bigr).
\]
Thus, for $v\in N_y$, $T_v$ is indeed a linear operator on $P_y$ independent of the choice of extensions. Self-adjointness of $T_v$ has not been established.

Since $A$ is a smooth tensor linear in $v$, $S=T+A$ has the same tensoriality. No derivatives of a chosen extension of $v$ are prescribed.
\end{proof}

\begin{proposition}[Evolution of the quartic form along radical directions]\label{prop:root-quartic-flow}
If $v\in\Gamma(N)$ and $a,i,j,k\in\Gamma(P)$, then
\begin{equation}\label{eq:Bflow}
 (\nabla^h_vB)(a,i,j,k)=B((T_v-A_v)a,i,j,k).
\end{equation}
\end{proposition}
\begin{proof}
In \eqref{eq:twisted}, take $e=v\in\Gamma(N)$ and the other arguments $a,i,j,k$ in $\Gamma(P)$. The first term on the right contains the radical argument $v$. Writing the second in terms of $A_v$ gives
\[
 (\nabla^h_vB)(a,i,j,k)-(\nabla^h_aB)(v,i,j,k)
 =-B(a,A_vi,j,k).
\]
Expand the covariant derivative in the middle term:
\[
 \begin{aligned}
 (\nabla^h_aB)(v,i,j,k) ={}&a\bigl(B(v,i,j,k)\bigr) -B(\nabla^h_av,i,j,k)-B(v,\nabla^h_ai,j,k)\ -B(v,i,\nabla^h_aj,k)-B(v,i,j,\nabla^h_ak)\\
 ={}&-B(\nabla^h_av,i,j,k) ={}-B((\nabla^h_av)_P,i,j,k) ={}B(T_va,i,j,k).
 \end{aligned}
\]
As in Proposition \ref{prop:root-autoparallel}, the scalar derivative is zero and the $i,j,k$ correction terms retain $v$. The last equality uses $T_va=-(\nabla^h_av)_P$.

Substituting into the twisted Codazzi equation gives
\[
 (\nabla^h_vB)(a,i,j,k)
 =B(T_va,i,j,k)-B(a,A_vi,j,k).
\]
Since $A_v\in\Cent(b)$, it can be moved from the second argument slot to the first:
\[
 B(a,A_vi,j,k)=B(A_va,i,j,k).
\]
This proves \eqref{eq:Bflow}.
\end{proof}

\begin{proposition}[Centroaffine splitting and commutation of centroid operators]\label{prop:root-centro-splitting}
At the same point, all operators $A_v,T_v,S_v$ ($v\in N_y$) belong to $\Cent(b_y)$ and hence commute with one another. Moreover,
\[
 S_vX=-(\nabla_Xv)_P,\qquad \nabla_NN\subset N,\qquad \nabla_NP\subset P.
\]
Thus $\nabla$-parallel translation along local radical leaves preserves the decomposition $N\oplus P$. This commutativity does not assert that $T_v$ or $S_v$ is self-adjoint.
\end{proposition}
\begin{proof}
Since $B$ is fully symmetric in its four arguments, so is $\nabla^h_vB$. Interchanging $a$ and $i$ in \eqref{eq:Bflow} therefore gives
\[
 b((T_v-A_v)a,i,j,k) =(\nabla^h_vB)(a,i,j,k) =(\nabla^h_vB)(i,a,j,k) =b((T_v-A_v)i,a,j,k) =b(a,(T_v-A_v)i,j,k).
\]
This says precisely that $T_v-A_v\in\Cent(b)$. Proposition \ref{prop:root-cartan-centroid} gives $A_v\in\Cent(b)$, and closure of the centroid under linear combinations then yields
\[
 T_v\in\Cent(b),\qquad
 S_v:=T_v+A_v\in\Cent(b).
\]
Commutativity of $\Cent(b)$ gives commutation of $A_v,T_v,S_v$ with $A_w,T_w,S_w$ at the same point. It does not imply self-adjointness of $T_v$ or $S_v$.

The geometric meaning of $S_v=T_v+A_v$ comes from the centroaffine connection. By \eqref{eq:centro},
\[
 \nabla_Xv=\nabla^h_Xv-C^\sharp(X,v).
\]
Taking the $P$-component and changing the sign gives
\[
 -(\nabla_Xv)_P =-(\nabla^h_Xv)_P+(C^\sharp(X,v))_P =T_vX+A_vX =S_vX.
\]
Thus the difference $A_v$ between the two splitting operators in Definition \ref{def:root-splitting-operators} comes from the difference $\nabla-\nabla^h=-C^\sharp$ between the connections.

Moreover, Lemma \ref{lem:root-cartan-mixed} and Proposition \ref{prop:root-autoparallel} have shown that
\[
 \nabla^h_NN\subset N,\quad C^\sharp(N,N)\subset N,
 \qquad
 \nabla^h_NP\subset P,\quad C^\sharp(N,P)\subset P.
\]
Subtracting the two pairs of relations gives
\[
 \nabla_NN\subset N,\qquad \nabla_NP\subset P.
\]
Thus $\nabla$ also preserves $N,P$ along radical leaves. Its splitting operator $S_v$ will enter the curvature calculation in Section \ref{sec:riccati}.

If $P=0$, all the transverse operators vanish and commutativity is immediate. If $P\ne0$, we use Proposition \ref{prop:quartic-centroid-algebra}. Thus the conclusions also include the trivial cases of zero rank and full rank.
\end{proof}

Autoparallelness keeps a geodesic with radical initial velocity in a radical leaf locally. Theorem \ref{thm:centro-periodic} extends the ambient centroaffine geodesic to all real times; extension inside the maximal-rank set is still to be proved in Section \ref{sec:extension}.

\subsection{The Riccati equation and cancellation after raising an index}\label{subsec:paper-7-4}\label{sec:riccati}
We derive the Riccati equation for $S_v$ from centroaffine curvature and the evolution equation for the raised tensor $\mathcal B$ from \eqref{eq:Bflow} and $\nabla h=2C$. The Cartan terms from changing connections and raising an index will cancel.

Choose $y_0\in U$ and $0\ne v_0\in N_{y_0}$. The corresponding $\nabla$-geodesic $y(t)$ is complete and periodic by Subsection \ref{sec:centro-periodic}. Work on an interval $I$ containing the initial time where it stays in $U$; autoparallelness and geodesic uniqueness keep its velocity in $N$. Put
\[
 v(t)=\dot y(t)\in N_{y(t)},\qquad
 \nabla_vv=0,\qquad a(t)=h(v,v)>0.
\]
Along radical directions, both $\nabla^h$ and $\nabla$ preserve $N$ and $P$. Denote their induced connections on $P$ along radical directions by $\nabla^{h,P}$ and $\nabla^P$, respectively:
\[
 \nabla^{h,P}_vX=\nabla^h_vX,\qquad
 \nabla^P_vX=\nabla_vX,
 \qquad v\in\Gamma(N),\quad X\in\Gamma(P).
\]
These connections induce derivatives on operators and tensors on $P$ only in directions in $N$; elsewhere we retain $\nabla^h,\nabla$. We use the operators $A_v,T_v,S_v$ from Definitions \ref{def:root-cartan-contraction} and \ref{def:root-splitting-operators}, together with their geometric expressions in Proposition \ref{prop:root-centro-splitting}. All three operators belong to $\Cent(b)$, and $A_v$ is self-adjoint. No self-adjointness assumption on $S_v$ will be made.

\begin{definition}[Splitting operator and connection along radical directions on the quotient bundle]\label{def:root-quotient-connection}
Let $Q:=(T\Sigma|_U)/N$ be the quotient bundle on $U$, and let $\overline X$ denote the quotient class of a local vector field $X$. For $v\in\Gamma(N)$, put
\[
 S_v\overline X:=-\overline{\nabla_Xv},\qquad
 \nabla^{Q}_v\overline X:=\overline{\nabla_vX}.
\]
The second expression defines the connection $\nabla^{Q}$ along radical directions induced on $Q$ by the centroaffine connection $\nabla$; it is defined here only for $v\in N$.
\end{definition}
We first check that both expressions are independent of the representative.

\begin{lemma}[Well-definedness of the quotient-bundle expressions]\label{lem:root-quotient-well-defined}
The two expressions in Definition \ref{def:root-quotient-connection} are independent of the representative of $X$. Under the identification of $Q$ with $P$, they correspond respectively to $S_v$ on $P$ and the induced connection $\nabla^P$ along radical directions.
\end{lemma}
\begin{proof}
Replacing $X$ by $X+u$, $u\in\Gamma(N)$, changes the expressions by $-\overline{\nabla_uv}$ and $\overline{\nabla_vu}$, both zero since $\nabla_NN\subset N$. Orthogonal representatives in $P$ give $S_vX=-(\nabla_Xv)_P$ and $\nabla^P_v$.

Linearity over functions and the product rule for the connection also give
\[
 \nabla^{Q}_{fv}\overline X=f\nabla^{Q}_v\overline X,
 \qquad
 \nabla^{Q}_v(f\overline X)=v(f)\overline X+f\nabla^{Q}_v\overline X.
\]
These are the two required rules for a connection along radical directions.
\end{proof}

To set $\nabla_X(\nabla_vv)=0$ in the curvature calculation, we need $\nabla_vv=0$ on a neighbourhood, not merely along the curve. The following extension supplies this.

\begin{lemma}[Local geodesic extension of the velocity of a radical geodesic]\label{lem:root-geodesic-field}
Near each point of a radical $\nabla$-geodesic with nonzero velocity, its velocity extends to a smooth nonzero radical vector field $v$ such that $\nabla_vv=0$ throughout the neighbourhood and the original curve within that neighbourhood is an integral curve of $v$.
\end{lemma}
\begin{proof}
The extension may be constructed as follows. Near a point of the curve, take a local hypersurface $H$ transverse to the velocity there. Choose a smooth nonzero radical vector field along $H$ whose value at that point is the prescribed velocity. Starting from each point of $H$, take the $\nabla$-geodesic with the chosen radical initial velocity for a short time. Since the radical distribution is integrable and autoparallel for $\nabla$, these geodesics lie locally in their respective radical leaves. Denote the resulting map by $\Phi(t,z)$, where $z\in H$. At the initial point, $\partial_t\Phi$ and $T H$ span the whole tangent space. The inverse function theorem therefore shows that $\Phi$ gives local coordinates after the neighbourhood is reduced. Set
\[
 v(\Phi(t,z))=\partial_t\Phi(t,z).
\]
This gives the required geodesic vector field in radical directions. The original geodesic is the member with the prescribed initial data.

\end{proof}
We carry out the following calculation in such a neighbourhood. Repeating it at each point of the curve gives the equation on the whole interval $I$.

\begin{proposition}[Riccati equation for the centroaffine splitting operator]\label{prop:root-riccati}
Along the radical $\nabla$-geodesic within the current constant-rank interval,
\begin{equation}\label{eq:riccati}
 \nabla^P_vS_v=S_v^2+a(t)\Id.
\end{equation}
The left-hand side is the covariant derivative along the curve of an operator on $P$, induced by $\nabla^P$, rather than an ordinary derivative in a frame chosen in advance.
\end{proposition}
\begin{proof}
Choose $v$ as in Lemma \ref{lem:root-geodesic-field}. In its local flow coordinates, choose a vector field $X$ whose components are independent of the flow parameter. Then
\[
 [v,X]=0.
\]
The value of $X$, hence its quotient class, is arbitrary. It need not lie in $P$, since the calculation is in the quotient bundle.

Since $\nabla$ is torsion-free,
\[
 \nabla_vX-\nabla_Xv=[v,X]=0,
\]
and hence
\[
 \nabla^{Q}_v\overline X
 =\overline{\nabla_vX}
 =\overline{\nabla_Xv}
 =-S_v\overline X.
\]
With our curvature convention,
\[
 R^{\nabla}(v,X)v
 =\nabla_v\nabla_Xv-\nabla_X\nabla_vv-\nabla_{[v,X]}v.
\]
The last two terms vanish because $\nabla_vv$ vanishes throughout the neighbourhood and $[v,X]=0$. Passing to quotient classes and using the product rule for the covariant derivative of an operator gives
\[
 \overline{R^{\nabla}(v,X)v} =\nabla^{Q}_v\bigl(\overline{\nabla_Xv}\bigr) =\nabla^{Q}_v(-S_v\overline X) =-(\nabla^{Q}_vS_v)\overline X -S_v(\nabla^{Q}_v\overline X) =-(\nabla^{Q}_vS_v)\overline X+S_v^2\overline X.
\]
On the other hand, $R^{\nabla}=R_0$ in \eqref{eq:centro} gives
\[
 R^{\nabla}(v,X)v=h(X,v)v-h(v,v)X.
\]
The first term is a multiple of $v$, so it belongs to $N$ and vanishes in the quotient bundle. Thus
\[
 \overline{R^{\nabla}(v,X)v}=-a(t)\overline X.
\]
Comparing the two expressions and using the arbitrariness of $\overline X$, we obtain
\[
 -\nabla^{Q}_vS_v+S_v^2=-a(t)\Id.
\]
Representing the quotient bundle by $P$ identifies $\nabla^{Q}$ with $\nabla^P$ and gives \eqref{eq:riccati}.
The left-hand side is the covariant derivative of an operator on $P$, meaning explicitly
\[
 (\nabla^P_vS_v)X=\nabla^P_v(S_vX)-S_v(\nabla^P_vX),\qquad X\in\Gamma(P).
\]
The square term comes from $-S_v(\nabla^{Q}_v\overline X)=S_v^2\overline X$; moving the curvature term $-a(t)\Id$ gives the positive sign of $a(t)\Id$.
\end{proof}

\begin{lemma}[Derivative of the quartic form under a change of connection]\label{lem:root-covariant-quartic-change}
If $X,Y,Z,W\in\Gamma(P)$, then
\[
 (\nabla^P_vb)(X,Y,Z,W)=b((S_v+2A_v)X,Y,Z,W).
\]
The term $2A_v$ includes contributions from both the original evolution equation and the connection corrections in the four covariant argument slots.
\end{lemma}
\begin{proof}
All vector arguments $X,Y,Z,W$ below belong to $P$. By \eqref{eq:Bflow},
\[
 (\nabla^{h,P}_vb)(X,Y,Z,W)
 =(\nabla^h_vB)(X,Y,Z,W)
 =b((T_v-A_v)X,Y,Z,W).
\]
Moreover, along radical directions,
\[
 \nabla^P_vX=\nabla^{h,P}_vX-A_vX.
\]
The definition of the derivative of a covariant tensor of order four gives
\[
 \begin{aligned}
 (\nabla^P_vb)(X,Y,Z,W) ={}&(\nabla^{h,P}_vb)(X,Y,Z,W) +b(A_vX,Y,Z,W)\\
 &+b(X,A_vY,Z,W) +b(X,Y,A_vZ,W) +b(X,Y,Z,A_vW).
 \end{aligned}
\]
The four correction terms are positive because the covariant-slot minus sign combines with $\nabla^P-\nabla^{h,P}=-A_v$ in direction $v$.

Since $A_v\in\Cent(b)$, in each of the last four terms $A_v$ can be moved to the first argument slot. Hence
\[
 \begin{aligned}
 (\nabla^P_vb)(X,Y,Z,W) &=b((T_v-A_v)X,Y,Z,W) +4b(A_vX,Y,Z,W)\\
 &=b((T_v+3A_v)X,Y,Z,W) =b((S_v+2A_v)X,Y,Z,W).
 \end{aligned}
\]
The last step uses $S_v=T_v+A_v$. When an index is raised below, the remaining $2A_v$ will cancel the term arising from differentiation of the metric.
\end{proof}

\begin{definition}[The quartic form with an index raised]\label{def:root-raised-quartic}
Raise the last index of $b$ with $h|_P$ to define the $P$-valued symmetric cubic tensor $\mathcal B$:
\[
 h(\mathcal B(X,Y,Z),W)=b(X,Y,Z,W).
\]
Positive definiteness of $h|_P$ determines the output vector uniquely. The tensor $\mathcal B$ is a vector-valued cubic tensor, whereas $b$ is a real-valued quartic form. Differentiation of the former must include the connection on its output vector.
\end{definition}

\begin{proposition}[Exact cancellation of the Cartan terms after raising an index]\label{prop:root-raised-cancellation}
Along a radical $\nabla$-geodesic, the tensor with an index raised satisfies
\begin{equation}\label{eq:raised}
 (\nabla^P_v\mathcal B)(X,Y,Z)=\mathcal B(S_vX,Y,Z).
\end{equation}
\end{proposition}
\begin{proof}
The proof uses $\nabla h=2C$, self-adjointness of $A_v$, and $A_v\in\Cent(b)$.
By the definition for a vector-valued tensor, its covariant derivative is
\[
 (\nabla^P_v\mathcal B)(X,Y,Z) ={}\nabla^P_v\bigl(\mathcal B(X,Y,Z)\bigr) -\mathcal B(\nabla^P_vX,Y,Z) -\mathcal B(X,\nabla^P_vY,Z) -\mathcal B(X,Y,\nabla^P_vZ).
\]
Since $\nabla^P$ does not preserve $h|_P$, the tensor product rule includes its derivative:
\[
 (\nabla^P_vb)(X,Y,Z,W) ={}(\nabla^P_v(h|_P))(\mathcal B(X,Y,Z),W) +h((\nabla^P_v\mathcal B)(X,Y,Z),W).
\]
The slot corrections have already been subtracted, so no parallelness of $X,Y,Z,W$ is assumed.

Since $\nabla h=2C$,
\[
 (\nabla^P_v(h|_P))(U,W)
 =(\nabla_vh)(U,W)
 =2C(v,U,W)=2h(A_vU,W),\qquad U,W\in P.
\]
Combining this with the formula in Lemma \ref{lem:root-covariant-quartic-change}, we obtain
\[
 h((\nabla^P_v\mathcal B)(X,Y,Z),W) ={}b((S_v+2A_v)X,Y,Z,W) -2h(A_v\mathcal B(X,Y,Z),W).
\]
By the definition of $\mathcal B$ and nondegeneracy of $h$, this is equivalent to
\[
 (\nabla^P_v\mathcal B)(X,Y,Z)
 =\mathcal B((S_v+2A_v)X,Y,Z)
  -2A_v\mathcal B(X,Y,Z).
\]

It remains to show that applying $A_v$ to the output vector is the same as applying it to the first argument. For any $W\in P$, we calculate
\[
 h(A_v\mathcal B(X,Y,Z),W) =h(\mathcal B(X,Y,Z),A_vW) =b(X,Y,Z,A_vW) =b(A_vX,Y,Z,W) =h(\mathcal B(A_vX,Y,Z),W).
\]
The first equality uses self-adjointness of $A_v$ with respect to $h$, and the third uses $A_v\in\Cent(b)$. Nondegeneracy of $h$ gives
\[
 A_v\mathcal B(X,Y,Z)=\mathcal B(A_vX,Y,Z).
\]
Consequently,
\[
 (\nabla^P_v\mathcal B)(X,Y,Z) =\mathcal B(S_vX,Y,Z) +2\mathcal B(A_vX,Y,Z) -2A_v\mathcal B(X,Y,Z) =\mathcal B(S_vX,Y,Z).
\]
This proves \eqref{eq:raised}.
Thus the $+2A_v$ from the connection change and the original evolution equation cancels the $-2A_v$ from differentiating the metric, using self-adjointness and the centroid property.
\end{proof}

\begin{lemma}[Parallel frames of the active space are fixed ambient vectors]\label{lem:root-fixed-active-frame}
Along the current radical geodesic interval, every $\nabla^P$-parallel frame $E_1,\ldots,E_q$ consists of constant vectors in the ambient vector space $V$. Thus
\[
 P_{y(t)}=P_0:=P_{y(0)}\subset V
\]
is one fixed subspace throughout the interval; the matrix of the restricted metric may still vary with $t$.
\end{lemma}
\begin{proof}
Put $q=\dim P$, and choose along the curve a $\nabla^P$-parallel frame in $P$,
\[
 E_1(t),\ldots,E_q(t),\qquad \nabla^P_vE_\alpha=0.
\]
Such a frame exists because $\nabla$ preserves $P$ along radical directions. Since $v\in N$ and $E_\alpha\in P$, we also have $h(v,E_\alpha)=0$. The centroaffine Gauss formula therefore gives
\[
 D_vE_\alpha
 =\nabla_vE_\alpha-h(v,E_\alpha)y
 =\nabla^P_vE_\alpha=0.
\]
Thus $E_\alpha(t)$ remains constant as a vector in the fixed vector space $V$:
\[
 E_\alpha(t)=E_\alpha(0).
\]
Consequently, on the current interval $I$, all $P_{y(t)}$, regarded as linear subspaces of $V$, equal
\[
 P_0:=P_{y_0}=\operatorname{span}\{E_1(0),\ldots,E_q(0)\}.
\]
This gives an explicit identification of $P$ along the curve with the fixed vector space $P_0$.

Constancy of the ambient vectors does not imply constancy of the metric matrix. Indeed,
\[
 \frac{\dd}{\dd t}h(E_\alpha,E_\beta) =(\nabla^P_v(h|_P))(E_\alpha,E_\beta) =2C(v,E_\alpha,E_\beta),
\]
need not vanish; the metric derivative in Proposition \ref{prop:root-raised-cancellation} is therefore necessary.
\end{proof}

\begin{definition}[Insertion map of the tensor with an index raised]\label{def:root-matrix-insertion}
Identify the active spaces with $P_0$ by the fixed ambient frame, and define
\[
 M(t):P_0\longrightarrow\operatorname{Hom}(\Sym^2P_0,P_0),\qquad
 (M(t)X)(Y,Z):=\mathcal B(t)(X,Y,Z).
\]
Here $\operatorname{Hom}(\Sym^2P_0,P_0)$ denotes the space of symmetric bilinear maps with values in $P_0$. Let $S(t)$ denote the matrix of $S_{v(t)}$ in the same fixed frame.
\end{definition}

\begin{proposition}[Equation for the injective map in a fixed frame]\label{prop:root-matrix-system}
In the fixed frame above, $M(t)$ is injective and
\begin{equation}\label{eq:matrixode}
 S'=S^2+a(t)\Id,\qquad M'=MS,
\end{equation}
where $MS$ means applying $S$ first and then $M$.
\end{proposition}
\begin{proof}
Write $S_v$ as the matrix $S(t)$ in the fixed frame above. If
\[
 S_vE_\beta=\sum_{\alpha=1}^q S^\alpha{}_\beta(t)E_\alpha,
\]
then $\nabla^P_vE_\alpha=0$ gives
\[
 (\nabla^P_vS_v)E_\beta
 =\sum_{\alpha=1}^q\frac{\dd S^\alpha{}_\beta}{\dd t}E_\alpha.
\]
Thus the covariant derivative of the operator is the ordinary matrix derivative in this frame.

We next prove that $M(t)$ is injective. Indeed,
\[
 \begin{aligned}
 X\in\ker M(t)
 &\ \Longleftrightarrow\
 \mathcal B(t)(X,Y,Z)=0
 \quad\text{for all }Y,Z\in P_0\\
 &\ \Longleftrightarrow\
 b(t)(X,Y,Z,W)=0
 \quad\text{for all }Y,Z,W\in P_0
 \ \Longleftrightarrow\ X\in\operatorname{rad}b(t).
 \end{aligned}
\]
The second equivalence uses nondegeneracy of $h$. Since $b(t)$ has trivial radical, $\ker M(t)=0$.

For fixed vectors $X,Y,Z\in P_0$, the corresponding vector fields along the curve are $\nabla^P$-parallel. Hence \eqref{eq:raised} gives
\[
 (M'(t)X)(Y,Z) =(\nabla^P_v\mathcal B)(X,Y,Z) =\mathcal B(S_vX,Y,Z) =(M(t)(S(t)X))(Y,Z).
\]
Thus $M'=M\circ S=MS$, with $S$ acting first on the inserted argument. Together with \eqref{eq:riccati}, this gives \eqref{eq:matrixode}.
\end{proof}

These identities require only smoothness on a constant-rank set. By Subsection \ref{sec:centro-periodic}, $y(t)$ and the positive periodic coefficient $a(t)$ exist globally, with a uniform positive lower bound. The geometric $S,M$ remain defined only while the curve stays in $U$.

Injectivity of $M$ and maximality of rank will control this interval; the fixed ambient frame will allow comparison with the globally smooth tensor.

\subsection{Extension of radical geodesics within the maximal-rank region}\label{subsec:paper-7-5}\label{sec:extension}
We now prove extension within the maximal-rank set, which does not follow from completeness of the centroaffine geodesic alone.

The explicit solution remains geometric before its denominator first becomes singular: any finite endpoint reached before that time retains maximal rank.

Retain the maximal insertion rank $q$ and maximal-rank open set $U$ defined in \eqref{eq:maximal-rank-set}, and assume $0<q<m$. Choose $y_0\in U$ and $0\ne v_0\in N_{y_0}$, let $y(t)$ be the corresponding complete $\nabla$-geodesic, and put
\[
 v(t)=\dot y(t),\qquad a(t)=h(v(t),v(t)).
\]
Autoparallelness and uniqueness keep $v(t)$ radical while $y(t)\in U$. In the fixed ambient parallel frame of the preceding subsection, all active spaces equal $P_0=P_{y_0}$, and $S,M$ satisfy \eqref{eq:matrixode} with initial values
\[
 S(0)=S_0,\quad M(0)=M_0,
\]
where $M_0$ is injective.

\begin{definition}[Scalar fundamental solutions and denominator matrix]\label{def:root-denominator-matrix}
Let $f,g$ be the two \textbf{standard fundamental solutions} of $z''+a(t)z=0$ with
\[
 f(0)=1,\quad f'(0)=0,\qquad g(0)=0,\quad g'(0)=1.
\]
For the initial splitting operator $S_0$, define the \textbf{denominator matrix}
\begin{equation}\label{eq:root-denominator}
 Y(t):=f(t)\Id-g(t)S_0.
\end{equation}
The complete geodesic defines $a$ on the whole real line. Hence the solutions $f,g$ of the linear equation, and therefore $Y$, are also defined on the whole real line. This does not presuppose that the curve remains in $U$.
\end{definition}

\begin{lemma}[The Wronskian and the linear equation for the denominator]\label{lem:root-denominator-wronskian}
The fundamental solutions satisfy $fg'-f'g=1$, so they never vanish simultaneously. Moreover,
\[
 Y''+a(t)Y=0,\qquad Y(0)=\Id,\qquad Y'(0)=-S_0.
\]
\end{lemma}
\begin{proof}
Their Wronskian
\[
 W(t)=f(t)g'(t)-f'(t)g(t)
\]
satisfies
\[
 W' =f'g'+fg''-f''g-f'g' =f(-ag)-(-af)g=0.
\]
Since $W(0)=1$, we obtain
\[
 fg'-f'g=1.
\]
Thus $f,g$ never vanish simultaneously.

For the denominator matrix $Y$ defined in \eqref{eq:root-denominator}, the initial values of the fundamental solutions and the scalar equation give
\[
 Y(0)=\Id,\qquad Y'(0)=-S_0,\qquad Y''=f''\Id-g''S_0 =-a(t)f\Id+a(t)gS_0=-a(t)Y.
\]
Only the globally defined $a(t)$ and the initial $S_0$ enter this construction; no extension within $U$ has been assumed.
\end{proof}

\begin{proposition}[Denominator representation of the matrix system]\label{prop:root-mobius-solution}
On an interval containing $0$ where $Y$ is invertible, the solution of the matrix initial value problem \eqref{eq:matrixode} is
\begin{equation}\label{eq:mobius}
 S=-Y'Y^{-1},\qquad M=M_0Y^{-1}.
\end{equation}
Consequently, on the intersection of this interval with the interval of definition of the geometric matrices, the formula represents the original $S,M$. Neither self-adjointness nor diagonalizability of $S_0$ is required.
\end{proposition}
\begin{proof}
On an interval containing $0$ where $Y$ is invertible, define the candidate matrices $\widetilde S,\widetilde M$ by the respective right-hand sides of \eqref{eq:mobius}. We verify the equations and initial values before identifying these functions with the geometric matrices.

Differentiating $YY^{-1}=\Id$ gives
\[
 Y'Y^{-1}+Y(Y^{-1})'=0,
 \qquad
 (Y^{-1})'=-Y^{-1}Y'Y^{-1}.
\]
Hence
\[
 \widetilde S' =-Y''Y^{-1}-Y'(Y^{-1})' =-Y''Y^{-1}+Y'Y^{-1}Y'Y^{-1} =a(t)\Id+\widetilde S^2,
\]
and
\[
 \widetilde S(0)=-Y'(0)Y(0)^{-1}=S_0.
\]
We also calculate
\[
 \widetilde M' =M_0(Y^{-1})' =-M_0Y^{-1}Y'Y^{-1} =(M_0Y^{-1})(-Y'Y^{-1}) =\widetilde M\widetilde S,
\]
and $\widetilde M(0)=M_0$. The calculation preserves the order of matrix multiplication and requires neither self-adjointness nor diagonalizability of $S_0$.

Uniqueness for ordinary differential equations gives \eqref{eq:mobius} on the common interval where the geometric $S,M$ are defined and $Y$ is invertible.
\end{proof}

\begin{definition}[First singular time and maximal-rank residence interval]\label{def:root-first-singular-time}
Define the \textbf{first positive singular time} of the denominator matrix by
\[
 t_*:=\inf\{t>0:\det Y(t)=0\}.
\]
If there is no positive singular time, set $t_*=\infty$. Let the right endpoint of the \textbf{maximal-rank residence interval} from the initial data be
\[
 \tau:=\sup\{T>0:y([0,T))\subset U\}.
\]
Since $Y(0)=\Id$ and $U$ is open, $t_*,\tau>0$. If $t_*<\infty$, continuity of $\det Y$ gives $\det Y(t_*)=0$, and $Y$ is invertible throughout $[0,t_*)$.
\end{definition}

\begin{lemma}[The fixed space and tensor limit at a finite endpoint]\label{lem:root-finite-endpoint-limit}
Suppose $0<t_0<\infty$, and that $y(t)$ lies in $U$ with radical velocity for $0\le t<t_0$. Then the fixed active space satisfies $P_0\subset T_{y(t_0)}\Sigma$. The insertion map $\widehat M_{t_0}$ obtained by raising an index of $B|_{P_0^4}$ with the endpoint metric $h|_{P_0}$ is precisely the finite limit of $M(t)$ in the fixed frame as $t\uparrow t_0$. This conclusion requires neither maximal rank at the endpoint nor invertibility of $Y(t_0)$.
\end{lemma}
\begin{proof}
For $0\le t<t_0$, the parallel frame
\[
 E_1,\ldots,E_q
\]
consists of fixed ambient vectors spanning $P_{y(t)}=P_0$. These vectors remain independent at the endpoint. To check tangency there, note that for $t<t_0$,
\[
 \dd E_{y(t)}(E_\alpha)=0;
\]
By continuity,
\[
 \dd E_{y(t_0)}(E_\alpha)=0,
\]
and hence
\[
 P_0\subset T_{y(t_0)}\Sigma.
\]
Although the endpoint rank has not yet been determined, the metric and the quartic form can already be restricted to this fixed space:
\[
 \widehat h_{t_0}=h_{y(t_0)}|_{P_0},\qquad
 \widehat b_{t_0}=B_{y(t_0)}|_{P_0^4}.
\]
The form $\widehat h_{t_0}$ is positive definite, so we may use it to raise one index of $\widehat b_{t_0}$ and define
\[
 \widehat h_{t_0}\bigl(\widehat{\mathcal B}_{t_0}(X,Y,Z),W\bigr)
 =\widehat b_{t_0}(X,Y,Z,W),
 \qquad X,Y,Z,W\in P_0.
\]
This gives the endpoint map
\[
 \widehat M_{t_0}:
 P_0\longrightarrow\operatorname{Hom}(\Sym^2P_0,P_0),
 \qquad
 (\widehat M_{t_0} X)(Y,Z)=\widehat{\mathcal B}_{t_0}(X,Y,Z).
\]
No endpoint rank is assumed in this definition.

We must also show that $\widehat M_{t_0}$ is indeed the limit of $M(t)$. In the fixed frame, write
\[
 H_{\alpha\beta}(t)=h_{y(t)}(E_\alpha,E_\beta),\qquad
 b_{\alpha\beta\gamma\delta}(t)
 =B_{y(t)}(E_\alpha,E_\beta,E_\gamma,E_\delta).
\]
These components have continuous limits at $t=t_0$, and $H(t_0)$ remains positive definite. Thus the inverse matrix $H(t)^{-1}=(H^{\alpha\beta}(t))$ also converges to $H(t_0)^{-1}$. The component formula for raising an index is
\[
 \mathcal B(t)^\mu{}_{\alpha\beta\gamma}
 =\sum_{\nu=1}^q H^{\mu\nu}(t)b_{\alpha\beta\gamma\nu}(t).
\]
Every term converges continuously, so
\[
 \widehat M_{t_0}
 =\lim_{t\uparrow t_0}M(t).
\]
Only continuity, the fixed frame and positive definiteness were used, without invertibility of $Y(t_0)$.
\end{proof}

\begin{theorem}[Maximal-rank extension before the first singularity of the denominator]\label{thm:root-maximal-rank-extension}
In the geometric setting of this subsection, $\tau\ge t_*$. Consequently, for every $0\le t<t_*$,
\[
 y(t)\in U,\qquad \dot y(t)\in N_{y(t)},\qquad P_{y(t)}=P_0,
\]
and \eqref{eq:mobius} continues to represent the original geometric tensors in the same fixed ambient frame.
\end{theorem}
\begin{proof}
We first show that the velocity remains a radical direction on $[0,\tau)$. The set of times when the velocity belongs to $N$ is nonempty and relatively closed, since $B(\dot y,\cdot,\cdot,\cdot)=0$ is a closed condition by continuity. Since the radical distribution is autoparallel for $\nabla$ and geodesics are locally unique, from any such time the curve remains in a local radical leaf for a short time in both directions. The set is therefore relatively open as well. Connectedness of the interval implies that the velocity belongs to $N$ at every time.

Suppose, to the contrary, that $\tau<t_*$. Then $\tau$ is finite, and the endpoint $p=y(\tau)$ of the global geodesic exists. By Lemma \ref{lem:root-finite-endpoint-limit}, the endpoint map is the limit of the geometric matrices. Proposition \ref{prop:root-mobius-solution} and invertibility of $Y(\tau)$ then give
\[
 \widehat M_\tau=\lim_{t\uparrow\tau}M(t)=M_0Y(\tau)^{-1}.
\]
Since $M_0$ is injective, $\widehat M_\tau$ is injective as well. Indeed,
\[
 \widehat M_\tau X=0
 \ \Longrightarrow\
 M_0Y(\tau)^{-1}X=0
 \ \Longrightarrow\
 Y(\tau)^{-1}X=0
 \ \Longrightarrow\
 X=0.
\]
Injectivity of the limit follows from this explicit composition with invertible $Y(\tau)^{-1}$, not merely from injectivity before the endpoint.

Put $p=y(\tau)$, and consider the full insertion map at the endpoint,
\[
 \iota_p:T_p\Sigma\longrightarrow\Sym^3(T_p^*\Sigma),
 \qquad
 \iota_p(X)=B_p(X,\cdot,\cdot,\cdot).
\]
We claim that $\iota_p|_{P_0}$ is injective. If $X\in P_0$ and $\iota_p(X)=0$, restricting the other three arguments to $P_0$ still gives
\[
 \widehat b_\tau(X,Y,Z,W)=0
 \quad\text{for all }Y,Z,W\in P_0.
\]
By nondegeneracy of $\widehat h_\tau$, this is equivalent to $\widehat M_\tau X=0$. Injectivity of the endpoint map gives $X=0$. Thus $\iota_p$ is injective on a $q$-dimensional space, and hence
\[
 \rank B_p=\rank\iota_p\ge q.
\]
On the other hand, $q$ is the maximum rank on the whole indicatrix, so $\rank B_p\le q$. Therefore
\[
 \rank B_p=q,\qquad p=y(\tau)\in U.
\]
Maximality supplies the upper rank bound matching the lower bound from the restricted tensor. Since $U$ is open and the global curve $y$ is continuous at $\tau$, there is an $\varepsilon>0$ such that $y([0,\tau+\varepsilon))\subset U$. This already contradicts the definition of $\tau$. For completeness, still under the contrary assumption, we also check below how the radical velocity and the fixed ambient frame would continue across the endpoint.

Since $\rank\iota_p=q$, the matrix of $\iota_p$ in local frames has a nonzero $q\times q$ minor. By continuity, this minor remains nonzero near $p$, so the nearby rank is at least $q$. Maximality bounds it above by $q$. Thus $p$ has a neighbourhood where the rank is identically $q$, and the whole neighbourhood lies in $U$. By the constant-rank theorem, $N$ and $P$ are smooth subbundles there.

To check the endpoint velocity, choose arbitrary smooth tangent vector fields $X,Y,Z$ near $p$. For $t<\tau$ sufficiently close to $\tau$,
\[
 B_{y(t)}(\dot y(t),X,Y,Z)=0.
\]
Letting $t\uparrow\tau$ and using smoothness of the tensor and the curve gives
\[
 B_p(\dot y(\tau),X,Y,Z)=0.
\]
Since the endpoint values of $X,Y,Z$ are arbitrary, we obtain
\[
 \dot y(\tau)\in N_p.
\]
The inclusion $\nabla_NN\subset N$ implies that the local $\nabla$-geodesic with initial data $(p,\dot y(\tau))$ remains in a radical leaf. By uniqueness of geodesics, it is a segment of the original global curve $y(t)$. Thus the same curve continues past $\tau$ as a radical geodesic.

Moreover, $P$ is now smooth near the endpoint, and $P_{y(t)}=P_0$ for $t<\tau$. Passing to the limit gives $P_p=P_0$. The original fixed vectors $E_1,\ldots,E_q$ therefore still form a frame of $P$ at the endpoint. Continuing them by $\nabla^P$-parallel translation leaves them the same fixed ambient vectors. Hence the geometric definitions of $S,M$ also extend in the same frame.

The openness argument above rules out $\tau<t_*$ and proves $\tau\ge t_*$. For every $0\le t<t_*$, the curve consequently lies in $U$, its velocity belongs to the radical, and \eqref{eq:mobius} holds in the same fixed ambient frame.
\end{proof}

The extension uses the globally complete geodesic, the explicit solution, smoothness, and maximality of rank; it assumes neither analyticity nor completeness of maximal-rank radical leaves.

If $t_*=\infty$, the curve remains in $U$ for all nonnegative time. The next subsection uses $a(t)\ge a_0>0$ to force a finite singularity when $S_0$ has a real eigenvalue.

\subsection{Exclusion of real eigenvalues and the global rank gap}\label{subsec:paper-7-6}\label{sec:blowup}
We use the preceding extension result to exclude real eigenvalues of $S_v$ for $v\ne0$, and then use the centroid to constrain the radical dimension.

Continue to use $q,U$ from \eqref{eq:maximal-rank-set}, and assume $0<q<m$.
Fix $y_0\in U$ and $0\ne v_0\in N_{y_0}$, and work along the corresponding complete centroaffine geodesic $y(t)$. Put
\[
 S_0=S_{v_0},\qquad M_0=M(0),\qquad
 a(t)=h_{y(t)}(\dot y,\dot y)\ge a_0>0.
\]
Here $M_0$ is injective, and the positive lower bound $a_0$ is supplied by Section \ref{sec:centro-periodic}. We use the standard fundamental solutions $f,g$ from Definition \ref{def:root-denominator-matrix} and the denominator matrix $Y$ defined in \eqref{eq:root-denominator}. We first locate a finite singularity and then contradict boundedness of the original tensor in the fixed frame.

\begin{lemma}[Finite zeros for a positive-coefficient equation and Riccati blow-up]\label{lem:root-scalar-zero}
For any $\lambda\in\R$, put $u_\lambda=f-\lambda g$. Its first positive zero $\tau_\lambda$ exists, and
\[
 0<\tau_\lambda\le
 \frac{\frac\pi2-\arctan(\lambda/\sqrt{a_0})}{\sqrt{a_0}}
 <\frac\pi{\sqrt{a_0}}.
\]
On $0\le t<\tau_\lambda$, we have $u_\lambda>0$, while $s=-u_\lambda'/u_\lambda$ satisfies $s'=s^2+a(t)$ and $s(0)=\lambda$, and tends to $+\infty$ as $t\uparrow\tau_\lambda$.
\end{lemma}
\begin{proof}
The scalar solution $u_\lambda$ is the stated linear combination of $f,g$. Since $f,g$ both satisfy $z''+a(t)z=0$ with the initial values specified in the preceding subsection,
\[
 u_\lambda''+a(t)u_\lambda=0,\qquad
 u_\lambda(0)=1,\quad u_\lambda'(0)=-\lambda.
\]
On an interval where $u_\lambda\ne0$, define
\[
 s(t)=-\frac{u_\lambda'(t)}{u_\lambda(t)}.
\]
Direct differentiation gives
\[
 s' =-\frac{u_\lambda''}{u_\lambda} +\frac{(u_\lambda')^2}{u_\lambda^2} =a(t)+s^2,
\]
and $s(0)=\lambda$. Thus the function constructed from the second-order linear equation solves the scalar Riccati initial value problem
\[
 s'=s^2+a(t),\qquad s(0)=\lambda
\]
until its denominator $u_\lambda$ first vanishes.

Let $\tau_\lambda$ be the first positive zero of $u_\lambda$, temporarily setting $\tau_\lambda=\infty$ if no positive zero exists. Since $u_\lambda(0)=1$, we have $\tau_\lambda>0$ and $u_\lambda(t)>0$ for $0\le t<\tau_\lambda$, so $s(t)$ is smooth there. Using $a(t)\ge a_0>0$, we obtain
\[
 \frac{\dd}{\dd t}\arctan\frac{s(t)}{\sqrt{a_0}} =\frac{\sqrt{a_0}}{a_0+s(t)^2}\,s'(t) =\frac{\sqrt{a_0}(s(t)^2+a(t))}{a_0+s(t)^2} \ge\sqrt{a_0}.
\]
Integration from $0$ to $t$ gives
\[
 \arctan\frac{s(t)}{\sqrt{a_0}}
 \ge
 \arctan\frac{\lambda}{\sqrt{a_0}}+\sqrt{a_0}\,t.
\]
The left-hand side is strictly less than $\pi/2$. Hence $s$ cannot be smooth on the entire interval $[0,T_\lambda]$, where
\[
 T_\lambda:=
 \frac{\frac{\pi}{2}-\arctan(\lambda/\sqrt{a_0})}{\sqrt{a_0}}
 \in\left(0,\frac{\pi}{\sqrt{a_0}}\right).
\]
If $\tau_\lambda>T_\lambda$, then $u_\lambda$ would be nonzero and $s$ smooth on $[0,T_\lambda]$. Substitution of $t=T_\lambda$ into the preceding inequality would make the left-hand side at least $\pi/2$, a contradiction. Therefore
\[
 0<\tau_\lambda\le T_\lambda<\infty.
\]

We can also determine the direction of divergence of $s$. The nontrivial solution $u_\lambda$ of the linear equation cannot satisfy both
\[
 u_\lambda(\tau_\lambda)=u_\lambda'(\tau_\lambda)=0.
\]
Otherwise, uniqueness of the initial value problem would make it identically zero, contrary to $u_\lambda(0)=1$. Thus this zero is simple. Since $u_\lambda>0$ before the zero, we have $u_\lambda'(\tau_\lambda)<0$. Taylor expansion gives
\[
 u_\lambda(t)
 =u_\lambda'(\tau_\lambda)(t-\tau_\lambda)
  +O((t-\tau_\lambda)^2),
\]
and therefore
\[
 s(t)=\frac{1}{\tau_\lambda-t}+O(1)
 \longrightarrow+\infty
 \qquad(t\uparrow\tau_\lambda).
\]
\end{proof}

Only finiteness of the positive zero will be used below; the upper bound and blow-up asymptotic are additional conclusions.

\begin{lemma}[A real eigendirection at the first singularity of the denominator]\label{lem:root-first-singular-eigenvector}
If $S_0$ has a real eigenvalue, the first positive singular time $t_*$ of $Y$ is finite. At that time $g(t_*)\ne0$, and
\[
 \lambda_0=\frac{f(t_*)}{g(t_*)}\in\R
\]
is an eigenvalue of $S_0$. There is a fixed nonzero real vector $z\in P_0$ such that
\[
 Y(t)z=u(t)z,\qquad u(t)=f(t)-\lambda_0g(t),
\]
where $u>0$ on $[0,t_*)$ and $u(t_*)=0$. Diagonalizability of $S_0$ is not required.
\end{lemma}
\begin{proof}
Suppose $S_0$ has a real eigenvalue $\lambda$. By Lemma \ref{lem:root-scalar-zero}, the corresponding $u_\lambda$ vanishes at a finite positive time $\tau_\lambda$.

Choose a nonzero real eigenvector $z_\lambda$ for $\lambda$. Then
\[
 Y(t)z_\lambda
 =(f(t)-\lambda g(t))z_\lambda
 =u_\lambda(t)z_\lambda.
\]
Thus $Y(\tau_\lambda)$ is singular, and the first positive singular time $t_*$ of $Y$ satisfies
\[
 0<t_*\le\tau_\lambda<\infty.
\]
The chosen $\lambda$ need not cause the first singularity, so we select an eigenvalue corresponding to $t_*$.

The preceding subsection proved that
\[
 fg'-f'g=1,
\]
so $f,g$ never vanish simultaneously. If $g(t_*)=0$, then $f(t_*)\ne0$, and
\[
 Y(t_*)=f(t_*)\Id
\]
would be invertible, contrary to the definition of $t_*$. Hence $g(t_*)\ne0$. Since
\[
 Y(t_*)=g(t_*)\left(\frac{f(t_*)}{g(t_*)}\Id-S_0\right)
\]
is singular, the real number
\[
 \lambda_0:=\frac{f(t_*)}{g(t_*)}
\]
must be an eigenvalue of $S_0$. Choose
\[
 0\ne z\in\ker(S_0-\lambda_0\Id)
\]
to be a real vector; such a vector is a nonzero kernel vector of a real singular matrix. Put
\[
 u(t)=f(t)-\lambda_0g(t).
\]
Then $u(t_*)=0$, and
\[
 Y(t)z=u(t)z.
\]
For $0\le t<t_*$, invertibility of $Y(t)$ gives $u(t)\ne0$. Since $u(0)=1$, the function $u(t)$ is positive on this interval and tends to zero as $t\uparrow t_*$. We have used only an ordinary eigenvector, without requiring diagonalizability of $S_0$.
\end{proof}

\begin{theorem}[Absence of real eigenvalues in nonzero radical directions]\label{thm:root-no-real-spectrum}
In the branch $0<q<m$,
\begin{equation}\label{eq:no-real}
 y\in U,\quad 0\ne v\in N_y
 \quad\Longrightarrow\quad S_v\text{ has no real eigenvalues}.
\end{equation}
\end{theorem}
\begin{proof}
Suppose, to the contrary, that $S_0$ has a real eigenvalue. Choose the finite time $t_*$, real eigenvector $z$, and function $u$ from Lemma \ref{lem:root-first-singular-eigenvector}. Theorem \ref{thm:root-maximal-rank-extension} ensures that the geometric equations remain applicable up to, but not including, $t_*$.

Section \ref{sec:extension} shows that, for every $0\le t<t_*$, the curve lies in the maximal-rank open set and the denominator representation \eqref{eq:mobius} gives the original geometric tensor in the same fixed parallel frame. Since $Y(t)z=u(t)z$ and $u(t)\ne0$, that representation gives directly
\[
 M(t)z=\frac{M_0z}{u(t)}.
\]
Since $M_0$ is injective and $z\ne0$, we have $M_0z\ne0$. Choose any norm on the fixed space
$\operatorname{Hom}(\Sym^2P_0,P_0)$. Then
\[
 \|M(t)z\|
 =\frac{\|M_0z\|}{|u(t)|}
 \longrightarrow\infty
 \qquad(t\uparrow t_*).
\]
The blow-up is in $M$, obtained from the original tensor $B$, not merely in the splitting operator $S$.

On the other hand, Theorem \ref{thm:root-maximal-rank-extension} keeps the curve in $U$ with radical velocity for $0\le t<t_*$. Thus Lemma \ref{lem:root-finite-endpoint-limit} applies at the finite endpoint $t_*$ and gives a finite limit for $M(t)$ in the same fixed frame. For fixed $z$, the vector $M(t)z$ is therefore bounded, contrary to the blow-up above. This conclusion requires neither maximal rank at the endpoint nor invertibility of $Y(t_*)$.

The contradiction applies to arbitrary $y_0\in U$, $0\ne v_0\in N_{y_0}$, proving \eqref{eq:no-real} without self-adjointness or diagonalizability.
\end{proof}

We now use the common eigenvector theorem for commuting families of complex matrices \cite[Lemma 1.3.19]{horn-johnson2013}, together with the real linearity of $v\mapsto S_v$.
\begin{lemma}[Common complex eigenvector and eigenvalue functional]\label{lem:root-common-eigenvector}
Fix $y\in U$, and complexify the active space to $P_{\C}=P_y\otimes_{\R}\C$. There exist $0\ne z\in P_{\C}$ and a real linear map
\[
 \ell:N_y\longrightarrow\C,
\]
such that $S_vz=\ell(v)z$ for every $v\in N_y$. We call $\ell$ the \textup{\textbf{eigenvalue functional}} corresponding to this common complex eigenline.
\end{lemma}

\begin{theorem}[Rank gap on the maximal-rank region]\label{thm:root-rank-gap}
If the global maximal insertion rank satisfies $0<q<m$, then on the maximal-rank open set $U$,
\begin{equation}\label{eq:rankgap}
 \dim N=1,\qquad q=m-1\text{ is even}.
\end{equation}
This conclusion constrains only the maximal-rank region of the nonzero degenerate branch.
\end{theorem}
\begin{proof}
Choose $z$ and $\ell$ as in Lemma \ref{lem:root-common-eigenvector}, and put $r=\dim N_y$.

Consider the real linear functional
\[
 \operatorname{Im}\ell:N\longrightarrow\R.
\]
If $r\ge2$, the rank--nullity formula gives
\[
 \dim\ker(\operatorname{Im}\ell)
 =r-\rank(\operatorname{Im}\ell)\ge r-1\ge1.
\]
Hence there is $0\ne v\in N$ with $\operatorname{Im}\ell(v)=0$, that is, $\lambda:=\ell(v)\in\R$. Thus
\[
 S_vz=\lambda z.
\]
Although $z$ is complex, this still implies that the real operator $S_v$ has a real eigenvalue. Indeed, write
\[
 z=z_1+iz_2,\qquad z_1,z_2\in P.
\]
Since $S_v$ is the complex linear extension of a real operator, comparison of real and imaginary parts gives
\[
 S_vz_1=\lambda z_1,\qquad S_vz_2=\lambda z_2.
\]
As $z\ne0$, at least one of $z_1,z_2$ is nonzero and is therefore a real eigenvector. This contradicts \eqref{eq:no-real}.

Thus $r\le1$. Since $q<m$ gives $r\ge1$, we have $r=1$ and $q=m-1$. Choose any nonzero radical vector $v$. The operator $S_v$ is a real matrix of order $q$, so its characteristic polynomial has real coefficients, and its nonreal roots occur in conjugate pairs with their algebraic multiplicities. By \eqref{eq:no-real}, it has no real roots; hence its degree $q$ is even. This proves \eqref{eq:rankgap}.
\end{proof}

The rank gap concerns $U$ only, not lower-rank points. If $m$ is even, it excludes the branch $0<q<m$, since $q=m-1$ would be odd. Thus the remaining analysis of this branch is needed only when $m$ is odd (equivalently, when $n=m+1$ is even).

This smooth argument used the positive lower bound from Section \ref{sec:centro-periodic}, maximal-rank extension, injectivity of $M_0$, and commutativity of the centroid. We next construct a compatible complex structure for use in the Cartan equations.

\subsection{Complex structures in the centroid}\label{subsec:paper-7-7}\label{sec:complex}
Fix $y\in U$ in the degenerate branch and $0\ne v\in N_y$. Write $P=P_y$, $b=b_y$, $S=S_v\in\Cent(b_y)$; by Section \ref{sec:blowup}, $S$ has no real eigenvalues.
We need to construct from $S$ a \emph{complex structure} still belonging to $\Cent(b)$, that is, a real linear operator $J$ with $J^2=-\Id$, and to transfer the same real polynomial to other symmetric tensors. This is an algebraic construction on a fixed vector space. It requires neither real-analyticity nor diagonalizability of $S$.

For a general matrix, prescribing only the values of a polynomial at the eigenvalues does not determine its action on Jordan blocks. We first recall the meanings of a complex structure and the minimal polynomial, and then construct a real polynomial that eliminates all nilpotent terms in the Jordan blocks.

\begin{definition}[Linear complex structure and compatibility with a quartic form]\label{def:quartic-compatible-complex-structure}
A real linear operator $J$ on a real vector space $P$ satisfying $J^2=-\Id$ is called a \textbf{linear complex structure}. Through the rule
\[
 (a+ic)x=ax+cJx,\qquad a,c\in\R,
\]
it makes the same underlying real space into a complex vector space, denoted by $(P,J)$. If also $J\in\Cent(b)$, we say that $J$ is \textbf{compatible} with $b$. Compatibility here means movement between slots and includes no orthogonality condition with respect to a given inner product.
\end{definition}

Associativity of this complex scalar multiplication can be checked directly:
\[
 (a\Id+cJ)(d\Id+eJ)=(ad-ce)\Id+(ae+cd)J.
\]
This is precisely complex multiplication. The complex dimension of $(P,J)$ is $\dim_{\R}P/2$, whereas the complexification $P\otimes_{\R}\C$ has complex dimension $\dim_{\R}P$. The two spaces serve different purposes in this subsection.

\begin{definition}[Minimal polynomial and factor exponents]\label{def:operator-minimal-polynomial}
The \textbf{minimal polynomial} of an operator $S$ on a finite-dimensional real space is the monic real polynomial of least degree satisfying $\mu(S)=0$.
In its complex factorization, the exponent of $(t-\lambda)^{r_\lambda}$ records the size of the largest Jordan block for $\lambda$. It controls the action of a polynomial on the generalized eigenspace and need not equal the algebraic multiplicity of $\lambda$ in the characteristic polynomial.
\end{definition}
It exists because the Cayley--Hamilton theorem gives an annihilating polynomial, and uniqueness follows from division with remainder.

The following real interpolation statement is a direct application of the polynomial Chinese remainder theorem \cite[Proposition 1.10]{atiyah-macdonald}. For the roots and factor exponents of the minimal polynomial, see \cite[Corollary 3.3.4 and Theorem 3.3.6]{horn-johnson2013}.
\begin{lemma}[Real interpolation retaining minimal-polynomial multiplicities]\label{lem:real-complex-polynomial}
Let $P\ne0$ be a finite-dimensional real vector space, let $S\in\operatorname{End}_{\R}(P)$ have no real eigenvalues, and let $\mu(t)\in\R[t]$ be the monic minimal polynomial of $S$. If $\Lambda_+$ denotes the set of distinct eigenvalues of $S$ in the upper half-plane, then
\[
 \mu(t)=\prod_{\lambda\in\Lambda_+}
 (t-\lambda)^{r_\lambda}(t-\bar\lambda)^{r_\lambda},
 \qquad r_\lambda\ge1.
\]
There is a unique polynomial $p$ with complex coefficients and degree less than $\deg\mu$ satisfying
\[
 p(t)\equiv i\pmod{(t-\lambda)^{r_\lambda}},\qquad p(t)\equiv-i\pmod{(t-\bar\lambda)^{r_\lambda}} \qquad(\lambda\in\Lambda_+).
\]
In fact, $p$ belongs to $\R[t]$, and $\mu$ divides $p^2+1$.
\end{lemma}

\begin{proof}
Apply the Chinese remainder theorem to the pairwise relatively prime ideals $((t-\lambda)^{r_\lambda})$ and $((t-\bar\lambda)^{r_\lambda})$. It gives the unique representative of the prescribed degree. Conjugation of the coefficients preserves this system of congruences, so uniqueness gives $p\in\R[t]$. Each factor divides either $p-i$ or $p+i$; hence $\mu\mid(p^2+1)$.
\end{proof}

The congruence $p(t)\equiv i\pmod{(t-\lambda)^r}$ prescribes not only $p(\lambda)=i$, but also
\[
 p^{(j)}(\lambda)=0\qquad(1\le j\le r-1).
\]
On a Jordan block $S=\lambda\Id+H$ of size at most $r$, the Taylor formula for the polynomial gives
\[
 p(\lambda\Id+H)
 =\sum_{j=0}^{r-1}\frac{p^{(j)}(\lambda)}{j!}H^j
 =i\Id.
\]
The lower-half-plane blocks give $-i\Id$. Thus the congruences control generalized eigenspaces, including the terms omitted by interpolation at distinct eigenvalues alone.

\begin{theorem}[A complex structure in the centroid of a quartic form]\label{thm:quartic-complex-structure}
Let $P\ne0$ be a finite-dimensional real vector space, and let $b\in\Sym^4(P^*)$ have trivial radical. If $\Cent(b)$ contains a real operator $S$ with no real eigenvalues, then there is a real polynomial $p$ such that $J=p(S)\in\Cent(b)$ and $J^2=-\Id$.
This conclusion requires neither diagonalizability of $S$ nor orthogonality of $J$ with respect to a prescribed inner product.
\end{theorem}

\begin{proof}
Take the real polynomial $p$ in Lemma \ref{lem:real-complex-polynomial}. Its divisibility property and $\mu(S)=0$ give
\[
 p(S)^2+\Id=\mu(S)r(S)=0.
\]
Hence
\begin{equation}\label{eq:Jcomplex}
 J=p(S),\qquad J^2=-\Id.
\end{equation}
Proposition \ref{prop:quartic-centroid-algebra} showed that $\Cent(b)$ is a real algebra with identity. If $p(t)=\sum_{j=0}^d a_jt^j$, then $S^j\in\Cent(b)$ for every power, and therefore $J=\sum_{j=0}^d a_jS^j\in\Cent(b)$.
Real coefficients make $J$ an operator on $P$, and centroid membership gives its movement between all slots.

\end{proof}

\begin{proposition}[Transfer of the same real polynomial between slots]\label{prop:polynomial-slot-transfer}
Let $\mathcal T\in\Sym^k(P^*)$, $k\ge2$, and suppose that a real operator $S$ satisfies
\[
 \mathcal T(Sx_1,x_2,\ldots,x_k)=\mathcal T(x_1,Sx_2,\ldots,x_k)
\]
for all arguments. Then, for every $p\in\R[t]$, the operator $p(S)$ can also be moved between any two argument slots of $\mathcal T$. Consequently, if another symmetric tensor also permits movement of $S$ between slots, the \emph{same} $J=p(S)$ constructed above applies to both tensors.
\end{proposition}
\begin{proof}
The movement identity holds for the identity operator. For an integer $j\ge1$, repeated use of the movement relation for $S$ gives
\[
 \mathcal T(S^jx_1,x_2,\ldots) =\mathcal T(S^{j-1}x_1,Sx_2,\ldots) =\cdots=\mathcal T(x_1,S^jx_2,\ldots).
\]
If $p(t)=\sum_j a_jt^j$, multiply these equalities by $a_j$ and add. Multilinearity gives the movement of $p(S)$ between the first two slots. Full symmetry of $\mathcal T$ then permits interchange of arguments and gives movement between any two slots.
\end{proof}

On the degenerate maximal-rank set, $P_y\ne0$, $\rad b_y=0$, and $S_v\in\Cent(b_y)$ has no real eigenvalues for $v\ne0$. Theorem \ref{thm:quartic-complex-structure} therefore applies at each point. We use the polynomial form of $J$ below.

Section \ref{sec:no-cubic} will prove the additional movement relation for the Cartan cubic form, after which Proposition \ref{prop:polynomial-slot-transfer} supplies the same $J$ for both tensors.

Only pointwise existence and polynomial transfer are asserted; no smooth field or orthogonality of $J$ is required.

\subsection{The Cartan differential equation in radical directions}\label{subsec:paper-8-1}\label{sec:A-ode}
Continue in the degenerate branch, where $N$ is one-dimensional, $P$ has positive even rank $q=m-1$, and $S_v$ has no real eigenvalues for $v\ne0$. The Cartan Codazzi equation will force $A_v$ to be scalar, with scalar coefficient constant along unit-speed radical curves. Closed orbits will then eliminate this scalar.

On a sufficiently small open subset $O$ of $U$, choose a smooth unit radical vector field $v$ such that
\[
 N|_O=\R v,\qquad h(v,v)=1.
\]
This choice uses only a local trivialization of the line bundle and does not require orientability of the whole bundle $N$. Put
\[
 c=C(v,v,v),\qquad A=A_v,\qquad T=T_v,\qquad S=S_v=T+A.
\]
Here $A,T$ are the operators from Definitions \ref{def:root-cartan-contraction} and \ref{def:root-splitting-operators}, respectively, for the unit radical vector field $v$.
By Lemma \ref{lem:root-cartan-mixed} and Proposition \ref{prop:root-autoparallel}, the mixed Cartan components satisfy \eqref{eq:CNNP}, and $\nabla^h$ preserves $N,P$ along radical directions. Proposition \ref{prop:root-centro-splitting} also gives $A,T\in\Cent(b)$, hence $AT=TA$. The operator $A$ is self-adjoint by symmetry of $C$, whereas $S$ has no real eigenvalues. We also use full symmetry of $\nabla^h C$ from Subsections \ref{subsec:paper-2-1}--\ref{subsec:paper-2-4}.

Along an integral curve of $v$, use $\nabla^{h,P}$, induced on $P$ by $\nabla^h$, and write
\[
 A'=\nabla^{h,P}_vA,\qquad
 A'X=\nabla^{h,P}_v(AX)-A(\nabla^{h,P}_vX)
\]
for the covariant derivative of an operator on $P$. This is defined because $\nabla^h_v$ preserves $P$.

\begin{lemma}[Cartan derivative equation on a one-dimensional radical distribution]\label{lem:cartan-root-ode}
With the preceding notation, the unit radical vector field satisfies
\[
 \nabla^h_vv=0,\qquad
 C(v,v,\cdot)^\sharp=cv,
\]
and along its integral curves,
\begin{equation}\label{eq:Aode}
 A'=(2A-c\Id)T.
\end{equation}
\end{lemma}

\begin{proof}
The inclusion $\nabla^h_NN\subset N$ implies that $\nabla^h_vv$ is a multiple of $v$. On the other hand, $\nabla^h h=0$ and $h(v,v)=1$ give
\[
 0=v\bigl(h(v,v)\bigr)=2h(\nabla^h_vv,v).
\]
Thus $\nabla^h_vv=0$, and the integral curves of $v$ are unit-speed geodesics of $h$.

Since $C(N,N,P)=0$, for any $Z\in T\Sigma|_O$ write $Z=h(Z,v)v+Z_P$ to obtain
\[
 C(v,v,Z)=h(Z,v)C(v,v,v)+C(v,v,Z_P)
         =c\,h(v,Z).
\]
By the definition of raising an index, this is $C(v,v,\cdot)^\sharp=cv$.

We must also determine the transverse derivative $\nabla^h_Xv$. For $X\in P$, the unit-length condition gives
\[
 0=X\bigl(h(v,v)\bigr)=2h(\nabla^h_Xv,v).
\]
Since $N=\R v$, it follows that $\nabla^h_Xv\in P$. The projection in the definition of $T$ may therefore be omitted, giving
\begin{equation}\label{eq:unit-root-transverse-derivative}
 \nabla^h_Xv=-TX.
\end{equation}
This uses $\dim N=1$; orthogonality to $v$ alone would not remove higher-dimensional radical components.

We next express the derivative of $A$ in terms of the covariant derivative of $C$.
Choose $\nabla^{h,P}$-parallel $X,Y\in P$ along a unit-speed radical curve, with arbitrary transverse initial values. Since $\nabla^h_v$ preserves $P$, $A'$ is the matrix derivative in this frame.

Differentiate $h(AX,Y)=C(v,X,Y)$ along the curve, using $\nabla^h_vv=0$:
\[
 h(A'X,Y) =v\bigl(C(v,X,Y)\bigr) =(\nabla^h_vC)(v,X,Y) =(\nabla^h_XC)(v,v,Y).
\]
The last equality uses full symmetry of $\nabla^h C$ to interchange the differentiation direction $v$ and the argument $X$.

We expand the correction terms in the covariant derivative.
At the point under consideration, extend $X,Y$ to local sections of $P$. Their parallelness along the curve is used only to express the derivative of $A$ as a covariant derivative of $C$. Since $C(v,v,Y)$ vanishes identically on $O$,
\[
 \begin{aligned}
 (\nabla^h_XC)(v,v,Y) ={}&X\bigl(C(v,v,Y)\bigr)\ -C(\nabla^h_Xv,v,Y)-C(v,\nabla^h_Xv,Y) -C(v,v,\nabla^h_XY)\\
 ={}&C(TX,v,Y)+C(v,TX,Y) -C(v,v,\nabla^h_XY).
 \end{aligned}
\]
Using $\nabla^h_Xv=-TX$, each of the first two terms equals $h(ATX,Y)$. The last term becomes $-c\,h(v,\nabla^h_XY)$ by $C(v,v,\cdot)^\sharp=cv$. Differentiating $h(v,Y)=0$ in direction $X$ gives
\[
 0=h(\nabla^h_Xv,Y)+h(v,\nabla^h_XY)
   =-h(TX,Y)+h(v,\nabla^h_XY).
\]
Hence $h(v,\nabla^h_XY)=h(TX,Y)$, and therefore
\[
 h(A'X,Y) =2h(ATX,Y)-c\,h(TX,Y) =h\bigl((2A-c\Id)TX,Y\bigr).
\]
Since $X,Y$ are arbitrary and $h|_P$ is nondegenerate, \eqref{eq:Aode} follows. Its $2A$ comes from the two $v$-slots and its $-cT$ from the $Y$-slot.
\end{proof}

\begin{remark}
Although the scalar function $C(v,v,Y)$ vanishes identically, $(\nabla^h_XC)(v,v,Y)$ includes slot corrections and need not vanish. Here the curves use $h$-arc length; their relation to the centroaffine affine parameter will be computed next.
\end{remark}

Self-adjointness in this equation and the absence of real eigenvalues now force $A$ to be scalar.

\begin{proposition}[Scalar form of Cartan contraction in a radical direction]\label{prop:cartan-root-scalar}
In the positive non-full maximal-rank branch, every local unit radical vector field $v$ on $U$ satisfies
\[
 A=\frac c2\Id,\qquad A'=0,\qquad v(c)=0.
\]
Equivalently, for every $X,Y\in P$,
\[
 C(v,X,Y)=\frac c2h(X,Y),\qquad c=C(v,v,v),
\]
and $c$ is constant along unit-speed radical curves.
\end{proposition}

\begin{proof}
In an $h|_P$-orthonormal, $\nabla^{h,P}$-parallel frame, $A$ and its derivative are symmetric. Equivalently, differentiating $h(AX,Y)=h(X,AY)$ for parallel $X,Y$ gives
\[
 h(A'X,Y)=h(X,A'Y).
\]
Hence $(2A-c\Id)T$ is self-adjoint. This uses metric compatibility of $\nabla^{h,P}$ and would not hold by the same argument for a nonmetric connection.

Fix $y\in O$. The following is linear algebra on $P_y$. Since $A$ is self-adjoint, $P_y$ is the orthogonal direct sum of its real eigenspaces. Choose one eigenvalue $a$, and put
\[
 E_a=\ker(A-a\Id)\ne0.
\]
Since $AT=TA$, for $X\in E_a$ we have
\[
 A(TX)=T(AX)=aTX.
\]
Thus $T$ preserves $E_a$, and $S=T+A$ also preserves $E_a$.

Suppose $2a-c\ne0$. For any $X,Y\in E_a$, we have $TX,TY\in E_a$. Self-adjointness of $(2A-c\Id)T$ therefore gives
\[
 (2a-c)h(TX,Y) =h\bigl((2A-c\Id)TX,Y\bigr) =h\bigl(X,(2A-c\Id)TY\bigr) =(2a-c)h(X,TY).
\]
Dividing by the nonzero real number $2a-c$ shows that $T|_{E_a}$ is self-adjoint for the positive definite inner product $h|_{E_a}$.

This yields a real eigenvalue of $S$, contradicting the known conclusion.
On the nonzero finite-dimensional positive definite inner product space $E_a$, the self-adjoint operator $T|_{E_a}$ has a real eigenvalue. Thus there are $0\ne X\in E_a$ and $\tau\in\R$ such that $TX=\tau X$. Then
\[
 SX=(T+A)X=(\tau+a)X.
\]
This gives the real eigenvalue $\tau+a$ of $S$, contrary to \eqref{eq:no-real}. Hence $2a-c\ne0$ is impossible, and every eigenvalue of $A$ equals $c/2$. Since $A$ is self-adjoint and orthogonally diagonalizable,
\[
 A_y=\frac{c(y)}2\Id_{P_y}.
\]
Thus $A=(c/2)\Id$ throughout $O$, without assuming self-adjointness of $T$ on all of $P$ or smooth eigenspaces.

Substitution into \eqref{eq:Aode} gives $A'=0$. On the other hand, $A=(c/2)\Id$ is an identity of operator fields on $O$, so
\[
 A'=\nabla^{h,P}_v\!\left(\frac c2\Id\right)
    =\frac{v(c)}2\Id.
\]
Since $q=\dim P>0$, the identity operator is nonzero, and $v(c)=0$. Finally, $h(AX,Y)=C(v,X,Y)$ gives the stated expression for the Cartan components.
\end{proof}

By Lemma \ref{lem:cartan-root-ode} and Proposition \ref{prop:cartan-root-scalar}, for a local unit radical vector field $v$, all Cartan components with at least one radical argument are controlled by the same scalar $c=C(v,v,v)$.
So far, $c$ is known to be constant along radical curves; transverse constancy and $c=0$ have not been proved. Replacing $v$ by $-v$ changes $c,A,T$ simultaneously in sign, preserving the conclusion.

We next compare the two parameters on complete closed radical orbits to prove $c=0$; the transverse part $C|_{P^3}$ is treated in Section \ref{sec:no-cubic}.

\subsection{Closed orbits and vanishing of the remaining scalar}\label{subsec:paper-8-2}\label{sec:closed-orbit}
Retain the local unit radical vector field $v$ and scalar $c$ from Lemma \ref{lem:cartan-root-ode} and Proposition \ref{prop:cartan-root-scalar}. To eliminate $c$, we first prove that the complete periodic orbit stays in $U$, and then compare its affine and unit-speed parameters.

\begin{proposition}[Elimination of radical Cartan components by closed orbits]\label{prop:closed-root-cartan}
In the positive non-full maximal-rank branch, the Cartan tensor on the maximal-rank open set $U$ satisfies
\[
 C(N,T\Sigma,T\Sigma)=0.
\]
In particular, for every local unit radical vector field $v$, we have $c=C(v,v,v)=0$, $A_v=0$, and $S_v=T_v$.
\end{proposition}

\begin{proof}
Choose arbitrary $y_0\in U$ and $0\ne w_0\in N_{y_0}$, and let $y(t)$ be the centroaffine geodesic with these initial data. Put
\[
 y(0)=y_0,\qquad w(t)=\dot y(t),\qquad w(0)=w_0,
 \qquad \nabla_ww=0.
\]
Section \ref{sec:centro-periodic} ensures that the curve exists on the whole real line and has a period $\tau_{\mathrm{per}}>0$ satisfying
\[
 y(t+\tau_{\mathrm{per}})=y(t),\qquad w(t+\tau_{\mathrm{per}})=w(t).
\]
The velocity $w(t)$ is everywhere nonzero. Otherwise, at some time the initial velocity would be zero, and uniqueness of the geodesic initial value problem would make the curve constant, contrary to $w_0\ne0$.

It remains to show that this global periodic curve stays in $U$ with radical velocity.

We now use the absence of real eigenvalues to show that the orbit remains in the maximal-rank open set.
Put $S_0=S_{w_0}$. Section \ref{sec:blowup} proved that $S_0$ has no real eigenvalues. The function $a(t)=h(w,w)$ is already defined along the global curve. Take the standard fundamental solutions $f,g$ from Definition \ref{def:root-denominator-matrix}. By Lemma \ref{lem:root-denominator-wronskian}, they satisfy $fg'-f'g=1$, so they never vanish simultaneously. Consider the denominator matrix $Y$ defined in \eqref{eq:root-denominator}.
If $g(t)=0$, then $f(t)\ne0$, so $Y(t)=f(t)\Id$ is invertible. If $g(t)\ne0$, then
\[
 Y(t)=g(t)\left(\frac{f(t)}{g(t)}\Id-S_0\right).
\]
Since $f(t)/g(t)$ is real and $S_0$ has no real eigenvalues, the operator in parentheses is also invertible. Thus $Y$ is invertible for every real time.

The maximal-rank extension result of Section \ref{sec:extension} therefore applies to every finite positive time and gives
\[
 y(t)\in U,\qquad w(t)\in N_{y(t)}\qquad(t\ge0).
\]
For extension to negative time, apply the same argument to the reversed curve $\widehat y(t)=y(-t)$. Its initial velocity is $-w_0$, and its initial splitting operator $S_{-w_0}=-S_0$ still has no real eigenvalues. Thus the whole curve remains in the maximal-rank set and tangent to the radical distribution. The integral below requires only one positive period.

The local conclusions now apply along the closed orbit; this required maximal-rank extension and invertibility of $Y$, not completeness alone.

Along the curve, define
\[
 u(t)=|w(t)|_h>0,\qquad v(t)=\frac{w(t)}{u(t)},
 \qquad c(t)=C_{y(t)}(v(t),v(t),v(t)).
\]
The nonzero periodic velocity $w=uv$ trivializes $\gamma^*N$, where $\gamma:\R/\tau_{\mathrm{per}}\mathbb Z\to U$ denotes the closed orbit. Its normalization $v$ is a smooth periodic unit section; no orientation of $N$ on $U$ is required. Local unit radical fields agreeing with $v(t)$ give compatible conclusions on overlapping intervals.

Put
\[
 s(t)=\int_0^t u(\tau)\,\dd\tau.
\]
Since $u>0$, the function $s$ is strictly increasing, and
\[
 \frac{\dd s}{\dd t}=u,\qquad
 \frac{\dd y}{\dd s}=v,\qquad h(v,v)=1.
\]
Thus $s$ is an arc-length parameter, whereas $t$ is an affine parameter for $\nabla$. The preceding subsection gives $\nabla^h_vv=0$, so $s$ is also the unit-speed parameter of this radical curve regarded as an $h$-geodesic. Moreover, $v(c)=0$ gives
\[
 \frac{\dd c(t)}{\dd t}=u(t)\,v(c)=0.
\]
Hence $c(t)$ is one constant along the whole orbit; denote it by $c_0$.

Comparison of the two connections gives the rate of change of the speed.
Since $\nabla_XY=\nabla^h_XY-C(X,Y,\cdot)^\sharp$ and $\nabla_ww=0$,
\[
 \nabla^h_ww =C(w,w,\cdot)^\sharp =u^2C(v,v,\cdot)^\sharp =u^2c_0v.
\]
On the other hand, $w=uv$ and $\nabla^h_vv=0$ give
\[
 \nabla^h_ww =\nabla^h_w(uv) =\frac{\dd u}{\dd t}\,v+u\nabla^h_wv =\frac{\dd u}{\dd t}\,v+u^2\nabla^h_vv =\frac{\dd u}{\dd t}\,v.
\]
Since $v\ne0$, comparison yields $\dd u/\dd t=c_0u^2$. Using $\dd s/\dd t=u$, we obtain
\[
 \frac{\dd}{\dd s}\log u
 =\frac1u\frac{\dd}{\dd t}\log u
 =\frac1{u^2}\frac{\dd u}{\dd t}
 =c_0.
\]
The two factors $1/u$ arise respectively from changing parameter and differentiating $\log u$.

We finally integrate over a period to eliminate the constant $c_0$.
Periodicity of both position and velocity of the centroaffine geodesic gives
\[
 u(\tau_{\mathrm{per}})
 =\sqrt{h_{y(\tau_{\mathrm{per}})}(w(\tau_{\mathrm{per}}),w(\tau_{\mathrm{per}}))}
 =\sqrt{h_{y(0)}(w(0),w(0))}
 =u(0).
\]
Let the $h$-length of one period be
\[
 L=s(\tau_{\mathrm{per}})-s(0)=\int_0^{\tau_{\mathrm{per}}}u(t)\,\dd t>0.
\]
On this period, $u$ is continuous and strictly positive, so $L$ is finite and positive. Integrating $\dd(\log u)/\dd s=c_0$ with respect to arc length gives
\[
 0 =\log u(\tau_{\mathrm{per}})-\log u(0) =\int_{s(0)}^{s(\tau_{\mathrm{per}})}c_0\,\dd s =c_0L.
\]
Thus $c_0=0$, using return of both position and velocity, not constancy of the centroaffine geodesic's $h$-speed.

The initial point $y_0\in U$ was arbitrary. At $y_0$ on the orbit above, take $v=w_0/|w_0|_h$. Since $c_0=0$, Proposition \ref{prop:cartan-root-scalar} gives
\[
 C(v,v,v)=0,\qquad C(v,v,X)=0,\qquad
 C(v,X,Y)=0\quad(X,Y\in P_{y_0}).
\]
Write arbitrary $Z,W\in T_{y_0}\Sigma$ as $Z=\alpha v+X$ and $W=\beta v+Y$. Multilinearity and full symmetry give
\[
 C(v,Z,W) ={}\alpha\beta C(v,v,v)+\alpha C(v,v,Y) +\beta C(v,X,v)+C(v,X,Y)=0.
\]
Since $N_{y_0}=\R v$, insertion of any radical vector also annihilates $C$. Thus $C(N,T\Sigma,T\Sigma)=0$ on all of $U$, or equivalently, for every local unit radical vector field,
\begin{equation}\label{eq:C-null}
 C(v,\cdot,\cdot)=0\quad\text{throughout the maximal-rank region}.
\end{equation}
The identity $h(A_vX,Y)=C(v,X,Y)$ gives $A_v=0$, and $S_v=T_v+A_v$ then gives $S_v=T_v$. This proves the proposition.
\end{proof}

Only the transverse cubic form remains:
\begin{equation}\label{eq:transverse-cartan-cubic}
 c_P:=C|_{P^3}.
\end{equation}
Since $S_v=T_v$, the operator $T_v$ inherits the established spectral and centroid properties; we next prove movement of $T_v$ between the slots of $c_P$.

\subsection{A common complex structure and compatible cubic forms}\label{subsec:paper-8-3}\label{sec:no-cubic}
By Proposition \ref{prop:closed-root-cartan} and \eqref{eq:C-null}, it remains to study the transverse cubic form $c_P$ defined in \eqref{eq:transverse-cartan-cubic}. We use the Codazzi equation to obtain a common complex structure for $c_P$ and the restricted quartic form $b$, then \eqref{eq:F} and positivity of $h|_P$ to prove $c_P=0$.

\begin{lemma}[A common complex structure from the Codazzi equation]\label{lem:common-complex-cubic-quartic}
In the positive non-full maximal-rank branch under consideration, let $c_P$ be the transverse cubic form defined in \eqref{eq:transverse-cartan-cubic}. For a local unit radical vector field $v$, the operator $T=T_v=S_v$ satisfies
\[
 c_P(TX,Y,Z)=c_P(X,TY,Z)=c_P(X,Y,TZ).
\]
At each $y\in U$, there is a real polynomial $p$ such that $J=p(T_y)$ satisfies $J^2=-\Id$ and can be moved between any two slots of both $b_y$ and $c_{P,y}$.
\end{lemma}
\begin{proof}
Choose a unit radical vector field $v$ on a sufficiently small open subset of $U$, and put $T=T_v=S_v$. Here $\dim N=1$, so \eqref{eq:unit-root-transverse-derivative} applies to every transverse vector field $X\in\Gamma(P)$.
Since $\nabla^h_vP\subset P$, the induced derivative of $c_P$ satisfies
\[
 (\nabla^{h,P}_vc_P)(X,Y,Z)=(\nabla^h_vC)(X,Y,Z)
 \qquad(X,Y,Z\in P).
\]

Full symmetry of $\nabla^h C$ gives
\[
 (\nabla^h_vC)(X,Y,Z)=(\nabla^h_XC)(v,Y,Z).
\]
Expanding the right-hand side by definition, we obtain
\[
 (\nabla^h_XC)(v,Y,Z) ={}X\bigl(C(v,Y,Z)\bigr)-C(\nabla^h_Xv,Y,Z) -C(v,\nabla^h_XY,Z)-C(v,Y,\nabla^h_XZ).
\]
By \eqref{eq:C-null}, $C(v,\cdot,\cdot)$ vanishes identically on the open set under consideration. Hence the first and the last two terms vanish. Substituting $\nabla^h_Xv=-TX$ gives
\[
 (\nabla^{h,P}_vc_P)(X,Y,Z)=c_P(TX,Y,Z).
\]
The left-hand side is fully symmetric in $X,Y,Z$. Interchanging these arguments and using symmetry of $c_P$ itself, we obtain
\[
 c_P(TX,Y,Z)=c_P(X,TY,Z)=c_P(X,Y,TZ).
\]
This is precisely the movement property of $T$ between the slots of $c_P$.

At each $y\in U$, we have $T_y=S_y\in\Cent(b_y)$, and $T_y$ has no real eigenvalues. By Theorem \ref{thm:quartic-complex-structure}, with $T_y$ in place of $S$ in \eqref{eq:Jcomplex}, we obtain a real polynomial $p$ and a complex structure $J=p(T_y)\in\Cent(b_y)$. The movement property just proved and Proposition \ref{prop:polynomial-slot-transfer} show that this $J$ can also be moved between any two slots of $c_{P,y}$. Thus the two tensors have the \emph{same} compatible complex structure, as required.
\end{proof}

The construction of $J$ is only pointwise. The algebraic lemma below does not differentiate $J$ and requires no smooth choice of $J$ on all of $U$.

\begin{lemma}[Vanishing of a compatible cubic form under a common complex structure]\label{lem:no-cubic}
Let $(P,h)$ be a finite-dimensional positive definite real inner product space, and let
\[
 b\in\Sym^4(P^*),\qquad \gamma\in\Sym^3(P^*).
\]
Assume that $b$ has trivial radical, meaning that its insertion map $w\mapsto b(w,\cdot,\cdot,\cdot)$ is injective. Suppose that a real linear operator $J$ satisfies $J^2=-\Id$ and can be moved between any two argument slots of both $b$ and $\gamma$.

Define the raised-index vector $\gamma^\sharp(X,Y)$ using $h$:
\[
 h\bigl(\gamma^\sharp(X,Y),Z\bigr)=\gamma(X,Y,Z).
\]
If, in addition, for every $e,a,x,y,z\in P$ the compatibility identity
\[
 0={}b\bigl(a,\gamma^\sharp(e,x),y,z\bigr) -b\bigl(e,\gamma^\sharp(a,x),y,z\bigr) -b\bigl(a,\gamma^\sharp(e,y),x,z\bigr) +b\bigl(e,\gamma^\sharp(a,y),x,z\bigr),
\]
holds, then $\gamma=0$. Orthogonality of $J$ with respect to $h$ is not required.
\end{lemma}

\begin{proof}
Complexify $P$ to $P_{\C}=P\otimes_{\R}\C$, and denote the corresponding complex multilinear extensions by $b_{\C}$, $\gamma_{\C}$, and $h_{\C}$. The form $h_{\C}$ is nondegenerate and complex bilinear. Its inverse gives an isomorphism
\[
 h_{\C}^{-1}:P_{\C}^*\longrightarrow P_{\C}.
\]
Explicitly, $q=h_{\C}^{-1}\alpha$ means that $h_{\C}(q,w)=\alpha(w)$ for every $w\in P_{\C}$. Thus the complexification of the operation of raising an index is
\[
 \gamma_{\C}^\sharp(e,x)
 =h_{\C}^{-1}\bigl(\gamma_{\C}(e,x,\cdot)\bigr).
\]
The compatibility identity is real multilinear and therefore holds for arbitrary complex vectors as well. Its internal contractions are still performed by $h_{\C}^{-1}$.

Since $J^2=-\Id$, the minimal polynomial of the complexification of $J$ divides $t^2+1=(t-\mathrm i)(t+\mathrm i)$. The two roots are distinct, so there is an eigenspace direct sum
\[
 P_{\C}=P_+\oplus P_-,\qquad
 P_\pm=\ker(J\mp\mathrm i\Id).
\]
Because $J$ is real, complex conjugation interchanges $P_+$ and $P_-$. This direct sum is defined by eigenspaces; it is not an orthogonal decomposition assumed in advance.

Let $\tau$ denote either $b_{\C}$ or $\gamma_{\C}$. If two of its arguments are $X_+\in P_+$ and $Y_-\in P_-$, moving $J$ between these slots gives
\[
 \mathrm i\,\tau(X_+,Y_-,\ldots) =\tau(JX_+,Y_-,\ldots) =\tau(X_+,JY_-,\ldots) =-\mathrm i\,\tau(X_+,Y_-,\ldots).
\]
Thus mixed components vanish by symmetry, and each tensor has only its pure positive and pure negative parts. Write
\[
 b_\pm=b_{\C}|_{P_\pm^4},\qquad
 \gamma_\pm=\gamma_{\C}|_{P_\pm^3}.
\]

We must also show that the purely negative quartic form $b_-$ has trivial radical.
First, the full form $b_{\C}$ has trivial radical. Indeed, if $w=u+\mathrm i v$, with $u,v\in P$, and
\[
 b_{\C}(w,X,Y,Z)=0\qquad(X,Y,Z\in P_{\C}),
\]
then, taking real vectors $X,Y,Z\in P$ in particular, its real and imaginary parts give
\[
 b(u,X,Y,Z)=0,\qquad b(v,X,Y,Z)=0.
\]
Since the real quartic form $b$ has trivial radical, $u=v=0$, and hence $w=0$.

If $w\in P_-$ annihilates $b_-$, components with the other arguments in $P_-$ vanish by assumption, and those with any argument in $P_+$ vanish by the mixed-component identity. Multilinearity gives
\[
 b_{\C}(w,\cdot,\cdot,\cdot)=0.
\]
Thus $w=0$, proving that $b_-$ has trivial radical.

Choose arbitrary
\[
 e,x\in P_+,\qquad a,y,z\in P_-.
\]
The second and third terms contain the mixed covectors $\gamma_{\C}(a,x,\cdot)$ and $\gamma_{\C}(e,y,\cdot)$ and vanish. The fourth term has fixed arguments from both eigenspaces in $b_{\C}$ and also vanishes, leaving
\[
 b_{\C}\bigl(a,\gamma_{\C}^\sharp(e,x),y,z\bigr)=0.
\]

Fix $e,x$, and put
\[
 \alpha=\gamma_{\C}(e,x,\cdot)\in P_{\C}^*,\qquad
 q=h_{\C}^{-1}\alpha.
\]
Decompose $q$ according to the eigenspace direct sum as $q=q_++q_-$. Since $a,y,z\in P_-$, the term containing $q_+$ is mixed and vanishes. The preceding identity reduces to
\[
 b_-(a,q_-,y,z)=0\qquad(a,y,z\in P_-).
\]
Symmetry and triviality of the radical of $b_-$ give $q_-=0$. Consequently,
\[
 h_{\C}^{-1}\alpha=q\in P_+.
\]
This eliminates only the negative component of the raised vector; positivity is still needed to show $\alpha=0$.

We now use positivity of the conjugate pairing to prove $\alpha=0$.
Since $e,x\in P_+$, vanishing of mixed components gives $\alpha|_{P_-}=0$. Complex conjugation of a covector is defined by
\[
 \bar\alpha(w)=\overline{\alpha(\bar w)}.
\]
If $w\in P_+$, then $\bar w\in P_-$, so $\bar\alpha(w)=0$. Hence $\bar\alpha|_{P_+}=0$. Together with $q\in P_+$, this gives
\[
 h_{\C}^{-1}(\alpha,\bar\alpha)
 :=\bar\alpha\bigl(h_{\C}^{-1}\alpha\bigr)
 =\bar\alpha(q)=0.
\]
Write $\alpha=\mu+\mathrm i\nu$, where $\mu,\nu\in P^*$ are real covectors. Complex bilinear expansion and symmetry of $h^{-1}$ give
\[
 h_{\C}^{-1}(\alpha,\bar\alpha) =h^{-1}(\mu,\mu)+h^{-1}(\nu,\nu)+\mathrm i\bigl(h^{-1}(\nu,\mu)-h^{-1}(\mu,\nu)\bigr) =|\mu|_{h^{-1}}^2+|\nu|_{h^{-1}}^2.
\]
Positive definiteness of $h^{-1}$ forces $\mu=\nu=0$, hence $\alpha=0$. Positivity is used with the conjugate pairing, not with the complex bilinear expression $h_{\C}(w,w)$.

We can now eliminate the entire cubic form.
The identity $\alpha=0$ gives, for every $e,x\in P_+$,
\[
 \gamma_{\C}(e,x,\cdot)=0.
\]
Thus $\gamma_+=0$. Since $\gamma$ is real, for arbitrary $w_1,w_2,w_3\in P_-$ we also have
\[
 \gamma_{\C}(w_1,w_2,w_3)
 =\overline{\gamma_{\C}(\bar w_1,\bar w_2,\bar w_3)}=0,
\]
where $\bar w_j\in P_+$. Hence $\gamma_-=0$. The mixed components already vanish by the movement property, so $\gamma_{\C}=0$. Restriction to the real space gives $\gamma=0$.
\end{proof}

\begin{corollary}[Vanishing of the Cartan tensor on the degenerate branch]\label{cor:degenerate-branch-cartan-zero}
In the degenerate branch $0<q<m$, we have $C=0$ throughout the maximal-rank open set $U$.
\end{corollary}
\begin{proof}
Fix any $y\in U$, and apply the lemma with $P=P_y$, $h=h_y|_{P_y}$, $b=b_y$, and $\gamma=c_{P,y}$. Positive definiteness of $h$, triviality of the radical of $b$, and existence of a common complex structure have all been established. It remains to verify the compatibility identity.

For any $X,Y\in P_y$, use the inner product on the full tangent space to define $C^\sharp(X,Y)$. By \eqref{eq:C-null} and symmetry of $C$,
\[
 h\bigl(C^\sharp(X,Y),v\bigr)=C(X,Y,v)=0.
\]
Thus $C^\sharp(X,Y)\in P_y$ and equals the vector obtained by raising an index of $c_{P,y}$ with $h_y|_{P_y}$. Taking all external arguments in $P_y$ in \eqref{eq:F} therefore gives the required compatibility identity for $b_y,c_{P,y}$; in an adapted orthonormal basis, internal radical terms vanish by \eqref{eq:C-null} and the radical property of $B$.

The lemma gives $c_P=0$ on $U$. Expanding three tangent vectors into radical and transverse parts, terms with a radical argument vanish by \eqref{eq:C-null}, and the all-transverse term vanishes by $c_P=0$. Thus
\[
 C=0\qquad\text{on the entire maximal-rank open set }U.
\]
\end{proof}

We have proved $C=0$ on the nonempty maximal-rank set in the degenerate branch, without analyticity; this neither excludes the branch nor implies $B=0$. Later, analyticity of $F$ will extend this vanishing. The next section treats the nonempty full-rank set independently.

\section{Generic irreducibility on the full-rank set}\label{sec:paper-7}
We consider the quartic form $B$ of Definition \ref{def:twisted-quartic} on the nonempty set $\Omega$ where its radical is trivial. This is the full-rank branch $q=m$; no conclusion from the degenerate branch is used. Excluding smooth orthogonal pure-block decompositions, then applying constant-rank kernels and spectral projections, will give an open dense set of irreducible points.

\subsection{Exclusion of orthogonal pure-block decompositions}\label{subsec:paper-9-1}\label{sec:no-splitting}
Throughout this subsection we assume that the full-rank set
\begin{equation}\label{eq:full-rank-locus}
 \Omega:=\{y\in\Sigma:\rad B_y=0\}
\end{equation}
is nonempty. It is open because full rank is an open condition. We first exclude a smooth orthogonal decomposition of $B$ into two nonzero pure blocks on an open subset of $\Omega$.

Here $P,Q$ denote candidate block subbundles; $P$ no longer denotes specifically the orthogonal complement of the radical used in Section \ref{sec:paper-6}.

\begin{definition}[Smooth orthogonal pure-block decomposition]\label{def:orthogonal-pure-block}
Two smooth subbundles $P,Q$ of positive rank on an open set are said to give an \textbf{orthogonal pure-block decomposition} of $B$ if
\[
 T\Sigma=P\oplus^\perp Q,\qquad B=B_P+B_Q,
\]
where $B_P=B|_{P^4}$ and $B_Q=B|_{Q^4}$ are each extended by zero along the other subbundle. Equivalently, $B$ vanishes whenever its four arguments include vectors from both $P$ and $Q$. This includes both the case of three arguments in one block and one in the other, and the case of two arguments in each block.
\end{definition}

\begin{proposition}[No nontrivial orthogonal pure-block decomposition on the set of trivial radical]\label{prop:no-pure-block}
Let $\mathcal V\subset\Omega$ be a nonempty open set. There are no smooth subbundles $P,Q$ of positive rank on $\mathcal V$ giving a smooth orthogonal pure-block decomposition in the sense of Definition \ref{def:orthogonal-pure-block}.
\end{proposition}

\begin{proof}
Suppose such a decomposition exists on $\mathcal V$, to which all identities and fields below are restricted. Use $C^\sharp$, obtained by raising an index with $h$, as in Definition \ref{def:twisted-quartic}.

We first prove that each restricted quartic form has trivial radical.
Fix any $p\in \mathcal V$. If $w\in P_p$ satisfies
\[
 (B_P)_p(w,X_1,X_2,X_3)=0
 \qquad(X_1,X_2,X_3\in P_p),
\]
then, for arbitrary $U_1,U_2,U_3\in T_p\Sigma$, decomposition of each argument according to $P_p\oplus Q_p$ and the vanishing of the mixed components give
\[
 B_p(w,U_1,U_2,U_3)
 =(B_P)_p\bigl(w,(U_1)_P,(U_2)_P,(U_3)_P\bigr)=0.
\]
Hence $w\in\rad B_p=0$. Interchanging $P,Q$ shows that both restricted forms have trivial radical throughout $\mathcal V$.

We next use the compatibility identity to prove that the Cartan cubic form also has only pure-block components.
In the vector form \eqref{eq:quartic-compatibility-vector} of algebraic compatibility, take $e,j\in Q$ and $a,i,k\in P$. In the first, second and fourth terms, $B$ already has fixed arguments from both blocks. These three terms therefore vanish, regardless of the decomposition of the contracted vector. Only the third term remains, so
\[
 0=-B_P\bigl(a,(C^\sharp(e,j))_P,i,k\bigr).
\]
Here $a,i,k$ range arbitrarily over $P$. Since $B_P$ has trivial radical and is completely symmetric,
\[
 (C^\sharp(e,j))_P=0.
\]
Thus, for every $X\in P$,
\[
 C(e,j,X)=h\bigl(C^\sharp(e,j),X\bigr)=0.
\]
This proves $C(Q,Q,P)=0$. Interchanging $P,Q$ gives $C(P,P,Q)=0$ in the same way. These are the only two types of mixed arguments for a cubic form, and hence
\[
 C=C_P+C_Q,\qquad C_P=C|_{P^3},\quad C_Q=C|_{Q^3}.
\]

Take arbitrary smooth local vector fields
\[
 X_1,X_2,X_3\in\Gamma(P),\qquad Y,Z\in\Gamma(Q).
\]
Substituting $(e,a,i,j,k)=(Y,X_1,X_2,X_3,Z)$ in the vector form \eqref{eq:twisted-vector} of the twisted Codazzi equation, we obtain
\[
 \begin{aligned}
 &(\nabla^h_YB)(X_1,X_2,X_3,Z)
       -(\nabla^h_{X_1}B)(Y,X_2,X_3,Z)\\
 &\quad=B\bigl(Y,C^\sharp(X_1,X_2),X_3,Z\bigr)
       -B\bigl(X_1,C^\sharp(Y,X_2),X_3,Z\bigr)=0.
 \end{aligned}
\]
The first term on the right already has the fixed arguments $Y\in Q$, $X_3\in P$, and the second has $X_1\in P$, $Z\in Q$. Both are therefore mixed components and vanish. The two covariant derivatives are consequently equal.

We first compute the covariant derivative with two arguments in each block. By definition,
\[
 \begin{aligned}
 &(\nabla^h_{X_1}B)(Y,X_2,X_3,Z)=X_1\bigl(B(Y,X_2,X_3,Z)\bigr)-B(\nabla^h_{X_1}Y,X_2,X_3,Z)\\
 &\quad-B(Y,\nabla^h_{X_1}X_2,X_3,Z)-B(Y,X_2,\nabla^h_{X_1}X_3,Z)-B(Y,X_2,X_3,\nabla^h_{X_1}Z).
 \end{aligned}
\]
The first term vanishes because the scalar function $B(Y,X_2,X_3,Z)$ is identically zero on $\mathcal V$. Each of the remaining four terms changes only one argument; the other three still include vectors from both $P$ and $Q$. Every such term therefore vanishes as well. It follows that
\[
 (\nabla^h_{X_1}B)(Y,X_2,X_3,Z)=0,
 \qquad (\nabla^h_YB)(X_1,X_2,X_3,Z)=0.
\]
We have not yet assumed that $\nabla^h$ preserves the two subbundles. The correction terms vanish solely because the unchanged arguments remain mixed.

We now expand $(\nabla^h_YB)(X_1,X_2,X_3,Z)$ by definition:
\[
 \begin{aligned}
 &(\nabla^h_YB)(X_1,X_2,X_3,Z)=Y\bigl(B(X_1,X_2,X_3,Z)\bigr)-B(\nabla^h_YX_1,X_2,X_3,Z)\\
 &\quad-B(X_1,\nabla^h_YX_2,X_3,Z)-B(X_1,X_2,\nabla^h_YX_3,Z)-B(X_1,X_2,X_3,\nabla^h_YZ).
 \end{aligned}
\]
The first term is again the ordinary derivative of an identically zero scalar function. In the three correction terms involving $\nabla^h_YX_r$, the unchanged arguments include two vectors in $P$ and the vector $Z\in Q$, so these terms also vanish. Only the last term can retain a pure-block component: after decomposing $\nabla^h_YZ$ according to $P\oplus Q$, its $Q$ component still gives mixed arguments, whereas its $P$ component puts all four arguments in $P$. Thus
\[
 (\nabla^h_YB)(X_1,X_2,X_3,Z)
 =-B_P\bigl(X_1,X_2,X_3,(\nabla^h_YZ)_P\bigr).
\]
Combining this with $(\nabla^h_YB)(X_1,X_2,X_3,Z)=0$, already proved, gives
\[
 B_P\bigl(X_1,X_2,X_3,(\nabla^h_YZ)_P\bigr)=0.
\]
At each point, the values of $X_1,X_2,X_3$ may be prescribed arbitrarily in $P$. Since $B_P$ has trivial radical, $(\nabla^h_YZ)_P=0$. Therefore
\[
 \nabla^h_QQ\subset Q.
\]

We now use metric compatibility to prove that both subbundles are parallel in every direction.
For $Y,Z\in\Gamma(Q)$ and $X\in\Gamma(P)$, orthogonality gives $h(X,Z)=0$, and hence
\[
 0=Y\bigl(h(X,Z)\bigr)
   =h(\nabla^h_YX,Z)+h(X,\nabla^h_YZ).
\]
We already know that $\nabla^h_YZ\in Q$, so the second term is zero. Since $Z$ is arbitrary, we obtain $\nabla^h_YX\in P$, that is, $\nabla^h_QP\subset P$. Interchanging $P,Q$ and repeating the mixed covariant derivative computation and the metric compatibility argument gives the other two inclusions. Thus
\[
 \nabla^h_PP\subset P,\qquad \nabla^h_QP\subset P,\qquad \nabla^h_PQ\subset Q,\qquad \nabla^h_QQ\subset Q.
\]
Every tangent vector field $W$ decomposes as $W_P+W_Q$, and the connection is $C^\infty$-linear in the differentiating direction. Consequently,
\[
 \nabla^h_W\Gamma(P)\subset\Gamma(P),\qquad
 \nabla^h_W\Gamma(Q)\subset\Gamma(Q).
\]
Thus $P,Q$ are parallel in all tangent directions.

Fix any $p\in \mathcal V$. Since both blocks are nonzero, we may choose unit vectors $X\in P_p$, $Y\in Q_p$ and extend them to local sections of $P,Q$, respectively. Since $Q$ is parallel, the connection preserves $Q$ in every direction, and therefore
\[
 R^h(X,Y)Y
 =\nabla^h_X\nabla^h_YY-\nabla^h_Y\nabla^h_XY-\nabla^h_{[X,Y]}Y
 \in Q_p.
\]
Moreover, $X\perp Q_p$ and $X,Y$ are orthonormal. The mixed plane they span thus satisfies
\[
 K(X,Y)=h\bigl(R^h(X,Y)Y,X\bigr)=0.
\]

On the other hand, the Gauss equation gives curvature one for the same mixed plane.
The identity $C=C_P+C_Q$ proved above has the following consequence: for $X\in P_p$, the Cartan contraction operator $C_X$ preserves $P_p$ and vanishes on $Q_p$; for $Y\in Q_p$, the operator $C_Y$ preserves $Q_p$ and vanishes on $P_p$. For example, if $U\in Q_p$, then $C(X,U,W)=0$ for every $W$, whence $C_XU=0$; if $U\in P_p$, then $C(X,U,Q_p)=0$, whence $C_XU\in P_p$. It follows that
\[
 C_XC_Y=C_YC_X=0,\qquad [C_X,C_Y]=0.
\]
By the Gauss equation \eqref{eq:gauss},
\[
 R^h(X,Y)Y=R_0(X,Y)Y-[C_X,C_Y]Y=R_0(X,Y)Y,
\]
where $R_0(X,Y)Y=h(Y,Y)X-h(X,Y)Y=X$, since $X,Y$ are orthonormal. Therefore
\[
 K(X,Y)=h(X,X)=1.
\]
This contradicts the zero curvature of the same mixed plane, excluding the assumed decomposition.
\end{proof}

Under a pure-block decomposition, triviality of the radical of each restriction also implies triviality of the radical of $B$. Indeed, if $w=w_P+w_Q$ annihilates $B$, taking the other three arguments in $P$ shows that $w_P$ annihilates $B_P$, so $w_P=0$; taking them in $Q$ similarly gives $w_Q=0$. Thus, in the presence of a pure-block decomposition, triviality of the radicals of both blocks is equivalent to $\mathcal V\subset\Omega$. If the two blocks are initially known to have trivial radical only at one point, a nonzero minor of each insertion map remains nonzero nearby; shrinking the neighbourhood on which the decomposition is defined then keeps both radicals trivial.

The decomposition must be smooth on an open set: the mixed scalar components must vanish identically before differentiation. A splitting at a single point is insufficient. The next subsection constructs the required local decomposition from smooth spectral projections; irreducibility at every point of $\Omega$ is not asserted.

The contradiction uses the curvature of a single mixed plane, not constant curvature of the entire indicatrix. Parallelness gives zero, while the vanishing mixed Cartan terms leave the positive model term $R_0$ in the Gauss equation.

\subsection{The kernel of the commutator map and generic irreducibility}\label{subsec:paper-9-2}\label{sec:commutant}
To deduce generic irreducibility, we must turn pointwise invariant subspaces into smooth subbundles before differentiating. We obtain smooth kernel sections on a constant-rank neighbourhood, use a spectral gap to construct smooth orthogonal projections, and apply Section \ref{sec:no-splitting}.

Continue to work on the fixed Minkowski indicatrix $(\Sigma,h)$, put $m=\dim\Sigma$, and consider the nonempty full-rank open set $\Omega$ defined in \eqref{eq:full-rank-locus}. For $u\in\Omega$, define the contraction operators of $B_u$ by
\[
 h_u\bigl(B_{XY}Z,W\bigr)=B_u(X,Y,Z,W),
 \qquad X,Y,Z,W\in T_u\Sigma.
\]
Complete symmetry gives $B_{XY}=B_{XY}^*$. Here and below adjoints are taken with respect to $h_u$.

\begin{definition}[Common invariant subspaces, irreducibility, and the self-adjoint common commutant]\label{def:family-irreducibility}
Fix $u\in\Omega$. A subspace $E\subset T_u\Sigma$ is called a \textbf{common invariant subspace} of the contraction family $\{B_{XY}\}$ if $B_{XY}E\subset E$ for every $B_{XY}$. The family is called \textbf{irreducible} if its only common invariant subspaces are the zero space and the whole space. Define its \textbf{self-adjoint common commutant} by
\begin{equation}\label{eq:selfadjoint-commutant}
 \mathcal Z_B(u)
 :=\{S\in\operatorname{End}(T_u\Sigma):
       S=S^*,\ [S,B_{XY}]=0\text{ for all }X,Y\in T_u\Sigma\}.
\end{equation}
\end{definition}

This is a real vector space containing $\R\Id$. Membership in the common commutant requires $S$ to commute with every $B_{XY}$; different operators $B_{XY}$ need not commute with one another. By contrast, $\Cent(b)$ in Subsections \ref{subsec:paper-7-1}--\ref{subsec:paper-7-7} is defined by moving an operator between tensor slots. These two spaces will be used according to their respective definitions.

The following criterion is an immediate consequence of the real self-adjoint spectral theorem \cite[Corollary 2.5.11(a)]{horn-johnson2013}.
\begin{proposition}[Pointwise criterion for irreducibility]\label{prop:commutant-irreducibility}
Fix $u$. The only common invariant subspaces of the contraction family $\{B_{XY}\}$ are $0$ and $T_u\Sigma$ if and only if
\[
 \mathcal Z_B(u)=\R\Id.
\]
\end{proposition}

We next express the commutation conditions as the kernel of a smooth linear map.
Choose a smooth $h$-orthonormal frame $e_1,\ldots,e_m$ on a coordinate neighbourhood, and write $B_{ij}=B_{e_i e_j}$. If $X=\sum_i x^ie_i$, $Y=\sum_j y^je_j$, then
\[
 B_{XY}=\sum_{i,j}x^iy^jB_{ij}.
\]
Hence commuting with all $B_{XY}$ is equivalent to commuting with the finitely many $B_{ij}$.

Denote by $\mathcal E_u$ the space of self-adjoint operators on $T_u\Sigma$. Its dimension is
\[
 d=\frac{m(m+1)}2.
\]
In the chosen frame define the linear map
\[
 \mathcal L_u:\mathcal E_u\longrightarrow
       \bigoplus_{1\le i\le j\le m}\mathfrak{so}(T_u\Sigma,h_u),
 \qquad S\longmapsto\bigl([S,B_{ij}]\bigr)_{i\le j}.
\]
Here $\mathfrak{so}(T_u\Sigma,h_u)$ denotes the space of skew-adjoint operators. The commutator of two self-adjoint operators is indeed skew-adjoint, since
\[
 (SB_{ij}-B_{ij}S)^*=B_{ij}S-SB_{ij}.
\]
The matrix coefficients of $\mathcal L_u$ depend smoothly on $u$, and
\[
 \ker\mathcal L_u=\mathcal Z_B(u),\qquad
 \dim\ker\mathcal L_u=d-\rank\mathcal L_u.
\]
Changing the smooth orthonormal frame does not change this kernel, and hence changes neither its dimension nor the rank of $\mathcal L_u$. Since $\Id\in\ker\mathcal L_u$, we always have
\[
 \dim\mathcal Z_B(u)\ge1,\qquad
 \rank\mathcal L_u\le d-1.
\]
Thus the set of points under consideration is precisely
\begin{equation}\label{eq:irreducible-locus}
 G:=\{u\in\Omega:\mathcal Z_B(u)=\R\Id\}
   =\{u\in\Omega:\rank\mathcal L_u=d-1\}.
\end{equation}

The next two local lemmas explain how a pointwise element of the kernel gives the smooth decomposition that we need.

We use the constant-rank theorem for vector bundle homomorphisms and the extension lemma for sections of vector bundles \cite[Theorem 10.34 and Lemma 10.12]{lee2013}.
\begin{lemma}[Constant-rank kernels and smooth extension of kernel elements]\label{lem:constant-rank-kernel}
Suppose that a smooth vector bundle map has constant rank $r$ on a neighbourhood. Then its kernel is a smooth vector subbundle. Every vector in the kernel at any point extends, near that point, to a smooth section of the kernel.
\end{lemma}

We use the construction of spectral projections from finite-dimensional perturbation theory; see \cite[Chapter II, Section 1.4]{kato}.
\begin{lemma}[A spectral gap gives a smooth orthogonal spectral projection]\label{lem:spectral-gap-projection}
Let $S(u)$ be a smooth field of self-adjoint operators, and suppose that $S(p)$ is not scalar. After shrinking the neighbourhood of $p$, there exists a smooth field of orthogonal projections $\Pi$ satisfying
\[
 \Pi^2=\Pi=\Pi^*,\qquad 0<\rank\Pi<m,
\]
whose image is the sum of a group of eigenspaces of $S$. Every operator commuting with $S(u)$ also commutes with $\Pi(u)$.
\end{lemma}

Only a gap between two spectral groups is needed, not constant multiplicities or smooth eigenvectors. In a smooth local $h$-orthonormal frame, choose a positively oriented closed contour $\Gamma$ separating two nonempty spectral groups of $S(p)$, and put $\delta=\operatorname{dist}(\Gamma,\sigma(S(p)))>0$. Shrink a neighbourhood $U$ of $p$ so that $\|S(u)-S(p)\|<\delta/2$ in the operator norm. The self-adjoint eigenvalue perturbation bound gives $\operatorname{dist}(\Gamma,\sigma(S(u)))>\delta/2$ throughout $U$, so the fixed contour remains disjoint from the spectrum and separates these groups. Hence
\[
 \Pi(u)=\frac{1}{2\pi\mathrm i}
 \int_\Gamma(z\Id-S(u))^{-1}\,\dd z.
\]
The resolvent along the fixed contour depends smoothly on $u$, giving the required smooth projection. By the constant-rank theorem for vector bundles, its image and kernel are smooth subbundles; no parallelness of these subbundles has yet been asserted.

We now prove generic irreducibility by means of the constant-rank kernel and spectral projections.
\begin{theorem}[Generic irreducibility]\label{thm:irreducible}
Let $B$ be a smooth completely symmetric quartic form on a fixed Minkowski indicatrix, satisfying the twisted Codazzi equation of Definition \ref{def:twisted-quartic}, and suppose that the full-rank open set defined in \eqref{eq:full-rank-locus} satisfies $\Omega\ne\varnothing$. Then $G$, defined in \eqref{eq:irreducible-locus}, is an open dense subset of $\Omega$. On $G$, the family of contractions of the quartic form is irreducible.
\end{theorem}
\begin{proof}
Using the characterization of $G$ by the rank of $\mathcal L_u$ in \eqref{eq:irreducible-locus}, we first prove openness and then relative density by contradiction.

If $p\in G$, then $\rank\mathcal L_p=d-1$, so some minor of order $(d-1)$ is nonzero at $p$. By continuity it remains nonzero nearby, and hence
\[
 \rank\mathcal L_u\ge d-1.
\]
On the other hand, the identity operator always belongs to the kernel, so $\rank\mathcal L_u\le d-1$. The rank is therefore constantly $d-1$ nearby, which means that all nearby points belong to $G$.

Suppose, to the contrary, that $G$ is not dense in $\Omega$. Then $\Omega\setminus\overline G$ is nonempty. Choose a nonempty coordinate open set $\mathcal V$ in this complement. Thus
\[
 \dim\ker\mathcal L_u\ge2,\qquad
 \rank\mathcal L_u\le d-2
 \qquad(u\in \mathcal V).
\]
In $\mathcal V$, choose a point $p$ at which $\mathcal L$ attains its maximal rank, and denote that rank by $r$. The existence of such a point does not require compactness of $\mathcal V$: rank takes only finitely many integer values, so its values on the nonempty set $\mathcal V$ have a maximum, attained at some point.

If $r>0$, a minor of order $r$ nonzero at $p$ gives a neighbourhood $U\subset \mathcal V$ of $p$ on which $\rank\mathcal L_u\ge r$. Maximality gives the reverse inequality $\rank\mathcal L_u\le r$. Hence the rank is constantly $r$ on $U$. If $r=0$, the rank is already constantly zero on $\mathcal V$, and we may simply take a small neighbourhood $U$. By the constant-rank kernel lemma,
\[
 \ker\mathcal L|_U\longrightarrow U
\]
is a smooth subbundle of rank $d-r\ge2$.

Since $\dim\ker\mathcal L_p\ge2$, choose a nonscalar self-adjoint operator $S_p$ in this kernel. By the constant-rank kernel lemma, extend it to a smooth kernel section $S$ on $U$. Then
\[
 S(u)=S(u)^*,\qquad [S(u),B_{XY}(u)]=0.
\]
To ensure explicitly that the field remains nonscalar nearby, take its trace-free part
\[
 \widehat S(u)=S(u)-\frac{\tr S(u)}m\Id.
\]
Since $\Id$ always belongs to the kernel, $\widehat S$ is still a smooth kernel section. Moreover, $\widehat S(p)\ne0$, and after shrinking the neighbourhood we have $\widehat S(u)\ne0$. A nonzero trace-free operator cannot be scalar, so $\widehat S$ is nonscalar throughout this neighbourhood.

Apply the spectral projection lemma and shrink $U$ once more to obtain a smooth nontrivial orthogonal projection $\Pi$. Every $B_{XY}$ commutes with $\widehat S$, hence also with $\Pi$. Put
\[
 P=\operatorname{im}\Pi,\qquad Q=\ker\Pi.
\]
Then $P,Q$ are smooth subbundles of positive rank, and
\[
 T\Sigma|_U=P\oplus^\perp Q.
\]
The commutation relation ensures that every $B_{XY}$ preserves both subbundles. For example, if $Z\in P$, then
\[
 \Pi B_{XY}Z=B_{XY}\Pi Z=B_{XY}Z,
\]
so $B_{XY}Z\in P$; the argument for $Q$ is the same.

For arbitrary $X,Y$ and $Z\in P$, $W\in Q$, invariance and orthogonality give
\[
 B(X,Y,Z,W)=h(B_{XY}Z,W)=0.
\]
Whenever the four arguments include vectors from both $P$ and $Q$, complete symmetry allows us to move two such vectors to the last two slots. All mixed components therefore vanish, and $P,Q$ give a smooth orthogonal pure-block decomposition on all of $U$ in the sense of Definition \ref{def:orthogonal-pure-block}. We continue to denote the corresponding blocks by $B_P,B_Q$.

Since $U\subset\Omega$, both restrictions have trivial radical. For $B_P$, for example, if $v\in P_u$ annihilates $B_P$, then for arbitrary $Z_1,Z_2,Z_3\in T_u\Sigma$,
\[
 B_u(v,Z_1,Z_2,Z_3)
 =(B_P)_u\bigl(v,(Z_1)_P,(Z_2)_P,(Z_3)_P\bigr)=0,
\]
whence $v\in\rad B_u=0$. The argument for $B_Q$ is identical.

The resulting smooth orthogonal pure-block decomposition on the nonempty open set $U\subset\Omega$ contradicts Proposition \ref{prop:no-pure-block}.

Thus $G$ is open and dense in $\Omega$. Since $\mathcal Z_B=\R\Id$ on $G$, the pointwise criterion gives irreducibility there.
\end{proof}

Openness and density are relative to $\Omega$: the exceptional set has empty interior there, but need not be empty, and no dimension or measure estimate is obtained. Since $\Omega$ is nonempty and open, so is $G$; density in the entire indicatrix is not asserted.

The later rigidity argument needs only a nonempty open set of irreducible points. The same contradiction on any exceptional open set gives the stronger density conclusion.

Subsections \ref{subsec:paper-11-1}--\ref{subsec:paper-11-3} give a curvature identity at these points; Section \ref{sec:quartic-rigidity} extends that identity analytically. Neither the smooth spectral projections nor irreducibility itself is extended. The present argument uses only smoothness, finite-dimensional spectral theory, and the preceding geometric proposition.

\section{Curvature-compatible quartic forms and Jordan rigidity}\label{sec:paper-8}
We pass from angular compatibility to curvature isotropy by a contraction identity and Jordan closure. A quartic extremum supplies a rank-one spectral projection; irreducibility then gives all self-adjoint operators and hence $K=sR_0$.

\subsection{Commutator curvature and compatible operators}\label{subsec:paper-10-1}\label{sec:curvature-compatible}
We work in a fixed finite-dimensional real Euclidean space $(V,h)$, put $m=\dim V$, and prescribe completely symmetric tensors
\[
 C\in\Sym^3V^*,\qquad B\in\Sym^4V^*.
\]
In the geometric application, $V$ is an indicatrix tangent space and $h,C,B$ are evaluated at that point. The present constructions are pointwise; the next subsection uses $\mathcal F_B(C)=0$ to prove compatibility of the quartic contractions.

Since $h$ is positive definite, for every $X,Y\in V$ there is a unique vector $C_XY$ satisfying
\[
 h(C_XY,Z)=C(X,Y,Z)\qquad(Z\in V).
\]
For fixed $X$, the map $C_X:V\to V$ is linear, and it depends linearly on $X$. Complete symmetry gives two distinct properties, both of which will be used:
\[
 C_XY=C_YX,\qquad C_X^*=C_X.
\]
The first follows because the two sides have the same inner product with every $Z$; the second follows from
\[
 h(C_XY,Z)=C(X,Y,Z)=C(X,Z,Y)=h(Y,C_XZ)
\]
Adjoints are always taken with respect to $h$.

Similarly, for fixed $U,V\in V$, define the contraction operator $B_{UV}$ of the quartic form by
\[
 h(B_{UV}X,Y)=B(U,V,X,Y).
\]
It is bilinear in $U,V$, satisfies $B_{UV}=B_{VU}$, and obeys
\[
 h(B_{UV}X,Y)=B(U,V,Y,X)=h(X,B_{UV}Y).
\]
Thus the contractions studied in Section \ref{sec:paper-7} are self-adjoint.

\begin{definition}[Commutator curvature of a cubic form]\label{def:commutator-curvature}
With the commutator convention $[A,D]=AD-DA$, the tensor $K$ defined below is called the \textbf{commutator curvature} of $C$:
\[
 K(X,Y):=[C_X,C_Y],\qquad
 K(X,Y)Z=C_X(C_YZ)-C_Y(C_XZ).
\]
\end{definition}

The tensor $K$ is bilinear in $X,Y$. By the Gauss equation \eqref{eq:gauss}, its relation to the indicatrix curvature is
\[
 R^h(X,Y)Z=R_0(X,Y)Z-K(X,Y)Z,
 \qquad R_0(X,Y)Z=h(Y,Z)X-h(X,Z)Y.
\]
Thus $K$ occurs with a minus sign in our Gauss convention.

To fix the order of the four indices, we first write, for vectors,
\[
 K(x,y,w,p):=h(K(x,y)p,w).
\]
The last slot $p$ is the vector acted upon, and the third slot $w$ is paired with the output. After choosing an orthonormal basis $e_1,\ldots,e_m$, write
\[
 C_{ijk}:=C(e_i,e_j,e_k),\qquad
 K_{xywp}:=h(K(e_x,e_y)e_p,e_w).
\]
Here lower-case letters in subscripts index basis vectors. By self-adjointness of $C_{e_x}$,
\[
 h(C_{e_x}C_{e_y}e_p,e_w) =h(C_{e_y}e_p,C_{e_x}e_w) =\sum_{\ell=1}^m h(C_{e_y}e_p,e_\ell)h(C_{e_x}e_w,e_\ell) =\sum_{\ell=1}^m C_{yp\ell}C_{xw\ell}.
\]
Interchanging $x,y$ and subtracting gives
\begin{equation}\label{eq:Kcomponents}
 K_{xywp}=\sum_{\ell=1}^m
       \bigl(C_{xw\ell}C_{yp\ell}-C_{yw\ell}C_{xp\ell}\bigr).
\end{equation}
This also shows that the contraction over the internal index $\ell$ uses $h$. We shall often omit such summation signs in an orthonormal basis.

\begin{proposition}[Algebraic symmetries of commutator curvature]\label{prop:commutator-curvature-symmetries}
The tensor $K$ above is skew in the first two slots and in the last two slots, is symmetric under interchange of the two pairs, and satisfies the first Bianchi identity:
\[
 \begin{gathered}
 K_{xywp}=-K_{yxwp},\qquad
 K_{xywp}=-K_{xypw},\qquad
 K_{xywp}=K_{wpxy},\\
 K(X,Y)Z+K(Y,Z)X+K(Z,X)Y=0.
 \end{gathered}
\]
It is therefore an algebraic curvature tensor.
\end{proposition}
\begin{proof}
Skew-symmetry in the first two slots follows directly from $[C_Y,C_X]=-[C_X,C_Y]$. Since each $C_X$ is self-adjoint,
\[
 K(X,Y)^*
 =(C_XC_Y-C_YC_X)^*
 =C_YC_X-C_XC_Y=-K(X,Y),
\]
and hence
\[
 h(K(X,Y)P,W)=-h(K(X,Y)W,P),
\]
which gives skew-symmetry in the last two slots.

Symmetry under interchange of the two pairs may be checked directly from \eqref{eq:Kcomponents}. Omitting the internal summation sign and using complete symmetry of $C$ and commutativity of scalar multiplication, we obtain
\[
 K_{wpxy} =C_{wx\ell}C_{py\ell}-C_{px\ell}C_{wy\ell} =C_{xw\ell}C_{yp\ell}-C_{yw\ell}C_{xp\ell} =K_{xywp}.
\]

Finally, expand all six terms in the Bianchi cyclic sum and group them according to the outermost contraction operator:
\[
 \begin{aligned}
 &K(X,Y)Z+K(Y,Z)X+K(Z,X)Y=C_XC_YZ-C_YC_XZ+C_YC_ZX-C_ZC_YX+C_ZC_XY-C_XC_ZY\\
 &\quad=C_X(C_YZ-C_ZY)+C_Y(C_ZX-C_XZ)+C_Z(C_XY-C_YX)=0.
 \end{aligned}
\]
The last step uses $C_XY=C_YX$. All these symmetries are pointwise algebraic identities.
\end{proof}

We use the definition of curvature compatibility of Mantica--Molinari \cite[arXiv version, Section 4, equation (17)]{mm}. Raising one index of their symmetric tensor of order two by $h$ gives exactly the self-adjoint operator used here; the vector and component formulas below express the same condition.
\begin{definition}[Curvature-compatible operator]\label{def:K-compatible}
A self-adjoint operator $A:V\to V$ is called \textbf{$K$-compatible} if, for arbitrary $x,y,z\in V$,
\[
 K(x,y)Az+K(y,z)Ax+K(z,x)Ay=0.
\]
Equivalently, in any orthonormal basis, writing $A_{zp}=h(Ae_z,e_p)$, we require
\begin{equation}\label{eq:Kcompatible}
 K_{xywp}A_{zp}+K_{yzwp}A_{xp}+K_{zxwp}A_{yp}=0
\end{equation}
for all free indices. Here the repeated index $p$ is summed, while $x,y,z,w$ are free indices.
\end{definition}

Compatibility inserts the same operator $A$ into the acted-upon slot of each Bianchi term. Self-adjointness imposes no restriction on the signs or multiplicities of its eigenvalues.

For the following properties of the compatible space, see \cite[arXiv version, Section 4 and the beginning of Section 4.1]{mm}.
\begin{proposition}[Linear structure and identity of the compatible space]\label{prop:compatible-space-vector}
For fixed $K$, the $K$-compatible self-adjoint operators form a real vector space containing $\R\Id$. This is the basic property used in Section \ref{sec:jordan} to obtain a unital Jordan algebra.
\end{proposition}

The next contraction identity will show that $\mathcal F_B(C)=0$ implies, for every $U,V$,
\[
 K(x,y)B_{UV}z+K(y,z)B_{UV}x+K(z,x)B_{UV}y=0.
\]
Complete symmetry alone does not give this conclusion.

\subsection{Angular compatibility and curvature compatibility}\label{subsec:paper-10-2}\label{sec:CF}
We prove curvature compatibility by a contraction identity for arbitrary completely symmetric $C,B$, imposing the angular identity $\mathcal F_B(C)=0$ from Lemma \ref{lem:quartic-algebraic-compatibility} only at the end.

Continue to take an orthonormal basis $e_1,\ldots,e_m$ in the fixed Euclidean space $(V,h)$, and let $C\in\Sym^3V^*$, $B\in\Sym^4V^*$. Following \eqref{eq:F}, but without setting it equal to zero, define
\[
 \mathcal F_{eaijk}:=\sum_{\ell=1}^m\bigl( C_{ei\ell}B_{a\ell jk}-C_{ai\ell}B_{e\ell jk} -C_{ej\ell}B_{a\ell ik}+C_{aj\ell}B_{e\ell ik}\bigr).
\]
The indices have already been raised in the orthonormal basis; all internal contractions use $h$.

Write $B_{uv}=B_{e_u e_v}$ and define the tensor with six indices
\[
 \mathcal T_{uvxyzw}:=\sum_p\bigl(
 K_{xywp}B_{uvzp}+K_{yzwp}B_{uvxp}+K_{zxwp}B_{uvyp}\bigr).
\]
Since $B_{uv}e_z=\sum_pB_{uvzp}e_p$, we have
\[
 \mathcal T_{uvxyzw} =h\bigl(K(e_x,e_y)B_{uv}e_z +K(e_y,e_z)B_{uv}e_x +K(e_z,e_x)B_{uv}e_y,\ e_w\bigr).
\]
Thus $\mathcal T=0$ says exactly that all contractions $B_{uv}$ by basis vectors satisfy the compatibility condition of the preceding subsection. Every $B_{UV}$ is a linear combination of these operators, and the compatibility condition is linear in the operator. Hence this is also equivalent to $K$-compatibility of every $B_{UV}$.

Below, the notation $\sum_{\operatorname{cyc}(x,y,z)}$ always denotes only the three choices
\[
 (x,y,z),\qquad(y,z,x),\qquad(z,x,y).
\]
The six free indices are $u,v,x,y,z,w$. The indices $u,v,w$ do not take part in the cyclic sum, and the internal summation indices $r,\ell$ are independent.

\begin{proposition}[Curvature contraction identity for the angular compatibility tensor]\label{prop:curvature-contraction-identity}
For arbitrary completely symmetric $C,B$ as above, without assuming $\mathcal F=0$, one always has
\begin{equation}\label{eq:CF}
2\mathcal T_{uvxyzw}=\sum_{\operatorname{cyc}(x,y,z)}\sum_r\big[C_{xwr}(2\mathcal F_{urvyz}-2\mathcal F_{urvzy} +\mathcal F_{uyvzr}-\mathcal F_{uzvyr}) +C_{uxr}\mathcal F_{vwyzr}+C_{vxr}\mathcal F_{uwyzr}\big].
\end{equation}
\end{proposition}
\begin{proof}
Each $\mathcal F$ on the right contains four terms. Expanding the three cyclic groups therefore gives $3\times6\times4=72$ products with their coefficients. To exhibit the cancellations, denote by $\mathscr A$ the cyclic sum of the first four $\mathcal F$ terms on the right, and by $\mathscr D$ that of the last two:
\[
 \begin{aligned}
 \mathscr A&:=\sum_{\operatorname{cyc}(x,y,z)}\sum_r C_{xwr}
   \bigl(2\mathcal F_{urvyz}-2\mathcal F_{urvzy}+\mathcal F_{uyvzr}-\mathcal F_{uzvyr}\bigr),\\
 \mathscr D&:=\sum_{\operatorname{cyc}(x,y,z)}\sum_r
   \bigl(C_{uxr}\mathcal F_{vwyzr}
         +C_{vxr}\mathcal F_{uwyzr}\bigr).
 \end{aligned}
\]
We shall prove $\mathscr A=2\mathcal T+\mathcal E$ and $\mathscr D=-\mathcal E$, where $\mathcal E$ is an intermediate cross-contraction.

We first expand the first four terms and separate the principal terms from the cross terms.
Fix one choice of $(x,y,z)$. By the definition of $\mathcal F$, the two terms containing $C_{uv\ell}$ in the first difference cancel because $B_{r\ell yz}=B_{r\ell zy}$; the two terms containing $C_{rv\ell}$ cancel in the same way. Thus
\[
 \mathcal F_{urvyz}-\mathcal F_{urvzy} =\sum_\ell\bigl(-C_{uy\ell}B_{r\ell vz} +C_{ry\ell}B_{u\ell vz} +C_{uz\ell}B_{r\ell vy} -C_{rz\ell}B_{u\ell vy}\bigr).
\]
In the other difference, the terms containing $C_{uv\ell}$ cancel by complete symmetry of $B$, and those containing $C_{yz\ell}$ cancel because $C_{yz\ell}=C_{zy\ell}$. What remains is
\[
 \mathcal F_{uyvzr}-\mathcal F_{uzvyr} =\sum_\ell\bigl(-C_{yv\ell}B_{u\ell zr} +C_{zv\ell}B_{u\ell yr} -C_{uz\ell}B_{y\ell vr} +C_{uy\ell}B_{z\ell vr}\bigr).
\]
Multiply the first identity by $2$, add the second, and use complete symmetry of $C,B$ to arrange the indices within each factor. This gives
\[
 \begin{aligned}
 &2\mathcal F_{urvyz}-2\mathcal F_{urvzy}
       +\mathcal F_{uyvzr}-\mathcal F_{uzvyr}\\
 &\quad=\sum_\ell\bigl[
       2C_{yr\ell}B_{uvz\ell}-2C_{zr\ell}B_{uvy\ell}        +C_{uz\ell}B_{vr\ell y}-C_{uy\ell}B_{vr\ell z}        +C_{vz\ell}B_{ur\ell y}-C_{vy\ell}B_{ur\ell z}\bigr].
 \end{aligned}
\]
The first two terms inside the brackets form the principal part that will give $2\mathcal T$; the remaining four terms define the cross-contraction
\[
 \mathcal E:=\sum_{\operatorname{cyc}(x,y,z)}\sum_{r,\ell} C_{xwr}\bigl(C_{uz\ell}B_{vr\ell y}-C_{uy\ell}B_{vr\ell z} +C_{vz\ell}B_{ur\ell y}-C_{vy\ell}B_{ur\ell z}\bigr).
\]
The quantities $\mathscr A,\mathscr D,\mathcal E$ retain the original six free indices, whose subscripts are temporarily suppressed in the proof.

The separation into principal terms and cross-contractions gives
\[
 \mathscr A-\mathcal E
 =2\sum_{\operatorname{cyc}(x,y,z)}\sum_{r,\ell}
 C_{xwr}\bigl(C_{yr\ell}B_{uvz\ell}
             -C_{zr\ell}B_{uvy\ell}\bigr).
\]
Within the existing three-term cyclic sum, the second term may be regrouped as
\[
 \sum_{\operatorname{cyc}(x,y,z)}
 C_{xwr}C_{zr\ell}B_{uvy\ell}
 =\sum_{\operatorname{cyc}(x,y,z)}
 C_{ywr}C_{xr\ell}B_{uvz\ell}.
\]
This merely writes the same three cyclic terms in a different order; it does not change the values of the fixed free indices. Hence
\[
 \mathscr A-\mathcal E =2\sum_{\operatorname{cyc}(x,y,z)}\sum_{r,\ell} \bigl(C_{xwr}C_{yr\ell}-C_{ywr}C_{xr\ell}\bigr)B_{uvz\ell} =2\sum_{\operatorname{cyc}(x,y,z)}\sum_\ell K_{xyw\ell}B_{uvz\ell} =2\mathcal T_{uvxyzw}.
\]
The second equality uses \eqref{eq:Kcomponents}. In the last equality the summation index $p$ in the definition of $\mathcal T$ has simply been renamed $\ell$. Therefore $\mathscr A=2\mathcal T+\mathcal E$.

We next expand the last two terms and identify their cross-contractions.
Again fixing one choice of $(x,y,z)$, expand the last two $\mathcal F$ terms and group them according to whether $B$ contains $w$ or $u,v$:
\[
 \begin{aligned}
 &C_{uxr}\mathcal F_{vwyzr}+C_{vxr}\mathcal F_{uwyzr}=\sum_\ell\bigl[(C_{uxr}C_{vy\ell}+C_{vxr}C_{uy\ell})B_{wr\ell z}-(C_{uxr}C_{vz\ell}+C_{vxr}C_{uz\ell})B_{wr\ell y}\\
 &\quad-C_{uxr}C_{wy\ell}B_{vr\ell z}+C_{uxr}C_{wz\ell}B_{vr\ell y}-C_{vxr}C_{wy\ell}B_{ur\ell z}+C_{vxr}C_{wz\ell}B_{ur\ell y}\bigr].
 \end{aligned}
\]
The index $r$ in this identity is temporarily fixed and remains to be summed in the outer sum. Denote the double sum and cyclic sum of the last four terms by
\[
 \mathcal J:=\sum_{\operatorname{cyc}(x,y,z)}\sum_{r,\ell} \bigl(-C_{uxr}C_{wy\ell}B_{vr\ell z} +C_{uxr}C_{wz\ell}B_{vr\ell y} -C_{vxr}C_{wy\ell}B_{ur\ell z} +C_{vxr}C_{wz\ell}B_{ur\ell y}\bigr).
\]
Interchanging the two summation indices $r,\ell$ throughout these terms and using symmetry of $B$, we obtain
\[
 \mathcal J=\sum_{\operatorname{cyc}(x,y,z)}\sum_{r,\ell} \bigl[-C_{ywr}(C_{ux\ell}B_{vr\ell z} +C_{vx\ell}B_{ur\ell z}) +C_{zwr}(C_{ux\ell}B_{vr\ell y} +C_{vx\ell}B_{ur\ell y})\bigr].
\]
Collect the three cyclic terms containing $C_{ywr}$ and those containing $C_{zwr}$ into expressions containing $C_{xwr}$. Using the same cyclic rearrangement as in the identification of $\mathscr A-\mathcal E$, we obtain
\[
 \mathcal J=\sum_{\operatorname{cyc}(x,y,z)}\sum_{r,\ell} C_{xwr}\bigl(-C_{uz\ell}B_{vr\ell y}+C_{uy\ell}B_{vr\ell z} -C_{vz\ell}B_{ur\ell y}+C_{vy\ell}B_{ur\ell z}\bigr) =-\mathcal E.
\]

Finally, we use symmetry to cancel the remaining cyclic terms.
Fix $u,v,w$ and put
\[
 H_{xyz}:=\sum_{r,\ell}
       (C_{uxr}C_{vy\ell}+C_{vxr}C_{uy\ell})B_{wr\ell z}.
\]
After summation over $r$, the first two terms inside the brackets in the expansion of $C_{uxr}\mathcal F_{vwyzr}+C_{vxr}\mathcal F_{uwyzr}$ give exactly $H_{xyz}-H_{xzy}$. Interchanging the summation indices $r,\ell$ in $H_{yxz}$, and then using $B_{w\ell rz}=B_{wr\ell z}$ and commutativity of the scalar factors, gives
\[
 H_{yxz} =\sum_{r,\ell} (C_{uy\ell}C_{vxr}+C_{vy\ell}C_{uxr})B_{wr\ell z} =H_{xyz}.
\]
Consequently the cyclic sum of this part is
\[
 \sum_{\operatorname{cyc}(x,y,z)}(H_{xyz}-H_{xzy}) =H_{xyz}+H_{yzx}+H_{zxy} -H_{xzy}-H_{yxz}-H_{zyx}=0.
\]
It follows that $\mathscr D=\mathcal J=-\mathcal E$. Together with $\mathscr A=2\mathcal T+\mathcal E$, already proved, this shows that the right side of \eqref{eq:CF} is
\[
 \mathscr A+\mathscr D
 =(2\mathcal T+\mathcal E)-\mathcal E=2\mathcal T.
\]
This proves the identity using only complete symmetry and rearrangement of the indicated sums.
\end{proof}

We now impose angular compatibility.
\begin{corollary}[Curvature compatibility of quartic contraction operators]\label{cor:quartic-contractions-compatible}
If $C\in\Sym^3V^*$, $B\in\Sym^4V^*$ satisfy $\mathcal F_B(C)=0$, then for every pair $U,V\in V$ the self-adjoint operator $B_{UV}$ is $K$-compatible, where $K(X,Y)=[C_X,C_Y]$. That is, for arbitrary $x,y,z\in V$,
\[
 K(x,y)B_{UV}z+K(y,z)B_{UV}x+K(z,x)B_{UV}y=0.
\]
\end{corollary}
\begin{proof}
In \eqref{eq:CF}, the condition $\mathcal F=0$ makes every term on the right vanish. Thus $2\mathcal T=0$, and hence $\mathcal T=0$. By the interpretation of $\mathcal T$ at the start of this subsection, for fixed basis vectors $e_u,e_v$ the corresponding vector-valued cyclic sum has zero inner product with every $e_w$. That vector-valued cyclic sum is therefore zero.

Now write $U=\sum_uU^ue_u$, $V=\sum_vV^ve_v$, so that $B_{UV}=\sum_{u,v}U^uV^vB_{uv}$. Taking this linear combination of the identities for contractions by basis vectors proves the assertion. Self-adjointness already follows from complete symmetry of $B$.
\end{proof}

\begin{remark}[Scope of the contraction identity and the compatibility conclusion]
Only the corollary uses $\mathcal F=0$. The calculation requires no triviality of the radical, irreducibility, or analyticity: it is pointwise in a fixed inner product space for completely symmetric $C,B$. It proves neither the converse implication nor commutation of $B_{UV}$ with each $K(x,y)$.
\end{remark}

\subsection{Jordan closure of compatible operators}\label{subsec:paper-10-3}\label{sec:jordan}
We apply the Jordan closure theorem of Mantica--Molinari \cite[arXiv version, Theorem 4.2]{mm} to the compatible contractions obtained in Subsection \ref{sec:CF}, recording its polynomial and spectral consequences.

Continue to fix the inner product space $(V,h)$ and the algebraic curvature tensor $K$, and write
\[
 \mathcal J_K=\{A=A^*:A\text{ is }K\text{-compatible}\}.
\]
The vector space $\mathcal J_K$ contains $\R\Id$. Compatibility retains the cyclic meaning of Definition \ref{def:K-compatible}, not commutation with curvature operators.

\begin{definition}[Unital real Jordan subalgebra]\label{def:unital-jordan-subalgebra}
On the space $\Sym(V,h)$ of self-adjoint operators define
\[
 A\circ D:=\frac{AD+DA}{2}.
\]
A real linear subspace $\mathcal A\subset\Sym(V,h)$ containing $\Id$ and closed under $\circ$ is called a \textbf{unital real Jordan subalgebra}.
\end{definition}
Unlike the ordinary product $AD$, the Jordan product of self-adjoint operators is always self-adjoint.

The following statement is \cite[arXiv version, Theorem 4.2]{mm}. For a self-adjoint operator $A$, set $a(X,Y)=h(AX,Y)$; raising one index of $a$ recovers $A$. With $K_{xywp}=h(K(x,y)p,w)$, the authors' compatibility equation (17) is \eqref{eq:Kcompatible}, and their symmetrized product of symmetric tensors corresponds to $A\circ D=(AD+DA)/2$. No commutation hypothesis on $A,D$ is imposed. Compatibility is linear in $K$, so an overall change of curvature sign does not affect this identification.
\begin{theorem}[Jordan closure of compatible self-adjoint operators]\label{thm:jordan-closure}
If $A,D\in\mathcal J_K$, then $A\circ D\in\mathcal J_K$. Hence $\mathcal J_K$ is a unital Jordan subalgebra of the space of self-adjoint operators.
\end{theorem}

\begin{definition}[Jordan subalgebra generated by a family of operators]\label{def:generated-jordan-algebra}
Given a family $\mathcal S$ of self-adjoint operators, the intersection of all real Jordan subalgebras containing $\Id$ and $\mathcal S$ is called the \textbf{unital Jordan subalgebra generated by} $\Id,\mathcal S$.
\end{definition}

Jordan closure preserves all powers of a single operator. Combined with the self-adjoint spectral theorem and Lagrange interpolation, it gives the following standard conclusion; see \cite[arXiv version, Theorem 4.2 and the paragraph following its proof]{mm} and \cite[Section 0.9.11 and Corollary 2.5.11(a)]{horn-johnson2013}.
\begin{corollary}[Polynomials and spectral projections of a single compatible operator]\label{cor:compatible-polynomial-projector}
Let $A\in\mathcal J_K$. Then every real polynomial $p(A)$ belongs to $\mathcal J_K$; the orthogonal projection onto any eigenspace of $A$ also belongs to $\mathcal J_K$. These operators also belong to the unital Jordan subalgebra generated by $\Id,A$.
\end{corollary}
Explicitly, if the distinct eigenvalues of $A$ are $\lambda_1,\ldots,\lambda_t$, the orthogonal projection onto the $\lambda_a$ eigenspace is
\[
 \Pi_a=\prod_{b\ne a}\frac{A-\lambda_b\Id}{\lambda_a-\lambda_b}.
\]
When $t=1$, the empty product is understood as $\Id$. The rank of the projection equals the dimension of the corresponding eigenspace, so $\Pi_a$ has rank one only when $\lambda_a$ is simple.

These are statements about powers of one operator and its pointwise spectral projections. They require neither pairwise commutation nor smooth eigenvectors, and assert no closure under ordinary products of different operators.

Denote by $\mathcal A$ the smallest Jordan subalgebra containing $\Id$ and all $B_{UV}$. Angular compatibility and the closure theorem give
\[
 \mathcal A\subset\mathcal J_K.
\]
Polynomials and spectral projections of $A=B_{vv}$ therefore belong to $\mathcal A$ and remain compatible. A quartic extremum will supply a simple extreme eigenvalue, and hence a rank-one projection, before irreducibility is used.

\subsection{Extrema of symmetric quartic forms and rank-one spectral projections}\label{subsec:paper-11-1}\label{sec:extreme}

Let $V$ be a finite-dimensional real vector space, $d\ge1$, and $T\in\Sym^dV^*$. Its diagonal polynomial is
\[
 q(x):=T(\underbrace{x,\ldots,x}_{d\text{ arguments}}).
\]
Multilinearity makes $q$ a homogeneous polynomial of degree $d$. The following elementary coefficient argument shows that its diagonal values determine $T$.

\begin{lemma}[Diagonal values determine a symmetric tensor]\label{lem:paper-polarization}
We have $q\equiv0$ if and only if $T=0$. If $V$ is also equipped with an inner product, then vanishing of $q$ on the unit sphere already implies $T=0$.
\end{lemma}
\begin{proof}
If $T=0$, then $q\equiv0$. Conversely, suppose that $q\equiv0$. For arbitrary $v_1,\ldots,v_d\in V$, consider the polynomial
\[
 P(t_1,\ldots,t_d):=q(t_1v_1+\cdots+t_dv_d).
\]
By multilinearity, the coefficient of $t_1\cdots t_d$ in $P$ is
\[
 \sum_{\sigma\in S_d}T(v_{\sigma(1)},\ldots,v_{\sigma(d)})
 =d!\,T(v_1,\ldots,v_d),
\]
where $S_d$ is the permutation group of $\{1,\ldots,d\}$ and the equality uses complete symmetry of $T$. Since $P\equiv0$, this coefficient is zero. Thus $T(v_1,\ldots,v_d)=0$ for arbitrary arguments, and $T=0$.

Finally, if $q$ vanishes on the unit sphere, homogeneity gives
\[
 q(x)=|x|^d q\left(\frac{x}{|x|}\right)=0
 \qquad(x\ne0).
\]
Since $d\ge1$, we also have $q(0)=0$. Hence $q\equiv0$, and the first assertion gives $T=0$.
\end{proof}

We seek a rank-one projection in the generated Jordan algebra by finding a simple extreme eigenvalue of a quartic contraction and isolating it with an interpolation polynomial.

Continue to fix a finite-dimensional real Euclidean space $(V,h)$, let $0\ne B\in\Sym^4V^*$, and retain the notation
\[
 h(B_{UV}x,y)=B(U,V,x,y),\qquad |x|^2=h(x,x).
\]
These assumptions suffice for the extremum argument; $C$ and angular compatibility enter only in its curvature interpretation.

\begin{proposition}[Quartic extrema and a spectral gap for a contraction operator]\label{prop:quartic-extreme-spectrum}
There exist a sign $\varepsilon\in\{1,-1\}$, a unit vector $v$, and a real number $\beta>0$ such that, on putting
\[
 \widetilde B=\varepsilon B,\qquad
 q(x)=\widetilde B(x,x,x,x),\qquad
 A=\widetilde B_{vv}=\varepsilon B_{vv},
\]
we have
\[
 \beta=q(v)=\max_{|x|=1}q(x),\qquad Av=\beta v,
\]
and
\[
 h(Aw,w)\le\frac{\beta}{3}|w|^2
 \qquad(w\in v^\perp).
\]
Thus $\beta$ is the simple largest eigenvalue of $A$, and all other eigenvalues are at most $\beta/3$.
\end{proposition}
\begin{proof}
We first choose a positive spherical maximum.
Write $q_B(x)=B(x,x,x,x)$. By Lemma \ref{lem:paper-polarization}, $B\ne0$ implies that $q_B$ is not identically zero. Choose $x_0$ with $q_B(x_0)\ne0$; then $x_0\ne0$. Choose $\varepsilon=\pm1$ so that $\varepsilon q_B(x_0)>0$. By homogeneity of degree four,
\[
 q\left(\frac{x_0}{|x_0|}\right)
 =\frac{\varepsilon q_B(x_0)}{|x_0|^4}>0.
\]
The unit sphere in a finite-dimensional space is compact and $q$ is continuous, so it attains a maximum $\beta$ on the sphere. The preceding inequality ensures $\beta>0$. Take a unit vector $v$ at which this maximum is attained.

We next use the first variation to determine an eigenvector.
Complete symmetry of $\widetilde B$ gives
\[
 h(Ax,y)=\widetilde B(v,v,x,y)
       =\widetilde B(v,v,y,x)=h(x,Ay).
\]
Thus $A$ is self-adjoint. For an arbitrary nonzero $w\perp v$, take the curve
\[
 x(t)=\cos(|w|t)v+\sin(|w|t)\frac{w}{|w|}.
\]
Since $|v|=1$ and $w\perp v$, we have $|x(t)|=1$ and
\[
 x(0)=v,\qquad x'(0)=w,\qquad x''(0)=-|w|^2v.
\]
Put $f(t)=q(x(t))$. The function $f$ has a maximum at $t=0$, so $f'(0)=0$. Differentiating in each of the four argument slots and using complete symmetry gives
\[
 f'(t)=4\widetilde B(x'(t),x(t),x(t),x(t)).
\]
Hence
\[
 0=f'(0)=4\widetilde B(w,v,v,v)=4h(Av,w).
\]
This holds for all $w\perp v$, including $w=0$, and therefore $Av$ is parallel to $v$. Moreover,
\[
 h(Av,v)=\widetilde B(v,v,v,v)=\beta
\]
together with $|v|=1$, gives $Av=\beta v$.

Finally, the second variation gives a strict spectral gap.
Differentiate the preceding expression for $f'(t)$ once more. Differentiation in the first slot produces one second-variation term, while differentiation in the remaining three slots produces three equal first-variation terms. Thus
\[
 f''(t) =4\widetilde B(x''(t),x(t),x(t),x(t)) +12\widetilde B(x'(t),x'(t),x(t),x(t)).
\]
The coefficients here are $4$ and $4\cdot3=12$, respectively. Setting $t=0$ and using $f''(0)\le0$, we obtain
\[
 0\ge f''(0) =-4|w|^2\widetilde B(v,v,v,v) +12\widetilde B(w,w,v,v) =-4\beta|w|^2+12h(Aw,w).
\]
Hence $h(Aw,w)\le(\beta/3)|w|^2$. The inequality also holds for $w=0$, without using a curve involving $w/|w|$.

Self-adjointness and $Av=\beta v$ also give
\[
 h(Aw,v)=h(w,Av)=\beta h(w,v)=0
 \qquad(w\perp v).
\]
Thus $v^\perp$ is invariant under $A$. By the self-adjoint spectral theorem, $A|_{v^\perp}$ has an orthonormal eigenbasis. If $w\in v^\perp$ is a nonzero eigenvector with eigenvalue $\lambda$, then
\[
 \lambda=\frac{h(Aw,w)}{|w|^2}\le\frac{\beta}{3}<\beta.
\]
The last strict inequality uses precisely $\beta>0$. Thus, in the orthogonal invariant decomposition $V=\R v\oplus v^\perp$, the eigenvalue $\beta$ occurs only on the one-dimensional subspace $\R v$, and the rest of the spectrum is at most $\beta/3$. Since a self-adjoint operator is diagonalizable, the algebraic and geometric multiplicities of $\beta$ are both one. If other eigenvalues exist, their distance from $\beta$ is at least $2\beta/3$. This proves the proposition.
\end{proof}

If $\varepsilon=1$, the original contraction $B_{vv}$ has simple largest eigenvalue $\beta$; if $\varepsilon=-1$, the original contraction $B_{vv}=-A$ has simple smallest eigenvalue $-\beta$. Thus, when retaining the original $B$, the precise conclusion is that some contraction has a simple \emph{extreme} eigenvalue.

We extract the projection onto this simple eigenspace by interpolation.
\begin{corollary}[A rank-one orthogonal projection in the generated algebra]\label{cor:rank-one-projector}
Denote by $\mathcal A$ the smallest real Jordan subalgebra containing $\Id$ and all $B_{UV}$. Then $\mathcal A$ contains a rank-one orthogonal projection. More precisely, for $v,A$ above, the orthogonal projection $P_v$ onto $\R v$ is a real polynomial in $A$, and also a real polynomial in the original contraction $B_{vv}$.
\end{corollary}
\begin{proof}
Let $\lambda_1,\ldots,\lambda_t$ be the distinct eigenvalues of $A$ other than $\beta$. By Proposition \ref{prop:quartic-extreme-spectrum}, each $\lambda_j\le\beta/3<\beta$, so all denominators below are nonzero. Define
\begin{equation}\label{eq:projector}
 p(s)=\prod_{j=1}^t\frac{s-\lambda_j}{\beta-\lambda_j},\qquad
 p(A)=P_v=v\otimes v^\flat.
\end{equation}
Here $v^\flat(z)=h(v,z)$, and hence $(v\otimes v^\flat)z=h(v,z)v$. We verify the operator identity in this formula.

The polynomial satisfies $p(\beta)=1$, $p(\lambda_j)=0$. If $Az=\lambda z$, then $A^rz=\lambda^rz$ for every integer $r\ge0$, and therefore
\[
 p(A)z=p(\lambda)z.
\]
By the proposition, the $\beta$ eigenspace is exactly $\R v$. By the self-adjoint spectral theorem, the other eigenspaces have orthogonal direct sum $v^\perp$. Thus $p(A)$ is the identity on $\R v$ and zero on $v^\perp$, and is precisely $P_v$. It satisfies
\[
 P_v^2=P_v,\qquad P_v^*=P_v,\qquad\rank P_v=1.
\]
If there are no other eigenvalues, then $V=\R v$ and $t=0$. The empty product is taken to be $1$, and again $p(A)=\Id=P_v$.

We now explain why $P_v$ belongs to $\mathcal A$. Since $B_{vv}\in\mathcal A$ and $\mathcal A$ is a real linear space, $A=\varepsilon B_{vv}\in\mathcal A$. From $\Id\in\mathcal A$ and
\[
 A\circ A^r=\frac{AA^r+A^rA}{2}=A^{r+1},
\]
induction places all ordinary powers of $A$ in $\mathcal A$. Taking real linear combinations then gives $p(A)=P_v\in\mathcal A$. This is exactly the polynomial closure for a single operator in Section \ref{sec:jordan}; it uses only powers of the same operator, not closure under ordinary products of different operators.

Finally, put $r(s)=p(\varepsilon s)$. Then
\[
 P_v=p(A)=p(\varepsilon B_{vv})=r(B_{vv}).
\]
Thus $P_v$ is also a polynomial in $B_{vv}$. Since $\widetilde B_{UV}=\varepsilon B_{UV}$, both contraction families generate the same unital real Jordan algebra.
\end{proof}

If, in addition, $C\in\Sym^3V^*$ is given and $\mathcal F_B(C)=0$, Section \ref{sec:CF} makes every $B_{UV}$ compatible with $K=[C,C]$, while Jordan closure in Section \ref{sec:jordan} gives
\[
 \mathcal A\subset\mathcal J_K.
\]
Hence the projection $P_v\in\mathcal A$ constructed here is also $K$-compatible. The angular compatibility identity is linear in $B$, so $\mathcal F_{\widetilde B}(C)=\varepsilon\mathcal F_B(C)=0$; curvature compatibility of an operator is linear as well. The sign choice above therefore preserves both kinds of compatibility.

Only $B\ne0$, complete symmetry, and positive definiteness of $h$ are needed here. No sign condition on $B$ is imposed, and triviality of the radical, irreducibility, and analyticity are not assumed. The construction gives one suitable direction $v$, not a simple eigenvalue for every contraction, and makes no smooth choice as the point varies.

\subsection{Rank-one generation of an irreducible Jordan subalgebra}\label{subsec:paper-11-2}\label{sec:generation}
We prove that a rank-one projection and irreducibility force a unital real Jordan subalgebra to contain all self-adjoint operators. In the application, Section \ref{sec:paper-7} supplies irreducibility, while Subsections \ref{subsec:paper-10-1}--\ref{subsec:paper-10-3} preserve curvature compatibility.

Let $(V,h)$ be a finite-dimensional real Euclidean space, and let $\mathcal A\subset\Sym(V,h)$ be a unital real Jordan subalgebra. Thus it is a real linear subspace of the self-adjoint operators, contains $\Id$, and is closed under
\[
 D\circ H=\frac{DH+HD}{2}
\]
The action of $\mathcal A$ on $V$ is called \textbf{irreducible} if every linear subspace $W\subset V$ satisfying
\[
 D(W)\subset W\qquad\text{for every }D\in\mathcal A,
\]
is either $W=\{0\}$ or $W=V$.

For $a,b\in V$, write $b^\flat(x)=h(b,x)$ and define
\[
 H_{ab}:=a\otimes b^\flat+b\otimes a^\flat,
 \qquad H_{ab}x=h(b,x)a+h(a,x)b.
\]
This is bilinear in $a,b$ and satisfies $H_{ab}=H_{ba}$ and $H_{ab}^*=H_{ab}$. For a fixed unit vector $v$, the orthogonal projection onto $\R v$ is
\[
 P_v=v\otimes v^\flat,\qquad H_{vv}=2P_v.
\]

\begin{theorem}[An irreducible Jordan subalgebra containing a rank-one projection]\label{thm:irreducible-jordan-generation}
Under the preceding assumptions, if $\mathcal A$ acts irreducibly on $V$ and $P_v\in\mathcal A$, then
\begin{equation}\label{eq:allSym}
 \mathcal A=\Sym(V,h).
\end{equation}
\end{theorem}
\begin{proof}
We first construct a nonzero subspace from the projection.
Define
\[
 W:=\{w\in V:H_{vw}\in\mathcal A\}.
\]
Since $w\mapsto H_{vw}$ is linear and $\mathcal A$ is a real linear space, $W$ is a linear subspace. More explicitly, if $w,z\in W$ and $\alpha,\gamma\in\R$, then
\[
 H_{v,\alpha w+\gamma z}=\alpha H_{vw}+\gamma H_{vz}\in\mathcal A.
\]
Moreover, $H_{vv}=2P_v\in\mathcal A$, so $v\in W$ and $W\ne\{0\}$. It remains to prove invariance.

We next prove that $H_{ab}\in\mathcal A$ whenever $a,b\in W$.
Put $E:=W\cap v^\perp$ and take $w\in E$. Since $|v|=1$ and $h(v,w)=0$,
\[
 H_{vw}v=w,\qquad H_{vw}w=|w|^2v.
\]
Hence, for any $x\in V$,
\[
 H_{vw}^2x =H_{vw}\bigl(h(w,x)v+h(v,x)w\bigr) =h(w,x)w+|w|^2h(v,x)v,
\]
that is,
\[
 H_{vw}^2=|w|^2P_v+w\otimes w^\flat.
\]
The condition $w\in W$ gives $H_{vw}\in\mathcal A$, and Jordan closure gives $H_{vw}^2=H_{vw}\circ H_{vw}\in\mathcal A$. Thus
\[
 Q_w:=w\otimes w^\flat=H_{vw}^2-|w|^2P_v\in\mathcal A.
\]
Here $Q_w$ is defined for every $w\in E$ and is zero when $w=0$. For $w\ne0$, it is $|w|^2$ times the orthogonal projection onto $\R w$.

If $w,z\in E$, then $w+z\in E$, so $Q_{w+z},Q_w,Q_z$ all belong to $\mathcal A$. Expansion and subtraction give
\[
 Q_{w+z}-Q_w-Q_z =w\otimes z^\flat+z\otimes w^\flat =H_{wz}\in\mathcal A.
\]
At this stage the polarization applies to two vectors in $E$.

Now take arbitrary $a,b\in W$. Since $v\in W$ and $|v|=1$, we may write
\[
 a=\alpha v+w,\qquad b=\gamma v+z,\qquad w,z\in E.
\]
Indeed, take $\alpha=h(a,v)$ and $\gamma=h(b,v)$. By bilinearity,
\[
 H_{ab}=2\alpha\gamma P_v+\alpha H_{vz}+\gamma H_{vw}+H_{wz}.
\]
Every term on the right belongs to $\mathcal A$, so $H_{ab}\in\mathcal A$ for all $a,b\in W$.

We now prove that $W$ is a common invariant subspace.
Take any $D\in\mathcal A$. First we check $Dv\in W$. Since $D$ is self-adjoint, for any $x\in V$,
\[
 (DP_v+P_vD)x =h(v,x)Dv+h(v,Dx)v =h(v,x)Dv+h(Dv,x)v =H_{v,Dv}x.
\]
Therefore
\[
 H_{v,Dv}=DP_v+P_vD=2D\circ P_v\in\mathcal A,
\]
and the definition of $W$ gives $Dv\in W$.

Next take any $w\in W$. Again using self-adjointness of $D$, termwise expansion gives
\[
 (DH_{vw}+H_{vw}D)x =h(w,x)Dv+h(v,x)Dw+h(Dw,x)v+h(Dv,x)w =(H_{Dv,w}+H_{v,Dw})x.
\]
Hence
\[
 H_{v,Dw}=2D\circ H_{vw}-H_{Dv,w}.
\]
The first term belongs to $\mathcal A$ because $D,H_{vw}\in\mathcal A$. The second also belongs to $\mathcal A$: we have just proved $Dv\in W$, and the assertion $H_{ab}\in\mathcal A$ proved above applies to $a=Dv,b=w\in W$. Consequently $H_{v,Dw}\in\mathcal A$, that is, $Dw\in W$.

Thus $W\ne\{0\}$ is invariant under every $D\in\mathcal A$, and irreducibility gives $W=V$. Only squares and symmetrized products have been used; ordinary products $DH$ need not belong to $\mathcal A$.

Finally, we write a spanning formula for all self-adjoint operators.
Combining $H_{ab}\in\mathcal A$ for $a,b\in W$ with $W=V$ gives $H_{ab}\in\mathcal A$ for all $a,b\in V$. Choose an orthonormal basis $e_1,\ldots,e_m$. If $E_{ij}$ denotes the matrix unit sending $e_j$ to $e_i$ and all other basis vectors to zero, then
\[
 H_{e_i e_i}=2E_{ii},\qquad
 H_{e_i e_j}=E_{ij}+E_{ji}\quad(i\ne j).
\]
The matrix $(d_{ij})$ of an arbitrary self-adjoint operator $D$ is symmetric, and hence
\[
 D=\frac12\sum_i d_{ii}H_{e_i e_i}
     +\sum_{i<j}d_{ij}H_{e_i e_j}\in\mathcal A.
\]
Thus $\Sym(V,h)\subset\mathcal A$. Together with the original reverse inclusion this gives \eqref{eq:allSym}. No classification of Jordan algebras has been used.
\end{proof}

For the algebra generated by $\Id$ and all $B_{UV}$, irreducibility follows from that of the contraction family: any invariant subspace for the algebra is invariant for every generator.

The equivalent formulation obtained in Section \ref{sec:commutant} at generic points of the open set of trivial radical is
\[
 \mathcal Z_B:=\{Q=Q^*:QB_{UV}=B_{UV}Q\text{ for all }U,V\}
 =\R\Id.
\]
Proposition \ref{prop:commutant-irreducibility} therefore gives irreducibility of the contraction family and its generated algebra at these generic points.

\subsection{The algebraic curvature isotropy theorem}\label{subsec:paper-11-3}\label{sec:isotropy}
We combine the preceding algebraic results. Fix a finite-dimensional Euclidean space $(V,h)$, with $m=\dim V\ge2$, and $C\in\Sym^3V^*$. As before, by
\[
 h(C_XY,Z)=C(X,Y,Z)
\]
we define the self-adjoint operator $C_X$, and then put
\[
 K(X,Y):=[C_X,C_Y].
\]
Complete symmetry gives $C_XY=C_YX$, and Subsection \ref{subsec:paper-10-1} establishes the algebraic curvature symmetries of $K$. Write
\[
 R_0(X,Y)Z=h(Y,Z)X-h(X,Z)Y.
\]

Recall the definition from Subsection \ref{subsec:paper-10-1}: a self-adjoint operator $A$ is $K$-compatible if, for arbitrary $x,y,z\in V$,
\[
 K(x,y)Az+K(y,z)Ax+K(z,x)Ay=0.
\]
This cyclic compatibility identity, rather than commutation with each $K(x,y)$, is used below.

We use the characterization of an algebraic curvature tensor compatible with every symmetric tensor in \cite[arXiv version, Proposition 4.7]{mm}. Under the convention $K_{xywp}=h(K(x,y)p,w)$, its equation (22) is precisely the form \eqref{eq:isotropy}, with its scalar coefficient denoted here by $s$.
\begin{proposition}[Curvature isotropy when all self-adjoint operators are compatible]\label{prop:all-compatible-isotropic}
With $(V,h)$, $m\ge2$, and $K=[C,C]$ as in this subsection, if every self-adjoint operator is $K$-compatible, then there exists $s\in\R$ such that
\begin{equation}\label{eq:isotropy}
 K(X,Y)Z=s\{h(Y,Z)X-h(X,Z)Y\}
 \qquad(X,Y,Z\in V).
\end{equation}
\end{proposition}

\begin{theorem}[Algebraic curvature isotropy theorem]\label{thm:algebra}
Let $B\in\Sym^4V^*$ be nonzero and satisfy $\mathcal F_B(C)=0$. Suppose that the contraction family
\[
 \{B_{uv}:u,v\in V\}
\]
has only scalar operators in its self-adjoint common commutant, that is,
\[
 \{Q=Q^*:QB_{uv}=B_{uv}Q\ \text{for all }u,v\}
 =\mathbb R\Id.
\]
Then there exists $s\in\mathbb R$ for which \eqref{eq:isotropy} holds.
\end{theorem}
\begin{proof}
Denote by $\mathcal A$ the smallest real Jordan subalgebra generated by $\Id$ and all $B_{uv}$.

The contraction identity of Subsection \ref{subsec:paper-10-2} shows that, when $\mathcal F_B(C)=0$, every $B_{uv}$ satisfies the compatibility identity, that is, $B_{uv}\in\mathcal J_K$. By the closure theorem of Mantica--Molinari \cite[arXiv version, Theorem 4.2]{mm}, the $K$-compatible self-adjoint operators are closed under the Jordan product
\[
 D\circ H=\frac{DH+HD}{2}
\]
and $\Id$ is compatible as well. Hence
\[
 \mathcal A\subset\mathcal J_K.
\]

The pointwise criterion of Subsection \ref{subsec:paper-9-2} says that the assumption that the self-adjoint common commutant consists only of scalars is equivalent to irreducibility of the contraction family $B_{uv}$ on $V$. If $W\subset V$ is preserved by every operator in $\mathcal A$, it is in particular preserved by every generator $B_{uv}$. Thus $\mathcal A$ is irreducible as well.

Since $B\ne0$, the quartic extremum argument of Subsection \ref{subsec:paper-11-1} produces a unit vector $v$ for which the rank-one orthogonal projection
\[
 P_v=v\otimes v^\flat
\]
belongs to $\mathcal A$. Apply Theorem \ref{thm:irreducible-jordan-generation}: an irreducible unital Jordan subalgebra containing a rank-one projection must be
\[
 \mathcal A=\Sym(V,h).
\]
It follows that every $P_z$ belongs to $\mathcal A$ and is therefore $K$-compatible. The hypotheses of Proposition \ref{prop:all-compatible-isotropic} are consequently satisfied, and \eqref{eq:isotropy} follows.
\end{proof}

No sign condition on $B$ is required: nonvanishing is used in the extremum construction, and replacing $B$ by $-B$ preserves both the generated real Jordan algebra and $\mathcal F_B(C)=0$.

In the Hessian geometry considered above,
\[
 R^h=R_0-K,
\]
so \eqref{eq:isotropy} is equivalent to
\[
 R^h=(1-s)R_0.
\]
This is pointwise isotropy; $s$ may vary along the indicatrix. Section \ref{sec:quartic-rigidity} extends the identity analytically and then proves constancy.

\section{Analytic rigidity of twisted Codazzi quartic forms on the indicatrix}\label{sec:quartic-rigidity}\label{sec:schur}\label{subsec:paper-12-2}
We prove Theorem \ref{thm:quartic-rigidity-intro} on a fixed Minkowski space. The nonfull-rank branch gives $C=0$ on an open set, and the full-rank branch gives curvature isotropy there. Analytic continuation, followed in the latter case by Schur's lemma and constant-curvature normalization, gives $R^h=R_0$; the cone formula then gives Hessian flatness.

\begin{definition}[Curvature isotropy and its defect tensor]\label{def:curvature-isotropy-defect}
Let $(S,h)$ be a Riemannian manifold of dimension $m\ge2$, and put
\[
 R_0(X,Y)Z=h(Y,Z)X-h(X,Z)Y.
\]
Define the normalized scalar curvature $\kappa_0$ and the \textbf{isotropy defect tensor} $\mathcal Q$ by
\begin{equation}\label{eq:isotropy-defect}
 \kappa_0:=\frac{\Scal(h)}{m(m-1)},\qquad
 \mathcal Q:=R^h-\kappa_0R_0.
\end{equation}
We say that the curvature of $h$ is \textbf{isotropic} at a point $p$ if there is a real number $\kappa(p)$ such that $R^h_p=\kappa(p)R_{0,p}$. We say that $h$ has \textbf{constant sectional curvature} if there is a constant $k$, independent of the point, such that $R^h=kR_0$.
\end{definition}

If $R^h_p=\kappa(p)R_{0,p}$, contraction gives
\[
 \Ric_p=(m-1)\kappa(p)h_p,\qquad
 \Scal_p=m(m-1)\kappa(p),
\]
and hence $\kappa(p)=\kappa_0(p)$. Thus $\mathcal Q_p=0$ characterizes isotropy at that point; it does not yet require $\kappa_0$ to be independent of the point.

The following is the tensor form of the identity theorem for real-analytic functions; see \cite[Section 4.1, p.~83]{krantz-parks2002}. Here $\mathcal Q$ is obtained from the analytic metric and its inverse matrix by finitely many differentiations, multiplications, and contractions. It is therefore an analytic tensor, and the identity theorem applies to its local coordinate components.
\begin{lemma}[Analytic continuation of the isotropy identity]\label{lem:analytic-isotropy-extension}
Let $(S,h)$ be a connected real-analytic Riemannian manifold, with $m=\dim S\ge2$. If the curvature of $h$ is isotropic on a nonempty open set $U\subset S$, then
\[
 \mathcal Q\equiv0,\qquad R^h=\kappa_0R_0
 \quad\text{on all of }S.
\]
\end{lemma}

We use the classical Schur lemma below; one may also take the Riemannian case of \cite[Lemma 7.1.1]{chern-shen-2005}.
\begin{lemma}[Schur's constancy conclusion and its dimension hypothesis]\label{lem:schur-constant-curvature}
Let $(S,h)$ be a connected smooth Riemannian manifold, with $m=\dim S\ge3$. If, on all of $S$, we have $R^h=\kappa R_0$, then $\kappa$ is constant.
\end{lemma}

When $m=2$, the dimension hypothesis of Schur's lemma is not satisfied, and one cannot use it to conclude that $\kappa$ is constant. The curvature tensor of a two-dimensional surface can always be written as $\kappa_0(p)R_0$. Lemma \ref{lem:analytic-isotropy-extension} on analytic continuation still applies, but does not imply constant Gaussian curvature.

\begin{proof}[Proof of Theorem \ref{thm:quartic-rigidity-intro}]
Let $V,F$ satisfy the hypotheses of the theorem, and put
\[
 n=\dim V\ge4,\qquad E=\tfrac12F^2,\qquad g=D^2E,\qquad
 \Sigma=\{y\in V:F(y)=1\},\qquad m=n-1\ge3.
\]
Let $h=g|_{T\Sigma}$, and let $C$ be the angular Cartan cubic tensor. Suppose that
\[
 B\in C^\infty(\Sigma,\Sym^4T^*\Sigma),\qquad B\not\equiv0,
\]
and that $B$ satisfies the global twisted Codazzi identity in the theorem.

Lemma \ref{lem:quartic-algebraic-compatibility} gives $\mathcal F_B(C)=0$ without any additional hypothesis on $B$. All index operations and contractions use $h$.

\textbf{Analyticity of the indicatrix and the maximal-rank open set.}
The function $F$ is real-analytic on the entire $V\setminus\{0\}$. For any $y\in\Sigma$, Euler's identity gives
\[
 \dd F_y(y)=F(y)=1.
\]
Thus $\dd F_y\ne0$, and $\Sigma$ is a real-analytic hypersurface. Since $g=D^2E$ and the full Cartan tensor $\mathcal C=\tfrac12D^3E$ are real-analytic, their tangential restrictions $h,C$ are also real-analytic. On the other hand, radial projection gives a diffeomorphism between $\Sigma$ and the standard topological sphere $S^{n-1}$, so that $\Sigma$ is compact and connected. The analytic continuation below is applied only to $C$ or to curvature tensors constructed from $h$; it does not require $B$ to be real-analytic.

We use the insertion rank of Definition \ref{def:quartic-radical-rank}, together with the maximal insertion rank $q$ and maximal-rank set $U$ defined in \eqref{eq:maximal-rank-set}. Since $B\not\equiv0$, Proposition \ref{prop:maximal-rank-bundles} gives $0<q\le m$, and $U$ is nonempty and open.

\textbf{First case: $0<q<m$.}
If $m$ is even, this case is already excluded by Theorem \ref{thm:root-rank-gap}. We may therefore assume that $m$ is odd. The entire strongly convex indicatrix has the centroaffine geodesic properties required in the nonfull-rank branch. The smooth symmetric tensor $B$ satisfies both twisted Codazzi and angular compatibility, so Corollary \ref{cor:degenerate-branch-cartan-zero} gives
\[
 C=0\qquad\text{on the nonempty open set }U.
\]
Since $C$ is a real-analytic tensor on the connected real-analytic manifold $\Sigma$, application of the real-analytic identity theorem to its local coordinate components gives
\[
 C\equiv0\qquad\text{on all of }\Sigma.
\]
Here we use vanishing on an open set, rather than at a single point. The Gauss equation \eqref{eq:gauss} of the indicatrix then gives
\[
 R^h=R_0\qquad\text{on all of }\Sigma.
\]
This settles the branch without excluding it or requiring a point of full rank.

\textbf{Second case: $q=m$.}
Use the full-rank set $\Omega$ defined in \eqref{eq:full-rank-locus}, the self-adjoint common commutant $\mathcal Z_B$ defined in \eqref{eq:selfadjoint-commutant}, and the set $G$ defined in \eqref{eq:irreducible-locus}. Since $q=m$, we have $\Omega=U$, so $\Omega$ is nonempty. The smoothness, complete symmetry, and global twisted Codazzi condition established above verify the hypotheses of the generic irreducibility theorem \ref{thm:irreducible}. It follows that $G$ is open and dense in $\Omega$.
Since $\Omega$ is nonempty and $G$ is dense in it, $G$ is also a nonempty open set. For each $y\in G$, triviality of the radical ensures that $B_y\ne0$. We have already proved $\mathcal F_{B_y}(C_y)=0$, while $\mathcal Z_B(y)=\mathbb R\Id$ is precisely the irreducibility condition in the algebraic theorem on curvature isotropy, Theorem \ref{thm:algebra}. Applying that theorem, we obtain a real number $s(y)$ such that
\[
 [C_X,C_Y]Z
 =s(y)\{h(Y,Z)X-h(X,Z)Y\}
 \qquad(X,Y,Z\in T_y\Sigma).
\]
The Gauss equation now gives
\[
 R^h_y=(1-s(y))R_{0,y}\qquad(y\in G).
\]
The coefficient is not yet known to be constant. No regularity assertion about the pointwise choice $s(y)$ is needed: taking the scalar curvature trace gives $1-s(y)=\kappa_0(y)$ on $G$, so the analytic isotropy defect tensor $\mathcal Q$ defined in \eqref{eq:isotropy-defect} vanishes on $G$.

Since $h$ is real-analytic and $\Sigma$ is connected, Lemma \ref{lem:analytic-isotropy-extension} extends the preceding identity to
\[
 R^h=\kappa_0R_0\qquad\text{on all of }\Sigma.
\]
We now have $m\ge3$, so Schur's lemma \ref{lem:schur-constant-curvature} applies and shows that $\kappa_0$ is constant. Only at this stage is the constant-curvature hypothesis of Proposition \ref{prop:compact-indicatrix-normalization} satisfied. That proposition gives $\kappa_0=1$, and hence in this branch also we obtain
\[
 R^h=R_0\qquad\text{on all of }\Sigma.
\]

\textbf{Determination of the Hessian curvature from the indicatrix curvature.}
The two cases exhaust $0<q\le m$, and both show that the sectional curvature of $h$ is identically $1$. By Proposition \ref{prop:hessian-cone} and the cone structure \eqref{eq:cone}, we identify $g$ with the cone metric through the radial diffeomorphism given there. Here $\rho$ denotes the radial variable, to distinguish it from the maximal rank $q$. For angular lifts independent of $\rho$, denoted by $X,Y,Z$, the curvature of this cone metric is
\[
 R^g(X,Y)Z
 =R^h(X,Y)Z-\{h(Y,Z)X-h(X,Z)Y\}=0,
\]
and the curvature components involving a radial direction also vanish. Since the radial and angular directions span the entire tangent space, it follows that $R^g\equiv0$ on $V\setminus\{0\}$. This proves all the conclusions of the theorem.
\end{proof}

Analyticity extends $C=0$ in the nonfull-rank branch and the vanishing isotropy defect in the full-rank branch; $B$ need only be smooth. The restriction $n\ge4$ enters through Schur's lemma, while constant-curvature normalization requires only $m\ge2$.

The theorem proves Hessian flatness on a single Minkowski space independently of analytic Laugwitz rigidity. Combining the two later gives quadraticity of $F^2$.

The uses of regularity in the two rigidity arguments and their application in Part II are collected below. The geometric hypotheses at each cited step remain in force.
\begingroup
\small
\renewcommand{\arraystretch}{1.18}
\begin{longtable}{@{}>{\raggedright\arraybackslash}p{0.27\linewidth}>{\raggedright\arraybackslash}p{0.69\linewidth}@{}}
\toprule
Step & Regularity used and its role \\
\midrule
\endfirsthead
\toprule
Step & Regularity used and its role \\
\midrule
\endhead
Local induction and ideal division (Section \ref{sec:paper-3}) &
The real-analytic potential $f$ has a convergent Taylor expansion. A nonzero remaining germ has a finite leading degree; convergence turns coefficientwise membership in $\J_t$ into membership of the analytic germ. \\[3pt]
Continuation along a great circle (Section \ref{sec:periodic}) &
The analytic cubic $C$ and a fixed ambient normal define an analytic scalar $\widehat A$. The identity theorem extends its division-free differential equation around the circle. The final integration over a period requires only $C^2$ regularity of that scalar. \\[3pt]
Continuation from an open set (Section \ref{sec:quartic-rigidity}) &
Analyticity of the fibre norm makes $C$ and $h$ analytic. The identity theorem extends $C=0$ in the deficient-rank branch, and the vanishing $\mathcal Q=0$ of the isotropy defect tensor defined in \eqref{eq:isotropy-defect} in the full-rank branch, over the connected indicatrix. \\[3pt]
Quartic equations and transport (Sections \ref{sec:paper-6}, \ref{sec:paper-7}, and \ref{sec:paper-5}) &
Smoothness of $h,C,B$ suffices for the angular equations, the local bundle constructions and the maximal-rank continuation argument. Analyticity of the quartic $B$ is not used. \\[3pt]
The final open-set argument (Section \ref{sec:main}) &
Smoothness of $\mathcal B$ makes $\{x:\mathcal B(x,\cdot)\not\equiv0\}$ open. Once each fibre there is quadratic, smoothness of $g$ gives smooth Riemannian coefficients by evaluation along a local nonzero section. \\
\bottomrule
\end{longtable}
\endgroup
In particular, the proof requires neither an analytic structure on the base manifold nor fibrewise analyticity of mixed derivatives of $F$ involving base variables.

\part{Berwald rigidity of Landsberg metrics}\label{part:finsler}
We now allow the base point to vary, assuming only smoothness in the base variables. Section \ref{sec:finsler-foundations} collects the foundations, and Section \ref{sec:paper-5} verifies the twisted Codazzi equation for the angular Berwald quartic form. Section \ref{sec:paper-9} applies the two fibrewise rigidity theorems and concludes on the open set where the Berwald curvature is nonzero.

\section{Preliminaries on Finsler geometry}\label{sec:finsler-foundations}
We fix the normalizations and classical identities needed from Finsler geometry. For sprays, Berwald curvature, and Landsberg curvature, see \cite[Sections 4.2, 6.1, and 6.2]{shen-2001} and \cite[Sections 2.1 and 4.3]{chern-shen-2005}.

\begin{definition}[Finsler metric]\label{def:finsler-metric}
Let $M$ be a smooth manifold. A Finsler metric in this paper is a continuous function $F:TM\to[0,\infty)$ such that $F(x,0)=0$, and such that on $TM\setminus\{0\}$ it is smooth, strictly positive, and positively homogeneous of degree one in the fibre variable. The Hessian of the fibre energy $E=F^2/2$ is positive definite at each nonzero tangent vector. We do not require $F(x,-y)=F(x,y)$. At a fixed base point, $F_x=F(x,\cdot)$ is a strongly convex Minkowski norm in the sense of Section \ref{sec:paper-2}.
\end{definition}

This agrees with the definition of regularity on the entire slit tangent bundle in \cite[Section 1.1]{bcs} and \cite[Definition 1.2.1]{chern-shen-2005}. This section uses only smoothness; the fibrewise analyticity hypothesis of the main theorem is not needed for these classical results.

\subsection{The canonical spray and the Berwald and Landsberg tensors}\label{subsec:paper-1-3}
In local coordinates $x^1,\ldots,x^n$ on $M$, write $y=y^i\partial_{x^i}$, obtaining tangent bundle coordinates $(x^i,y^i)$. The Hessian in $g_{ij}(x,y)=E_{y^iy^j}(x,y)$ is taken only in the fibre variables.
The classical relation between energy variation and the canonical spray is as follows; see \cite[equations (4.27)--(4.31)]{shen-2001} and \cite[Section 3.1]{li-minkowski}.
\begin{proposition}[Energy variation and the canonical spray]\label{prop:canonical-spray}
Let $x(t)$ be a smooth curve with nonzero velocity. It is stationary for the energy functional with fixed endpoints if and only if
\begin{equation}\label{eq:spray}
 \ddot x^a+2G^a(x,\dot x)=0,\qquad
 2G^a=g^{ab}\bigl(E_{x^ky^b}y^k-E_{x^b}\bigr).
\end{equation}
The coefficients $G^a$ are smooth and positively homogeneous of degree two in $y$, and they define a spray on the slit tangent bundle by
\[
 \mathscr S=y^i\partial_{x^i}-2G^i\partial_{y^i}.
\]
\end{proposition}

\begin{definition}[Berwald curvature and Landsberg curvature]\label{def:berwald-landsberg-curvatures}
Indices following a comma denote fibre partial derivatives; for example, $G^a{}_{,i}=\partial G^a/\partial y^i$. Define
\[
 \mathcal B^a{}_{ijk}:=G^a{}_{,ijk}.
\]
This is the Berwald curvature. Put $E_a=E_{y^a}=g_{ab}y^b$, and define
\begin{equation}\label{eq:landsberg}
 L_{ijk}:=-\tfrac12E_a\mathcal B^a{}_{ijk}.
\end{equation}
This is the Landsberg curvature. A metric satisfying $\mathcal B\equiv0$ is called a Berwald metric; a metric satisfying $L\equiv0$ is called a Landsberg metric.
\end{definition}
For this normalization, see \cite[Definitions 6.1.1 and 6.2.1; equations (6.25)--(6.27)]{shen-2001}. The quadratic energy function in that book is denoted by $L$ and corresponds to our $F^2$. Here the letter $L$ denotes only the Landsberg tensor; the two uses should be distinguished.

For the equivalence between vanishing Berwald curvature and a quadratic spray, see \cite[Section 6.1]{shen-2001} and \cite[Section 4.3]{chern-shen-2005}.
\begin{proposition}[Tensoriality and the Berwald criterion]\label{prop:berwald-quadratic-spray}
The tensor $\mathcal B$ is a covariant three-tensor on the slit tangent bundle with values in the pullback tangent bundle, and is completely symmetric in its three lower indices; $L$ is a symmetric covariant three-tensor. If $n\ge2$, the following conditions are equivalent:
\begin{enumerate}[label=\textup{(\roman*)},leftmargin=2.3em,itemsep=2pt]
\item $\mathcal B\equiv0$;
\item On every base coordinate domain there are smooth functions depending only on $x$, with $\Gamma^a_{ij}=\Gamma^a_{ji}$, such that $G^a(x,y)=\tfrac12\Gamma^a_{ij}(x)y^iy^j$.
\end{enumerate}
In this case the coefficients define a torsion-free linear connection on $M$.
\end{proposition}

\begin{remark}[The Berwald and Landsberg conditions]
Equation \eqref{eq:landsberg} immediately gives $\mathcal B\equiv0\Rightarrow L\equiv0$. For this classical implication, see \cite[Proposition 2.1.3]{chern-shen-2005}. The converse is the Unicorn problem considered here. Our dimensional and analyticity assumptions were stated in Section \ref{sec:introduction}.

The Berwald property does not require the norm to arise from an inner product. For example, if in local coordinates $F(x,y)$ is independent of $x$, then \eqref{eq:spray} gives $G^a=0$, so the metric is Berwald even when $F(y)^2$ is not quadratic. Also, $\mathcal B$ and the Riemannian curvature $R^g$ of the slit fibre in Sections \ref{subsec:paper-2-1}--\ref{subsec:paper-2-4} are different tensors.
\end{remark}

\subsection{The Landsberg condition and symmetry of the Berwald tensor}\label{subsec:paper-1-4}
Lower the upper index of the Berwald curvature with the fundamental tensor, and put
\begin{equation}\label{eq:lowered-berwald}
 B_{\ell ijk}:=g_{\ell a}\mathcal B^a{}_{ijk}
              =g_{\ell a}G^a{}_{,ijk}.
\end{equation}
This tensor is always symmetric in its last three indices. The following characterization and its consequences follow from the classical contraction formula \eqref{eq:landsberg} and Euler's homogeneity identity; see \cite[Sections 6.1 and 6.2]{shen-2001}. The symmetry characterization is stated explicitly in \cite[Proposition 2.4]{crampin2011}.

\begin{proposition}[Symmetry characterization of the Landsberg condition]\label{prop:landsberg-symmetric-quartic}
Let $F$ be a smooth strongly convex Finsler metric. The following conditions are equivalent:
\begin{enumerate}[label=\textup{(\roman*)},leftmargin=2.3em,itemsep=2pt]
\item $F$ is Landsberg, that is, $L\equiv0$;
\item The full tensor $B$ defined in \eqref{eq:lowered-berwald} is completely symmetric on $TM\setminus\{0\}$.
\end{enumerate}
When these conditions hold, insertion of the reference vector $y$ into any slot of $B$ gives zero, and
\[
 B_{\ell ijk}=-E_aG^a{}_{,\ell ijk},\qquad
 B_{\ell ijk}(x,\lambda y)=\lambda^{-1}B_{\ell ijk}(x,y)\quad(\lambda>0).
\]
Consequently, on each entire slit fibre, $B$ is uniquely determined by its angular restriction to $\Sigma_x$.
\end{proposition}

The equivalence concerns the full tensor on the slit tangent bundle. Part I uses complete symmetry of its angular restriction, while radial vanishing permits recovery of the full tensor from that restriction.

\section{Landsberg angular equations}\label{sec:paper-5}\label{subsec:paper-6-1}\label{sec:twisted}
We verify the twisted Codazzi equation required by Theorem \ref{thm:quartic-rigidity-intro} and the algebraic compatibility identity it implies. Both identities use only smoothness.

Fix a base point $x$. We work with a general Landsberg metric, assuming neither flatness of the fibre Hessian metric nor the sphere commutativity condition of Sections \ref{sec:paper-3} and \ref{sec:paper-4}.

Let $F$ be a smooth strongly convex Landsberg metric, and put $V=T_xM$ and $E=F^2/2$. On $V\setminus\{0\}$, denote the fibre Hessian metric by $g$, the full Cartan tensor by $\mathcal C$, and the Levi-Civita connection by $\nabla^g$, corresponding to $g$. The Berwald tensor $B$ with its upper index lowered is defined as in \eqref{eq:lowered-berwald}.
On the indicatrix $\Sigma=\{y\in V:F(x,y)=1\}$, put $h=g|_{T\Sigma}$ and $C=\mathcal C|_{T\Sigma}$. Denote by $\nabla^h$ the Levi-Civita connection of $h$; the tangential restriction of $B$ is again denoted by $B$.

\begin{lemma}[Two angular identities for the Landsberg quartic form]\label{lem:B-angular-identities}
Under the hypotheses and notation above, the completely symmetric quartic form $B$ on the indicatrix satisfies the twisted Codazzi identity \eqref{eq:twisted}, and hence also the algebraic compatibility identity \eqref{eq:F}. Both identities hold on the entire slit fibre as well, with $C,\nabla^h,h$ replaced respectively by $\mathcal C,\nabla^g,g$. Thus the angular Berwald quartic form on each fixed fibre is a smooth twisted Codazzi quartic form in the sense of Definition \ref{def:twisted-quartic}.
\end{lemma}

\begin{proof}
Fix linear coordinates $y^1,\ldots,y^n$ on $V$. Write $\partial_i=\partial/\partial y^i$, use commas for fibre partial derivatives, and suppress the fixed base point $x$. Then
\[
 g_{ij}=E_{ij},\qquad
 \mathcal C_{ijk}=\tfrac12E_{ijk},\qquad
 \partial_e g_{ij}=2\mathcal C_{eij}.
\]
For the connection $\nabla^g$, Section \ref{sec:hess-curv} gives
\[
 \widetilde\Gamma^\ell_{ei}
 =g^{\ell p}\mathcal C_{eip}
 =\mathcal C_{ei}{}^\ell.
\]

By Proposition \ref{prop:landsberg-symmetric-quartic}, the full lowered Berwald tensor $B$ defined in \eqref{eq:lowered-berwald} is completely symmetric and satisfies $B_y(y,\cdot,\cdot,\cdot)=0$.

To eliminate the fourth fibre partial derivatives of the spray, interchange the first two slots using complete symmetry:
\[
 B_{aijk}=B_{iajk}=g_{i\ell}G^\ell{}_{,ajk}.
\]
The product rule gives
\[
 \begin{aligned}
 \partial_eB_{aijk}
 &=\partial_e\bigl(g_{i\ell}G^\ell{}_{,ajk}\bigr)=2\mathcal C_{ei\ell}G^\ell{}_{,ajk}
       +g_{i\ell}G^\ell{}_{,eajk},\\
 \partial_aB_{eijk}
 &=\partial_a\bigl(g_{i\ell}G^\ell{}_{,ejk}\bigr)=2\mathcal C_{ai\ell}G^\ell{}_{,ejk}
       +g_{i\ell}G^\ell{}_{,aejk}.
 \end{aligned}
\]
Since $G^\ell{}_{,eajk}=G^\ell{}_{,aejk}$, subtraction gives
\[
 \partial_eB_{aijk}-\partial_aB_{eijk}
 =2\mathcal C_{ei\ell}G^\ell{}_{,ajk}
  -2\mathcal C_{ai\ell}G^\ell{}_{,ejk}.
\]
Moreover, $B_{apjk}=g_{p\ell}G^\ell{}_{,ajk}$ implies
\[
 G^\ell{}_{,ajk}=g^{\ell p}B_{apjk}.
\]
Substitution and renaming of the contracted indices give
\[
 \partial_eB_{aijk}-\partial_aB_{eijk}
 =2\bigl(\mathcal C_{ei}{}^\ell B_{a\ell jk}
          -\mathcal C_{ai}{}^\ell B_{e\ell jk}\bigr).
\]
The factor $2$ comes from $Dg=2\mathcal C$.

For the connection calculation, put temporarily
\[
 S_i=\mathcal C_{ei}{}^\ell B_{a\ell jk}-\mathcal C_{ai}{}^\ell B_{e\ell jk},\qquad
 S_j=\mathcal C_{ej}{}^\ell B_{ai\ell k}-\mathcal C_{aj}{}^\ell B_{ei\ell k},\qquad
 S_k=\mathcal C_{ek}{}^\ell B_{aij\ell}-\mathcal C_{ak}{}^\ell B_{eij\ell}.
\]
These symbols denote only three expressions in the present calculation. Each still has five free indices $e,a,i,j,k$; they are not three components of a one-form.

The partial derivative difference $\partial_eB_{aijk}-\partial_aB_{eijk}$ obtained above is symmetric in $i,j,k$. Interchanging $i,j$ and using the complete symmetry of $B$ turns the right-hand side into $2S_j$; interchanging $i,k$ turns it into $2S_k$. Hence
\[
 \partial_eB_{aijk}-\partial_aB_{eijk}
 =2S_i=2S_j=2S_k,
 \qquad S_i=S_j=S_k.
\]

When covariantly differentiating a covariant four-tensor, one subtracts a connection term for each tensor slot. Thus
\[
 \begin{aligned}
 (\nabla^g_eB)_{aijk}
 ={}&\partial_eB_{aijk}
       -\mathcal C_{ea}{}^\ell B_{\ell ijk}-\mathcal C_{ei}{}^\ell B_{a\ell jk}
       -\mathcal C_{ej}{}^\ell B_{ai\ell k}
       -\mathcal C_{ek}{}^\ell B_{aij\ell},\\
 (\nabla^g_aB)_{eijk}
 ={}&\partial_aB_{eijk}
       -\mathcal C_{ae}{}^\ell B_{\ell ijk}-\mathcal C_{ai}{}^\ell B_{e\ell jk}
       -\mathcal C_{aj}{}^\ell B_{ei\ell k}
       -\mathcal C_{ak}{}^\ell B_{eij\ell}.
 \end{aligned}
\]
On subtraction, the connection terms from the first tensor slot give
\[
 -\mathcal C_{ea}{}^\ell B_{\ell ijk}
 +\mathcal C_{ae}{}^\ell B_{\ell ijk}=0,
\]
because $\mathcal C_{ea}{}^\ell=\mathcal C_{ae}{}^\ell$. The remaining three slots $i,j,k$ contribute $-S_i,-S_j,-S_k$, respectively, while the partial derivative terms give $2S_i$. Therefore
\[
 (\nabla^g_eB)_{aijk} -(\nabla^g_aB)_{eijk} =2S_i-S_i-S_j-S_k =-S_i.
\]
Expanding $S_i$, we obtain the differential equation on the full fibre:
\[
 (\nabla^g_eB)_{aijk} -(\nabla^g_aB)_{eijk} =\mathcal C_{ai}{}^\ell B_{e\ell jk} -\mathcal C_{ei}{}^\ell B_{a\ell jk}.
\]
This is the twisted Codazzi equation on the full fibre.

The proof of Lemma \ref{lem:quartic-algebraic-compatibility} uses only complete symmetry of $B$ and the twisted Codazzi equation. It therefore applies equally to $(g,\mathcal C,\nabla^g)$ and gives the algebraic compatibility identity on the entire slit fibre: equation \eqref{eq:F} with $C$ replaced by $\mathcal C$ and contractions taken with $g$. This is also the relation $S_i=S_j$ already obtained, so algebraic compatibility is not an additional assumption.

All internal contractions in these formulas use $g$. For example,
\[
 \mathcal C_{ei}{}^\ell B_{a\ell jk}
 =g^{pq}\mathcal C_{eip}B_{aqjk}.
\]
In general linear coordinates, two repeated lower indices cannot simply be summed as in Euclidean coordinates. Writing one of the indices raised keeps the inverse metric explicit.

We now restrict to the indicatrix, with the notation fixed above. By Proposition \ref{prop:hessian-cone}, for tangent vector fields $X,Y$ we have
\[
 \nabla^g_XY=\nabla^h_XY-h(X,Y)y.
\]
Substitute this identity in the definition of the covariant derivative of a quartic tensor. For arbitrary tangent vector fields $Y_1,\ldots,Y_4$, we obtain
\[
 (\nabla^g_XB)(Y_1,Y_2,Y_3,Y_4) =(\nabla^h_XB)(Y_1,Y_2,Y_3,Y_4) +\sum_{r=1}^4h(X,Y_r) B(Y_1,\ldots,\underbrace{y}_{r\text{-th slot}},\ldots,Y_4).
\]
Since $B(y,\cdot,\cdot,\cdot)=0$, every term in the sum vanishes. Thus taking the covariant derivative on the full fibre and then restricting gives the same result as restricting first to the indicatrix and then taking its covariant derivative.

The internal contractions also restrict exactly to angular contractions. At $y\in\Sigma$, choose a $g_y$-orthonormal basis of the full fibre, written $y,e_1,\ldots,e_m$, with each $e_\alpha\in T_y\Sigma$. For tangent vectors $X,Z,A,U,V$, for example,
\[
 \begin{aligned}
 \langle\mathcal C(X,Z,\cdot),B(A,\cdot,U,V)\rangle_{g^{-1}}
 &=\mathcal C(X,Z,y)B(A,y,U,V)+\sum_{\alpha=1}^m
       C(X,Z,e_\alpha)B(A,e_\alpha,U,V)\\
 &=\sum_{\alpha=1}^m
       C(X,Z,e_\alpha)B(A,e_\alpha,U,V).
 \end{aligned}
\]

The radial term vanishes, and the last line is contraction with $h^{-1}$. The full-fibre twisted Codazzi and algebraic compatibility identities therefore restrict to \eqref{eq:twisted} and \eqref{eq:F}, respectively.

\end{proof}

For the coordinate-free form of the twisted Codazzi equation, see Definition \ref{def:twisted-quartic}.

\begin{remark}[Scope of the angular equations]
Lemma \ref{lem:B-angular-identities} is a local smooth result within a fixed fibre. It requires neither fibre analyticity, flatness of the Hessian metric, compactness of the indicatrix, nor periodicity. The Landsberg condition ensures complete symmetry of $B$, while commutation of the fourth fibre partial derivatives of the spray supplies the initial differential relation. Equation \eqref{eq:twisted} controls the angular variation of $B$, and \eqref{eq:F} gives algebraic compatibility between $C$ and $B$ at a point. If $C=0$ at that point, the algebraic condition \eqref{eq:F} reduces to $0=0$, whereas the differential equation still says that the corresponding difference of covariant derivatives vanishes there. The analytic rigidity theorem of Part I uses the twisted Codazzi equation on the entire indicatrix directly; its degenerate branch concludes by vanishing of the Cartan tensor on a nonempty open set.
\end{remark}

\section{Proof of the main theorem}\label{sec:paper-9}
We first apply the two rigidity theorems of Part I to a fixed fibre. We then show that the open set where the Berwald curvature is nonzero must be empty.

\subsection{Application of the two rigidity theorems to a fixed fibre}\label{sec:fixed-fiber-application}
\begin{proposition}[A nonzero Berwald quartic form forces the fixed fibre to be quadratic]\label{prop:nonzero-fiber-quadratic}
Let $F$ be a smooth strongly convex Landsberg metric, regular on the entire slit tangent bundle of a smooth manifold of dimension $n\ge4$. Fix a base point $x$, and suppose that $F_x$ is real-analytic on the whole of $T_xM\setminus\{0\}$. If the angular Berwald quartic tensor of this fibre satisfies $B_x\not\equiv0$, then there is a positive definite symmetric bilinear form $q_x$ such that
\[
 F_x(y)^2=q_x(y,y)\qquad(y\ne0).
\]
\end{proposition}
\begin{proof}
Fix $V=T_xM$ and its indicatrix $\Sigma_x=\{y\in V:F_x(y)=1\}$. By Proposition \ref{prop:landsberg-symmetric-quartic}, $B_x$ is a smooth completely symmetric quartic form on the entire indicatrix; Lemma \ref{lem:B-angular-identities} shows that it satisfies the twisted Codazzi identity \eqref{eq:twisted}. Thus regularity of $F_x$ on the entire slit space, its strong convexity and real analyticity, together with $n\ge4$ and $B_x\not\equiv0$, verify all hypotheses of Theorem \ref{thm:quartic-rigidity-intro}. That theorem gives
\[
 R^{h_x}=R_{0,x},\qquad R^{g_x}=0,
 \qquad g_x=D^2(F_x^2/2).
\]
The latter curvature identity holds on the whole of $T_xM\setminus\{0\}$. Since $n\ge4$ also satisfies the condition $n\ge3$ required for analytic Laugwitz rigidity, Theorem \ref{thm:laugwitz} applies. We obtain a positive definite symmetric bilinear form $q_x$ such that $F_x(y)^2=q_x(y,y)$.
\end{proof}

These applications require analyticity only of the fixed-fibre norm, not of the base manifold or the Berwald quartic form. The conclusion under $B_x\not\equiv0$ does not yet eliminate Berwald curvature, which involves variation of the base point.

\subsection{Completion of the proof of the main theorem}\label{subsec:paper-12-4}\label{sec:main}

\begin{proof}[Proof of Theorem \ref{thm:main}]
Let
\[
 \mathscr U:=\{x\in M:\mathcal B(x,\cdot)\not\equiv0\}.
\]
This set is open. Indeed, if a component of $\mathcal B$ is nonzero at $(x_0,y_0)$, extend $y_0$ locally by keeping its fibre coordinates fixed. Smoothness ensures that the same component remains nonzero near $x_0$.

Since $g$ is invertible, lowering the upper index preserves nonvanishing. By Proposition \ref{prop:landsberg-symmetric-quartic}, radial vanishing and homogeneity then show that $\mathcal B(x,\cdot)\not\equiv0$ is equivalent to nonvanishing of the angular quartic form $B_x$. Thus Proposition \ref{prop:nonzero-fiber-quadratic} makes every fibre over $\mathscr U$ quadratic.

Suppose that $\mathscr U$ is nonempty. Choose a nonempty coordinate neighbourhood $W\subset\mathscr U$ and a nowhere vanishing smooth section $v_*$ of $TW$. Put
\[
 q_{ij}(x):=g_{ij}(x,v_*(x)).
\]
Since every fibre over $W$ is quadratic, its fundamental tensor is independent of the fibre variable. Hence
\[
 g_{ij}(x,y)=q_{ij}(x),\qquad F(x,y)^2=q_{ij}(x)y^iy^j\qquad(x\in W,\ y\ne0).
\]
The coefficients $q_{ij}$ are smooth, symmetric, and positive definite. Therefore $F|_W$ is Riemannian, and hence Berwald \cite[Section 4.3]{chern-shen-2005}. This contradicts $W\subset\mathscr U$. Thus $\mathscr U=\varnothing$, proving $\mathcal B\equiv0$.
\end{proof}

\begin{remark}[Propagation of quadratic fibres]\label{rem:quadratic-component-propagation}
Theorem \ref{thm:laugwitz} and classical Landsberg parallel translation also give a propagation statement in dimension $n\ge3$. Let $F$ be smooth, strongly convex and Landsberg on a connected component $M_0$, regular on its entire slit tangent bundle and real-analytic on every slit fibre. If $F_{x_0}^2$ is a positive-definite quadratic form at one point $x_0\in M_0$, then $F|_{M_0}$ is Riemannian.

Indeed, canonical parallel translation along a finite piecewise smooth path is a diffeomorphism of the entire slit fibres preserving their Hessian metrics \cite[Section 4.2, Proposition 4.4.1]{chern-shen-2005}; see also \cite[Section 4]{li-minkowski}. Since every point of $M_0$ can be joined to $x_0$ by such a path, flatness of the initial Hessian metric implies flatness on every target fibre. Theorem \ref{thm:laugwitz} makes each target norm quadratic. Evaluation of its fundamental tensor along local nonzero sections, as in the preceding proof, gives smooth Riemannian metrics; they agree on overlaps because their squared norms are all $F^2$. No completeness of the base manifold is required.

This also yields another proof of the main theorem. Proposition \ref{prop:nonzero-fiber-quadratic} would turn a fibre with $B_x\not\equiv0$ into a quadratic fibre. Its connected component would then be Riemannian and hence Berwald, a contradiction.
\end{remark}

\medskip
\noindent\textbf{Acknowledgments:} This study was funded by the National Key R \& D Program of China (Grant No. 2024YFA1013201) and the Science and Technology Research Program of Chongqing Municipal Education Commission (Grant No. KJZD-K202301110). The authors used ChatGPT 6 Astra to draft and revise text, to develop and check proofs, and to compare the paper's statements about the cited literature with the sources. The authors reviewed and verified all such output and take full responsibility for the content of this paper.

\end{document}